%% file: main.tex
\documentclass[review,onefignum,onetabnum]{macros} % 
\nolinenumbers
\title{An Adaptive CUR Algorithm and its Application to Reduced-Order Modeling of Random PDEs}

\author{%
  Grishma Palkar\thanks{Department of Mechanical Engineering, University of Pittsburgh, Pennsylvania, USA.}
  \and
  Hessam Babaee\footnotemark[1]\thanks{Author for correspondence: \email{h.babaee@pitt.edu}.}
}

\begin{document}
\maketitle

\begin{abstract}
Certain classes of CUR algorithms, also referred to as cross or pseudoskeleton algorithms, are widely used for low-rank matrix approximation when direct access to all matrix entries is costly. Their key advantage lies in constructing a rank-r approximation by sampling only r columns and r rows of the target matrix. This property makes them particularly attractive for reduced-order modeling of nonlinear matrix differential equations, where nonlinear operations on low-rank matrices can otherwise produce high-rank or even full-rank intermediates that must subsequently be truncated to rank $r$. CUR cross algorithms bypass the intermediate step and directly form the rank-$r$ matrix. However, standard cross algorithms may suffer from loss of accuracy in some settings, limiting their robustness and broad applicability. In this work, we propose a cross oversampling algorithm that augments the intersection with additional sampled columns and rows. We provide an error analysis demonstrating that the proposed oversampling improves robustness. We also present an algorithm that adaptively selects the number of oversampling entries based on efficiently computable indicators. We demonstrate the performance of the proposed CUR algorithm for time integration of several nonlinear stochastic PDEs on low-rank matrix manifolds.

  % The CUR decomposition is a powerful tool for solving large matrix differential equations (MDEs) on low-rank manifolds. CUR achieves computational efficiency by sampling a subset of rows and columns from the solution matrix.  However, as the approximation rank increases, CUR decomposition can become less accurate, leading to numerical instabilities. This instability can be mitigated by oversampling, which involves increasing the number of rows or columns sampled to enhance the numerical stability. Despite its benefits, oversampling introduces significant computational challenges, as selecting and processing additional rows and columns can greatly increase computational cost, particularly in large-scale applications. Striking the right balance between oversampling and computational efficiency is crucial to ensuring scalability and cost effectiveness. In this work, we address these challenges by proposing a scalable CUR decomposition framework. Our approach introduces a cost-effective oversampling strategy that selectively focuses on cross entries instead of entire rows and columns, significantly reducing computational overhead. By overcoming the limitations of conventional CUR oversampling methods, our approach effectively balances accuracy, stability, and computational efficiency, making it a practical and reliable solution for large-scale, high-dimensional applications.
\end{abstract}

% REQUIRED

% % REQUIRED
% \begin{AMS}

% \end{AMS}

% Site to get AMS subject classification: https://mathscinet.ams.org/msnhtml/msc2020.pdf 

\section{Introduction}\label{sec:intro}
Matrix low-rank approximations are essential in a wide range of science and engineering problems, from data compression and large-scale data analysis to reduced-order modeling and scientific computing \cite{UT19}. They also serve as building blocks for many tensor network low-rank approximations \cite{B23}. Efficiently constructing such approximations is therefore of great interest to numerous downstream applications.

% A key example, and the focus of this work, is the dynamical low-rank approximation (DLRA) of nonlinear partial differential equations, in which the semi-discrete PDE is recast as a matrix or tensor differential equation (MDE or TDE) and solved in low-rank form. Numerical efficiency in these applications hinges on performing all operations—linear or nonlinear—directly in compressed form. For instance, the PDE may involve fractional nonlinearities, requiring the computation at each time step of terms such as $1/\mathbf{A}$ (with division taken entrywise), where $\mathbf{A}$ is stored in low-rank form but $1/\mathbf{A}$ is full-rank. Computing the low-rank approximation of$1/\mathbf{A}$ via SVD would require access to all entries of the matrix—i.e., decompression—thereby negating the computational gains that DLRA is designed to deliver.

The singular value decomposition (SVD) provides the optimal low-rank approximation, and scalable algorithms exist for computing it even for very large matrices. In many applications, however, the dominant cost is not in computing the SVD itself but in acquiring the matrix entries needed for the SVD. 

One such case, and the focus of this work, is reduced-order modeling of partial differential equations (PDEs) via the dynamical low-rank approximation (DLRA) \cite{koch2007dynamical}. In DLRA, the semi-discrete PDE is recast as a matrix differential equation (MDE) in the form $\dot{\mathbf A} = \mathcal{F}(\mathbf A)$ and then solved in a compressed form, where $\mathbf A \in \mathbb{R}^{n\times s}$ is the solution matrix and $ \mathcal{F}(\mathbf A)$ is the right-hand side matrix.  Achieving computational efficiency in DLRA requires performing all linear and nonlinear operations involved in $\mathcal{F}(\mathbf A)$ without decompression by exploiting the fact that matrix $\mathbf A$ is available in a low-rank compressed form. However, $\mathcal{F}(\mathbf A)$  can be a very high-rank or even a full-rank matrix even when $\mathbf A$ is low-rank.  For example, fractional nonlinearities may require computing $(1/\mathbf A)_{ij} = 1/\mathbf A_{ij}$  (entrywise) when $\mathbf{A}$ is stored in a low-rank form. However, $1/\mathbf A$  is full-rank despite $\mathbf A$ being low rank (if $\mbox{rank}(\mathbf A)>1$). Thus, obtaining a low-rank representation of $\mathcal{F}(\mathbf A)$ via SVD in such cases demands full access to the matrix entries, after which the SVD can be computed. This access to all entries of $\mathbf A$ would erode the very efficiency that the reduced-order modeling seeks to achieve.

In this setting, the target matrix is $\mathcal{F}(\mathbf{A})$, and the objective is to construct an accurate low-rank approximation of $\mathcal{F}(\mathbf{A})$ while evaluating as few of its entries as possible. This requires robust data-efficient techniques where the low-rank approximation of $\mathcal{F}(\mathbf{A})$ is constructed by minimally sampling it while satisfying the stringent error tolerances often required in the numerical solution of PDEs. As we explain in \S \ref{sec:Intro_CUR} and \S \ref{sec:Intro_DLRA}, certain variants of CUR-based low-rank approximations are particularly well-suited to this task and have recently been applied to solving PDEs in low-rank form.

% This, in turn, requires constructing accurate low-rank approximations under the tight error tolerances necessary for efficient PDE solvers.        \\

\subsection{CUR Decompositions}\label{sec:Intro_CUR}
 The CUR decompositions became well-known after the work of Goreĭnov et al. \cite{goreinov1997theory,GT01,GO10}, although they have existed since the 1950s. See \cite{HH20} for historical notes. The CUR decompositions have received renewed attention in recent years, particularly through the work of Drineas et al. \cite{doi:10.1073/pnas.0803205106,doi:10.1137/07070471X} and others \cite{OZ18,HH21}. One of the key reasons for the popularity of CUR decompositions is their interpretability. A CUR decomposition represents the low-rank approximation using actual rows and columns of the target matrix, rather than the eigenrows and eigencolumns produced by SVD. This property is particularly valuable in data analysis, where retaining the original data features is often desirable because they are interpretable. The CUR low-rank approximations also inherit certain properties of the target matrix, for example, sparsity, non-negativity, and integer-valued entries. These structures are lost in the SVD low-rank approximation as the singular vectors are the linear combination of all columns and rows. 

CUR decompositions have different types with different accuracies and different data dependencies. A rank-$r$ CUR decomposition of matrix $\mathbf A \in \mathbb{R}^{n \times s}$ can be represented as 
\begin{equation}
\mathbf A \approx \mathbf C \mathbf Z \mathbf R,
\end{equation}
where  $\mathbf C \in \mathbb{R}^{n\times r}$ and $\mathbf R \in \mathbb{R}^{r\times s}$ are actual subsets of columns and rows of matrix $\mathbf A$ and $\mathbf Z \in \mathbb{R}^{r\times r}$ matrix. Different CUR decompositions can be obtained depending on (i) how $\mathbf Z$ is calculated, and (ii) which columns and rows are selected. We first discuss the common choices for calculating $\mathbf Z$. 

 The first choice is to take $\mathbf{Z}$ as the inverse of the column-row intersection matrix as follows: 
\begin{equation}\label{eq:Z_cross}
\mathbf{Z} = \mathbf{A}(\mathbf{p},\mathbf{s})^{-1},
\end{equation}
  where $\mathbf{p} = [p_1, p_2, \dots, p_r]$ is the vector of selected row indices and $\mathbf{s} = [s_1, s_2, \dots, s_r]$ is the vector of selected column indices. This form of CUR is also known as a \emph{(pseudo-)skeleton} or \emph{cross} approximation and was analyzed by Goreĭnov et al.~\cite{goreinov1997theory}. We refer to this CUR scheme as CUR-Cross or \texttt{CUR-CR}.  The second common choice for $\mathbf{Z}$ is obtained by orthogonal projection of $\mathbf{A}$ onto the selected columns and rows as follows
\begin{equation}\label{eq:Z_opt}
\mathbf{Z} = \mathbf{C}^{\dagger} \mathbf{A} \mathbf{R}^{\dagger}.
\end{equation} 
For given $\mathbf{C}$ and $\mathbf{R}$, this choice of $\mathbf{Z}$ yields the optimal CUR approximation in the Frobenius norm \cite{S99}, i.e., it minimizes $\| \mathbf{A} - \mathbf{C} \mathbf{Z} \mathbf{R} \|_F$. We refer to this CUR scheme as \texttt{CUR-opt}.

Beyond the choice of how $\mathbf{Z}$ is computed, another key factor influencing the quality of a CUR approximation is the selection of the row and column index sets, $\mathbf{p}$ and $\mathbf{s}$. There are two broad categories of sampling techniques: pivot-based methods and probabilistic sampling methods. Some of the prominent pivot-based methods include maximum volume (maxvol) \cite{GO10}, discrete empirical interpolation method (DEIM)\cite{sorensen2016deim}, and QDEIM \cite{DG16}.
 As demonstrated in \cite{goreinov1997theory}, the accuracy of \texttt{CUR-CR} is dependent on the determinant of the intersection submatrix ($\mathbf Z$). This determinant is often referred to as the matrix intersection volume. For better accuracy, rows and columns should be chosen to maximize the matrix volume \cite{GO10}.  Sorensen et al.~\cite{sorensen2016deim} proposed using the DEIM algorithm \cite{CS10} for column-row selection. They showed that the low-rank approximation error of \texttt{CUR-opt} depends on the norm of the inverse of the selected subcolumn and subrow matrices. This observation motivates the use of DEIM, since its greedy strategy explicitly seeks to minimize the norm of submatrices.

Sampling-based methods select indices from probability distributions derived from structural properties of the matrix. A common choice is leverage scores based on the row norms of the dominant singular vectors \cite{doi:10.1137/07070471X,doi:10.1073/pnas.0803205106}. Other approaches include uniform sampling \cite{doi:10.1137/110852310}, volume sampling \cite{cortinovis2020low,deshpande2006matrix,deshpande2006adaptive,goreinov1997theory}, determinantal point process (DPP) sampling \cite{derezinski2021determinantal}, and Batson-Spielman-Srivastava (BSS) sampling \cite{batson2009twice,boutsidis2014near}.
 Hybrid strategies such as L-DEIM \cite{gidisu2021hybrid} and deterministic leverage score methods \cite{papailiopoulos2014provable} have been proposed to combine the benefits of both pivoting and sampling.

% Among these, volume sampling has been shown to yield CUR approximations with near-optimal error guarantees \cite{cortinovis2020low,zamarashkin2018existence}. To reduce cost, sampling can also be performed on sketch matrices \cite{drineas2012fast}.

Although the \texttt{CUR-opt} scheme is optimal in terms of accuracy,  it is not data efficient because it requires access to all entries of $\mathbf A$. Therefore, for problems where the cost of computing the entries of the target matrix dominates the overall cost, this scheme is not efficient.  On the other hand, the \texttt{CUR-CR} scheme is data efficient as it can produce a low-rank approximation of $\mathbf{A}$ by accessing only  $r(n+s)-r^2$ entries---those in the selected columns and rows. However,  \texttt{CUR-CR} can become unstable, especially as the rank increases \cite{PN25}. This is due to the fact that the matrix $\mathbf A(\mathbf p, \mathbf s)$ becomes ill-conditioned as the rank increases. To this end, both the pseudoinverse \cite{dong2023simpler, martinsson2020randomized, sorensen2016deim} and the $\epsilon$-pseudoinverse \cite{PN25}---where singular values of $\mathbf{A}(\mathbf{p}, \mathbf{s})$ smaller than $\epsilon$ are truncated prior to inversion---have been proposed. It was shown in \cite{TDBCUR} that performing a QR decomposition of the selected columns (or rows) and inverting the submatrix of the orthonormalized columns can partially alleviate ill-conditioning even when $\mathbf{A}(\mathbf{p}, \mathbf{s})$ is singular,  such as in over-rank approximation. However, as we demonstrate in this paper, numerical stability in \texttt{CUR-CR} is not fully resolved even with QR, since the condition number of the submatrix of the orthonormalized columns can still be large---though no longer unbounded. This can be particularly problematic when CUR is used for reduced-order modeling of PDEs, which require strict accuracy thresholds.

A reasonable compromise between \texttt{CUR-CR} and \texttt{CUR-opt}  is to enhance the stability of \texttt{CUR-CR} via oversampling.  The existing approaches involve sampling larger than rank ($r$) columns and rows. The number of oversampled columns and rows might be different, i.e., $\mathbf C \in \mathbb{R}^{n \times (r+m_c)}$ and $\mathbf R \in \mathbb{R}^{ (r+m_r)\times  s}$, where $m_c$ and $m_r$ are the number of oversampled columns and rows, respectively. This approach has been used in conjunction with randomized algorithms. See for example, \cite{doi:10.1137/07070471X,doi:10.1073/pnas.0803205106}. Oversampling of rows and columns results in $(r+m_c)\times (r+m_r)$ intersection matrices. One approach to obtain a rank-$r$ approximation is to perform SVD on this matrix and truncate at rank-$r$ \cite{OZ18}. Another approach is to first reduce the dimensions of the selected columns and rows to $r$ via SVD and then obtain $\mathbf Z$ via the orthogonal projection of the target matrix onto these rank-$r$ spaces \cite[Section 4.2.2]{HH21}. However, this approach would require access to all entries of the target matrix.  We refer the reader to \cite{PN25} for a comprehensive analysis of various oversampling and regularization CUR techniques.  

 \subsection{Dynamical Low-Rank Approximation}\label{sec:Intro_DLRA}
The DLRA is a powerful computational method for reduced order modeling of generic matrix and tensor differential equations \cite{koch2007dynamical, koch2010dynamical}.  The DLRA is closely related to the Dirac-Frankel variational principle introduced in quantum chemistry for solving the Schrödinger equation \cite{beck2000multiconfiguration}. The DLRA reduces the dimensionality of MDEs/TDEs by explicitly formulating evolution equations in the low-rank factorized form, thereby reducing the computational complexity, where the low-rank structures are represented by time-dependent bases (TDBs).

The key advantage of the DLRA is its ability to adapt in real time to changes in system dynamics. Another significant benefit is that it does not require a separate offline phase to collect data for constructing the low-rank subspace, unlike many existing ROM techniques, such as those based on proper orthogonal decoposition (POD). Instead, in the DLRA formulation, explicit evolution equations for the TDBs are directly derived from the full-order model (FOM). This approach eliminates the need for precomputed snapshots, enabling DLRA to function as an adaptive, real-time ROM that adjusts to dynamic changes on the fly. This also makes DLRA more computationally costly to solve than POD-based or similar data-driven ROMS. However, DLRA can tackle a much wider range of problems, for convection-dominated problems and turbulent flows.  There are other dimension reduction methods and ROMs based on TDB that differ from the DLRA framework (e.g.,\cite{Towne_2018,peherstorfer2020model, huang2023predictive, PADOVAN2024112597}). %However, these approaches are not further addressed in this paper, as the primary focus is on DLRA-based ROMs.

 Recently, the DLRA and it's closely related variations have been employed to solve a growing list of diverse problems including random PDEs \cite{sapsis2009dynamically,cheng2013dynamically,musharbash2018dual,patil2020real}, Boltzmann transport and Vlasov equations \cite{kusch2023robust,ceruti2022unconventional,einkemmer2019quasi,hu2022adaptive}, turbulent combustion \cite{ramezanian2021fly,JLBC25}, shallow water equations \cite{koellermeier2024macro}, control \cite{blanchard2019control}, linear sensitivity analysis \cite{donello2022computing}, chemical kinetics \cite{nouri2022skeletal,nouri2024skeletal}, hydrodynamic stability analysis \cite{babaee2016minimization,kern2021transient,kern2024onset}, and deep learning \cite{savostianova2024robust,schotthofer2022low}.

While solving MDEs in the low-rank form can result in significant computational savings for linear or quadratically nonlinear problems, the cost of solving these equations can increase substantially for MDEs with high-order polynomial nonlinearity or non-polynomial nonlinearities, such as exponential or fractional nonlinearities. If these types of nonlinearities exist in the MDEs,  the computational cost of solving these equations in the low-rank form can become as expensive as solving the FOM. This is because in solving the MDE $\dot{\mathbf A} = \mathcal F(\mathbf A)$, when matrix $\mathbf A \in \mathbb{R}^{n\times s}$ is approximated with a rank-$r$ matrix $\hat{\mathbf A} \in \mathbb{R}^{n\times s}$, the rank of matrix $\mathcal{F}(\hat{\mathbf A}) \in \mathbb{R}^{n\times s}$ increases with $\mathcal{O}(r^p)$ when $\mathcal F$ has polynomial nonlinearity of order $p$ and $\mathcal{F}(\hat{\mathbf A})$ becomes full-rank if $\mathcal F$ has non-polynomial nonlinearity. When $\mathcal{F}(\hat{\mathbf A})$ is full-rank, all the entries of this matrix must be explicitly computed, leading to  floating-point operation (FLOP) costs comparable to those of the FOM. This issue is analogous to reduced-order modeling of nonlinear PDEs with POD, where handling the nonlinear terms in the reduced form can become as costly as solving the FOM \cite{CS10}.

% Recently, a new sparse interpolation algorithm was introduced to solve TDB ROM evolution equations. This algorithm significantly reduces computational complexity for general nonlinear stochastic partial differential equations (SPDEs) \cite{NADERI2023115813}. The key to this method is the use of an interpolatory CUR decomposition, which constructs a low-rank approximation of $\mathcal F(\hat{\mathbf A})$ by selecting specific columns and rows. The selection of these rows and columns is performed using the discrete empirical interpolation method (DEIM) \cite{CS10}, which ensures that the approximation remains computationally efficient while preserving the essential dynamics of the system. 

To address the computational cost issues, a CUR-based approach was presented \cite{NADERI2023115813}, where a cross approximation of $\mathcal{F}(\hat{\mathbf A})$ is constructed. More recently,  in \cite{TDBCUR}, a new algorithm was introduced that applied CUR decomposition to time-discrete MDEs. This algorithm uses \texttt{CUR-CR} low-rank approximation based on DEIM in conjunction with row oversampling and is particularly well-suited for solving nonlinear MDEs on low-rank manifolds. One of the key advantages of this algorithm is its ability to handle arbitrarily nonlinear terms in a cost-effective manner. This approach has been extended to higher-order implicit time integration schemes \cite{NAB25} as well as solving nonlinear TDEs in Tucker tensor \cite{GB24} and tensor train \cite{GBTT24,D25} forms. These techniques have recently been used for solving the Boltzmann equation \cite{AE25} and turbulent combustion \cite{JLBC25}. 

% The algorithm avoids the inversion of singular values, which is a source of instability when standard time integration schemes are used for solving DLRA equations. Doing so ensures robustness even in the presence of small or zero singular values. Additionally, the algorithm achieves high-order temporal accuracy. To further improve the accuracy and efficiency of the CUR approximation, \cite{TDBCUR} also incorporated a rank-adaptive approach and a row oversampling strategy. These techniques allow the algorithm to dynamically adjust the rank of the approximation during computation, ensuring that the low-rank model remains accurate as the system evolves. Row oversampling improves the accuracy of CUR low-rank approximation.

The CUR decomposition in \cite{TDBCUR} requires sampling entire rows of the matrix, which causes the cost of data acquisition to scale with the number of columns. Moreover, the method in \cite{TDBCUR} employs only row oversampling. In general, however, oversampling both rows and columns is necessary to ensure the stability of the CUR algorithm. While it is possible to incorporate column oversampling in a manner similar to \cite{TDBCUR} - by sampling additional rows and columns followed by truncation to a rank-$r$ approximation-this approach effectively incurs the computational cost of solving a significantly higher-rank problem. Such overhead becomes particularly burdensome for high-dimensional systems.
% \begin{enumerate}
%     \item The method accounts for row oversampling but does not address column oversampling.
    
%     \item The proposed row oversampling strategy requires sampling entire rows of the matrix, causing the cost of data acquisition to scale with the size of columns. This becomes particularly expensive for high-dimensional systems.
    
%     \item The number of oversampling points significantly affects the stability of the CUR decomposition. In Donello et al.~\cite{TDBCUR}, it is treated as a fixed hyperparameter that must be predetermined. Allowing the number of oversampling points to adapt dynamically during computation could improve the method’s accuracy.
% \end{enumerate}

\subsection{Contributions}
As demonstrated in this paper, it has become clear that oversampling is essential for enhancing the stability of the CUR decomposition. It should be considered a fundamental component of a robust cross algorithm, not merely an optional enhancement. However, oversampling necessitates access to additional entries of the target matrix, which, when solving MDEs on low-rank manifolds, results in increased computational cost.  Achieving a balance between robustness and computational efficiency is essential for scalability. Our work addresses these challenges by introducing an adaptive cost-efficient oversampling strategy integrated into the CUR decomposition, featuring the following elements:
\begin{enumerate}
    \item We propose a cross CUR scheme that oversamples only the row-column intersection entries, rather than entire rows and columns, thereby requiring fewer matrix entries to achieve a given error threshold. The scheme is agnostic to the choice of sampling strategy, accommodating both probabilistic and deterministic pivot-based methods. %The proposed CUR scheme shows good experimental performance and, in some cases, outperforms the existing oversampling algorithms.
    \item We present a deterministic error bound in the spectral norm. Building on this bound, we propose an adaptive oversampling strategy in which the number of oversampling points is determined dynamically, eliminating the need for the practitioner to specify it in advance.
    \item We demonstrate the utility of the cross CUR scheme in the time integration of random PDEs on low-rank matrix manifolds for a diverse set of PDEs.
    % Cost-effective oversampling strategy: We propose a novel approach that focuses on sampling cross entries rather than entire rows and columns. This significantly reduces computational overhead while preserving the quality of the approximation.

% \item Rank-adaptivity mechanism: We develop a dynamic rank-adaptive scheme that minimizes regression error to adjust the rank of the low-rank approximation during computation. This ensures that the rank evolves optimally, enhancing accuracy, particularly in high-dimensional or nonlinear systems.

% \item Hyperparameter sensitivity analysis: To address the reliance on finely tuned hyperparameters in oversampling method, we perform a detailed sensitivity analysis. This minimizes dependence on hyperparameters, making the framework more robust and easier to implement.
\end{enumerate}

% \hb{check this synopsis after revision:} 
The remainder of this paper is structured as follows: \S \ref{sec:meth} provides an overview of DLRA  and CUR decomposition for time-integration of MDEs on a low-rank manifold. We present an error bound for the resulting low-rank approximation. We then introduce a cross-oversampling algorithm, with rank adaptivity and adaptive sampling strategies to solve MDEs on a low-rank manifold. \S \ref{sec:demo} demonstrates the effectiveness of our approach through a series of numerical experiments. These include a toy problem, applications to several highly nonlinear stochastic PDEs such as the Burgers equation, the Allen-Cahn equation, and the Korteweg-de Vries equation, as well as applications to heat conduction and radiation problems on a rectangular geometry with a finite element mesh. Finally, \S \ref{sec:con} concludes the paper with a summary of our findings and a discussion of their broader implications.

\section{Methodology}\label{sec:meth}
\subsection{Preliminaries}

We begin by introducing the notation used throughout this paper. Let $\mathbf{p} = [p_1, p_2, \dots, p_r]$ denote the indices of the selected rows, and $\mathbf{s} = [s_1, s_2, \dots, s_r]$ the indices of the selected columns. Moreover, let $\overline{\mathbf{p}} = [p_1, p_2, \dots, p_r,$ $ p_{r+1}, \dots, p_{r+m_r}]$ denote the indices of the selected rows together with $m_r$ oversampling rows, and similarly $\overline{\mathbf{s}} = [s_1, s_2, \dots, s_r, s_{r+1}, \dots, s_{r+m_c}]$ denote the selected columns together with $m_c$ oversampling columns. For a matrix $\mathbf{A} \in \mathbb{R}^{n \times s}$, we adopt MATLAB-style indexing: $\mathbf{A}(\mathbf{p}, :)$ denotes the submatrix consisting of rows indexed by $\mathbf{p}$, and $\mathbf{A}(:, \mathbf{s})$ denotes the submatrix consisting of columns indexed by $\mathbf{s}$. Vectors are denoted in lowercase boldface (e.g., $\mathbf{a} \in \mathbb{R}^n$), and matrices in uppercase boldface (e.g., $\mathbf{A} \in \mathbb{R}^{n \times s}$). %Low-rank approximations are denoted with a hat symbol (e.g., $\hat{\mathbf{V}}$).
We denote the Moore-Penrose pseudoinverse of $\mathbf A$ with $\mathbf A^{\dagger}$. The identity matrix of size $n \times n $ is shown with $\mathbf I_n$. The indexing matrices can be defined: $\mathbf P = \mathbf I_n(:,\mathbf p) \in \mathbb{R}^{n \times r}$, $\ol{\mathbf P} = \mathbf I_n(:,\ol{\mathbf p})\in \mathbb{R}^{n \times (r+m_r)} $, $\mathbf S = \mathbf I_s(:,\mathbf s) \in \mathbb{R}^{s \times r}$, and $\ol{\mathbf S} = \mathbf I_s(:,\ol{\mathbf s}) \in \mathbb{R}^{s \times (r+m_c)}$.

Now, we present the definition of a low-rank matrix manifold.
\begin{definition}[Low-rank matrix manifold]
 The low-rank matrix manifold $\mathcal{M}_r$ is defined as the set of all rank $r$ matrices. 
\begin{equation*}
\mathcal{M}_r = \{\hat{\mathbf{A}} \in \mathbb{R}^{n \times s}: \ \mbox{rank}(\hat{\mathbf{A}}) = r \}, 
\end{equation*}
 Any member of the set $\mathcal{M}_r$ will be denoted by a hat symbol $( \hat{ \ \ } )$, e.g. $\hat{ \mathbf{A}}$.
\end{definition}

\subsection{Problem setup}
We consider nonlinear SPDEs in the form of:
\begin{equation}\label{eq:FOM_Cont}
\frac{\partial v}{\partial t} = f(v;x,\bs \xi),
\end{equation}
with appropriate initial and boundary conditions. In the above equation, $v=v(x;t,\bs\xi)$ where $x$ represents the spatial coordinate, $\bs \xi \in \mathbb{R}^d$ denotes the vector of random parameters, $t$ is time, and $f(v;x,\bs\xi)$ is the nonlinear differential operator.

Discretizing  Eq. \eqref{eq:FOM_Cont} in $x$ and $\bs\xi$ results in the following nonlinear MDE:
\begin{equation}\label{eq:FOM}
\frac{\mathrm d \mathbf A}{\mathrm d t} = \mathcal F(\mathbf{A}), \quad t\in I
\end{equation}
where the time interval is given by $I=[0,T_f]$, $\mathbf {A}(t) \in \mathbb{R}^{n \times s}$ is a time-dependent matrix, where $n$ is the spatial degrees of freedom and $s$ is the number of random samples. The map $\mathcal F(\mathbf{A}):  \mathbb{R}^{n\times s} \rightarrow \mathbb{R}^{n\times s}$ is obtained by discretizing  $f(v;x,\bs\xi)$ in $x$ and $\bs\xi$. Eq. \eqref{eq:FOM} is supplemented with appropriate initial conditions, namely $\mathbf{A}(t_0)=\mathbf{A}_0$. Furthermore, we assume that boundary conditions have already been incorporated into Eq. \eqref{eq:FOM}. We will refer to Eq. \eqref{eq:FOM} as the FOM. Both $n$ and $s$ are typically large, and as a result, solving the MDE Eq. \eqref{eq:FOM} scales at least as $\mathcal{O}(ns)$ both in terms of memory and FLOPs.

\subsection{Dynamical low-rank approximation for time integration of MDEs on low-rank manifolds}

 The DLRA enables approximating generic MDEs on low-rank matrix manifolds~\cite{koch2007dynamical}. The DLRA formulation is equivalent to the dynamically orthogonal (DO) decomposition \cite{SL09}. In the  DLRA formulation,  $\mathbf{A}(t)$ is approximated at each instant by a low-rank matrix
 $\hat{\mathbf{A}}(t)$ of the form
\begin{equation}
    \hat{\mathbf A}(t)=\mathbf{U}(t) {\mathbf \Sigma}(t) \mathbf Y(t)^{\mathrm{T}},
\end{equation}
where $\mathbf{U}(t) \in \mathbb{R}^{n \times r}$ and $\mathbf{Y}(t) \in \mathbb{R}^{s \times r}$ are time-dependent orthonormal spatial and parametric bases, respectively, $\mathbf{\Sigma}(t) \in \mathbb{R}^{r \times r}$ is a full matrix, and $r \ll \min(n, s)$ is the rank of the approximation. For brevity, the explicit dependence on time is omitted where convenient.

The central idea is to minimize
the instantaneous residual by optimally updating $\mathbf{U}$ and $\mathbf{Y}$ in time. We substituted the low-rank approximation into the FOM Eq. \eqref{eq:FOM}, the residual is defined as
\begin{equation}
\mathbf R(t)=\frac{\mathrm{d}}{\mathrm{~d} t}\left(\mathbf U \mathbf \Sigma \mathbf {Y}^{\mathrm T}\right)-\mathcal{F}\left(\mathbf {U} \mathbf {\Sigma} \mathbf {Y}^{\mathrm T}\right).
\end{equation}

% where $\mathcal{F}$ denotes the right-hand side of the MDE Eq. (\ref{eq:FOM}).

The evolution equations for $\mathbf{U}$, $\mathbf{\Sigma}$, and $\mathbf{Y}$ are obtained by minimizing the residual norm
\begin{equation}\label{eq:ResMin}
\mathcal{J}\left( \dot{\mathbf{U}}, \dot{\boldsymbol{\Sigma}}, \dot{\mathbf{Y}} \right) = \left\| \frac{\mathrm{d}}{\mathrm{d}t} \left( \mathbf{U} \mathbf{\Sigma} \mathbf{Y}^{\mathrm T} \right) - \mathcal{F}\left(  \mathbf{U} \mathbf{\Sigma} \mathbf{Y}^{\mathrm T} \right) \right\|_F^2,
\end{equation}
subject to the orthonormality constraints on $\mathbf{U}$ and $\mathbf{Y}$. Riemannian optimization techniques~\cite{koch2007dynamical} or Lagrange multipliers~\cite{ramezanian2021fly} can be used to solve this constrained minimization problem. The resulting evolution equations for the low-rank factors  are given by:
\begin{subequations}\label{eq:DLRA_evolution}
\begin{align}
\dot{\mathbf{\Sigma}} & =\mathbf{U}^{\mathrm T} \mathbf{F} \mathbf{Y},
\label{eq:Sigma_evolution} \\
\dot{\mathbf{U}} & =\left(\mathbf{I}_n-\mathbf{U} \mathbf{U}^{\mathrm T}\right) \mathbf{F} \mathbf{Y} \mathbf{\Sigma}^{-1}, \label{eq:U_evolution} \\
\dot{\mathbf{Y}} & =\left(\mathbf{I}_s-\mathbf{Y} \mathbf{Y}^{\mathrm T}\right) \mathbf{F}^{\mathrm T} \mathbf{U} \mathbf{\Sigma}^{-\mathrm T}, \label{eq:Y_evolution}
\end{align}
\end{subequations}
 where $\mathbf{F}=\mathcal{F}\left(\mathbf{U} \mathbf{\Sigma} \mathbf{Y}^{\mathrm T}\right)$. % This is equivalent to orthogonal projection of $\mathcal{F}\left(\mathbf{U} \boldsymbol{\Sigma} \mathbf{Y}^{\mathrm{T}}\right)$ onto the tangent space of the low-rank matrix manifold at $\hat{\mathbf{V}}=\mathbf{U} \boldsymbol{\Sigma} \mathbf{Y}^{\mathrm{T}}$.  
  From a geometric perspective, the DLRA~\cite{koch2007dynamical} evolution equations are obtained via the orthogonal projection of $\mathcal{F}(\hat{\mathbf A})$ onto the tangent plane of the low-rank matrix manifold at $\hat{\mathbf A}$. 
  
  % utilizes the properties of the tangent space $\mathcal{T}_{\hat{\mathbf{A}}} \mathcal{M}_r$ [definition (\ref{def:Tangent})] to solve the constrained residual minimization problem given by Eq. (\ref{eq:ResMin}). In this framework, the residual is minimized subject to the constraint that $\dot{\hat{\mathbf{A}}}(t) \in \mathcal{T}_{\hat{\mathbf{A}}} \mathcal{M}_r$. %The solution to this constrained minimization problem is given by 
% \begin{equation}\label{eq:Min_tan}
% \dot{\hat{\mathbf{V}}} = \mathcal{P}_{\mathcal{T}_{\hat{\mathbf{V}}}}\left( \mathcal{F}(t, \hat{\mathbf{V}}) \right),
% \end{equation}
% where $\mathcal{P}_{\mathcal{T}_{\hat{\mathbf{V}}}}$ [definition (\ref{def:Ortho_Proj})] denotes the orthogonal projection onto the tangent space $\mathcal{T}_{\hat{\mathbf{V}}}$ at $\hat{\mathbf{V}} = \mathbf{U} \boldsymbol{\Sigma} \mathbf{Y}^{\mathrm T}$. 
% The equations \ref{eq:Sigma_evolution}-\ref{eq:Y_evolution} can easily be derived from equation \ref{eq:Min_tan}.

The time integration of  Eqs. (\ref{eq:U_evolution}-\ref{eq:Y_evolution}) can become unstable or very stiff when $\mathbf \Sigma$ is singular or near-singular. This is particularly problematic because obtaining more accurate results requires resolving smaller singular values. This issue has been addressed with new time integration schemes that are robust in the presence of zero or arbitrarily small singular values \cite{ceruti2022unconventional,ceruti2022rank,ceruti2023rank,ceruti2024robust,KieriVandereycken2019,rodgers2022adaptive,doi:10.1137/21M1431229}.  

The computational savings from solving the above MDE on low-rank manifolds diminish when the exact rank of $\mathcal{F}(\hat{\mathbf{A}})$ is high. In the cases where $\mathcal{F}(\hat{\mathbf{A}})$ is full-rank---even though $\hat{\mathbf{A}}$ itself is rank-$r$---the savings vanish because all entries of the matrix $\mathbf F = \mathcal{F}(\hat{\mathbf{A}})$ need to be computed. The rank of $\mathcal{F}(\hat{\mathbf{A}})$ can substantially exceed $r$, for instance, when many terms are added or when $\mathcal{F}$ contains polynomial nonlinearities; indeed, the higher the polynomial degree, the larger the resulting rank. When $\mathcal{F}$ involves non-polynomial nonlinearities, $\mathcal{F}(\hat{\mathbf{A}})$ is full-rank, apart from rare special cases. Consequently, any exact-rank treatment of $\mathbf{F}$ may be computationally prohibitive, particularly for nonlinear MDEs.

\subsection{CUR for time-integration of MDEs on low-rank manifolds}
Recently, a new time-integration methodology based on a CUR low-rank approximation was introduced, providing a cost-effective way to address the nonlinearity of $\mathcal{F}$ \cite{TDBCUR}. Numerical examples further demonstrated that the method remains robust even in the presence of small or vanishing singular values, and that standard high-order schemes—such as those from the Runge–Kutta family—can be employed. 

We begin by briefly outlining the methodology proposed in \cite{TDBCUR}, as it serves as the foundation for the approach developed in this paper.

We discretize Eq. (\ref{eq:FOM})  in time using the time integration scheme of choice. For simplicity in the exposition, we use explicit Euler temporal discretization as follows:
\begin{equation}\label{eq:FOM_euler} 
\mathbf{A}^k=\hat{\mathbf{A}}^{k-1}+\Delta t \overline{\mathbf{F}},    
\end{equation}
where $\Delta t$ is the step size and $\overline{\mathbf{F}}$ in this case is given by: $\overline{\mathbf{F}}=\mathcal{F}\left(\hat{\mathbf{A}}^{k-1}\right)$. Higher-order time integration schemes may be used with some modification, as explained in \cite[Section 2(f)]{TDBCUR}.   We assume that the solution at the previous time step, $\hat{\mathbf{A}}^{k-1}$, lies on the rank-$r$ manifold, i.e., it is a rank-$r$ matrix. However, the solution at the current time step, $\mathbf{A}^k$, is generally not of rank $r$ due to the nonlinear nature of the operator $\mathcal{F}$ and/or the addition of matrices, both of which typically increase the rank. Consequently, to prevent the rank from growing unbounded over time, we must bring $\mathbf{A}^k$ back onto the rank-$r$ manifold $\mathcal{M}_r$ at each time step—ideally in an optimal manner.

To illustrate this more clearly, consider the rank-$r$ approximation of $\mathbf{A}^k$
\begin{equation}
    \mathbf A^k = \hat{\mathbf A}^k + \mathbf R^k,
\end{equation}
where $\mathbf R^k$ is the residual due to the low-rank approximation. The best low-rank approximation in the Frobenius norm  can be obtained via SVD as follows:
\begin{equation}
  % \text{min}_{\hat{\mathbf V}^k \in \mathcal{M}_r} \|\hat{\mathbf V}^k-\mathbf V^k\|_F
  \hat{\mathbf A}^k_{\text{best}} = \mbox{\texttt{SVD}}(\mathbf A^k),
\end{equation}
%as follows:
%  \begin{equation}\label{eq:Ortho}
% \hat{\mathbf{V}}^k_{\text{best}}  = \mathbf{P}^{\Langle}_{\mathbf{U}^k}\mathbf V^k \mathbf{P}^{\Langle}_{\mathbf{Y}^k}
% \end{equation}
% where $\mathbf{P}^{\Langle}_{\mathbf{U}^k} = \mathbf U^k \mathbf {U^k}^{\mathrm T}$ and $\mathbf{P}^{\Langle}_{\mathbf{Y}^k} = \mathbf Y^k \mathbf {Y^k}^{\mathrm T}$ are the orthogonal projections [definition (\ref{def:Ortho_Proj})] onto the row and column spaces $ \mathbf U^k$ and $ \mathbf Y^k$ of $\hat{\mathbf{V}}^k_{\text{best}}$.
where \texttt{SVD}( $\mathbf A^k$) denotes the rank-$r$ truncated SVD of $\mathbf A^k$.  The above approach, also known as \emph{step truncation} \cite{rodgers2022adaptive, KieriVandereycken2019}, offers several important advantages: it is robust in the presence of small or zero singular values, simple to implement, and preserves the nominal temporal accuracy of the underlying time integration scheme—whether explicit or implicit \cite{RV23}.

However, the issue of computational cost for nonlinear MDEs is not mitigated. For non-polynomial nonlinearities, all entries of $\overline{\mathbf F}$ must be explicitly computed. In the case of a polynomial nonlinearity of degree $p$, the rank of $\overline{\mathbf F}$ scales as $r^p$, which can become prohibitive. Additionally, this approach requires computing the singular value decomposition (SVD). However, this cost can be reduced using randomized SVD algorithms.

In \cite{TDBCUR}, this challenge was addressed by employing a collocation scheme based on a \texttt{CUR-CR} algorithm, where $\mathbf A^k$ is approximated via:
$$
\hat{\mathbf A}^k = \texttt{CUR-CR}(\mathbf A^k).
$$
The above scheme sets the residual equal to zero at the CUR-selected columns and rows, i.e., $\mathbf R^k (\mathbf p,:) = \mathbf 0$ and  $\mathbf R^k (:,\mathbf s) = \mathbf 0$. The columns and rows are selected based on applying DEIM to the left and right singular vectors of $\hat{\mathbf A}^{k-1}$, i.e.,  the low-rank solution at the previous time step. Therefore, $\mathbf p$ and $\mathbf s$ are updated at each time step, resulting in a time adaptive CUR scheme.    To enhance the robustness and the stability of the \texttt{CUR-CR} algorithm, row oversampling was proposed in \cite{TDBCUR}, where an additional fixed number of rows are sampled to improve the condition number of subcolumn matrices.

Although the CUR algorithm used in \cite{TDBCUR} offers computational efficiency, there are aspects that can be improved. The method applies row oversampling to enhance stability; it does not employ a similar strategy for column selection.  The row oversampling strategy requires access to entire rows of the matrix, which can be computationally expensive when $s$ is large. Sweeping the entire row means the cost of oversampling is similar to that of row sampling for a higher-rank approximation.  The number of oversampled rows is treated as a fixed hyperparameter. However, the stability and accuracy of the decomposition are sensitive to this value. It is also difficult to know a priori how many rows should be oversampled.

In the next section, we introduce an adaptive cross-oversampling algorithm, informed by the error analysis in Section~\ref{sec:err_anl}. The algorithm dynamically determines the number of oversampling entries required, thereby tailoring the computational effort to the problem at hand. When applied to the solution of MDEs on low-rank manifolds, this strategy improves both robustness and overall performance.

% In the next section, we present an adaptive cross oversampling algorithm that oversamples the intersection of additional rows and columns rather than the entire row. 

% We denote the oversampling algorithm used in \cite{TDBCUR} as \texttt{R-OS}, where “\texttt{R-OS}” stands for row oversampling. In the next section, we present a cross oversampling algorithm that addresses all of the above limitations.  

\subsection{Adaptive cross oversampling}
We present an adaptive cross oversampling strategy that addresses the above limitations and efficiently enhances the accuracy and robustness of the \texttt{CUR-CR} algorithm.
% The method requires oversampling at the intersection of additional rows and columns.
% Rather than sampling entire rows and columns of the matrix $\mathbf{A}$, we strategically sample selected cross elements, thus reducing the cost of data acquisition.
We refer to this method as \texttt{CUR-CR-OS}, where \texttt{OS} indicates oversampling.

The proposed algorithm can be used as a generic CUR matrix low-rank approximation and is not restricted to the solution of MDEs. Its design is independent of the specific row and column selection strategy, making it compatible with a variety of pivot-based or probabilistic sampling techniques. To this end, we first present the methodology for a general matrix $\mathbf{A}$. In Section~\ref{sec:TDB-CUR}, we then illustrate how the algorithm can be combined with DEIM for column and row selection in the context of solving random PDEs on low-rank matrix manifolds.
% We aim to approximate ${\mathbf V}^{k}$ using a CUR decomposition such that
% \begin{equation}
% \hat{\mathbf V}^{k} = \texttt{CUR}({\mathbf V}^{k}).
% \end{equation}

The first step is to perform QR decompositions on the selected columns and rows:
\begin{subequations}
\begin{align}
\mathbf Q_c \mathbf R_c &= \texttt{qr}({\mathbf A}(:, \mathbf s)), \\
\mathbf Q_r \mathbf R_r &= \texttt{qr}({\mathbf A}(\mathbf p, :)^{\mathrm T}),
\end{align}
\end{subequations}
where $\mathbf Q_c \in \mathbb{R}^{n \times r}$ and $\mathbf Q_r \in \mathbb{R}^{s \times r}$ are orthonormal matrices, and $\mathbf R_c, \mathbf R_r \in \mathbb{R}^{r \times r}$ are upper triangular matrices and $\mathbf p = [p_1, p_2, \dots, p_r]$ and $\mathbf s =[s_1, s_2, \dots, s_r]$ are the selected rows and columns of $\mathbf A$, respectively. 

% Since the CUR components of ${\mathbf V}^{k}$ are not available at the current time step $k$, we leverage the approximation from the previous time step, namely $\hat{\mathbf V}^{k-1} = \mathbf U^{k-1} \mathbf \Sigma^{k-1} {\mathbf Y^{k-1}}^{\mathrm T}$, to compute the DEIM sampling points. Specifically, we set
% \begin{equation}
% \mathbf I = \texttt{DEIM}(\mathbf U^{k-1}), \quad \mathbf J = \texttt{DEIM}(\mathbf Y^{k-1}).
% \end{equation}

% Remark : The initial approximation is ideally obtained from the FOM initial condition via the rank-$r$ singular value decomposition, i.e., $\hat{\mathbf V}^0 = \texttt{SVD}_r(\mathbf V_0)$.

Next, we consider a CUR decomposition in the form of:
\begin{equation}
\hat{\mathbf A} = \mathbf Q_c \mathbf Z \mathbf Q_r^{\mathrm T},
\end{equation}
where $\mathbf Z \in \mathbb{R}^{r \times r}$ is obtained as follows:
\begin{equation}
\mathbf Z =  \mathbf Q_c^{\dagger}(\bar{\mathbf p}, :)  \mathbf A^k(\bar{\mathbf p}, \bar{\mathbf s}) \left( \mathbf Q_r^{\dagger}(\bar{\mathbf s}, :) \right)^{\mathrm T},
\end{equation}
where $\bar{\mathbf p}$ and $\bar{\mathbf s}$ are the oversampled row and column index sets, respectively. Specifically, $\bar{\mathbf p} \in \mathbb{R}^{r + m_r}$ and $\bar{\mathbf s} \in \mathbb{R}^{r + m_c}$, where $m_r$ and $m_c$ denote the number of additional rows and columns used for oversampling. Both $m_r$ and $m_c$ are treated as prescribed parameters (or hyperparameters) in this section. In Section~\ref{subsec:Adap_m}, we introduce an algorithm that determines suitable values of $m_r$ and $m_c$ adaptively. %These indices are selected using the GappyPOD+E algorithm \cite{peherstorfer2020stability}.

The oversampling indices $\bar{\mathbf p}$ and $\bar{\mathbf s}$ can be selected using a variety of strategies such as leverage score sampling \cite{doi:10.1137/07070471X,doi:10.1073/pnas.0803205106}, volume sampling \cite{cortinovis2020low,goreinov1997theory}, the BSS framework \cite{batson2009twice,boutsidis2014near}, deterministic leverage-score approaches \cite{papailiopoulos2014provable}, and GappyPOD+E \cite{peherstorfer2020stability}. The presented methodology is agnostic to the choice, but in practice, we use GappyPOD+E due to its simplicity and robustness.

%Deterministic pivoting-based methods such as GappyPOD+E \cite{peherstorfer2020stability} or stabilized variants of QDEIM improve the conditioning of the sampled submatrix. Randomized approaches include leverage-score sampling \cite{doi:10.1137/07070471X}, volume sampling \cite{cortinovis2020low}, and determinantal point processes (DPPs) \cite{derezinski2021determinantal}, which provide strong probabilistic guarantees. Hybrid techniques, such as L-DEIM \cite{gidisu2021hybrid}, combine deterministic pivoting with randomized oversampling to balance stability and accuracy. Our framework is agnostic to the choice, but in practice we use GappyPOD+E due to its simplicity and robustness.

\begin{figure}[t]
    \centering
    \includegraphics[width=0.8\textwidth]{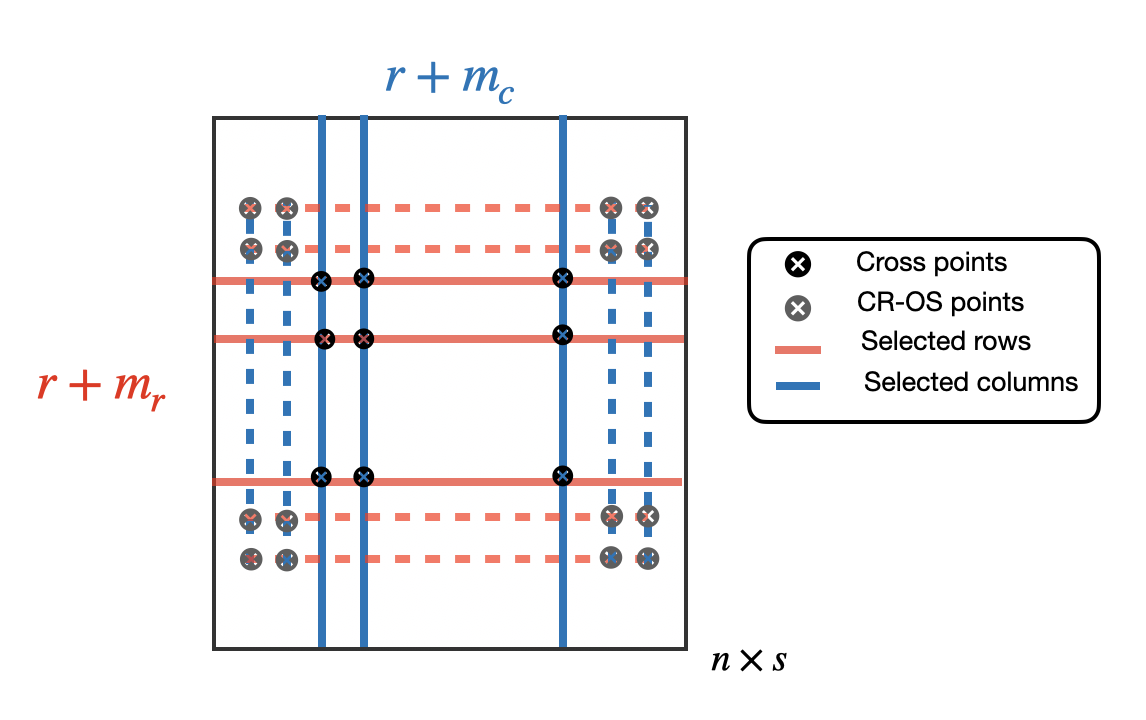}
     \caption{Schematic of CUR oversampling methodology using cross oversampling (\texttt{CR-OS}), to oversampling $m_r$ rows and $m_c$ columns of matrix $\mathbf A \in  \mathbb R^{n\times s}$. We observe the only extra cost resulting from oversampling is $m_rm_c$, which are the total number of cross points used for oversampling.}
    \label{fig:Compare_CostofData}
\end{figure}

The CUR approximation \(\hat{\mathbf A} = \mathbf Q_c \mathbf Z \mathbf Q_r^{\mathrm T}\) can be converted into an SVD form by computing the SVD of the small matrix \(\mathbf Z\) and rotating $\mathbf Q_c$ and $\mathbf Q_r$ according to their singular values. First, we compute the SVD of  \(\mathbf Z\):
\[
\mathbf Z = \mathbf \Psi_c \mathbf \Sigma^k \mathbf \Psi_r^{\mathrm T},
\]
where $\mathbf \Psi_c \in \mathbb{R}^{r \times r}$ and $\mathbf \Psi_r \in \mathbb{R}^{r \times r}$ are orthonormal matrices. Replacing $\mathbf Z$ with its SVD results in:
\[
\hat{\mathbf A} =  \mathbf Q_c \mathbf Z \mathbf Q_r^{\mathrm T} = \mathbf Q_c \mathbf \Psi_c \mathbf \Sigma \mathbf \Psi_r^{\mathrm T} \mathbf Q_r^{\mathrm T} = \mathbf U \mathbf \Sigma  {\mathbf Y}^{\mathrm T},
\]
where
\[
\mathbf U = \mathbf Q_c \mathbf \Psi_c \quad \mbox{and}  \quad \mathbf Y = \mathbf Q_r \mathbf \Psi_r.
\]
 We observe that in the proposed cross oversampling algorithm, only an additional $m_r m_c$ elements of $\mathbf A$ need to be computed for oversampling beyond the standard $r$ rows and $r$ columns used in the CUR approximation. Therefore, the total number of matrix entries required is:
\begin{equation}
n_d = nr + rs - r^2 + m_r m_c,
\end{equation}
where $r^2$ is subtracted since they have been accounted for twice in $r(n+s)$. Therefore, the additional number of entries required due to the proposed oversampling is $m_r m_s$.   
As a result, the number of additional entries of the target matrix required for oversampling does not scale with $n$ or $s$. A schematic of the \texttt{CUR-CR-OS} algorithm is shown in Figure (\ref{fig:Compare_CostofData}).

\subsection{Error analysis}\label{sec:err_anl}
In this section, we provide an error analysis of the proposed CUR algorithm. We use the results of the error analysis to develop an adaptive oversampling algorithm, where $m_r$ and $m_c$ are determined on the fly. 

The \texttt{CUR-CR-OS} algorithm can be expressed as: 
\begin{equation}\label{eq:obliq_proj}
 \hat{\mathbf A}= \mathcal P \mathbf A \mathcal  S,   
\end{equation}
where $\mathcal P \in \mathbb{R}^{n \times n}$ and $\mathcal S \in \mathbb{R}^{s \times s}$ are \emph{oblique projectors} onto the space spanned by the selected columns and rows of $\mathbf A$ as follows:
\begin{subequations}
\begin{align}
    \mathcal P &= \mathbf Q_c  (\ol{\mathbf P}^{\mathrm T} \mathbf Q_c)^{\dagger} \ol{\mathbf P}^{\mathrm T}, \label{eq:proj_P}\\
    \mathcal S &= \ol{\mathbf S} ( \mathbf Q_r^{\mathrm T} \ol{\mathbf S})^{\dagger} \mathbf Q_r^{\mathrm T}, \label{eq:proj_S}
\end{align}
\end{subequations}
In the above definition of the projectors,   $\mathcal P$ and $\mathcal S$ act from the left and right of  $\mathbf A$, respectively. It is straightforward to verify that  Eq. \ref{eq:obliq_proj} with the $\mathcal P$ and $\mathcal S$ definitions given by Eqs. (\ref{eq:proj_P})-(\ref{eq:proj_S}) is equivalent to the \texttt{CR-OS}.

The derivation of the error bound for the presented CUR method follows a similar approach to that presented  \cite{sorensen2016deim}. However, unlike in \cite{sorensen2016deim}, $\mathcal P $ and $\mathcal S$ are \emph{oblique} projectors as opposed \emph{orthogonal} projectors used in \cite{sorensen2016deim}. 

\begin{lemma}
The projectors $\mathcal P$ and $\mathcal S$ given by Eqs. (\ref{eq:proj_P})-(\ref{eq:proj_S}) can be expressed in terms of the elements of $\mathbf A$ and are equivalent to \begin{subequations}
\begin{align}
    \mathcal P &= \mathbf A \mathbf S  (\ol{\mathbf P}^{\mathrm T} \mathbf A \mathbf S)^{\dagger} \ol{\mathbf P}^{\mathrm T}, \label{eq:proj_PV}\\
    \mathcal S &= \ol{\mathbf S} ( \mathbf P^{\mathrm T}\mathbf A \ol{\mathbf S})^{\dagger} \mathbf P^{\mathrm T} \mathbf A\label{eq:proj_SV}. 
\end{align}
\end{subequations}
\end{lemma}
\begin{proof}
Using the QR factorization, we have: $\mathbf A \mathbf S=\mathbf Q_c\mathbf R_c$. Therefore, $\mathbf Q_c=\mathbf A \mathbf S \mathbf R_c^{-1}$ where $\mathbf R_c$ is invertible. Then, 
substitute $\mathbf Q_c=\mathbf A \mathbf S\mathbf R_c^{-1}$ into the definition of $\mathcal P$ given  by Eq. (\ref{eq:proj_P}):
\begin{align*}
\mathcal P
&= \mathbf A \mathbf S\mathbf R_c^{-1}\,\big(\overline{\mathbf P}^{\mathrm T}(\mathbf A \mathbf S)\mathbf R_c^{-1}\big)^\dagger\,\overline{\mathbf P}^{\mathrm T}\\
&= \mathbf A \mathbf S\mathbf R_c^{-1}\,\mathbf R_c\,(\overline{\mathbf P}^{\mathrm T}\mathbf A \mathbf S)^\dagger\,\overline{\mathbf P}^{\mathrm T}\\
&= \mathbf A \mathbf S\,(\overline{\mathbf P}^{\mathrm T}\mathbf A \mathbf S)^\dagger\,\overline{\mathbf P}^{\mathrm T}.
\end{align*}
% For any full-column-rank matrix $X$ and any invertible square matrix $R$,
% \[
% (XR)_L^\dagger \;=\; R^{-1}X_L^\dagger,
% \]
% because \((R^{-1}X_L^\dagger)(XR)=R^{-1}(X_L^\dagger X)R=I\).
% Apply this with \(X=\overline{\mathbf P}^{\mathrm T}AS\) and \(R=\mathbf R_c^{-1}\):
% \[
% \big(\overline{\mathbf P}^{\mathrm T}(AS)\mathbf R_c^{-1}\big)_L^\dagger
% = (\mathbf R_c^{-1})^{-1}(\overline{\mathbf P}^{\mathrm T}AS)_L^\dagger
% = \mathbf R_c\,(\overline{\mathbf P}^{\mathrm T}AS)_L^\dagger.
% \]
% Hence
% \[
% \mathcal P
% = (AS)\mathbf R_c^{-1}\,\mathbf R_c\,(\overline{\mathbf P}^{\mathrm T}AS)_L^\dagger\,\overline{\mathbf P}^{\mathrm T}
% = AS\,(\overline{\mathbf P}^{\mathrm T}AS)_L^\dagger\,\overline{\mathbf P}^{\mathrm T}.

% \]
Similarly, using QR factorization, we have $ (\mathbf P^{\mathrm T} \mathbf A)^{\mathrm T}= \mathbf  Q_r \mathbf  R_r$ with $\mathbf  R_r$ invertible, so $\mathbf  Q_r=\left(\mathbf  P^{\mathrm T} \mathbf  A\right)^{\mathrm T} \mathbf  R_r^{-1}$.
 Therefore, Eq. (\ref{eq:proj_S}) can be similarly  expressed as
\begin{align*}
    \mathcal  S &=\bar{\mathbf  S}\left(\left(  (\mathbf P^T \mathbf  A)^{\mathrm T}\mathbf  R_r^{-1} \right)^{\mathrm T}\bar{\mathbf S}\right)^{\dagger}\left( (\mathbf P^{\mathrm T} \bm A)^{\mathrm T} \mathbf  R_r^{-1} \right)^{\mathrm T}\\
    & = \bar{\mathbf  S}\left(\mathbf   R_r^{-\mathrm T} \mathbf P^T \mathbf  A \bar{\mathbf S}\right)^{\dagger}\left( {\mathbf  R_r^{-\mathrm T}} \mathbf P^{\mathrm T} \mathbf A \right)\\
    & = \bar{\mathbf  S}\left( \mathbf P^T \mathbf  A \bar{\mathbf S}\right)^{\dagger}\left( \mathbf  R_r^{\mathrm T} \mathbf  R_r^{-\mathrm T} \right) \mathbf P^{\mathrm T} \mathbf A \\
    & = \bar{\mathbf  S}\left( \mathbf P^T \mathbf  A \bar{\mathbf S}\right)^{\dagger}\mathbf P^{\mathrm T} \mathbf A. \\
\end{align*}

% Similarly, substituting $\mathbf Q_r = (\mathbf P^{\mathrm T}\mathbf A) \mathbf R_r^{-1}$ into Eqs. (\ref{eq:proj_S}) and applying the right-inverse property, we obtain $\mathcal S = \overline{\mathbf S}\,(\mathbf P^{\mathrm T}\mathbf A\,\overline{\mathbf S})_R^{\dagger}\,\mathbf P^{\mathrm T}\mathbf A$. \hb{show the details of the S projectors.}This completes the proof. 
\end{proof}

%% This should be a lemma ....
% First, we explicitly express the projectors $\mathcal P$ and $\mathcal S$ in terms of the elements of $\mathbf A$. To this end, we show that the CR-OS projectors can be expressed as:
% \begin{subequations}
% \begin{align}
%     \mathcal P &= \mathbf A \mathbf S  (\ol{\mathbf P}^{\mathrm T} \mathbf A \mathbf S)_L^{\dagger} \ol{\mathbf P}^{\mathrm T}, \label{eq:proj_PV}\\
%     \mathcal S &= \ol{\mathbf S} ( \mathbf P^{\mathrm T}\mathbf A \ol{\mathbf S})_R^{\dagger} \mathbf P^{\mathrm T} \mathbf A\label{eq:proj_SV}. 
% \end{align}
% \end{subequations}
% \hb{Prove the above by replacing the QR decomposition in Eqs. \ref{eq:proj_P}-\ref{eq:proj_S} ...}
%% ---------------------------------------------------------- End of Lemma -------------

We now derive an error bound for the \texttt{CUR-CR-OS} where the error is bounded by the error of the best rank-$r$ approximation obtained via SVD. To this end, let $\mathbf A \approx \mathbf U \mathbf \Sigma \mathbf Y^{\mathrm T}$ be the rank-$r$ SVD  of $\mathbf A$, where $\mathbf U \in \mathbb{R}^{n \times r}$ and $\mathbf Y \in \mathbb{R}^{s \times r}$ are the matrices of the first $r$ left and right singular vectors of $\mathbf A$, respectively. Also, let $\tilde{\mathcal P}$ and $\tilde{\mathcal S}$ be the oblique projectors onto the space spanned by  $\mathbf U$ and $\mathbf Y$, respectively as follows:
\begin{subequations}
\begin{align*}
    \tilde{\mathcal P} &= \mathbf U  (\mathbf P^{\mathrm T} \mathbf U)^{-1} \mathbf P^{\mathrm T},\\
    \tilde{\mathcal S} &= \mathbf S ( \mathbf Y^{\mathrm T} \mathbf S)^{-1}\mathbf Y^{\mathrm T}. 
\end{align*}
\end{subequations}

According to \cite[Lemma 4.1]{sorensen2016deim}, the error of the oblique projection of $\mathbf A$ onto $\mathbf U$ and $\mathbf Y$ is bounded by:
\begin{subequations}
\begin{align*}
    \| \mathbf A - \tilde{\mathcal P} \mathbf A\|  &\leq \eta_{p} \sigma_{r+1}, \\
    \| \mathbf A -  \mathbf A \tilde{\mathcal S}\| &\leq \eta_{s} \sigma_{r+1}, 
\end{align*}
\end{subequations}
where 
\[
\eta_{p} = \| (\mathbf P^{\mathrm T} \mathbf U)^{-1} \| \quad \mbox{and} \quad \eta_{s} = \| (\mathbf Y^{\mathrm T} \mathbf S)^{-1} \|.
\]
It is important to note that $\eta_p$ and $\eta_s$ require access to exact singular vectors of $\mathbf A$. It was shown in \cite[Theorem 4.1]{sorensen2016deim} that the \texttt{CUR-opt} error is bounded by 
\begin{equation}\label{err_bnd_opt}
    \|\mathbf A - \texttt{CUR-opt}(\mathbf A) \| \le (\eta_p + \eta_s) \sigma_{r+1},
\end{equation}
where $\texttt{CUR-opt}(\mathbf A) = \mathbf C \mathbf C^{\dagger} \mathbf A \mathbf R^{\dagger} \mathbf R$ and $\mathbf C = \mathbf A \mathbf S $ and $\mathbf R = \mathbf P^{\mathrm T} \mathbf A$. 
%% This should be a lemma -----------------------------------
\begin{lemma}\label{lem:obliq_proj}
The following bounds hold for the projection errors of $\mathbf A$:
\begin{align}
    \| (\mathbf I_n - \mathcal{P}) \mathbf A \| &\le \overline{\eta}_p \, \eta_s \, \sigma_{r+1}, \label{eq:proj_err1}\\
    \| \mathbf A (\mathbf I_s - \mathcal{S}) \| &\le \eta_p \, \overline{\eta}_s \, \sigma_{r+1}, \label{eq:proj_err2}
\end{align}
where $\mathcal{P}$ and $\mathcal{S}$ are the projection operators defined in
Eqs.~\eqref{eq:proj_P}--\eqref{eq:proj_S}, $\sigma_{r+1}$ is the $(r+1)$th singular value of $\mathbf A$, and
$\eta_p=\| (\mathbf P^{\mathrm T} \mathbf U)^{-1} \|$, $\eta_s=\| (\mathbf Y^{\mathrm T} \mathbf S)^{-1} \|$ are defined based on the exact singular vectors of $\mathbf A$, and  $\overline{\eta}_p = \| (\ol{\mathbf P}^{\mathrm T} \mathbf Q_c)^{\dagger} \|$, and $\overline{\eta}_s=\| ( \mathbf Q_r^{\mathrm T} \ol{\mathbf S})^{\dagger}\|$.
% are computed using the acutal subcolumns and subrows of $\mathbf A$.
\end{lemma}

% Now we show that 
% \begin{equation}
%     (\mathbf I_n - \mathcal P) \mathbf A = (\mathbf I_n - \mathcal P) \mathbf A (\mathbf I_s - \tilde {\mathcal S}).  
% \end{equation}
\begin{proof}
From Eq. (\ref{eq:proj_PV}), we have:
\begin{align*}
(\mathbf I_n - \mathcal P) \mathbf A &=(\mathbf I_n - \mathbf A \mathbf S  (\ol{\mathbf P}^{\mathrm T} \mathbf A \mathbf S)^{\dagger} \ol{\mathbf P}^{\mathrm T} ) \mathbf A\\
&= \mathbf A(\mathbf I_s- \mathbf S  (\ol{\mathbf P}^{\mathrm T} \mathbf A \mathbf S)^{\dagger} \ol{\mathbf P}^{\mathrm T} \mathbf A)\\
&=\mathbf A (\mathbf I_s - \mathbf \Phi),
\end{align*}
where $\mathbf \Phi = \mathbf S  (\ol{\mathbf P}^{\mathrm T} \mathbf A \mathbf S)^{\dagger} \ol{\mathbf P}^{\mathrm T} \mathbf A$. It is straightforward to verify that $\mathbf \Phi$ is an oblique projector onto the span of the columns of $\mathbf S$, since $\mathbf \Phi^2 = \mathbf \Phi$. Therefore, 

\[
\mathbf \Phi \mathbf S = \mathbf S  (\ol{\mathbf P}^{\mathrm T} \mathbf A \mathbf S)^{\dagger} \ol{\mathbf P}^{\mathrm T} \mathbf A \mathbf S = \mathbf S.
\]
Now we can show that $\mathbf \Phi \tilde{\mathcal S} = \tilde{\mathcal S}$, since:
\[
 \mathbf \Phi \tilde{\mathcal S} = \mathbf \Phi \mathbf S ( \mathbf Y^{\mathrm T} \mathbf S)^{\dagger}\mathbf Y^{\mathrm T} = \mathbf S ( \mathbf Y^{\mathrm T} \mathbf S)^{\dagger}\mathbf Y^{\mathrm T} = \tilde{\mathcal S}. 
\]
Therefore,
\[
\mathbf I_s - \mathbf \Phi = (\mathbf I_s - \mathbf \Phi)(\mathbf I_s - \tilde{\mathcal S}).
\]
As a result:
\begin{equation}
  (\mathbf I_n - \mathcal P) \mathbf A =   \mathbf A (\mathbf I_s - \mathbf \Phi)=(\mathbf I_n - \mathcal P) \mathbf A (\mathbf I_s - \tilde{\mathcal S}).
\end{equation}
Therefore,
\begin{align*}
\| (\mathbf I_n - \mathcal P) \mathbf A \| &= \| (\mathbf I_n - \mathcal P) \mathbf A (\mathbf I_s - \tilde{\mathcal S}) \|\\
& \leq   \| (\mathbf I_n - \mathcal P) \|  \|\mathbf A (\mathbf I_s - \tilde{\mathcal S}) \|\\
& = \ol{\eta}_p \eta_s \sigma_{r+1}.
\end{align*}
where $\ol{\eta}_p = \| \mathcal P\| = \| \mathbf I_p - \mathcal P\| = \| (\ol{\mathbf P}^{\mathrm T} \mathbf Q_c)^{\dagger} \| $.\\
Similarly, it can be shown that:
\begin{align*}
\| \mathbf A (\mathbf I_s - \mathcal S)  \| &= \| (\mathbf I_n - \tilde{\mathcal P}) \mathbf A (\mathbf I_s - \mathcal S) \|\\
& \leq   \| (\mathbf I_n - \tilde{\mathcal P})\mathbf A \|  \| (\mathbf I_s - \mathcal S) \|\\
& = \eta_p \ol{\eta}_s \sigma_{r+1},
\end{align*}
where $\ol{\eta}_s = \| \mathcal S\| = \| \mathbf I_s - \mathcal S\| = \|(\mathbf Q_r^{\mathrm T} \ol{\mathbf S})^{\dagger} \| $. This completes the proof.
\end{proof}

% \hb{show the details of this line. For the norm of projector, there is a paper cited in Mike's paper. }
%--------------------------------------------- End of Lemma ---------------------

% -------------- This should be the main Theorem -----------------
Now we prove an error bound for $\|\mathbf A -\mathcal P \mathbf A \mathcal S \|$.
\begin{theorem}\label{Err_bd}
The following bound holds for the projection error:
\begin{equation}\label{eq:err_bnd_CR_OS}
\| \mathbf A - \mathcal{P} \mathbf A \mathcal{S} \| 
\le \min \left\{
    \overline{\eta}_p \big( \eta_s + \eta_p \overline{\eta}_s \big) , \;
    \overline{\eta}_s \big( \eta_p + \overline{\eta}_p \eta_s \big) 
\right\} \sigma_{r+1}.
\end{equation}
Here, $\mathcal{P}$ and $\mathcal{S}$ are the projection operators defined in 
Eqs.~\eqref{eq:proj_P}--\eqref{eq:proj_S}, $\sigma_{r+1}$ is the $(r+1)$th singular value of $\mathbf A$, and $\eta_p$, $\eta_s$, $\overline{\eta}_p$, $\overline{\eta}_s$ are defined in Lemma \ref{lem:obliq_proj}.
\end{theorem}
\begin{proof}
We first split
\begin{align*}
\mathbf A - \mathcal P \mathbf A \mathcal S
&= \mathbf A - \mathcal P \mathbf A + \mathcal P \mathbf A - \mathcal P \mathbf A \mathcal S \\
&= (\mathbf I_n - \mathcal P)\,\mathbf A \;+\; \mathcal P \mathbf A\,(\mathbf I_s - \mathcal S).
\end{align*}
Taking norms and using sub-multiplicativity,
\begin{align*}
\| \mathbf A - \mathcal P \mathbf A \mathcal S \|
&\le \|(\mathbf I_n - \mathcal P)\mathbf A\| + \|\mathcal P\|\,\|\mathbf A(\mathbf I_s - \mathcal S)\| \\
&\le \ol{\eta}_p\,\eta_s\,\sigma_{r+1} \;+\; \ol{\eta}_p\,\eta_p\,\ol{\eta}_s\,\sigma_{r+1}
= \ol{\eta}_p\bigl(\eta_s + \eta_p \ol{\eta}_s\bigr)\sigma_{r+1},
\end{align*}
where inequalities \ref{eq:proj_err1} and \ref{eq:proj_err2} from Lemma \ref{lem:obliq_proj} are used.
Similarly,
\begin{align*}
\| \mathbf A - \mathcal P \mathbf A \mathcal S \|
&\le \|\mathbf A(\mathbf I_s - \mathcal S)\| + \|(\mathbf I_n - \mathcal P)\mathbf A\|\,\|\mathcal S\| \\
&\le \eta_p\,\ol{\eta}_s\,\sigma_{r+1} \;+\; \ol{\eta}_p\,\eta_s\,\ol{\eta}_s\,\sigma_{r+1}
= \ol{\eta}_s\bigl(\eta_p + \ol{\eta}_p \eta_s\bigr)\sigma_{r+1}.
\end{align*}
With $\ol{\eta}_p = \|\mathcal P\| = \|\mathbf I_n - \mathcal P\|$ and
$\ol{\eta}_s = \|\mathcal S\| = \|\mathbf I_s - \mathcal S\|$, we conclude
\[
\| \mathbf A - \mathcal P \mathbf A \mathcal S \|
\le \min
\left\{
\ol{\eta}_p\bigl(\eta_s + \eta_p \ol{\eta}_s\bigr),
\ol{\eta}_s\bigl(\eta_p + \ol{\eta}_p \eta_s\bigr)
\right\}\sigma_{r+1}.
\]
This completes the proof. 
\end{proof}
We make a few observations about the above error bound. First, it is a deterministic error bound in the spectral norm that is agnostic with respect to how rows and columns are selected. For example, it is valid for probabilistic sampling based methods or deterministic pivot-based techniques. These error bounds are not really complete without bounding $\eta_p$, $\eta_s$, $\ol{\eta}_p$, and $\ol{\eta}_s$.  Depending on the choice of index selection, there are some available bounds. For example, if the points are selected via the DEIM algorithm, the $\eta$'s are bounded as shown in \cite[Lemma 3.2]{CS10}. However, the bound in \cite{CS10} is pessimistic as it exponentially grows with the rank. In practice, when indices are selected based on DEIM, these  $\eta$'s are much smaller as the DEIM point selection algorithm is specifically designed to minimize these $\eta$'s using a greedy algorithm. When the indices are selected based on the maximum volume algorithm, much tighter bounds are presented in \cite[Proposition 5.1]{HH21}. As highlighted in \cite{HH21}, obtaining general bounds for the inverses of submatrices remains a challenging problem.  As we show, these error bounds can be very useful for devising an adaptive oversampling algorithm.

Comparing the error bound of \texttt{CUR-CR-OS} with that of the optimal \texttt{CUR-opt} can be both revealing and instructive. If all possible rows and columns are selected, \texttt{CUR-CR-OS}  exactly recovers \texttt{CUR-opt}. In this case, $\overline{\eta}_p = \overline{\eta}_s = 1$, and the error bound in inequality~\ref{eq:err_bnd_CR_OS} also recovers the optimal bound in inequality~\ref{err_bnd_opt}. Viewed this way, the \texttt{CUR-CR-OS} algorithm naturally fills the gap between \texttt{CUR-CR} on one end and \texttt{CUR-opt} on the other. Specifically, with no oversampling ($m_r = m_c = 0$), \texttt{CUR-CR-OS} coincides with \texttt{CUR-CR}, whereas when all points are included in oversampling, it becomes identical to \texttt{CUR-opt}.

In some practical applications, such as solving MDE in low-rank form, the exact singular vectors of $\mathbf A$ are typically not available, because calculating the exact singular values requires access to all entries of $\mathbf A$, which we want to avoid in the first place. As a result, $\eta_p$ and $\eta_s$ cannot be calculated for any posterior error analysis. On the other hand, $\ol{\eta}_p$ and $\ol{\eta}_s$  can be calculated a posteriori because they are based on the actual selected columns and rows of $\mathbf A$. 
We use this feature to propose an adaptive oversampling algorithm in the next section. 
% %%
% \begin{proof}

% \begin{align*}
% \mathbf A -\mathcal P \mathbf A \mathcal S &= \mathbf A - \mathcal P  \mathbf A +  \mathcal P \mathbf A - \mathcal P \mathbf A \mathcal S \\
% &=(\mathbf I_n - \mathcal P)  \mathbf A +  \mathcal P \mathbf A (\mathbf I_s -  \mathcal S)
% \end{align*}
% Therefore,
% \begin{align*}
% \| \mathbf A -\mathcal P \mathbf A \mathcal S \| &\leq \| (\mathbf I_n - \mathcal P)  \mathbf A \| +  \|\mathcal P \mathbf A (\mathbf I_s -  \mathcal S) \| \\
% & \leq \| (\mathbf I_n - \mathcal P)  \mathbf A \| +  \|\mathcal P \| \| \mathbf A (\mathbf I_s -  \mathcal S) \|\\
% & \leq  \ol{\eta}_p \eta_s \sigma_{r+1} +   \ol{\eta}_p \eta_p \ol{\eta}_s \sigma_{r+1}\\
% & = \ol{\eta}_p (\eta_s + \eta_p \ol{\eta}_s) \sigma_{r+1}.
% \end{align*}
% where $\ol{\eta}_p = \| \mathcal P\| = \| \mathbf I_n - \mathcal P\| = \|(\ol{\mathbf P}^{\mathrm T} \mathbf Q_c)_L^{\dagger} \| $.
% Similarly, it can be shown that:
% \begin{align*}
% \| \mathbf A -\mathcal P \mathbf A \mathcal S \| &\leq \|   \mathbf A (\mathbf I_s - \mathcal S) \| +  \|(\mathbf I_n -  \mathcal P) \mathbf A \mathcal S   \| \\
% & \leq \|   \mathbf A (\mathbf I_s - \mathcal S) \| +  \|(\mathbf I_n -  \mathcal P) \mathbf A\| \|\mathcal S   \|\\
% & \leq \eta_p  \ol{\eta}_s \sigma_{r+1} +   \ol{\eta}_p \eta_s \ol{\eta}_s \sigma_{r+1}\\
% & = \ol{\eta}_s (\eta_p + \ol{\eta}_p \eta_s) \sigma_{r+1}.
% \end{align*}
% Therefore, $\| \mathbf A -\mathcal P \mathbf A \mathcal S \|$ is bounded by the minimum of these two bounds.

% \end{proof}
% \hb{show the minimum explicitly.}

\subsubsection{Adaptive oversampling}\label{subsec:Adap_m}
In the presented oversampling algorithm, the number of oversampled rows ($m_r$) and columns ($m_c$) serve as hyperparameters. However, rather than relying on fixed pre-selected values for $m_r$ and $m_c$, we propose an adaptive strategy inspired by the above error analysis. Specifically, reducing  $\ol{\eta}_p = \| (\mathbf Q_c(\ol{\mathbf p},:))^{\dagger}\|$, $\ol{\eta}_s = \| (\mathbf Q_r(\ol{\mathbf s},:))^{\dagger}\|$  makes the \texttt{CUR-CUR-OS} arbitrarily closer to the optimal CUR approximation, which is obtained via the orthogonal projection of the target matrix onto the space spanned by selected columns and rows \cite{sorensen2016deim}. 
% \begin{align*}
%     \epsilon_f = \text{min}\{\eta_I(1+\eta_J),\eta_J(1+\eta_I)\}
% \end{align*}
% where $\eta_I = \| (\mathbf P^{\mathrm T} \mathbf Q_c)^{\dagger}\|$, $\eta_J = \| (\mathbf S^{\mathrm T} \mathbf Q_r)^{\dagger}\|$, as defined in Theorem (\ref{thm:errbd}).
We note that the additional computational cost of evaluating $\ol{\eta}_p$ and $\ol{\eta}_s$ is negligible, since the associated submatrices have typically small sizes and their inversion is already required as part of the CUR computation.
To this end, we dynamically adjust $m_r$ and $m_c$ such that the condition numbers $\ol{\eta}_p$ and $\ol{\eta}_s$ remain below a user-specified threshold $\epsilon_{\text{os}}$. This adaptive mechanism proceeds as follows:
% \begin{itemize}
% \item {Row Oversampling ($m_r$):}
\begin{itemize}
\item If $\ol{\eta}_p > \epsilon_{\text{os}}$, increase $m_r$  until $\ol{\eta}_p \leq \epsilon_{\text{os}}$.
\item If $\ol{\eta}_p < \epsilon_{\text{os}}$, decrease $m_r$  until the smallest value of $m_r$ is found such that $\ol{\eta}_p \leq \epsilon_{\text{os}}$.
\end{itemize}
The analogous approach is applied to column oversampling, yielding $m_c$. This approach allows for different values of  $m_r$  and  $m_c$, enabling independent control over row and column cross oversampling. 

% \item {Column Oversampling ($m_c$):}
% \begin{itemize}
%     \item If $\eta_J > \epsilon_{\text{os}}$, increment $m_c$ by $1$ (i.e., $m_c \leftarrow m_c + 1$) until $\eta_J \leq \epsilon_{\text{os}}$.
%     \item If $\eta_J < \epsilon_{\text{os}}$, decrement $m_c$ by $1$ (i.e., $m_c \leftarrow m_c - 1$) until the smallest value of $m_c$ is found such that $\eta_J \leq \epsilon_{\text{os}}$.
% \end{itemize}
% \end{itemize}
 % The presented approach allows the oversampling level to adapt dynamically, avoiding insufficient oversampling that could lead to numerical instability, as well as excessive oversampling that would incur unnecessary computational cost.

 The details of \texttt{CUR-CR-OS} are presented in Algorithm (\ref{alg:CR-OS}).

% offers several benefits. It ensures that the selected rows and columns capture sufficient variability in $\mathbf V^k$, avoiding both under-sampling (which may degrade accuracy) and over-sampling (which increases computational cost unnecessarily). It maintains the amplification factor $\epsilon_f$ within a bounded range, ensuring the numerical stability of the CUR decomposition. It enables real-time adjustment as the system evolves, maintaining high approximation fidelity over time. It effectively balances the trade-off between computational efficiency and approximation quality.

\begin{algorithm}[htbp]
\caption{Adaptive CUR oversampling (\texttt{CUR-CR-OS}) algorithm}
\label{alg:CR-OS}
\begin{algorithmic}
\State{\hspace*{\algorithmicindent} \textbf{Input}: $\mathbf A \in \mathbb{R}^{n \times {s}}$, $\mathbf p$, $\mathbf s$ ${m}_r$, ${m}_c$, $r$, $\epsilon_{os}$}
\State{\hspace*{\algorithmicindent} \textbf{Output}: $\mathbf Q_c$, $\mathbf Q_r$, $\mathbf Z$}

\State{$[\mathbf{Q}_c ,\sim] = \texttt{qr}(\mathbf A(:,\mathbf s),\texttt{`econ'})$} %\Comment{Compute the economy QR of solution at ${{\mathbf s}}$.}

\State{$[\mathbf{Q}_r ,\sim] = \texttt{qr}(\mathbf A(\mathbf p,:)^{\mathrm T},\texttt{`econ'})$} %\Comment{Compute the economy QR of solution at ${{\mathbf p}}$.}

\State{$\bar{\mathbf{s}} \gets \texttt{GPODE}(\mathbf {Q}_r, {r}+{m_c})$, $\bar{\mathbf{p}} \gets \texttt{GPODE}(\mathbf Q_c, {r}+{m_r})$ } %\Comment{Compute the oversampled indices.}

\State{$\ol{\eta}_p = \|\mathbf {Q}_c(\bar{\mathbf p},:)^{\dagger}\|$, $\ol{\eta}_s = \|  \mathbf {Q}_r(\bar{\mathbf s},:)^{\dagger}\|$ } %\Comment{Compute condition number for $\mathbf {Q_r}$ and $\mathbf {Q_c}$.}

\State{\textbf{Update: $m_r$}}
\While {$\ol{\eta}_p > \epsilon_{\text{os}}$}
     \State {$m_r \leftarrow m_r + 1$}
     \State{Recompute $\ol{\eta}_p$ : $\bar{\mathbf{p}} \gets \texttt{GPODE}(\mathbf {Q}_c, {r}+m_r)$ ,  $\ol{\eta}_p = \| ( \mathbf {Q}_c(\bar{\mathbf p},:))^{\dagger}\|$}

\EndWhile
\While{$\ol{\eta}_p < \epsilon_{\text{os}}$}
      \State {$m_r \leftarrow m_r - 1$}
      \State {Recompute $\ol{\eta}_p$ : $\bar{\mathbf{p}} \gets \texttt{GPODE}(\mathbf {Q}_c, {r}+m_r)$ , $\ol{\eta}_p = \| ( \mathbf Q_c(\bar{\mathbf p},:))^{\dagger}\|$}
\EndWhile

\State{\textbf{Update: $m_c$}}
\While {$\ol{\eta}_s > \epsilon_{\text{os}}$}
     \State {$m_c \leftarrow m_c + 1$}
     \State{Recompute $\ol{\eta}_c$ : $\bar{\mathbf{s}} \gets \texttt{GPODE}(\mathbf {Q}_r, {r}+m_c)$ ,  $\ol{\eta}_s = \| ( \mathbf {Q}_r(\bar{\mathbf s},:))^{\dagger}\|$}
\EndWhile
\While{$\ol{\eta}_s < \epsilon_{\text{os}}$}
     \State {$m_c \leftarrow m_c - 1$}
      \State {Recompute $\ol{\eta}_s$ : $\bar{\mathbf{s}} \gets \texttt{GPODE}(\mathbf {Q}_r, {r}+m_c)$ , $\ol{\eta}_s = \| ( \mathbf Q_r(\bar{\mathbf s},:))^{\dagger}\|$}
\EndWhile
% \hb{add compute $\mathbf Z$}\\
\State{\textbf{Compute $\mathbf Z$}}: $\mathbf Z =  \mathbf Q_c^{\dagger}(\bar{\mathbf p}, :)  \mathbf A(\bar{\mathbf p}, \bar{\mathbf s}) \left( \mathbf Q_r^{\dagger}(\bar{\mathbf s}, :) \right)^{\mathrm T}$
% \STATE{Use $\eta_s$ to update ${m}_c$ following similar procedure.}
% \STATE Use $\eta_s$ to update $m_c$ following a similar procedure.
\end{algorithmic}
\end{algorithm}

\subsection{Adaptive CUR for reduced-order modeling of MDEs}\label{sec:TDB-CUR}
We demonstrate how \texttt{CUR-CR-OS} can be employed for the time integration of nonlinear MDEs on low-rank matrix manifolds. Our primary focus is on MDEs originating from the discretization of random PDEs, though the methodology is also applicable to deterministic PDEs. In our demonstration, we use DEIM \cite{CS10} for the selection of rows and columns and use GappyPOD+E \cite{peherstorfer2020stability} for the selection of oversampling points.  DEIM requires access to exact or approximate singular vectors of the target matrix. However, these vectors are not available. To this end, we apply DEIM to the singular vectors of the solution at the previous time step, i.e., $\hat{\mathbf A}^{k-1}=\mathbf U^{k-1}\boldsymbol\Sigma^{k-1}{\mathbf Y^{k-1}}^{\mathrm{T}}$. We build a rank-$r$ approximation of $\mathbf A^k$ using \texttt{CUR-CR-OS} as follows:

% \hb{Explain the application of the CUR algorithm here. }
% \begin{enumerate}
% \item Compute the row and column indices by applying the DEIM algorithm to the column and row bases from the previous time step, respectively: $\mathbf p = \texttt{DEIM}(\mathbf U^{k-1})$ and $\mathbf s = \texttt{DEIM}(\mathbf Y^{k-1})$, where $\mathbf U^{k-1}$ and $\mathbf Y^{k-1}$ are the left and right singular vectors of $\hat{\mathbf A}^{k-1}$, respectively.

% \item Compute the $r$ DEIM-selected columns and rows of the FOM matrix: $\mathbf A^k(:,\mathbf s) \in \mathbb{R}^{n \times r}$ and $\mathbf A^k(\mathbf p,:) \in \mathbb{R}^{r \times s}$.
% \item \hb{Call \textttt{CUR-CR-OS}(A(:,s), A(p,:), $m_r^{k-1}$ and $m_c^{k-1}$, and $r^{k-1}$.  }
% \item \hb{Explain the steps on how you change this to SVD}
% \end{enumerate}

% \hb{Explaintha if you use QDEIM you do not need to calculate SVD. }

\begin{enumerate}
\item Determine the row and column index sets by applying the DEIM algorithm to the column and row bases from the previous time step: 
\[
\mathbf{p} = \texttt{DEIM}(\mathbf{U}^{k-1}), 
\quad 
\mathbf{s} = \texttt{DEIM}(\mathbf{Y}^{k-1}),
\]
where $\mathbf{U}^{k-1}$ and $\mathbf{Y}^{k-1}$ are the left and right singular vectors of $\hat{\mathbf{A}}^{k-1}$.  
% \emph{Remark:} If QDEIM is used, the indices are obtained via pivoted QR decompositions of $\mathbf{U}^{k-1}$ and $\mathbf{Y}^{k-1}$, and thus no new SVD of a large matrix is required.

\item Extract the $r$ DEIM-selected columns and rows of the FOM matrix:
\[
\mathbf{A}^k(:,\mathbf{s}) \in \mathbb{R}^{n \times r}, 
\qquad 
\mathbf{A}^k(\mathbf{p},:) \in \mathbb{R}^{r \times s}.
\]

\item Apply the oversampling routine to the sampled blocks using the settings from the previous time step:
\[
\big(\mathbf{Q}_c,\mathbf{Q}_r,\mathbf{Z}\big) =
\texttt{CUR-CR-OS}\!\Big(\mathbf{A}^k(\mathbf{p},:),
\mathbf{A}^k(:,\mathbf{s}),\,
m_r^{k-1}, m_c^{k-1}, r^{k-1}, \epsilon_{\text{os}}
\Big)\]
Here, \texttt{CUR-CR-OS} utilizes $\mathbf{A}^k(\mathbf{p},:)$ and $\mathbf{A}^k(:,\mathbf{s})$. It also adaptively determines $\overline{\mathbf{p}}$ and $\overline{\mathbf{s}}$ and samples $\mathbf{A}^k(\ol{\mathbf{p}},\ol{\mathbf{s}})$  to compute $\mathbf Z$:
\[
\mathbf{Z}
= \mathbf{Q}_c^{\dagger}(\overline{\mathbf{p}},:)\,
\mathbf{A}^k(\overline{\mathbf{p}},\overline{\mathbf{s}})\,
\big(\mathbf{Q}_r^{\dagger}(\overline{\mathbf{s}},:)\big)^{\mathrm{T}},
\]
while adjusting $(m_r,m_c)$ until the conditioning indicators $\overline{\eta}_p$ and $\overline{\eta}_s$ fall below $\epsilon_{\text{os}}$.

\item Compute the SVD of the small $r \times r$ matrix,
\[
\mathbf{Z} = \mathbf{\Psi}_c\,\boldsymbol{\Sigma}^k\,\mathbf{\Psi}_r^{\top},
\]
and rotate the orthonormal factors to update the low-rank state:
\[
\mathbf{U}^k = \mathbf{Q}_c\mathbf{\Psi}_c,
\qquad 
\mathbf{Y}^k = \mathbf{Q}_r\mathbf{\Psi}_r.
\]
\end{enumerate}
In Step 4, we express $\mathbf{A}^k$ in the SVD form. The reason is that these singular vectors are needed in the next time step for DEIM. However, if, for example, QDEIM is used, there is no need to perform Step 4, and $\mathbf{A}^k$ can be stored as $\mathbf{A}^k= \mathbf Q_c \mathbf Z \mathbf Q_c^{\mathrm T}$ since QDEIM can be applied to  $\mathbf Q_c $ and $\mathbf Q_r $. 

% The advantage of the above algorithm is that it constructs a rank-$r$ approximation by evaluating $\mathbf{A}^k$ at only $r$ columns and rows. Thus, we only require access to a total of $(n+s)r - r^2$ entries in the matrix $\mathbf A^k$. This CUR-based approach addresses the key challenges associated with the time integration of nonlinear MDEs on low-rank manifolds, as outlined earlier in this section. Moreover, it is minimally intrusive, enabling integration with existing solvers with little modification. 

\subsubsection{Rank adaptivity}\label{subsec:Adap_r}
Rank is a critical hyperparameter that governs the approximation error in reduced-order models and directly impacts the accuracy of TDB-ROMs. However, selecting an appropriate rank a priori is challenging. The optimal rank is problem-dependent, as the decay of singular values varies from one PDE to another. Furthermore, even within a single PDE, the required rank must change over time to maintain a desired accuracy threshold dynamically. These considerations strongly motivate the development of rank-adaptive methods. Accordingly, rank adaptivity has garnered considerable attention in recent years, leading to the emergence of several excellent techniques for various time integrations schemes \cite{dektor2021rank, HP22,ceruti2022rank}.

In \cite{TDBCUR}, a rank adaptivity criterion based on the ratio of the last resolved singular values to the Frobenius norm of the low-rank matrix was proposed.  The underlying rationale is that the relative magnitude of the last resolved singular value serves as a proxy for the approximation error. In this work, we use a different criterion based on the actual error. 
To this end,  we propose an error proxy, denoted by $\epsilon$, that estimates the accuracy of the low-rank approximation at a subset of points:
\begin{equation}\label{eq:epsilon}
     \epsilon = \frac{\| \mathbf A(\bar{\mathbf p}_r,\bar{\mathbf s}_r) - \hat{\mathbf A}(\bar{\mathbf p}_r,\bar{\mathbf s}_r)\|}{(r+\bar{m}_r)(r+\bar{m}_c)}.
\end{equation}
Here, $\bar{m}_r$ and $\bar{m}_c$ represent $1\%$ of the total row and column indices, respectively. The corresponding oversampled index sets $\bar{\mathbf p}_r \in \mathbb{R}^{r + \bar{m}_r}$ and $\bar{\mathbf s}_r \in \mathbb{R}^{r + \bar{m}_c}$ are selected using the GappyPOD+E algorithm ~\cite{peherstorfer2020stability}. These oversampled rows and columns are used consistently, regardless of the values of $m_r$ and $m_c$ obtained via the adaptive oversampling strategy in Section \ref{subsec:Adap_m}, since the latter may be zero.

The error proxy $\epsilon$ measures the regression error between the true matrix and its low-rank approximation over the oversampled entries. Rank adaptivity aims to minimize this error by adjusting the rank $r$. Rather than enforcing a fixed error tolerance, we define an upper and lower threshold, $\epsilon_u$ and $\epsilon_l$, respectively. The rank is updated adaptively as follows: if $\epsilon > \epsilon_u$, increase the rank: $r \leftarrow r + 1$ and if $\epsilon < \epsilon_l$, decrease the rank: $r \leftarrow r - 1$.

This adaptive procedure ensures that the rank is neither too low (causing underfitting) nor unnecessarily high (wasting computational resources). Importantly, this rank adjustment strategy does not require setting a strict threshold for $\epsilon$, allowing flexibility across applications. We emphasize that this is just one possible criterion for rank adaptation. Other strategies may be employed depending on the nature and requirements of the problem. Nevertheless, in our numerical experiments, this approach has proven to be both simple and effective, closely tracking the trend of the true error. Furthermore, the rank can also be incremented by more than one if necessary, allowing rapid adaptation until $\epsilon$ falls below the upper threshold $\epsilon_u$.

Algorithm \texttt{TDB-CUR (CR-OS)}  (\ref{alg:TDB-CUR-OS}) presents the procedure for solving MDEs in low-rank form using the adaptive cross oversampling approach.

\begin{algorithm}[htbp]
\caption{\texttt{TDB-CUR (CR-OS)} algorithm}
\label{alg:TDB-CUR-OS}
\begin{algorithmic}
\State {\hspace*{\algorithmicindent} \textbf{Input}: ${\mathbf{U}^{k-1}} \in \mathbb{R}^{n \times {r}}$, ${\boldsymbol{\Sigma}^{k-1}} \in \mathbb{R}^{{r} \times {r}}$, ${\mathbf{Y}^{k-1}}\in \mathbb{R}^{s \times {r}}$, ${r}^{k-1}$, ${m}_r^{k-1}$, ${m}_c^{k-1}$, $\bar{m}_r^{k-1}$, $\bar{m}_c^{k-1}$, $\epsilon_{os}$, $\epsilon_u$, $\epsilon_l$} \\

\State {\hspace*{\algorithmicindent} \textbf{Output}: $\mathbf{U}^{k}\in\mathbb R^{n\times r}$, $\boldsymbol{\Sigma}^{k}\in\mathbb R^{r\times r}$, $\mathbf{Y}^{k}\in\mathbb R^{s\times r}$, ${r}^{k}$, ${m}_r^{k}$, ${m}_c^{k}$}

\State {$\mathbf{s} \gets \texttt{DEIM}({\mathbf{Y}^{k-1}}, {r}^{k-1})$, $\mathbf{p} \gets \texttt{DEIM}({\mathbf U^{k-1}}, {r}^{k-1})$} %\Comment{Compute DEIM indices.}

\State {$\bar{\mathbf{s}} \gets \texttt{GPODE}({\mathbf{Y}^{k-1}}, {r}^{k-1}+{m_c}^{k-1})$, $\bar{\mathbf{p}} \gets \texttt{GPODE}({\mathbf U^{k-1}}, {r}^{k-1}+{m_r}^{k-1})$}
% \Comment{Compute oversampled indices.}

\State {$\bar{\mathbf{s}}_r \gets \texttt{GPODE}({\mathbf{Y}^{k-1}}, {r}^{k-1}+\bar{m}_c^{k-1})$, $\bar{\mathbf{p}}_r \gets \texttt{GPODE}({\mathbf U}, {r}^{k-1}+\bar{m}_r^{k-1})$}
% \Comment{Compute oversampled indices for rank adaptivity.}

\State {$\hat{\mathbf{A}}(:, \mathbf{s}) = {\mathbf{U}^{k-1}} {\boldsymbol{\Sigma}^{k-1}} {\mathbf{Y}^{k-1}}(\mathbf{s}, :)^{\mathrm T}$, $\hat{\mathbf{A}}(\mathbf{p}, :) = {\mathbf{U}^{k-1}}(\mathbf{p}, :) {\boldsymbol{\Sigma}^{k-1}} {\mathbf{Y}^{k-1}}^{\mathrm T}$}
% \Comment{Construct low-rank approximation.}

\State {$\hat{\mathbf{A}}(\bar{\mathbf{p}}, \bar{\mathbf{s}}) = {\mathbf{U}^{k-1}}(\bar{\mathbf{p}}, :) {\boldsymbol{\Sigma}^{k-1}} {\mathbf{Y}^{k-1}}(\bar{\mathbf{s}}, :)^{\mathrm T}$ , $\hat{\mathbf{A}}(\bar{\mathbf{p}}_r, \bar{\mathbf{s}}_r) = {\mathbf{U}^{k-1}}(\bar{\mathbf{p}}_r, :) {\boldsymbol{\Sigma}^{k-1}} {\mathbf{Y}^{k-1}}(\bar{\mathbf{s}}_r, :)^{\mathrm T}$}
% \Comment{Construct cross approximation.}

\State {$\mathbf{A}(:, \mathbf{s}) = \hat{\mathbf{A}}(:, \mathbf{s}) + \Delta t \bar{\mathbf{F}}(:, \mathbf{s})$, $\mathbf{A}(\mathbf{p}, :) = \hat{\mathbf{A}}(\mathbf{p}, :) + \Delta t \bar{\mathbf{F}}(\mathbf{p}, :)$}
% \Comment{Update solution at indices $\mathbf p$ and $\mathbf J$.}
\State{ $\mathbf{A}(\bar{\mathbf{p}}, \bar{\mathbf{s}}) = \hat{\mathbf{A}}(\bar{\mathbf{p}}, \bar{\mathbf{s}}) + \Delta t \bar{\mathbf{F}}(\bar{\mathbf{p}}, \bar{\mathbf{s}})$}
% \Comment{Update solution at the oversampled cross entries $\bar{\mathbf{p}}$ and $\bar{\mathbf{s}}$.}

\State {$\mathbf{A}(\bar{\mathbf{p}}_r, \bar{\mathbf{s}}_r) = \hat{\mathbf{A}}(\bar{\mathbf{p}}_r, \bar{\mathbf{s}}_r) + \Delta t \bar{\mathbf{F}}(\bar{\mathbf{p}}_r, \bar{\mathbf{s}}_r)$}
% \Comment{Update solution at cross entries $\bar{\mathbf{p}}_r$ and $\bar{\mathbf{s}}_r$ for rank adaptivity.}

\State{$\mathbf{Q}_c\mathbf{R}_c = \texttt{QR}(\mathbf A(:,\mathbf s),\texttt{'econ'})$}%\Comment{Compute the economy QR of $\mathbf A(:,\mathbf s)$.}}
\State {$\mathbf{Q}_r\mathbf{R}_r = \texttt{QR}(\mathbf A(\mathbf p,:)^{\mathrm T},\texttt{'econ'})$}%$\Comment{Compute the economy QR of $\mathbf A(\mathbf p,:)^{\mathrm T}$.}}

\State{$\mathbf{Z}=\mathbf{Q}_c^{\dagger}(\overline{\mathbf{p}},:)\,
\mathbf{A}^k(\overline{\mathbf{p}},\overline{\mathbf{s}})\,
\big(\mathbf{Q}_r^{\dagger}(\overline{\mathbf{s}},:)\big)^{\mathrm{T}}$}

\State{$\mathbf{\Psi}_c \boldsymbol{\Sigma}^{k} \mathbf{\Psi}_r^{\mathrm T} = \texttt{SVD}(\mathbf{Z},\texttt{'econ'}$)}%\Comment{Compute the economy SVD of $\mathbf C$.}
\State{$\mathbf{U}^{k}=\mathbf{Q}_c \mathbf{\Psi}_c$ }%\Comment{In-subspace rotation of the orthonormal basis, $\mathbf Q_c$.}}

\State{$\mathbf{Y}^{k}=\mathbf{Q}_r \mathbf{\Psi}_r$}% \Comment{In-subspace rotation of the orthonormal basis, $\mathbf Q_r$.}}

\State{\textbf{Rank adaptivity:} $\displaystyle
\epsilon \gets 
\frac{\left\lVert 
\mathbf A(\bar{\mathbf p}_r,\bar{\mathbf s}_r)
- \mathbf U^k(\bar{\mathbf p}_r,:)\,\boldsymbol{\Sigma}^k\,{\mathbf Y^k(\bar{\mathbf s}_r,:)}^{\mathrm{T}}
\right\rVert}
{(r+\bar m_r)(r+\bar m_c)}$}

% \WHILE {$\eta_p > \epsilon_{\text{os}}$}
%      \STATE {$\hat{m}_r \leftarrow m_r + 1$}
%      \STATE{Recompute $\eta_p$ : $\bar{\mathbf{p}} \gets \texttt{GPODE}(\mathbf {Q_c}, {r}+{\hat{m}_r})$ ,  $\eta_p = \| ( \mathbf {Q_c}(\bar{\mathbf p},:))^{\dagger}\|$}

% \ENDWHILE

\If{$\epsilon > \epsilon_u$}
  \State{ $ r^k \gets r^{k-1} + 1$}
\ElsIf{$\epsilon < \epsilon_l$}
  \State{$ r^k \gets r^{k-1} - 1$}
  \State{$\mathbf U^k \gets \mathbf U^k(:,1{:}  r^k)$,  $\boldsymbol{\Sigma}^k \gets \boldsymbol{\Sigma}^k(1{:} r^k,\,1{:} r^k)$, $\mathbf Y^k \gets \mathbf Y^k(:,1{:} r^k)$}
\EndIf

\State{\textbf{Adaptive Oversampling:}}{\quad $[\hat{m}_r^k, \hat{m}_c^k] \leftarrow \text{\texttt{CUR-CR-OS}}$}
% \Comment{Algorithm \ref{alg:CR-OS}}
\end{algorithmic}
\end{algorithm}

\section{Demonstration Cases}\label{sec:demo}
\subsection{Matrix Low-Rank Approximation}
% \hb{Use Figure \ref{xxx} instead of fig \ref{xxx}.}
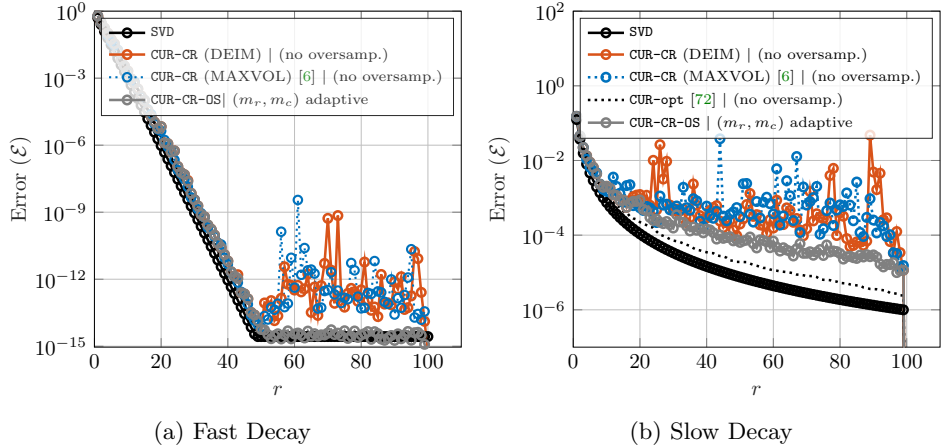
\begin{figure}[t]
  \centering

  \begin{subfigure}[b]{0.48\linewidth}
    \centering
    \input{figures/Err_r_fast}
    \caption{Fast Decay}
    \label{fig:Toy_problem1-a}
  \end{subfigure}
  \begin{subfigure}[b]{0.48\linewidth}
    \centering
    \input{figures/Err_r_slow}
    \caption{Slow Decay}
    \label{fig:Toy_problem1-b}
  \end{subfigure}

  % \begin{subfigure}[b]{0.48\linewidth}
  %   \centering
  %   \input{figures/Err_Cost_fastdecay}
  %   \caption{Fast Decay Cost}
  %   \label{fig:Toy_problem1-c}
  % \end{subfigure}
  % \begin{subfigure}[b]{0.48\linewidth}
  %   \centering
  %   \input{figures/Err_Cost_slowdecay}
  %   \caption{Slow Decay Cost}
  %   \label{fig:Toy_problem1-d}
  % \end{subfigure}

  \caption{Comparative analysis of fast and slow decay for matrix low-rank approximation: (a), (b) relative error $\mathcal{E}$ versus reduction order $r$.}
  % (c),(b) relative error $\mathcal{E}$ versus \% of matrix entries.}
  \label{fig:Toy_problem1}
\end{figure}

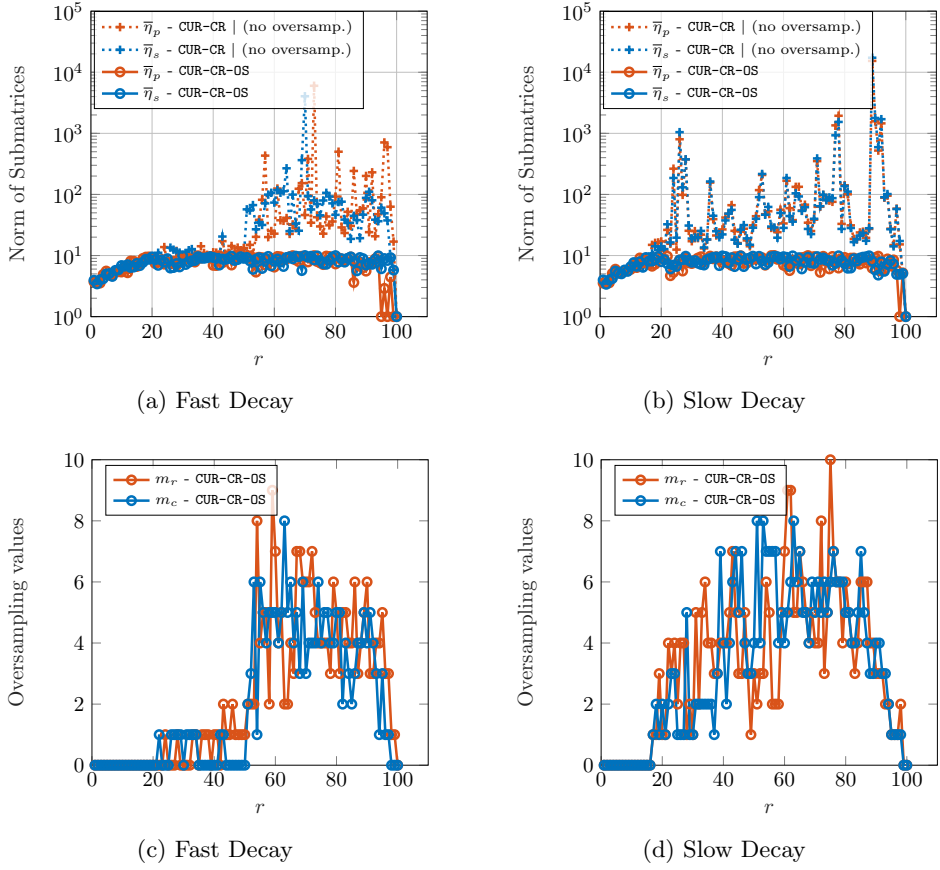
\begin{figure}
  \centering

  \begin{subfigure}[b]{0.48\linewidth}
    \centering
    \input{figures/Amp_I_J_r_fast}
    \caption{Fast Decay}
    \label{fig:Toy_problem2-a}
  \end{subfigure}\hfill
  \begin{subfigure}[b]{0.48\linewidth}
    \centering
    \input{figures/Amp_I_J_r_slow}
    \caption{Slow Decay}
    \label{fig:Toy_problem2-b}
  \end{subfigure}

  % \vspace{0.5em}

  \begin{subfigure}[b]{0.48\linewidth}
    \centering
    \input{figures/Toy_m_fast}
    \caption{Fast Decay}
    \label{fig:Toy_problem2-c}
  \end{subfigure}\hfill
  \begin{subfigure}[b]{0.48\linewidth}
    \centering
    \input{figures/Toy_m_slow}
    \caption{Slow Decay}
    \label{fig:Toy_problem2-d}
  \end{subfigure}

  % \vspace{0.5em}

  % \begin{subfigure}[b]{0.48\linewidth}
  %   \centering
  %   \input{figures/Amp_Err_f}
  %   \caption{Fast Decay}
  % \end{subfigure}\hfill
  % \begin{subfigure}[b]{0.48\linewidth}
  %   \centering
  %   \input{figures/Amp_Err_s}
  %   \caption{Slow Decay}
  % \end{subfigure}

  \caption{Comparative analysis of fast and slow decay for matrix low-rank approximation: (a), (b) condition number for row $\eta_p$ and column $\eta_s$ for \texttt{CUR-CR} ($m=0$) i.e. without oversampling and \texttt{CUR-CR-OS} i.e. with adaptive column and row oversampling versus reduction order $r$. (c),(d) adaptive row ($m_r$) and column ($m_c$) oversampling used for \texttt{CUR-CR-OS} versus reduction order $r$.}
  \label{fig:Toy_problem2}
  % (e), (f) Error Ratio Bound and True Error Ratio for reduction order $r = 10, 50$ versus oversampling value $m$.
\end{figure}

In this example, we demonstrate the performance of \texttt{CUR-OS} algorithm. We consider the matrix given by
$$
\mathbf A=\left(e^{ \mathbf W_1}\right) e \mathbf D \left(e^{ \mathbf W_2}\right)^T,
$$
where the matrices $\mathbf W_1 \in \mathbb{R}^{n \times n}$ and $\mathbf W_2 \in \mathbb{R}^{n \times n}$ are randomly generated skew-symmetric matrices as follows: $\mathbf{W}_1=\left(\tilde{\mathbf{W}}_1-\tilde{\mathbf{W}}_1^T\right) / 2$ and $\mathbf{W}_2=\left(\tilde{\mathbf{W}}_2-\tilde{\mathbf{W}}_2^T\right) / 2$, where $\tilde{\mathbf{W}}_1 \in \mathbb{R}^{n \times n}$ and $\tilde{\mathbf{W}}_2 \in \mathbb{R}^{n \times n}$ are uniformly distributed random matrices and $n=100$ for all the cases considered. We examine two scenarios: \emph{fast decay} and \emph{slow decay}. The decay rate is governed by a diagonal matrix $\mathbf{D} \in \mathbb{R}^{n \times n}$, whose diagonal entries are given by $d_i = 1 / 2^i$ for the fast decay case and $d_i = 1 / i^3$ for the slow decay case, with $i \in {1,2,\ldots,n}$.

We use the relative $\ell_2$ error defined as  
\[
\mathcal{E} = \frac{\| \hat{\mathbf A} - \mathbf A \|_{2}}{\| \mathbf A \|_{2}},
\]  
where $\hat{\mathbf A}$ denotes the rank-$r$ approximation of the matrix $\mathbf A$ obtained using a low-rank approximation method.

In Figure \ref{fig:Toy_problem1-a},  the relative errors $\mathcal{E}$ for ranks \(1 \leq r \leq n\) for different algorithms are shown. For \texttt{CUR-CR}, we consider two different sampling algorithms of DEIM and MAXVOL. For \texttt{CUR-CR-OS}, DEIM is used for sampling the $r$ columns and rows, and GappyPOD+E is used for cross oversampling.  While both CUR without oversampling (\texttt{CUR-CR})  and \texttt{CUR-CR-OS} with adaptive oversampling perform similarly for small ranks, significant differences emerge as \(r\) increases.  The \texttt{CUR-CR} low-rank matrices constructed with DEIM and MAXVOL exhibit similarly large fluctuations in the approximation error, whereas \texttt{CUR-CR-OS} closely follows the SVD errors. This instability is even more pronounced in the slow-decay case, where \texttt{CUR-CR} becomes unstable as early as \(r = 20\) as shown in Figure \ref{fig:Toy_problem1-b}. In this figure, we also show the error of the optimal CUR ($\texttt{CUR-opt}$) \cite{SE16}, which yields the optimal CUR rank-$r$ approximation in the Frobenius norm.

In Figures~\ref{fig:Toy_problem2-a} and \ref{fig:Toy_problem2-b}, we plot the norms of the row and column submatrices ($\overline{\eta}_p$ and $\overline{\eta}_s$) versus rank. Without oversampling, these norms can grow significantly, particularly at higher ranks. For the slow-decay matrix, this growth begins at lower ranks.  The onset of large norms coincides with the increase in error observed in Figures~\ref{fig:Toy_problem1-a} and \ref{fig:Toy_problem1-b}, indicating that these easily computable norms provide an effective a posteriori indicator of the \texttt{CUR-CR} instabilities.    In \texttt{CUR-CR-OS}, adaptive oversampling is used where the norm of these submatrices is maintained below $\epsilon_{os}=10$.   

Figures~\ref{fig:Toy_problem2-c} and \ref{fig:Toy_problem2-d} show the number of entries selected for cross oversampling in the rows and columns. When  \texttt{CUR-CR} is stable, little to no oversampling is required. The amount of oversampling increases with rank, but then diminishes as CUR approaches the full-rank case, i.e., $r \approx n$.

% The instability observed at higher ranks is attributed to the ill-conditioning of the matrices $(\mathbf P^T \mathbf Q_c)^{\dagger}$ and $(\mathbf S^T \mathbf Q_r)^{\dagger}$, as reflected in the growth of the condition numbers \(\eta_p\) and \(\eta_s\) respectively, shown in Figure \ref{fig:Toy_problem2}(a) and \ref{fig:Toy_problem2}(b). In contrast, the \texttt{CUR-CR-OS} method, which employs adaptive cross row and column oversampling, successfully controls these condition numbers, resulting in improved stability. In Figure \ref{fig:Toy_problem2}(c) and \ref{fig:Toy_problem2}(d), we can see the values of $m_r$ and $m_c$ used in \texttt{CUR-CR-OS} for adaptive row and column oversampling.

% In Figure \ref{fig:Toy_problem2}(e) and \ref{fig:Toy_problem2}(f), we plot the error ratio bound \(\epsilon_f\) and the true error ratio = \(\left( \frac{\| \mathbf{V} - \hat{\mathbf{V}} \|}{\hat{\sigma}_{r+1}} \right)\) as functions of the oversampling parameter \(m\), where $m = m_r = m_c$, for ranks \(r = 10\) and \(r = 50\). As expected from the theoretical inequality in Eq. \eqref{eq:err_bd},  
% $\| \mathbf{A} - \hat{\mathbf{A}} \| \leq \epsilon_f \hat{\sigma}_{r+1}
% $, the true error ratio remains below the error ratio bound across all values of \(m\). Moreover, we observe that as the oversampling parameter increases, the error ratio bound converges to the true error ratio, confirming the tightness of the bound for both ranks.

\subsection{Application to Stochastic PDEs}

We will apply \texttt{CUR-CR-OS}  to a diverse set of stochastic PDEs described in Table \ref{table:1}. In particular, we consider one-dimensional Burgers, Allen-Cahn, and Korteweg–De Vries (KdV) equations, subject to stochastic boundary conditions, stochastic initial conditions, or both. Each case involves a $ d$-dimensional random space of $\bs \xi =[\xi_1, \xi_2, \dots, \xi_d]$. The random dimension in each case is shown in  Table \ref{table:1}.  In all cases considered in this section, $\xi_i$ random variables are zero-mean Gaussian with the standard deviation of 1, i.e., $\xi_i \sim \mathcal N(0,1)$. For stochastic sampling, we use Monte Carlo sampling and draw $s$ samples from the joint distribution of $\bs \xi$.  For the discretization of the spatial domain, we use the second-order finite difference scheme on a uniform grid of size $n$. We use the fourth-order explicit Runge-Kutta method to integrate in time with time-step $\Delta t$. 

\begin{table}[thbp]
\centering
\setlength{\tabcolsep}{5pt} %Gap before text Starts
\renewcommand{\arraystretch}{1.5} %Cell Height scaling
\setlength{\arrayrulewidth}{0.2mm} %Table Border Thickness
\begin{tabular}{|c|c|c|}
    \hline
    \multicolumn{3}{|c|}{\textbf{Random PDEs}} \\ \hline
    \textbf{Burgers} & \textbf{Allen-Cahn} & \textbf{Korteweg–De Vries} \\\hline

     $\dfrac{\partial v}{\partial t}+ \dfrac{1}{2} \dfrac{\partial v^2}{\partial x} = \nu \dfrac{\partial^{2} v}{\partial x^{2}}$ & $\dfrac{\partial v}{\partial t} = \nu \dfrac{\partial^{2} v}{\partial x^{2}} -v^3 + v $& $\dfrac{\partial v}{\partial t} +v \dfrac{\partial v}{\partial x} = \gamma \frac{\partial^{3} v}{\partial x^{3}}$ \\ 
     $x \in[0,1], t \in[0,5] $ & $x \in[0,2\pi], t \in[0,50]$ & $x \in[0,10], t \in[0,1]$ \\
     $n = 401, s = 256$& $n = 256, s = 256$& $n = 1024, s = 256$\\
      $ \nu=2.5\times 10^{-3}, d = 17$&  $\nu=5\times 10^{-3}, d = 4$ & $\gamma=2\times 10^{-4}, d = 4$\\\hline
     % Stochastic Boundary& Periodic Boundary  & Periodic Boundary \\ \hline
      $\Delta t = 2.5 \times 10^{-4}$&$\Delta t = 10^{-2}$ & $\Delta t = 1 10^{-4}$\\ \hline
    %  $ \sigma \sum_{i=1}^{d} \frac{\sin(i\pi t) \xi_{i}}{i^2}$ & &  \\ 
    % $ v(1, t; \bs\xi) = 0$ & & \\ \hline

\end{tabular}
\caption{Problem setup for random PDEs.}
\label{table:1}
\end{table}

The problem setup for the Burgers equation is identical to the case considered in \cite{TDBCUR}.  The stochasticity is introduced via the Dirichlet boundary condition at $x=0$ according to:
\[
v(0, t ; \boldsymbol{\xi})=-\sin (2 \pi t)+\sigma \sum_{i=1}^d \frac{\sin (i \pi t) \xi_i t}{i^2},\]
as well as the initial condition 
\[ v(x, 0 ; \boldsymbol{\xi})=  \sin (2 \pi x)\left[0.5\left(e^{\cos (2 \pi x)}-1.5\right) +\sigma \sum_{i=1}^d \sqrt{\lambda_{x_i}} \psi_i(x) \xi_i\right],\]
where $\sigma  = 0.001$. The deterministic Dirichlet boundary condition is used at the right boundary ($x=1$). 
  For the stochastic initial condition, $\lambda_{x_{i}}$ and $\psi_{i}(x)$ are the eigenvalues and eigenvectors of the spatial squared-exponential kernel, respectively. 
 
 For Allen-Cahn and Korteweg-De-Vries equations, the stochasticity is introduced via the initial condition. We consider:
 \[
 v(x,0;\bs \xi)=\\
 e^{-27(x-4.2)^2}-e^{-23.5\left(x-\frac{\pi}{2}\right)^2} +e^{-38(x-5.4)^2}+\frac{\tanh (2 \sin (x))}{3} +\sigma \sum_{i=1}^d \lambda_i \psi_i(x) \xi_i\]
  for Allen-Cahn and 
 \[
 v(x,0;\bs \xi)=\frac{1}{40} \log \left(1+\frac{\cosh (20)^2}{\cosh (20(x-2))^2}\right) +\sigma \sum_{i=1}^d \lambda_i \psi_i(x) \xi_i\]
 for Korteweg-De-Vries. Here, $\lambda_i$ and $\psi_i $ are the eigenvalues and eigenvectors of the spatial squared-exponential kernel, respectively. Construct the spatial squared-exponential kernel matrix $K_{i j}=\exp \left(-\frac{|x_i-x_j|^2}{2 \ell^2}\right)$ on the spatial grid. Then compute it's eigen decomposition $\mathbf{K}=\mathbf{L} \mathbf\Lambda \mathbf{L}^{\mathrm{T}}$; the eigenvalues are $\lambda_i$ (diagonal of $\mathbf\Lambda$ ) and the eigenvectors are $\psi_i$ columns of $\mathbf L$ (optionally truncating small $\lambda_i$). Periodic boundary conditions are for both Allen-Cahn and Korteweg-De-Vries.  We select the first $r$ sampling indices with QDEIM and the oversampling points with GappyPOD+E.
In all cases, the reference solution is obtained by integrating the FOM in time using the same $\Delta t$ as employed for the TDB-CUR models.

\begin{figure}
\centering
  \begin{subfigure}[b]{0.48\linewidth}
    \centering
    \input{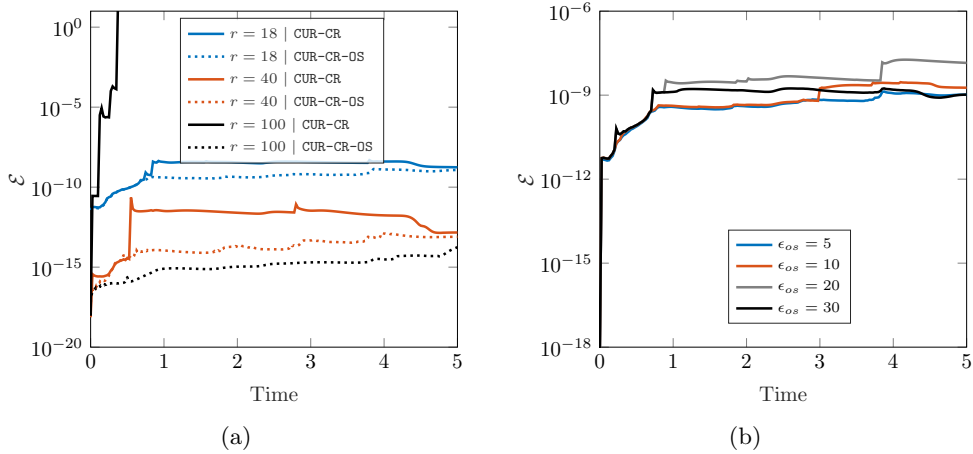}
     \caption{}
    \label{fig:Burger_Sensitivity_Cost-a}
  \end{subfigure}\hfill
  \begin{subfigure}[b]{0.48\linewidth}
    \centering
    \input{figures/Burger_Sensitivity}
    \caption{}
    \label{fig:Burger_Sensitivity_Cost-b}
  \end{subfigure}
  % \begin{subfigure}[b]{0.50\linewidth}
  %   \centering
  %   \input{figures/CpuCost_m}
  %   \caption{}
  %   \label{fig:Burger_Sensitivity_Cost-c}
  %   % \caption{Cost plot}   
  %   \end{subfigure}
\caption{Stochastic Burgers: (a) Analysis of the normalized \texttt{TDB-CUR-CR-OS} error ($\mathcal{E}$) to the rank ($r = 18, 40, 100$) without oversampling ($m=0$, solid line) and with adaptive oversampling ($m_r, m_c$  adaptive, dotted line). (b) Sensitivity study of the effect of chosen oversampling bound $(\epsilon_{os})$ on the normalized error, with constant rank $r=18$. }
% (c) CPU cost comparison for \texttt{CUR-CR} with only row oversampling \cite{TDBCUR} and \texttt{CUR-CR-OS} with sampling size $s=10^7$ for oversampling value $m = 20, 40, 60, 80, 100$ at constant rank $r=18$.
\label{fig:Burger_Sensitivity_Cost}
\end{figure}

% \subsubsection{Stochastic Burgers Equation}
In Figure \ref{fig:Burger_Sensitivity_Cost}, we evalaute different aspects of the performance of \texttt{CUR-CR-OS} for the stochastic Burgers case. 
In Figure \ref{fig:Burger_Sensitivity_Cost-a}, the average Frobenius error defined as  $\mathcal{E} = \| \hat{\mathbf A} - \mathbf A \|_{\mathcal{F}}/ns$  for CUR with no oversampling (\texttt{CUR-CR}) and CUR with cross oversampling  (\texttt{CUR-CR-OS}) is shown.  Three ranks of $r=18,40$ and $100$ are considered. Ranks r=40 and r=100 represent overapproximation cases, as they exceed the rank needed to capture singular values on the order of $\mathcal{O}(10^{-14}$). We observe that \texttt{CUR-CR} performs reasonably well for r=18, exhibits a sudden increase in error at r=40, and diverges at r=100. In contrast, \texttt{CUR-CR-OS} maintains good performance across all cases, yielding stable time integration.

In the \texttt{CUR-CR-OS}, $\epsilon_{os}$ controls the accuracy of the oversampling algorithm. In  Figure \ref{fig:Burger_Sensitivity_Cost-b}, the average Frobenius error for  $\epsilon_{os}=\{5, 10, 20, 30 \}$ for a fixed-rank approximation ($r=18$) is shown.  This demonstrates that the error remains relatively insensitive when $\epsilon_{os} \leq 10$.

\begin{figure}[t!]
\centering

\begin{minipage}[t]{1\textwidth}
    \centering
    \includegraphics[width=\textwidth]{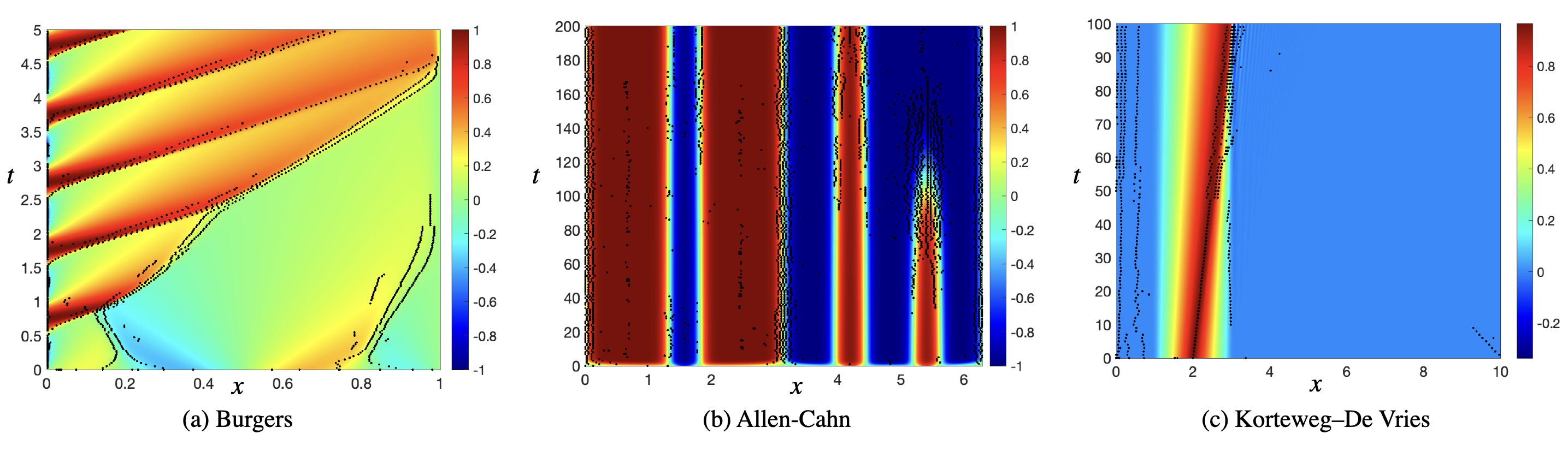}
\end{minipage}
\caption{Mean flow solution for Burgers, Allen-Cahn, and Korteweg-De Vries equations with the rank adaptive QDEIM points (black dots).}
\label{fig:mean_flow_spde}
\end{figure}

Figure \ref{fig:mean_flow_spde} shows the mean solution - obtained by averaging the columns of $\mathbf{A}^k$ - for the three stochastic PDEs, along with the rank-adaptive sampling points with tolerance $\epsilon_u = 10^{-8}$ and adaptive oversampling with $\epsilon_{os} = 10$. The oversampling points are omitted; only the $r$ selected row samples are displayed. The figure illustrates the evolution of row sampling over time, with many sampling points concentrated in regions of steep spatial gradients, where stochastic variations are most pronounced.

\begin{figure}[htbp]
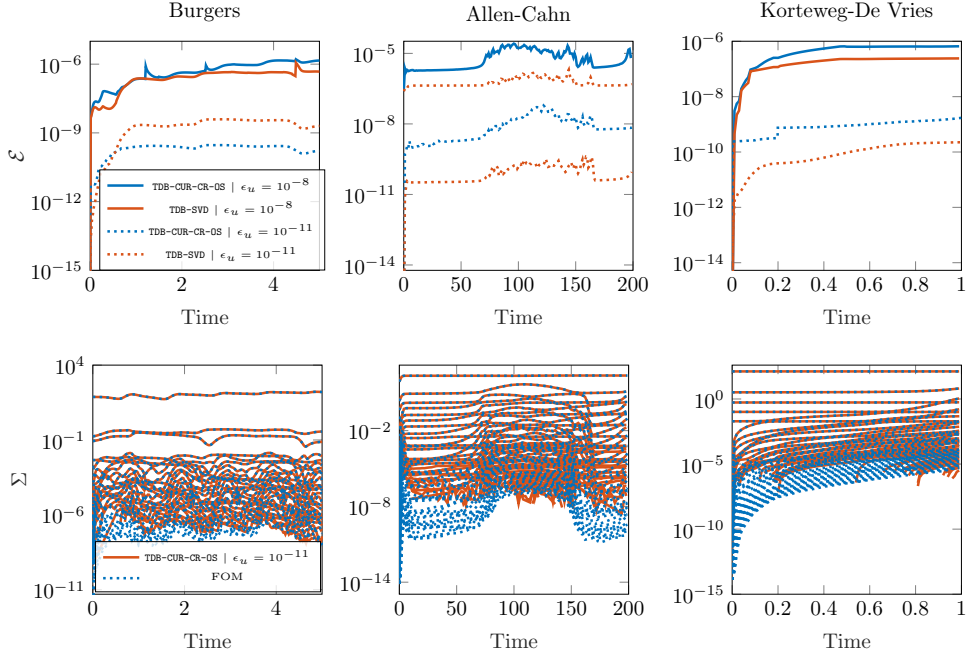

  \centering

  \begin{subfigure}[t]{0.34\linewidth}
    \centering
    \input{figures/Error_Burgers}
    % \caption{Burgers}
  \end{subfigure}\hfill
  \begin{subfigure}[t]{0.34\linewidth}
    \centering
    \input{figures/Error_AC}
    % \caption{Allen–Cahn}
  \end{subfigure}\hfill
  \begin{subfigure}[t]{0.32\linewidth}
    \centering
    \input{figures/Error_KdV}
    % \caption{KdV}
  \end{subfigure}

  \vspace{0.6em}

  \begin{subfigure}[t]{0.34\linewidth}
    \centering
    \input{figures/SingVal_Burgers}
    % \caption{Burgers singular values}
  \end{subfigure}\hfill
  \begin{subfigure}[t]{0.34\linewidth}
    \centering
    \input{figures/SingVal_AC}
    % \caption{Allen–Cahn singular values}
  \end{subfigure}\hfill
  \begin{subfigure}[t]{0.32\linewidth}
    \centering
    \input{figures/SingVal_KdV}
    % \caption{KdV singular values}
  \end{subfigure}

  \caption{Errors and singular values vs.\ time for Burgers, Allen–Cahn, and KdV with $\epsilon_u=10^{-8}$ and $\epsilon_u=10^{-10}$.}
  \label{fig:comparison_spde}
\end{figure}

In Figure~\ref{fig:comparison_spde}, we examine the errors and singular values over time for the  Burgers, Allen–Cahn, and KdV equations. We consider two values for the rank adaptivity threshold: $\epsilon_u = 10^{-8}$ and $10^{-11}$, which serve as upper bounds for the error proxy $\epsilon$ defined in Eq. (\ref{eq:epsilon}). For adaptive oversampling, we bound $\ol{\eta}_p$ and $\ol{\eta}_s$ by setting $\epsilon_{\text{os}} = 10$, as discussed in \S\ref{subsec:Adap_m}. In the first row of Figure~\ref{fig:comparison_spde}, we compare the \texttt{CUR-CR-OS} error with those obtained by performing SVD as a rank truncation method at each time step instead of CUR, i.e., $\hat{\mathbf A}_{\texttt{SVD}}^k = \texttt{SVD}(\hat{\mathbf A}_{\texttt{SVD}}^{k-1} + \Delta t \overline{\mathbf F} )$. 
 We observe that the \texttt{CUR-CR-OS} error closely follows the SVD error for all three PDEs. Moreover, reducing the threshold $\epsilon_u$ results in lower errors. The second row compares the singular values obtained from the \texttt{CUR-CR-OS} approximation (for $\epsilon_u = 10^{-11}$) with those from the FOM solution. We observe strong agreement between the two, demonstrating the accuracy and robustness of the proposed method.

\begin{figure}[h]
  \centering

  \begin{subfigure}[t]{0.34\linewidth}
    \centering
    \input{figures/rank_Burgers}
    % \caption{Rank (Burgers)}
  \end{subfigure}\hfill
  \begin{subfigure}[t]{0.34\linewidth}
    \centering
    \input{figures/rank_AC}
    % \caption{Rank (Allen–Cahn)}
  \end{subfigure}\hfill
  \begin{subfigure}[t]{0.32\linewidth}
    \centering
    \input{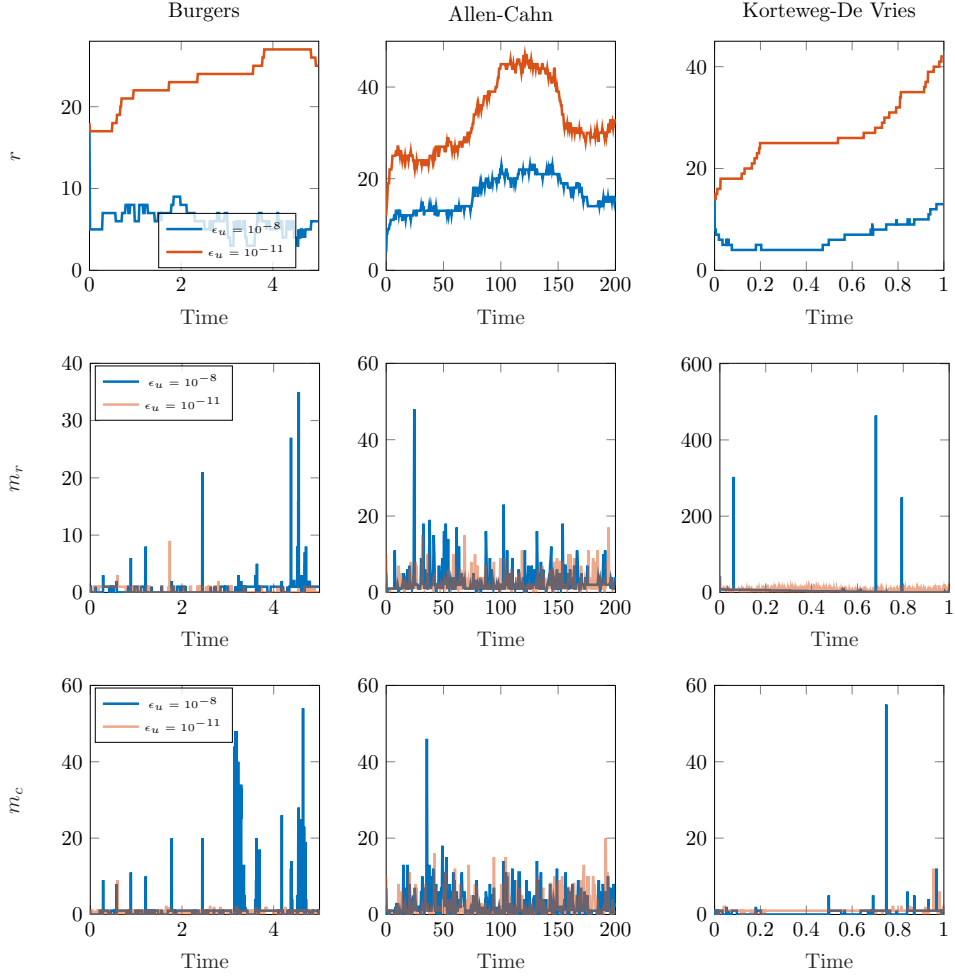}
    % \caption{Rank (KdV)}
  \end{subfigure}

  \vspace{0.6em}

  \begin{subfigure}[t]{0.34\linewidth}
    \centering
    \input{figures/mr_Burgers}
    % \caption{$m_r$ (Burgers)}
  \end{subfigure}\hfill
  \begin{subfigure}[t]{0.34\linewidth}
    \centering
    \input{figures/mr_AC}
    % \caption{$m_r$ (Allen–Cahn)}
  \end{subfigure}\hfill
  \begin{subfigure}[t]{0.32\linewidth}
    \centering
    \input{figures/mr_KdV}
    % \caption{$m_r$ (KdV)}
  \end{subfigure}

  \vspace{0.6em}

  \begin{subfigure}[t]{0.34\linewidth}
    \centering
    \input{figures/mc_Burgers}
    % \caption{$m_c$ (Burgers)}
  \end{subfigure}\hfill
  \begin{subfigure}[t]{0.34\linewidth}
    \centering
    \input{figures/mc_AC}
    % \caption{$m_c$ (Allen–Cahn)}
  \end{subfigure}\hfill
  \begin{subfigure}[t]{0.32\linewidth}
    \centering
    \input{figures/mc_KdV}
    % \caption{$m_c$ (KdV)}
  \end{subfigure}

  \caption{Rank and oversampling vs.\ time for Burgers, Allen–Cahn, and KdV with $\epsilon_u=10^{-8}$ and $\epsilon_u=10^{-10}$.}
  \label{fig:comparison_spde2}
\end{figure}

Figure~\ref{fig:comparison_spde2} presents the adaptive rank and oversampling values of \texttt{TDB-CUR-CR-OS} for the three SPDEs. The first row shows the rank evolution over time, governed by the rank-adaptivity mechanism. A stricter tolerance ($\epsilon_u = 10^{-11}$) leads to substantially higher ranks than a looser tolerance ($\epsilon_u = 10^{-8}$). The second and third rows display the temporal evolution of row and column oversampling values, respectively. The markedly different behaviors observed across the three SPDEs underscore the importance of an adaptive strategy for robustness and efficiency, as opposed to prescribing a fixed oversampling level.

\subsection{Heat Conduction-Radiation Problem}\label{HCR}
\begin{figure}[t!]
\centering

\begin{subfigure}[t]{0.48\textwidth}
    \centering
    \includegraphics[width=\textwidth]{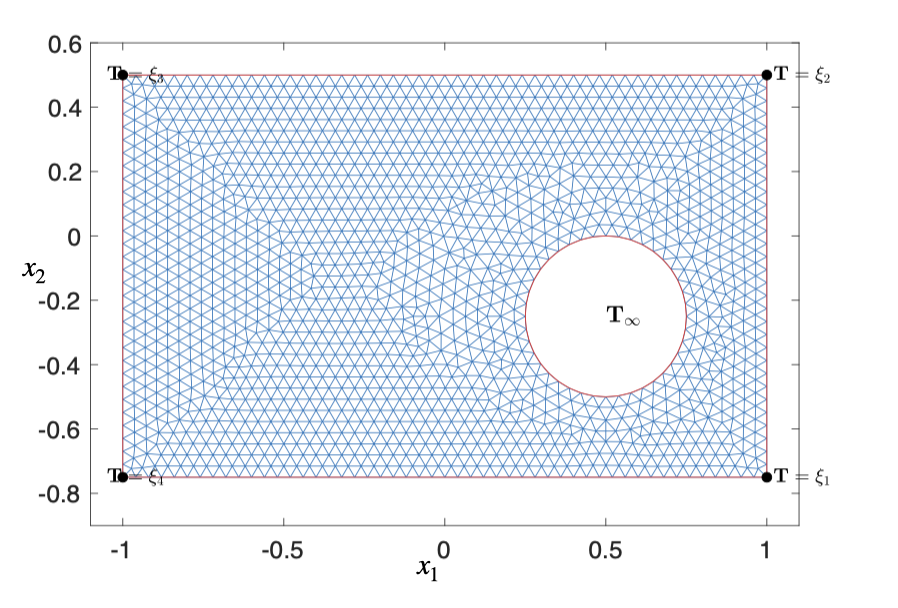}
    \caption{Geometry}
\end{subfigure}
\hfill
\begin{subfigure}[t]{0.48\textwidth}
    \centering
    \includegraphics[width=\textwidth]{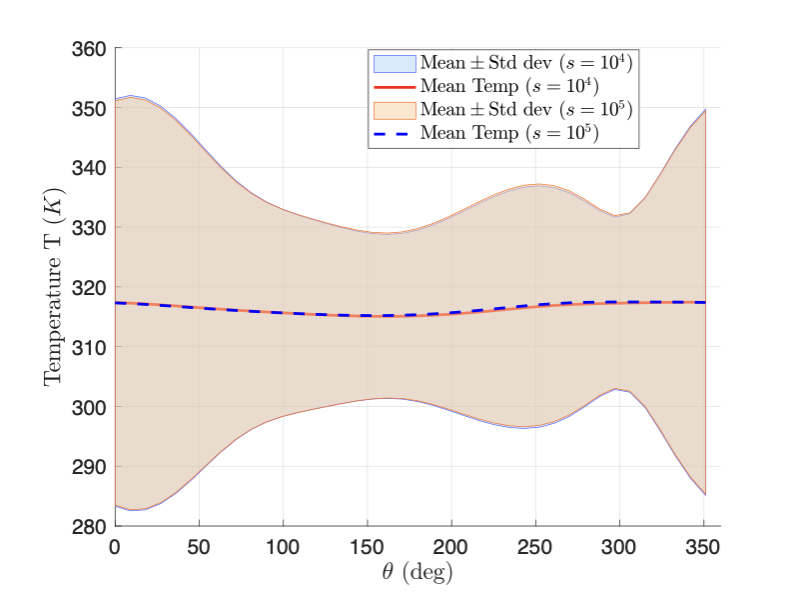}
    \caption{Cylinder Temperature}
\end{subfigure}

\caption{Heat Conduction Radiation Problem: (a) Finite element mesh of the transient heat conduction radiation problem. (b) Mean cylinder temperature.}
\label{Fig:HCR_geometry_Mean_Var}
\end{figure}

We solve a transient conduction radiation heat transfer problem.
The geometry of the problem, along with the finite element mesh, is shown in Figure \ref{Fig:HCR_geometry_Mean_Var}(a).
The cylinder is centered at $ (x_{1_0}, x_{2_0} )=(0.5,-0.25)$ and its radius is $r_0=0.25$. The governing equation is given by:
\begin{equation}\label{eq:HPC_0}
     \rho c_p \frac{\partial  T}{\partial t}=k\left(\frac{\partial^2  T}{\partial x_1^2}+\frac{\partial^2  T}{\partial x_2^2}\right),
\end{equation}
%\hb{use $(x_1,x_2)$ instead of $(x,y)$ to be consistent with Eq. (5). } 
where $(x_1, x_2)$ denote the spatial coordinates, $T$ is the temperature, $t$ is time, $k$ is the thermal conductivity, $\rho$ is the density and $c_p$ is the specific heat capacity. The plate is built with Carbon Steel with the following material properties: $k=50 \mathrm{~W} / \mathrm{m} \cdot \mathrm{K}$, $\rho=7850 \mathrm{~kg} / \mathrm{m}^3$, and $c_p=500 \mathrm{J} / \mathrm{kg} \cdot \mathrm{K}.$

 The temperature on the boundary of the outer edges of the rectangle is assumed to be of stochastic Dirichlet type. 
%\hb{do not use the edge ids to free to boundary conditions. Not a standard approach. You can refer to boundaries as cylinder and outer edges.} 
We consider linear variation across each edge. As a result, the entire boundary
profiles on the outer edges can be parameterized by the value of temperature at four corners. This
parameterization ensures that the temperature profile on the outer edges varies continuously from one
edge to another as the joint edges share the same value corner temperature value. The temperature at
the four corners is devoted by parameters $\bs \xi = \{\xi_1,\xi_2,\xi_3,\xi_4\}$, where
$\xi_1$ is the temperature at the $(x_1, x_2) = (1, -0.75)$, $\xi_2$ is the temperature at the $(x_1, x_2) = (1, 0.5)$, $\xi_3$ is the temperature at the $(x_1, x_2) = (-1, 0.5)$, and $\xi_4$ is the temperature at the $(x_1, x_2) = (-1, -0.75)$. The corner temperatures are considered to be parameters in the range of $273 \mathrm{~K} \leq \xi_i \leq 373\mathrm{~K}.$

 The cylinder surface loses or absorbs heat flux through radiation as follows:
\begin{equation}\label{eq:HCR_1}
    -k \frac{\partial T}{\partial \hat{n}} = \epsilon \sigma \left(T^4 - \T_{\infty}^4\right), 
    \quad (x_1, x_2) \in \text{cylindrical surface}
\end{equation}
where $\hat{n}$ is a unit normal vector to the surface of the insulated cylinder, $\epsilon=0.2$ is the emissivity, and $\sigma=5.67 \times 10^{-8} \mathrm{~W} / \mathrm{m}^2 \mathrm{~K}^4$ is the Stefan-Boltzmann constant, and $T_{\infty}$ is the ambient temperature. We assume $T_{\infty}=273 \mathrm{~K}$ in this case.

The FOM is solved using the finite element method (FEM) with linear (P1) basis functions. The FEM matrices are generated with the MATLAB FEM solver. We split the nodes into interior and Dirichlet boundary sets. As a result, the nodes on the cylinder are considred as interior points, since their values evolve with time.  The evolution equation for the interior nodes can be written as:  
\begin{equation}\label{eq:HCR_2}
    \mathbf{M}_{i i} \frac{d \mathbf{T}_i}{d t}=\mathbf{K}_{i i} \mathbf{T}_i+\mathbf{K}_{i B} \mathbf{T}_B-\epsilon \sigma \mathbf{g}_i \circ\left(\mathbf{T}_i^4-\mathbf{T}_{\infty}^4\right),
\end{equation}
where $\mathbf{M}_{i i}$ and $\mathbf{K}_{i i}$ represent the mass and stiffness matrices of the interior points, respectively, $\mathbf{T}_i$ and $\mathbf{T}_B$ are the temperature vectors at the interior and boundary points. The vector $\mathbf{g}_i$ is an indicator vector, which has the same size as $\mathbf{T}_i$ and is zero everywhere except at the points on the cylinder surface where its values are equal to 1. The symbol o denotes the element-wise product. Appendix (\ref{appendix:FEM}) shows the detailed derivation of Eq. (\ref{eq:HCR_2}).

We consider time in the interval of $0 \leq t \leq T_f$  uniformly distributed
time snapshots, where $T_f = \rho c_p/k = 78500$ (sec). We consider  $\Delta t = T_f/ 8000$. 

% Directly solving Eq. (\ref{eq:HCR_2}) requires inverting the mass matrix $\mathbf M_{ii}$. However, this poses a significant computational challenge: while $\mathbf M_{ii}$ is sparse, its inverse is dense, leading to high memory and computational costs. In such cases, applying the \texttt{CUR-CR-OS} algorithm would require sampling all elements of the inverse, resulting in a computational cost comparable to that of solving the FOM, thus negating any potential efficiency gains.

% To mitigate this problem we will diagonalize $\mathbf M_{ii}$ using diagonal scaling method described by O.C. Zienkiewicz in \cite{zienkiewicz2005finite}, section 16.2.4. The diagonal scaling method is a lumping technique used to approximate a consistent mass matrix $\mathbf{M}_{ii}$ by a diagonal matrix $\tilde{\mathbf{M}}_{ii}$, in a way that preserves total mass. The key idea is that you scale the diagonal entries of $\mathbf{M}_{ii}$ proportionally so that the total mass (i.e., sum of the diagonal entries) remains equal to that of the consistent matrix. $\tilde{\mathbf{M}}_{ii} = c \mathbf{M}_{ii}$, where $c $ is the ratio of all elements of $\mathbf M_{ii}$ divided by sum of all diagonal elements of $\mathbf M_{ii}$. This ensures that $\sum \mathbf{M}_{ii} = \tilde{\mathbf{M}}_{ii}$, thus preserving total mass.

We consider the diagonal scaling method for mass matrix lumping, as described in Ref. \cite[Section 16.2.4 ]{zienkiewicz2005finite}. This technique approximates the mass matrix $\mathbf{M}_{ii}$ with a diagonal matrix $\tilde{\mathbf{M}}_{ii}$ while preserving the total mass. The approach involves proportionally scaling the diagonal entries of $\mathbf{M}_{ii}$ as follows: set $\tilde{\mathbf{M}}_{ii} = c \cdot \text{diag}(\mathbf{M}_{ii})$, where the scaling factor $c$ is given by the ratio of the sum of all entries in $\mathbf{M}_{ii}$ to the sum of its diagonal entries. This ensures that $\sum \mathbf{M}_{ii} = \sum \tilde{\mathbf{M}}_{ii}$, thus preserving total mass. Hence instead of solving Eq. (\ref{eq:HCR_2}), we will solve,
\begin{equation}\label{eq:HCR_3}
    \frac{d \mathbf{T}_i}{d t}=\tilde{\mathbf{M}}_{ii} ^{-1}(\mathbf{K}_{i i} \mathbf{T}_i+\mathbf{K}_{i B} \mathbf{T}_B-\epsilon \sigma \mathbf{g}_i \circ\left(\mathbf{T}_i^4-\mathbf{T}_{\infty}^4\right)),
\end{equation}

% \gp{Need help rewriting this para well-}
% In the heat conduction–radiation problem, strong nonlinearities lead to high computational costs when using methods such as DLRA\cite{koch2007dynamical, koch2010dynamical} and the BUG integrators \cite{ceruti2022unconventional,ceruti2022rank,ceruti2023rank}, which require access to the full matrix. In contrast, \texttt{CUR-CR-OS}  significantly reduces computational complexity by sampling only the dominant rows and columns, making it more scalable for such problems.

\begin{figure}
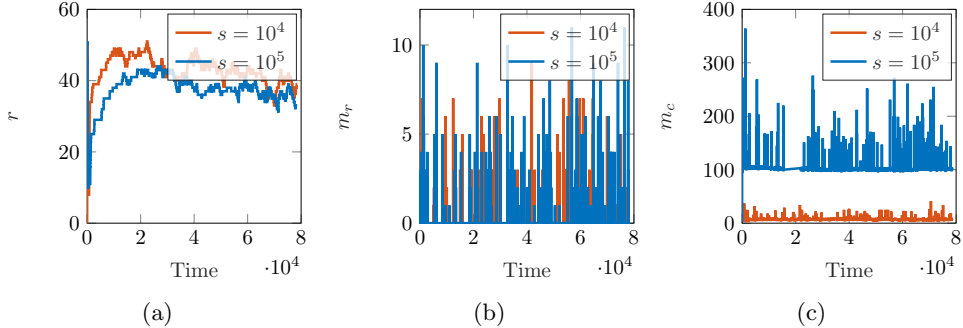

\centering

    \begin{subfigure}[t]{0.34\linewidth}
    \centering
    \input{figures/Rank_HCR}
    \caption{}
  \end{subfigure}
  \begin{subfigure}[t]{0.32\linewidth}
    \centering
    \input{figures/Oversamp_r_HCR}
    \caption{}
 \end{subfigure}
   \begin{subfigure}[t]{0.32\linewidth}
    \centering
    \input{figures/Oversamp_c_HCR}
    \caption{}
 \end{subfigure}
\caption{Heat Conduction Radiation Problem: Analysis for two sample sizes $s = 10 ^4$ and $s = 10 ^5$. (a) Adaptive rank $r$, (b) adaptive row oversampling $m_r$, and (c) adaptive column oversampling $m_c$.}
  \label{fig:HCR_r_m}
\end{figure}

We solve the problem for large sample sizes $s = 10^4$ and $s = 10^5$. In Figure~\ref{Fig:HCR_geometry_Mean_Var}(b), we plot the mean as well as the standard deviation of the temperature along the cylinder. The mean temperature profiles for both sample sizes are nearly identical, showing that the size of the samples is sufficient for computing the mean and variance. 
In Figure \ref{fig:HCR_r_m}(a) and (b), we compare the evolution of the rank and adaptive oversampling values over time for two different sample sizes, $s = 10^4$ and $s = 10^5$, in the TDB-CUR method. The required rank increases rapidly at early times, reflecting the need to capture transient dynamics. The rank then plateaus as the system approaches a quasi-steady regime. The adaptive oversampling value remains relatively constant throughout the simulation. For $s = 10^5$, the oversampling is significantly higher for the column oversampling, probably reflecting the need to account for a wider range of variability and more fine details in the data set.

% \begin{figure}
%   \centering

%   \begin{subfigure}[b]{0.48\linewidth}
%     \centering
%     \input{figures/Err_r_fast}
%     \caption{Fast Decay}
%   \end{subfigure}
%   \begin{subfigure}[b]{0.48\linewidth}
%     \centering
%     \input{figures/Err_r_slow}
%     \caption{Slow Decay}
%   \end{subfigure}

%   \begin{subfigure}[b]{0.48\linewidth}
%     \centering
%     \input{figures/Err_Cost_fastdecay}
%     \caption{Fast Decay Cost}
%   \end{subfigure}
%   \begin{subfigure}[b]{0.48\linewidth}
%     \centering
%     \input{figures/Err_Cost_slowdecay}
%     \caption{Slow Decay Cost}
%   \end{subfigure}

% % In your preamble (if not already):

\section{Conclusions}\label{sec:con}
\label{sec:conclusions}
In this work, we introduced a cost-efficient and adaptive oversampling framework for CUR decompositions aimed at solving large-scale matrix differential equations (MDEs) on low-rank manifolds. Recognizing that oversampling is essential for enhancing the stability of CUR approximations, we developed a novel cross-oversampling strategy that avoids the computational overhead of sampling entire rows and columns by focusing instead on a sparse set of cross entries. This significantly reduces the cost of data acquisition while preserving or even improving the quality of the approximation.

To further enhance the robustness and scalability of the method, we proposed two adaptive mechanisms:
An adaptive oversampling scheme based on the amplification factor and a rank-adaptive strategy based on a residual proxy error, enabling dynamic adjustment of the low-rank approximation in response to the evolving complexity of the system.

We derived an error bound for the resulting \texttt{CUR-CR-OS}  approximation and demonstrated both theoretically and empirically that our approach improves upon the row oversampling scheme of prior methods \cite{TDBCUR}. Numerical experiments on toy problems and non-linear stochastic PDEs (including Burgers, Allen–Cahn, and Korteweg–de Vries equations) confirmed that the proposed algorithm achieves high accuracy, stability, and computational efficiency across a range of scenarios. We further applied our method to high-dimensional heat conduction problems on complex geometries, showing that it scales effectively with the problem size and sample dimension. The adaptive oversampling and rank-adjustment strategies presented here make the method highly robust, reducing sensitivity to hyperparameter choices and allowing it to generalize across a wide class of problems.

\section*{Acknowledgments}
The authors would like to express their gratitude to their colleague, Sara Akhavan, for providing the solvers for the stochastic PDEs.

\section{Appendix}
% \subsection{Error bound}\label{appendix:ErrorBound}
% In the next theorem, we will show the error bound for the TDB-CUR (CR-OS) algorithm \ref{table:1}. The proof will follow a similar procedure to that used in the theorem 2.8 of \cite{TDBCUR}. Here, the indexing matrices $\bm P$ and $\bm S$ are of sizes $n\times r'$ and $s\times r'$, respectively and $r' \geq r$.
% \begin{theorem}\label{thm:errbd}
%     Let $\bm V \in  \mathbb R^{n\times s}$ and $\hat{\bm  V} $ be the rank-$r$ TDB-CUR (CR-OS) approximation to $\bm V$. Here, $\bm Q_r$ and $\bm Q_c$ are obtained from QR factorization of $\bm V$ at $\bm I$ rows and $\bm J$ columns, respectively. Let $\epsilon_f = \text{min}\{\eta_I(1+\eta_J),\eta_J(1+\eta_I)\}$, where  $\eta_I = \| (\bm Q_c(\bm I,:))^{\dagger}\|$, $\eta_J = \| (\bm Q_r(\bm J,:))^{\dagger}\|$ and $\hat{\sigma}_{r+1}= \text{max} \{ \| \bm V -\bm Q_c \bm Q_c^{\mathrm T} \bm V\| ,  \| \bm V -\bm V \bm Q_r \bm Q_r^{\mathrm T} \|\}$. Then the error of the oblique projection is bounded by
%     \begin{equation}\label{eq:err_bd}
%         \| \bm V - \hat{\bm  V}\| \leq \epsilon_f\hat{\sigma}_{r+1} \end{equation}
% \end{theorem}

\subsection{ Interior-Boundary Domain Splitting in FEM}\label{appendix:FEM}
The original equation we have is the following:
\[
\mathbf{M}\,\frac{\mathrm{d}\mathbf{T}}{\mathrm{d}t}
=
\mathbf{K}\mathbf{T}
\;-\;
\epsilon\sigma\mathbf{g}\circ\bigl(\mathbf{T}^4 - T_\infty^4\bigr).
\]
where \(T\in\mathbb{R}^N\) is the temperature matrix in all mesh nodes of \(N\). \(\mathbf M\in\mathbb{R}^{N\times N}\) is the mass matrix. \(\mathbf K\in\mathbb{R}^{N\times N}\) is the stiffness matrix. The symbol “\(\circ\)” denotes element‐wise multiplication. Splitting the temperature, mass, and stiffness matrix into interior and boundary points, we have
\[
\mathbf{T} = 
\begin{bmatrix}\mathbf{T}_i\\[6pt]\mathbf{T}_B
\end{bmatrix}
,\qquad
\mathbf{M} =
\begin{bmatrix}
\mathbf{M}_{ii} & \mathbf{M}_{iB} \\[6pt]
\mathbf{M}_{Bi} & \mathbf{M}_{BB}
\end{bmatrix}
,\qquad
\mathbf{K} =
\begin{bmatrix}
\mathbf{K}_{ii} & \mathbf{K}_{iB} \\[6pt]
\mathbf{K}_{Bi} & \mathbf{K}_{BB}
\end{bmatrix}
\]

$\mathbf{T}_i$ is the unknown temperature at the interior nodes, $\mathbf{T}_B$ is the known Dirichlet temperature at boundary nodes. $\mathbf{M}_{ii}$ is the mass matrix block that couples interior degrees of freedom (DOFs) with themselves. It arises from integrating shape functions over the interior of the domain and directly influences the time evolution of interior temperatures. $\mathbf{M}_{iB}$, $\mathbf{M}_{Bi}$ is the mass coupling between interior and boundary nodes. This block completes the symmetric structure of the mass matrix. $\mathbf{M}_{BB}$ is the mass matrix over boundary DOFs. $\mathbf{K}_{ii}$ is the stiffness matrix block for interactions between interior DOFs. 
$\mathbf{K}_{iB}$ is the stiffness coupling between the boundary and interior nodes, and $\mathbf{K}_{Bi}=\mathbf{K}_{iB}^{\mathrm T}$.
Now, by solving this system, we have two coupled equations:
\begin{align}
\mathbf{M}_{ii} \dot{\mathbf{T}}_i + \mathbf{M}_{iB} \dot{\mathbf{T}}_B
&=
\mathbf{K}_{ii} \mathbf{T}_i + \mathbf{K}_{iB} \mathbf{T}_B -\epsilon\sigma\mathbf{g}_i\circ(\mathbf{T}_i^4 - \mathbf{T}_\infty^4),
\label{eq:interior_block}\\[6pt]
\mathbf{M}_{Bi}\,\dot{\mathbf{T}}_i +\mathbf{M}_{BB}\,\dot{\mathbf{T}}_B
&=
\mathbf{K}_{Bi}\mathbf{T}_i + \mathbf{K}_{BB}\mathbf{T}_B -\epsilon\sigma\mathbf{g}_B\circ(\mathbf{T}_B^4 - \mathbf{T}_\infty^4).
\label{eq:boundary_block}
\end{align}

The temperature at the Dirichlet boundary is fixed, so
\[
\mathbf{\dot T}_B = 0.
\]
We can use equation \eqref{eq:boundary_block} as a consistency check and equation \eqref{eq:interior_block} as the main equation. The resulting ODE for the unknown interior temperatures $\mathbf{T}_i$ is
\[
\mathbf{M}_{ii}\frac{\mathrm{d}\mathbf{T}_i}{\mathrm{d}t}
=
\mathbf{K}_{ii}\,\mathbf{T}_i + \mathbf{K}_{iB}\,\mathbf{T}_B
-
\epsilon\sigma\mathbf{g}_i \circ (\mathbf{T}_i^4 - \mathbf{T}_\infty^4).
\]
% In the final equation, the term  $\mathbf{K}_{iB} \dot{\mathbf{T}}_B$ appears as a known forcing term that accounts for the influence of the prescribed boundary temperatures on the interior nodes.

% \bibliographystyle{siamplain}
\bibliographystyle{unsrt}
\bibliography{references,Hessam}
\end{document}

%% file: figures/Err_r_fast.tex
\definecolor{mycolor1}{rgb}{0.00000,0.44700,0.74100}%
\definecolor{mycolor2}{rgb}{0.85000,0.32500,0.09800}%
\definecolor{mycolor3}{rgb}{0.5,0.5,0.5}
\begin{tikzpicture}[scale=0.8]

\begin{axis}[%
width=2.4in,
height=2.2in,
at={(1.011in,0.642in)},
scale only axis,
xmin=0,
xmax=110,
xlabel style={font=\color{white!15!black}},
xlabel={$r$},
ymode=log,
ymin=9.10521593419642e-16,
ymax=1,
yminorticks=true,
ylabel style={font=\color{white!15!black}},
ylabel={Error ($\mathcal{E}$)},
axis background/.style={fill=white},
xmajorgrids,
ymajorgrids,
%yminorgrids,
legend style={fill opacity=0.8, legend cell align=left, font=\scriptsize,align=left, draw=white!15!black}
]
\addplot [color=black, line width=1.2pt,mark=o, mark options={solid, black}]
  table[row sep=crcr]{%
% 1	0.5\\
% 6	0.015625\\
% 11	0.000488281249999996\\
% 16	1.52587890625064e-05\\
% 21	4.76837158206147e-07\\
% 26	1.49011611965367e-08\\
% 31	4.65661284412446e-10\\
% 36	1.4551913714439e-11\\
% 41	4.54748817996246e-13\\
% 46	1.42165289089878e-14\\
% 51	2.66519969664132e-15\\
% 56	2.66594708660211e-15\\
% 61	2.66690036574752e-15\\
% 66	2.66728055057133e-15\\
% 71	2.66682516048327e-15\\
% 76	2.66733485939598e-15\\
% 81	2.66836210698518e-15\\
% 86	2.66853675687933e-15\\
% 91	2.66941595064921e-15\\
% 96	2.66892050457089e-15\\
1	0.5\\
2	0.25\\
3	0.125\\
4	0.0625\\
5	0.03125\\
6	0.015625\\
7	0.00781249999999999\\
8	0.00390624999999999\\
9	0.001953125\\
10	0.000976562499999994\\
11	0.000488281249999995\\
12	0.000244140625000001\\
13	0.000122070312500006\\
14	6.10351562499921e-05\\
15	3.05175781250067e-05\\
16	1.52587890625036e-05\\
17	7.6293945312459e-06\\
18	3.81469726562783e-06\\
19	1.90734863281508e-06\\
20	9.53674316404865e-07\\
21	4.76837158204211e-07\\
22	2.38418579104617e-07\\
23	1.19209289548903e-07\\
24	5.96046447702437e-08\\
25	2.98023223843118e-08\\
26	1.49011611998579e-08\\
27	7.45058059456053e-09\\
28	3.72529030296827e-09\\
29	1.86264515264495e-09\\
30	9.31322578113473e-10\\
31	4.65661293856696e-10\\
32	2.3283064632951e-10\\
33	1.1641532591305e-10\\
34	5.82076638873279e-11\\
35	2.91038218833091e-11\\
36	1.45519104097013e-11\\
37	7.2759620672413e-12\\
38	3.63798227363015e-12\\
39	1.8189881431935e-12\\
40	9.09497101172744e-13\\
41	4.54750425824627e-13\\
42	2.27370942120804e-13\\
43	1.13687027677109e-13\\
44	5.68469904301387e-14\\
45	2.84132401314214e-14\\
46	1.42120918244995e-14\\
47	7.10303222663371e-15\\
48	3.54791999252674e-15\\
49	2.70955471940679e-15\\
50	2.71036377298737e-15\\
51	2.71095604735791e-15\\
52	2.7105622267361e-15\\
53	2.71042257248127e-15\\
54	2.71088031538346e-15\\
55	2.71196137444729e-15\\
56	2.71217512219644e-15\\
57	2.71228988864786e-15\\
58	2.71265785288079e-15\\
59	2.71262949657756e-15\\
60	2.71218459824864e-15\\
61	2.71213308409764e-15\\
62	2.71252528968877e-15\\
63	2.71234043766083e-15\\
64	2.71290019734472e-15\\
65	2.71324456403862e-15\\
66	2.71356592370076e-15\\
67	2.71347167578088e-15\\
68	2.7139400583796e-15\\
69	2.7137563183285e-15\\
70	2.71357609581124e-15\\
71	2.71374053106012e-15\\
72	2.71408200738311e-15\\
73	2.71389690396055e-15\\
74	2.71317541349573e-15\\
75	2.71317824799264e-15\\
76	2.71276315868065e-15\\
77	2.71274872586547e-15\\
78	2.71197355196972e-15\\
79	2.71175297769799e-15\\
80	2.71151659845682e-15\\
81	2.71121736537427e-15\\
82	2.71104875424554e-15\\
83	2.71101049463029e-15\\
84	2.71123033907864e-15\\
85	2.71098469063338e-15\\
86	2.71092760269558e-15\\
87	2.710874344292e-15\\
88	2.7110165486418e-15\\
89	2.71071490838422e-15\\
90	2.71049327998484e-15\\
91	2.71064877835394e-15\\
92	2.7107131023673e-15\\
93	2.71055148513234e-15\\
94	2.71030522029895e-15\\
95	2.71035833249863e-15\\
96	2.71069195501787e-15\\
97	2.71050838302378e-15\\
98	2.71101422713943e-15\\
99	2.71116755528897e-15\\
100	2.71115259955627e-15\\
};
\addlegendentry{\texttt{SVD}}

\addplot [color=mycolor2, line width=1.2pt,mark=o, mark options={solid, mycolor2}]
  table[row sep=crcr]{%
% 1	0.601542638828667\\
% 6	0.0331198742100311\\
% 11	0.000974481781024973\\
% 16	3.72243621369948e-05\\
% 21	1.34398326499121e-06\\
% 26	5.98241953603287e-08\\
% 31	1.85825774037014e-09\\
% 36	9.19975782935928e-11\\
% 41	2.16811387475095e-12\\
% 46	1.20106340124163e-13\\
% 51	3.57875480972611e-14\\
% 56	7.96755847172213e-13\\
% 61	2.24152383140705e-12\\
% 66	8.51627337572838e-14\\
% 71	3.36509232390094e-12\\
% 76	1.62100366116228e-13\\
% 81	6.61008862300297e-13\\
% 86	7.36252660949446e-14\\
% 91	3.01136656428044e-13\\
% 96	2.37831242490086e-12\\
1	0.612017026941478\\
2	0.307492955625974\\
3	0.172597870753999\\
4	0.10315183594602\\
5	0.055238042882207\\
6	0.0347730455811607\\
7	0.0168155854116056\\
8	0.00663698807943946\\
9	0.00350182976668751\\
10	0.00209722904025338\\
11	0.000936163661808515\\
12	0.000588667909838891\\
13	0.000330926852474536\\
14	0.000157760331298599\\
15	0.000101192724614388\\
16	3.70729090709351e-05\\
17	1.99614132844045e-05\\
18	1.07540680489655e-05\\
19	5.6163752454094e-06\\
20	6.76389903492232e-06\\
21	1.26944155354616e-06\\
22	1.17387481777424e-06\\
23	4.92657320801356e-07\\
24	5.55922664023971e-07\\
25	1.1082205467812e-07\\
26	5.91535798822913e-08\\
27	4.78409374783443e-08\\
28	2.65825628025339e-08\\
29	7.01004042472748e-09\\
30	3.61183757287041e-09\\
31	1.91655855187787e-09\\
32	8.01103813623095e-10\\
33	6.25673740074713e-10\\
34	2.63499454773507e-10\\
35	1.73510618533137e-10\\
36	9.46052311717348e-11\\
37	4.07798869334737e-11\\
38	2.59833370005925e-11\\
39	1.25876144033457e-11\\
40	6.51345286969213e-12\\
41	2.15841102208393e-12\\
42	1.36623015161726e-12\\
43	1.59843022309598e-12\\
44	2.47960111249532e-13\\
45	2.18582922601466e-13\\
46	1.32669941239442e-13\\
47	3.25271476396699e-14\\
48	2.50168147980242e-14\\
49	1.33020816886694e-14\\
50	6.24906336590305e-15\\
51	8.66905366145931e-14\\
52	6.16907822860124e-14\\
53	1.10536533095182e-13\\
54	1.39548992104517e-14\\
55	2.29543680243729e-14\\
56	1.40980935604214e-13\\
57	3.61036892225008e-12\\
58	5.95779942047013e-13\\
59	5.33733236575714e-13\\
60	3.26897963762005e-13\\
61	5.31687566051212e-13\\
62	2.712575988893e-13\\
63	3.73714008566518e-13\\
64	1.71080850987282e-12\\
65	3.26865360366918e-14\\
66	4.50078485314146e-13\\
67	2.40244619279004e-14\\
68	5.2114317204297e-13\\
69	5.09896205731576e-12\\
70	5.30341274184355e-10\\
71	2.53590123698807e-12\\
72	8.36122179408212e-14\\
73	7.07644879989728e-10\\
74	1.71256470588935e-13\\
75	5.60892134887413e-14\\
76	2.34709734972889e-13\\
77	3.89760229965123e-13\\
78	2.78697742300831e-13\\
79	1.02426706829867e-13\\
80	3.35224882985759e-13\\
81	5.86158927955832e-12\\
82	1.16702094213222e-13\\
83	4.05744764064228e-14\\
84	5.71942420456249e-14\\
85	2.12175291327491e-14\\
86	1.63554821093861e-12\\
87	7.698942528711e-14\\
88	8.55240959151061e-14\\
89	7.1438782512577e-14\\
90	7.54503244018376e-13\\
91	3.6018826251802e-13\\
92	1.29655468308498e-12\\
93	4.43691279030281e-14\\
94	5.34653525545165e-14\\
95	1.78836760498874e-13\\
96	1.73501072383029e-11\\
97	8.37568421463267e-12\\
98	2.58201017324597e-13\\
99	1.3004096130707e-14\\
100	5.43201639289635e-16\\
};
\addlegendentry{\texttt{CUR-CR} (DEIM) $|$ (no oversamp.)}

\addplot [color=mycolor1, dotted, line width=1.0pt, mark=o, mark options={solid, mycolor1}]
  table[row sep=crcr]{%
1	0.612017026941478\\
2	0.307492955625974\\
3	0.172597870753999\\
4	0.0838431563268695\\
5	0.0429006609350698\\
6	0.0292544834401747\\
7	0.0184590490519216\\
8	0.00662496951488563\\
9	0.00336793203618044\\
10	0.00192288062933955\\
11	0.00081571582035379\\
12	0.000714251361103164\\
13	0.000318054996428819\\
14	0.00014666962056872\\
15	0.000105007901631428\\
16	3.96740189539395e-05\\
17	2.42946007270493e-05\\
18	9.13118645968891e-06\\
19	5.16476837135642e-06\\
20	4.32978482586316e-06\\
21	1.72310627350581e-06\\
22	8.45000558675833e-07\\
23	5.76361225785945e-07\\
24	2.20282121861064e-07\\
25	7.69230219737766e-08\\
26	5.16836279593637e-08\\
27	3.0334365142572e-08\\
28	1.36599040793446e-08\\
29	4.55574289512191e-09\\
30	3.61458613376336e-09\\
31	1.36744144562624e-09\\
32	8.92548167409696e-10\\
33	5.87202858044247e-10\\
34	2.05377126125133e-10\\
35	1.65300755409892e-10\\
36	8.1433410026216e-11\\
37	4.24058090865181e-11\\
38	3.57718482488519e-11\\
39	9.17746054498683e-12\\
40	4.69878958403529e-12\\
41	2.04291068205518e-12\\
42	1.30205446333521e-12\\
43	7.09424513290045e-13\\
44	2.7078922565359e-13\\
45	1.90449605792251e-13\\
46	7.51588342949524e-14\\
47	2.3476471807481e-14\\
48	1.40813603537642e-14\\
49	9.30149782280924e-15\\
50	6.89445330124789e-15\\
51	2.26026729997639e-14\\
52	1.57380598058645e-14\\
53	5.51563547743252e-14\\
54	4.29686544280736e-14\\
55	8.0509017630978e-14\\
56	1.34033073395701e-10\\
57	5.636404606361e-14\\
58	8.16765430109067e-14\\
59	9.11236452708004e-12\\
60	4.17986776501763e-13\\
61	3.52978721456077e-09\\
62	1.57563718495267e-12\\
63	2.4895316403619e-11\\
64	2.50080725693726e-12\\
65	2.70600531753787e-12\\
66	1.78830888770015e-12\\
67	4.85606453200322e-14\\
68	5.34286577796534e-14\\
69	2.96474397510717e-13\\
70	3.2828275973071e-13\\
71	2.06074307151522e-12\\
72	2.08643545553712e-13\\
73	3.10536965776619e-13\\
74	8.59449799391032e-14\\
75	1.47371720167717e-13\\
76	1.28031010264557e-13\\
77	5.32797553500253e-12\\
78	2.49931108159198e-12\\
79	1.65814493344619e-13\\
80	4.96944473007317e-14\\
81	5.24889353082987e-14\\
82	2.43988404021449e-13\\
83	2.19504929176118e-13\\
84	6.49885071563344e-12\\
85	2.85320094096777e-13\\
86	2.56241022223434e-13\\
87	9.74663452049781e-14\\
88	7.83351447218144e-14\\
89	9.83580228732949e-13\\
90	2.63249820746449e-14\\
91	1.42481346435462e-13\\
92	5.26584984099692e-14\\
93	1.96348196697036e-13\\
94	1.1906170818978e-12\\
95	2.20289383224309e-11\\
96	1.9379515791106e-14\\
97	2.77853104960562e-14\\
98	1.53458951370264e-14\\
99	3.61210252123422e-14\\
};
\addlegendentry{\texttt{CUR-CR} (MAXVOL) \cite{GO10} $|$ (no oversamp.)}

% \addplot [color=mycolor1, line width=1.2pt, mark=o, mark options={solid, mycolor1}]
%   table[row sep=crcr]{%
% 1	0.609457273253409\\
% 6	0.0294192981711891\\
% 11	0.000866536467998819\\
% 16	3.47190334856572e-05\\
% 21	1.38785151830582e-06\\
% 26	4.59326615019986e-08\\
% 31	1.71990049688393e-09\\
% 36	6.55987358678259e-11\\
% 41	1.97531033871732e-12\\
% 46	9.19934429133082e-14\\
% 51	4.24962122524889e-15\\
% 56	2.3827629418998e-15\\
% 61	3.12036286169722e-15\\
% 66	2.14930656842651e-15\\
% 71	2.28213866926914e-15\\
% 76	2.07909418047908e-15\\
% 81	1.68050976874271e-15\\
% 86	1.23336693881067e-15\\
% 91	0\\
% 96	0\\
% };
% \addlegendentry{CUR (CR-OS)$|$ $m = 10$}

% \addplot [color=mycolor3, line width=1.2pt, mark=o, mark options={solid, mycolor3}]
%   table[row sep=crcr]{%
% 1	0.601505186884026\\
% 6	0.0277149856664462\\
% 11	0.000846970826828632\\
% 16	3.29458363643862e-05\\
% 21	1.17579378308854e-06\\
% 26	3.95777684222712e-08\\
% 31	1.37570667211509e-09\\
% 36	5.43181002231016e-11\\
% 41	1.58472351420565e-12\\
% 46	6.32602101246281e-14\\
% 51	2.96097454682801e-15\\
% 56	1.24149056189841e-15\\
% 61	9.52032577178191e-16\\
% 66	9.10521593419642e-16\\
% 71	0\\
% 76	0\\
% 81	0\\
% 86	0\\
% 91	0\\
% 96	0\\
% };
% \addlegendentry{CUR (CR-OS)$|$ $m = 30$}

\addplot [color=mycolor3, line width=1.2pt,mark=o, mark options={solid, mycolor3}]
  table[row sep=crcr]{%
% 1	0.612017026941478\\
% 6	0.0347730455811607\\
% 11	0.000936163661808515\\
% 16	3.70729090709351e-05\\
% 21	1.26944155354616e-06\\
% 26	5.91535802003737e-08\\
% 31	1.91655846632844e-09\\
% 36	9.46045553100478e-11\\
% 41	2.15861250728194e-12\\
% 46	1.32785064257516e-13\\
% 51	5.36613008313402e-15\\
% 56	2.11420975744424e-15\\
% 61	3.911014807107e-15\\
% 66	3.65208287615038e-15\\
% 71	5.48671849353165e-15\\
% 76	2.23377422455477e-15\\
% 81	3.68409452992513e-15\\
% 86	2.09456947215038e-15\\
% 91	3.59433708928646e-15\\
% 96	4.04112081538649e-15\\
1	0.612017026941478\\
2	0.307492955625974\\
3	0.172597870753999\\
4	0.10315183594602\\
5	0.055238042882207\\
6	0.0347730455811607\\
7	0.0168155854116056\\
8	0.00663698807943946\\
9	0.00350182976668751\\
10	0.00209722904025338\\
11	0.000936163661808515\\
12	0.000588667909838891\\
13	0.000330926852474536\\
14	0.000157760331298599\\
15	0.000101192724614388\\
16	3.70729090709351e-05\\
17	1.99614132844045e-05\\
18	1.07540680489655e-05\\
19	5.6163752454094e-06\\
20	6.76389903492232e-06\\
21	1.26944155354616e-06\\
22	1.17387481785658e-06\\
23	4.92657320801356e-07\\
24	5.55922664002238e-07\\
25	1.1082205467812e-07\\
26	5.91535802003737e-08\\
27	4.78409378789888e-08\\
28	2.38135636524827e-08\\
29	7.01004049356125e-09\\
30	3.61183757287041e-09\\
31	1.91655846632844e-09\\
32	8.01103883142875e-10\\
33	6.18575892531085e-10\\
34	2.71139363155349e-10\\
35	1.73510618533137e-10\\
36	9.46045553100478e-11\\
37	4.07797991285971e-11\\
38	2.59834442240788e-11\\
39	1.25876144033457e-11\\
40	6.5135903666095e-12\\
41	2.15861250728194e-12\\
42	1.36621135650713e-12\\
43	1.08642908763004e-12\\
44	2.48036736980234e-13\\
45	2.185125245135e-13\\
46	1.32785064257516e-13\\
47	3.23007036534238e-14\\
48	2.49831013302146e-14\\
49	1.35739792977519e-14\\
50	5.63065994026826e-15\\
51	5.36613008313402e-15\\
52	3.34582029910389e-15\\
53	2.17059702407359e-15\\
54	3.00429891327904e-15\\
55	3.61493237842831e-15\\
56	2.11420975744424e-15\\
57	3.27348130959966e-15\\
58	4.88928801560728e-15\\
59	2.13700318364263e-15\\
60	3.51602375012891e-15\\
61	3.911014807107e-15\\
62	3.12905475916919e-15\\
63	2.3410328042944e-15\\
64	4.76566829684229e-15\\
65	2.84669538621056e-15\\
66	3.65208287615038e-15\\
67	2.83327513655906e-15\\
68	3.243808664898e-15\\
69	3.19979845373504e-15\\
70	3.91824421114526e-15\\
71	5.48671849353165e-15\\
72	3.34618563093815e-15\\
73	3.86188501492888e-15\\
74	2.91186833734804e-15\\
75	4.94667415684499e-15\\
76	2.23377422455477e-15\\
77	2.8135669816989e-15\\
78	4.8251585384013e-15\\
79	2.46308043938828e-15\\
80	3.21876927778426e-15\\
81	3.68409452992513e-15\\
82	4.51880714118766e-15\\
83	2.04359278179983e-15\\
84	4.49052442636433e-15\\
85	3.66764725621925e-15\\
86	2.09456947215038e-15\\
87	1.84162479583709e-15\\
88	2.31932605204523e-15\\
89	3.04264084093048e-15\\
90	1.82375618915461e-15\\
91	3.59433708928646e-15\\
92	2.15429234191871e-15\\
93	3.25267279961039e-15\\
94	4.57728435571205e-15\\
95	1.86085899921994e-15\\
96	4.04112081538649e-15\\
97	4.20380346824311e-15\\
98	3.11913863494676e-15\\
99	1.21383598841764e-15\\
100	5.43201639289635e-16\\
};
\addlegendentry{\texttt{CUR-CR-OS}$|$ $(m_r,m_c)$ adaptive}
\end{axis}

\end{tikzpicture}%

%% file: figures/Err_r_slow.tex
\definecolor{mycolor1}{rgb}{0.00000,0.44700,0.74100}%
\definecolor{mycolor2}{rgb}{0.85000,0.32500,0.09800}%
\definecolor{mycolor3}{rgb}{0.5,0.5,0.5}
\begin{tikzpicture}[scale=0.8]

\begin{axis}[%
width=2.4in,
height=2.2in,
at={(1.011in,0.642in)},
scale only axis,
xmin=0,
xmax=110,
xlabel style={font=\color{white!15!black}},
xlabel={$r$},
ymode=log,
ymin=1e-07,
ymax=100,
yminorticks=true,
ylabel style={font=\color{white!15!black}},
ylabel={Error ($\mathcal{E}$)},
axis background/.style={fill=white},
xmajorgrids,
ymajorgrids,
%yminorgrids,
% legend style={fill opacity=0.8,
% font=\scriptsize,
%   draw opacity=1,
%   text opacity=1,
%   at={(0.009,0.87)},
%   anchor=west}
legend style={fill opacity=0.8,draw opacity=1,
  text opacity=1,legend cell align=left, font=\scriptsize,align=left, draw=white!15!black}
]
\addplot [color=black, line width=1.2pt, mark=o, mark options={solid, black}]
  table[row sep=crcr]{%
% 1	0.125\\
% 6	0.00291545189504373\\
% 11	0.000578703703703705\\
% 16	0.000203541624262163\\
% 21	9.39143501126919e-05\\
% 26	5.0805263425289e-05\\
% 31	3.05175781250043e-05\\
% 36	1.97421672951351e-05\\
% 41	1.34974624770543e-05\\
% 46	9.63177715920624e-06\\
% 51	7.11197086936511e-06\\
% 56	5.39977212961448e-06\\
% 61	4.1958980900245e-06\\
% 66	3.32487706266881e-06\\
% 71	2.67918381343358e-06\\
% 76	2.19042216006086e-06\\
% 81	1.81367072444465e-06\\
% 86	1.51859596690575e-06\\
% 91	1.28421139146643e-06\\
% 96	1.09568268154392e-06\\
1	0.125\\
2	0.037037037037037\\
3	0.015625\\
4	0.00800000000000001\\
5	0.00462962962962963\\
6	0.00291545189504373\\
7	0.001953125\\
8	0.00137174211248286\\
9	0.00100000000000001\\
10	0.00075131480090159\\
11	0.000578703703703705\\
12	0.000455166135639524\\
13	0.000364431486880464\\
14	0.000296296296296298\\
15	0.000244140625\\
16	0.000203541624262163\\
17	0.000171467764060356\\
18	0.000145793847499641\\
19	0.000124999999999995\\
20	0.000107979699816442\\
21	9.39143501126902e-05\\
22	8.21895290540022e-05\\
23	7.23379629629701e-05\\
24	6.39999999999929e-05\\
25	5.68957669549374e-05\\
26	5.08052634252944e-05\\
27	4.55539358600594e-05\\
28	4.1002091106649e-05\\
29	3.7037037037037e-05\\
30	3.35671847202205e-05\\
31	3.05175781250023e-05\\
32	2.78264741074655e-05\\
33	2.54427030327714e-05\\
34	2.33236151603567e-05\\
35	2.14334705075398e-05\\
36	1.97421672951331e-05\\
37	1.82242309374525e-05\\
38	1.68580050236823e-05\\
39	1.56250000000029e-05\\
40	1.45093657956145e-05\\
41	1.34974624770521e-05\\
42	1.2577508898593e-05\\
43	1.17392937640889e-05\\
44	1.09739368998603e-05\\
45	1.02736911317448e-05\\
46	9.63177715920868e-06\\
47	9.04224537036644e-06\\
48	8.49985975230775e-06\\
49	7.99999999999841e-06\\
50	7.53857867637465e-06\\
51	7.11197086936688e-06\\
52	6.71695426425232e-06\\
53	6.35065792816454e-06\\
54	6.01051840720522e-06\\
55	5.69424198251149e-06\\
56	5.39977212961918e-06\\
57	5.12526138833261e-06\\
58	4.86904698142955e-06\\
59	4.62962962963342e-06\\
60	4.40565509888315e-06\\
61	4.19589809002701e-06\\
62	3.99924814134649e-06\\
63	3.81469726562269e-06\\
64	3.64132908510833e-06\\
65	3.47830926342361e-06\\
66	3.32487706266578e-06\\
67	3.1803378790944e-06\\
68	3.04405663163541e-06\\
69	2.91545189504581e-06\\
70	2.79399068483608e-06\\
71	2.67918381343948e-06\\
72	2.57058174834832e-06\\
73	2.46777091188057e-06\\
74	2.37037037037607e-06\\
75	2.27802886717332e-06\\
76	2.19042216006302e-06\\
77	2.10725062794766e-06\\
78	2.02823711713807e-06\\
79	1.95312500000619e-06\\
80	1.88167642316657e-06\\
81	1.81367072444287e-06\\
82	1.74890300059721e-06\\
83	1.68718280963133e-06\\
84	1.62833299408214e-06\\
85	1.57218861231788e-06\\
86	1.51859596690596e-06\\
87	1.46741172051584e-06\\
88	1.41850209017271e-06\\
89	1.37174211248596e-06\\
90	1.3270149727077e-06\\
91	1.28421139147173e-06\\
92	1.24322906371286e-06\\
93	1.20397214489954e-06\\
94	1.16635077999473e-06\\
95	1.13028067129172e-06\\
96	1.09568268153153e-06\\
97	1.06248246902918e-06\\
98	1.03061015213057e-06\\
99	1.00000000000296e-06\\
100	7.13000445905086e-15\\
};
\addlegendentry{\texttt{SVD}}

\addplot [color=mycolor2, line width=1.2pt, mark=o, mark options={solid, mycolor2}]
  table[row sep=crcr]{%
% 1	0.151824497546581\\
% 6	0.00657278835253934\\
% 11	0.00180859914841923\\
% 16	0.000555423865238335\\
% 21	0.000631305360333662\\
% 26	0.0264974104778594\\
% 31	0.000496321314204839\\
% 36	0.00234512668874502\\
% 41	0.0004498474327411\\
% 46	0.000215084390705511\\
% 51	0.000273860739927224\\
% 56	0.000188185844750728\\
% 61	0.000862807600267491\\
% 66	0.000221010230464426\\
% 71	0.00146826012462088\\
% 76	0.000259996123427099\\
% 81	0.000321148074809039\\
% 86	6.44753403345655e-05\\
% 91	0.00165677112751117\\
% 96	3.30289444029887e-05\\
1	0.151824497546581\\
2	0.0448962731242851\\
3	0.0222954125455349\\
4	0.013389754970706\\
5	0.00897173824581315\\
6	0.00657278835253934\\
7	0.0044010013174092\\
8	0.00229579109610078\\
9	0.00191979022881901\\
10	0.00203304402608858\\
11	0.00180859914841923\\
12	0.00120327911021725\\
13	0.00115438870411914\\
14	0.000835817116454422\\
15	0.000870391363669555\\
16	0.000555423865238335\\
17	0.000850169035798551\\
18	0.00108293549958555\\
19	0.00118319850669784\\
20	0.00122413088404352\\
21	0.000631305360333662\\
22	0.00116946540360112\\
23	0.000633663811223663\\
24	0.0100996545640765\\
25	0.00047411436728369\\
26	0.0264974104778594\\
27	0.00319092424415899\\
28	0.009471684068826\\
29	0.000637330214961093\\
30	0.00043898366690699\\
31	0.000496321314204839\\
32	0.000442155984959364\\
33	0.000496838569816199\\
34	0.000217429791511047\\
35	0.000287432526775576\\
36	0.00234512668874502\\
37	0.000632872986033364\\
38	0.000333543000949005\\
39	0.00027366911618855\\
40	0.00028622306932199\\
41	0.0004498474327411\\
42	0.000549350424980221\\
43	0.000616559000023944\\
44	0.000238183560300454\\
45	0.000159371145131107\\
46	0.000215084390705511\\
47	0.000246384916081667\\
48	0.000138010265738582\\
49	9.78723560712725e-05\\
50	0.000198515070975747\\
51	0.000273860739927224\\
52	0.000659584483091049\\
53	0.00130771325134394\\
54	0.000318094801548477\\
55	0.000354394534064828\\
56	0.000188185844750728\\
57	0.000179989976356583\\
58	0.00011039328457289\\
59	0.000259294399824494\\
60	0.000155025895709438\\
61	0.000862807600267491\\
62	0.000229917685984629\\
63	0.000191854481463341\\
64	0.000538516177597138\\
65	0.000459362518900841\\
66	0.000221010230464426\\
67	0.000279130380235772\\
68	0.000136593799783823\\
69	9.60655311117289e-05\\
70	0.000148231261544724\\
71	0.00146826012462088\\
72	0.000240800086339804\\
73	0.000374320515962041\\
74	0.000301788053557086\\
75	0.000274127703418759\\
76	0.000259996123427099\\
77	0.00396613760651091\\
78	0.00615941118818204\\
79	9.79720781525038e-05\\
80	0.000396685092402552\\
81	0.000321148074809039\\
82	8.07859514459273e-05\\
83	4.92019947881724e-05\\
84	5.30387626331064e-05\\
85	6.61174927715469e-05\\
86	6.44753403345655e-05\\
87	5.34121170067946e-05\\
88	6.9033442750291e-05\\
89	0.0477763857610969\\
90	0.00490159662503336\\
91	0.00165677112751117\\
92	0.00455147870056577\\
93	0.000259547059020914\\
94	0.000291301859077611\\
95	6.72810874164224e-05\\
96	3.30289444029887e-05\\
97	0.000124706913989429\\
98	3.41646045400636e-05\\
99	1.15230285254642e-05\\
100	6.40900860117206e-16\\
};
\addlegendentry{\texttt{CUR-CR} (DEIM) $|$ (no oversamp.)}

\addplot [color=mycolor1, dotted, line width=1.2pt, mark=o, mark options={solid, mycolor1}]
  table[row sep=crcr]{%
1	0.151824497546581\\
2	0.0448962731242851\\
3	0.0222954125455349\\
4	0.010858975236114\\
5	0.00659292657611463\\
6	0.00522719849038983\\
7	0.00467553434383177\\
8	0.00232512093319599\\
9	0.00168158564525236\\
10	0.00149057946110797\\
11	0.00141407020942598\\
12	0.00358668266493799\\
13	0.00201261774091629\\
14	0.00132316083433487\\
15	0.000762608438596459\\
16	0.00212415086211743\\
17	0.00224464152779355\\
18	0.000617340809766764\\
19	0.00068421738231352\\
20	0.000539979014124943\\
21	0.00045950704198698\\
22	0.000588679440056898\\
23	0.000537125329424212\\
24	0.00043755130409246\\
25	0.000217762218326199\\
26	0.00111566969162074\\
27	0.000638009349724891\\
28	0.000675914182900264\\
29	0.00047591132842664\\
30	0.000509436628581331\\
31	0.000604898040847869\\
32	0.000470875879564633\\
33	0.00190941493913423\\
34	0.000544718584594897\\
35	0.000465013438709032\\
36	0.000946493050759969\\
37	0.000829058146739204\\
38	0.000237032847048691\\
39	0.000170973935896922\\
40	0.000599387723690186\\
41	0.00033227480869118\\
42	0.000282900259135839\\
43	0.000431471841694446\\
44	0.0381093773441789\\
45	0.000323373909495608\\
46	0.000241670135994014\\
47	0.000360048158716408\\
48	0.000419717047977073\\
49	0.000546901563696224\\
50	0.000312399082892191\\
51	0.00237556593386196\\
52	0.00241868430177243\\
53	0.000364871743372032\\
54	0.000223639628743088\\
55	0.000795844713731672\\
56	0.000288999127116426\\
57	0.000763178130400215\\
58	0.00103733455672577\\
59	0.000696771384054931\\
60	0.000370962967522982\\
61	0.00599305221700581\\
62	0.000186395238996026\\
63	0.000327489928621637\\
64	0.00287702519517117\\
65	0.00113718199262757\\
66	0.000245573074147714\\
67	0.0127609158838066\\
68	0.000117375786478522\\
69	0.00194818827959874\\
70	0.000639324870455796\\
71	0.00079293206457987\\
72	0.00101402828833314\\
73	0.00201193301229296\\
74	0.000105719266469712\\
75	0.000148563219698154\\
76	0.000452651093689264\\
77	0.000727085982855632\\
78	0.000180954626029571\\
79	0.000470281284306141\\
80	0.00064901297916341\\
81	0.000586341390735243\\
82	0.00135414613821191\\
83	0.000637237900680933\\
84	0.000243202378864293\\
85	0.00264914550948931\\
86	0.000118011433385224\\
87	9.21626194860972e-05\\
88	0.000298877168016449\\
89	0.000946518485555822\\
90	0.000253493486923857\\
91	0.000114480631758651\\
92	0.000219831993414775\\
93	0.000176448203533017\\
94	9.16484213920818e-05\\
95	0.000102883816729997\\
96	3.36462862074815e-05\\
97	3.20619737677028e-05\\
98	1.42424544109718e-05\\
99	1.5329569945385e-05\\
};
\addlegendentry{\texttt{CUR-CR} (MAXVOL) \cite{GO10} $|$ (no oversamp.)}

\addplot [color=black, dotted, line width=1.2pt, mark options={solid, black}]
  table[row sep=crcr]{%
1	0.151483314856081\\
2	0.0439569266981717\\
3	0.021541137611937\\
4	0.0113615614481329\\
5	0.0063749539856963\\
6	0.00464055701775634\\
7	0.00311819659473502\\
8	0.00205895342848202\\
9	0.00153753449461765\\
10	0.00121255478957999\\
11	0.00088477055941628\\
12	0.000766109061304936\\
13	0.000601219637469348\\
14	0.000565347779946577\\
15	0.000432459161051577\\
16	0.000412528253726384\\
17	0.000316378810343578\\
18	0.000279971124461793\\
19	0.00026415452166394\\
20	0.000223293334106451\\
21	0.000205592228827659\\
22	0.000184633831604379\\
23	0.000165168284375729\\
24	0.000159736887296482\\
25	0.00012495837908351\\
26	0.000102859078311251\\
27	9.56541339471151e-05\\
28	8.78700509354766e-05\\
29	7.78657311473359e-05\\
30	7.15532406483345e-05\\
31	6.87697696854226e-05\\
32	6.54868965430471e-05\\
33	5.98949452730665e-05\\
34	5.52481276367253e-05\\
35	5.15425429742863e-05\\
36	4.50210324800904e-05\\
37	4.05575925194561e-05\\
38	3.90574769236754e-05\\
39	3.65085865743081e-05\\
40	3.48497124408701e-05\\
41	3.31456131630465e-05\\
42	3.28800649307134e-05\\
43	2.88532347049001e-05\\
44	2.76412973049099e-05\\
45	2.70278477767818e-05\\
46	2.48287179442624e-05\\
47	2.11144895212137e-05\\
48	2.05230867899088e-05\\
49	1.92615526212504e-05\\
50	1.91819556225216e-05\\
51	1.90612194713464e-05\\
52	1.89682648672917e-05\\
53	1.73425404188434e-05\\
54	1.59670382082986e-05\\
55	1.51866238138835e-05\\
56	1.42471786760995e-05\\
57	1.32574144371401e-05\\
58	1.31023838177521e-05\\
59	1.20553117217995e-05\\
60	1.14769788383935e-05\\
61	1.11944598579782e-05\\
62	1.06883440740646e-05\\
63	1.04739586734429e-05\\
64	1.02776172340591e-05\\
65	9.92251271526022e-06\\
66	9.5515505694601e-06\\
67	8.80312721373571e-06\\
68	8.67514017283176e-06\\
69	8.42188885646728e-06\\
70	8.02620025826618e-06\\
71	8.02410608154505e-06\\
72	8.02178308506479e-06\\
73	7.56576042585443e-06\\
74	6.57063727959302e-06\\
75	6.3705292688246e-06\\
76	6.331470458733e-06\\
77	6.22397409876052e-06\\
78	5.76265962723003e-06\\
79	5.70090780659611e-06\\
80	5.24044247992823e-06\\
81	5.15933533941712e-06\\
82	5.11228027661805e-06\\
83	4.7945131529598e-06\\
84	4.63912343723691e-06\\
85	4.52071640759011e-06\\
86	4.46183250585559e-06\\
87	3.9787999241184e-06\\
88	3.83614714748835e-06\\
89	3.62008408276913e-06\\
90	3.58127311922439e-06\\
91	3.55915938362859e-06\\
92	3.43615622941569e-06\\
93	3.34281965199913e-06\\
94	3.25020004261416e-06\\
95	2.81576687044002e-06\\
96	2.80016733707749e-06\\
97	2.69351020292039e-06\\
98	2.43790034018667e-06\\
99	2.43145956794867e-06\\
100	2.14200198935819e-11\\
};
\addlegendentry{\texttt{CUR-opt} \cite{SE16} $|$ (no oversamp.)}

% \addplot [color=mycolor1, line width=1.2pt, mark=o, mark options={solid, mycolor1}]
%   table[row sep=crcr]{%
% 1	0.158611784690574\\
% 6	0.00706025479089822\\
% 11	0.00184752882769932\\
% 16	0.000930233803666817\\
% 21	0.000580911170157419\\
% 26	0.000345635785656025\\
% 31	0.000256701860435389\\
% 36	0.000201788872021563\\
% 41	0.000142182569825807\\
% 46	0.00011611818134255\\
% 51	9.77755579110592e-05\\
% 56	7.80821502454374e-05\\
% 61	7.09168806960619e-05\\
% 66	5.53779067320353e-05\\
% 71	5.22201373624814e-05\\
% 76	4.44770494272362e-05\\
% 81	3.13510063623425e-05\\
% 86	2.08564224263083e-05\\
% 91	0\\
% 96	0\\
% };
% \addlegendentry{CUR (CR-OS)$|$ $m = 10$}

% \addplot [color=mycolor3, line width=1.2pt, mark=o, mark options={solid, mycolor3}]
%   table[row sep=crcr]{%
% 1	0.158060304655326\\
% 6	0.00663260081818688\\
% 11	0.00176569267972969\\
% 16	0.000854274049235245\\
% 21	0.000487630141651666\\
% 26	0.00029461052488542\\
% 31	0.000204788963031614\\
% 36	0.000149942905273339\\
% 41	0.000114815571327512\\
% 46	8.65832698599681e-05\\
% 51	6.90795137332028e-05\\
% 56	5.51548769450634e-05\\
% 61	4.3567264973104e-05\\
% 66	3.54231720562117e-05\\
% 71	0\\
% 76	0\\
% 81	0\\
% 86	0\\
% 91	0\\
% 96	0\\
% };
% \addlegendentry{CUR (CR-OS)$|$ $m = 30$}

\addplot [color=mycolor3, line width=1.2pt, mark=o, mark options={solid, mycolor3}]
  table[row sep=crcr]{%
% 1	0.151824497546581\\
% 6	0.00657278835253934\\
% 11	0.00180859914841923\\
% 16	0.000555423865238335\\
% 21	0.000486796183651411\\
% 26	0.00022152556693511\\
% 31	0.000350735215154105\\
% 36	0.000163708892502235\\
% 41	0.000138462440074253\\
% 46	8.64322623130167e-05\\
% 51	7.02075065989672e-05\\
% 56	7.4989266893873e-05\\
% 61	4.77825198730904e-05\\
% 66	2.76387472590046e-05\\
% 71	3.5945844886655e-05\\
% 76	5.4684453852651e-05\\
% 81	2.28415595322877e-05\\
% 86	2.92591109237976e-05\\
% 91	2.5969219952053e-05\\
% 96	1.83115742510543e-05\\
1	0.151824497546581\\
2	0.0448962731242851\\
3	0.0222954125455349\\
4	0.013389754970706\\
5	0.00897173824581315\\
6	0.00657278835253934\\
7	0.0044010013174092\\
8	0.00229579109610078\\
9	0.00191979022881901\\
10	0.00203304402608858\\
11	0.00180859914841923\\
12	0.00120327911021725\\
13	0.00115438870411914\\
14	0.000835817116454422\\
15	0.000870391363669555\\
16	0.000555423865238335\\
17	0.000721905478196333\\
18	0.000806461589298723\\
19	0.000799753080374223\\
20	0.000878118232312764\\
21	0.000486796183651411\\
22	0.000447492254372737\\
23	0.000298651104197731\\
24	0.000247264024433311\\
25	0.000267253609367835\\
26	0.00022152556693511\\
27	0.00020325513041636\\
28	0.000231490925225473\\
29	0.000190380832218622\\
30	0.0001550468885872\\
31	0.000350735215154105\\
32	0.000225634717886191\\
33	0.000221490392338811\\
34	0.000174711382652395\\
35	0.000129528342039578\\
36	0.000163708892502235\\
37	0.000135743772812866\\
38	0.0001500354407012\\
39	0.000115179859403323\\
40	0.000144929629082574\\
41	0.000138462440074253\\
42	0.000116666686060056\\
43	7.8462817181727e-05\\
44	0.000100044505632964\\
45	8.20557278654419e-05\\
46	8.64322623130167e-05\\
47	9.10133602495146e-05\\
48	7.02259482024012e-05\\
49	7.55889744531245e-05\\
50	6.39322332108398e-05\\
51	7.02075065989672e-05\\
52	7.35566574788427e-05\\
53	8.93182498416187e-05\\
54	5.9201650553333e-05\\
55	4.53591306297766e-05\\
56	7.4989266893873e-05\\
57	7.22143791527668e-05\\
58	6.62868921533771e-05\\
59	4.61288078533633e-05\\
60	6.9443625318493e-05\\
61	4.77825198730904e-05\\
62	3.73871434803081e-05\\
63	4.60119852372079e-05\\
64	4.0007327739845e-05\\
65	3.70488276600385e-05\\
66	2.76387472590046e-05\\
67	4.34953275569247e-05\\
68	5.93072959865962e-05\\
69	3.25722090595159e-05\\
70	3.26970378062589e-05\\
71	3.5945844886655e-05\\
72	2.89235117713108e-05\\
73	3.51981206336857e-05\\
74	4.10991592774909e-05\\
75	2.85401327853765e-05\\
76	5.4684453852651e-05\\
77	3.53414268770641e-05\\
78	3.76410602037046e-05\\
79	3.62952370616097e-05\\
80	3.22037851308122e-05\\
81	2.28415595322877e-05\\
82	2.64716631080322e-05\\
83	2.72942756910464e-05\\
84	2.49585815797731e-05\\
85	2.37795000687065e-05\\
86	2.92591109237976e-05\\
87	2.59722643552914e-05\\
88	3.00107722928067e-05\\
89	2.64735061521058e-05\\
90	2.96026193804e-05\\
91	2.5969219952053e-05\\
92	1.6895539269035e-05\\
93	1.40845667061451e-05\\
94	1.6450739223484e-05\\
95	1.62304508641586e-05\\
96	1.83115742510543e-05\\
97	9.94470738822751e-06\\
98	1.1220571678291e-05\\
99	1.15230285254642e-05\\
100	6.40900860117206e-16\\
};
\addlegendentry{\texttt{CUR-CR-OS} $|$ $(m_r,m_c)$ adaptive}

\end{axis}
\end{tikzpicture}%

%% file: figures/Amp_I_J_r_fast.tex
\definecolor{mycolor1}{rgb}{0.00000,0.44700,0.74100}%
\definecolor{mycolor2}{rgb}{0.85000,0.32500,0.09800}%
\definecolor{mycolor3}{rgb}{0.5,0.5,0.5}
\begin{tikzpicture}[scale=0.8]

\begin{axis}[%
width=2.2in,
height=2.0in,
at={(1.011in,0.642in)},
scale only axis,
xmin=0,
xmax=110,
xlabel style={font=\color{white!15!black}},
xlabel={$r$},
ymode=log,
ymin=1,
ymax=100000,
yminorticks=true,
ylabel style={font=\color{white!15!black}},
ylabel={Norm of Submatrices},
axis background/.style={fill=white},
xmajorgrids,
ymajorgrids,
%yminorgrids,
legend style={fill opacity=0.8,legend cell align=left,
font=\scriptsize,
  draw opacity=1,
  text opacity=1,
  at={(0.009,0.84)},
  anchor=west}
]

\addplot [color=mycolor2, dotted, line width=1.2pt, mark=+, mark options={solid, mycolor2}]
  table[row sep=crcr]{%
1	3.69276285917803\\
2	3.53833890253118\\
3	3.54889762765579\\
4	4.89254785677607\\
5	5.5677894075573\\
6	4.74824134873114\\
7	4.93225279761555\\
8	5.35891709588462\\
9	5.42548051082557\\
10	5.37565947954572\\
11	6.46992895118493\\
12	5.22896723367493\\
13	8.1258972464157\\
14	8.07570892147291\\
15	6.6153313180071\\
16	7.22322261904848\\
17	8.95868807612574\\
18	9.42748414324589\\
19	9.48582452756929\\
20	8.72803528438014\\
21	7.42176636388717\\
22	7.73209763341806\\
23	8.51491565487108\\
24	13.7549827449543\\
25	7.72760232737181\\
26	7.52368365193236\\
27	8.03685225228425\\
28	11.9033645583992\\
29	8.85784822124667\\
30	8.40598908916947\\
31	8.95864871002309\\
32	8.9321934832376\\
33	10.0446348363922\\
34	10.5267229726354\\
35	9.90950273585176\\
36	14.1961360556278\\
37	10.9998958422565\\
38	10.2644320760122\\
39	9.54733980195312\\
40	11.2986111419479\\
41	11.1592136120381\\
42	9.93653135126903\\
43	17.1993094902583\\
44	10.6411993721533\\
45	11.0702490937089\\
46	14.3263972910851\\
47	11.3539283969048\\
48	13.4826169925815\\
49	12.2139327273535\\
50	13.3020931728584\\
51	10.8528851254223\\
52	14.9388379342136\\
53	10.6820104520475\\
54	20.2233817658972\\
55	38.3820748805214\\
56	81.7653984112616\\
57	431.819794452004\\
58	20.0534648130049\\
59	25.1792480963844\\
60	22.7951458910906\\
61	124.897828116576\\
62	37.7761669280106\\
63	37.0117401122161\\
64	41.0804059296927\\
65	22.5856708302703\\
66	32.2778411579381\\
67	27.9678512013402\\
68	123.913149307159\\
69	153.103973409127\\
70	46.0362615825138\\
71	381.196502598822\\
72	34.4256171758556\\
73	6006.89508222115\\
74	42.3895914167328\\
75	29.9495323422234\\
76	22.793326008744\\
77	43.5411560844968\\
78	31.2072899620826\\
79	36.0917012581723\\
80	92.2979582853065\\
81	500.472548554441\\
82	28.8738448436962\\
83	23.381717234434\\
84	51.9409866757854\\
85	30.5281174365899\\
86	242.499747562909\\
87	47.0654697647335\\
88	55.1969219791658\\
89	58.9436593456484\\
90	200.755182857572\\
91	24.1363068147537\\
92	230.067117362982\\
93	20.6773747976729\\
94	47.4522982665814\\
95	90.040303110586\\
96	709.843113709358\\
97	593.712036316249\\
98	62.991382075409\\
99	16.8987583650246\\
100	1\\
};
\addlegendentry{$\ol{\eta}_p$ - \texttt{CUR-CR} $|$ (no oversamp.)}

\addplot [color=mycolor1, dotted, line width=1.2pt, mark=+, mark options={solid, mycolor1}]
  table[row sep=crcr]{%
1	3.93824962202478\\
2	3.43293627900168\\
3	4.00120166011026\\
4	4.16562021179195\\
5	4.96780272664394\\
6	5.05817443488834\\
7	4.55167593224573\\
8	5.7113700946033\\
9	5.84142651683713\\
10	6.85149844347051\\
11	6.70280319828595\\
12	7.06931604062937\\
13	6.38332609927467\\
14	6.47924823836974\\
15	7.41124567726841\\
16	6.84283493305121\\
17	7.71499754065457\\
18	8.67887161231186\\
19	8.69277506702064\\
20	9.49033280245568\\
21	8.55929067572811\\
22	11.0976121381551\\
23	8.4699802299633\\
24	7.82830497758762\\
25	8.88941874277892\\
26	13.4601361869387\\
27	12.4668507368485\\
28	11.6919923763648\\
29	11.1767463308672\\
30	9.92466456500508\\
31	11.0431078685123\\
32	12.2005265862072\\
33	12.7464133343894\\
34	10.5718941152652\\
35	9.56515479653703\\
36	9.06912141815402\\
37	9.01390655405411\\
38	9.20156140377768\\
39	9.59686242015802\\
40	9.2859511505521\\
41	9.00275876124246\\
42	10.0494268791739\\
43	20.4573715504375\\
44	8.65938746301282\\
45	8.95970159454677\\
46	9.77122438960001\\
47	8.3027749479553\\
48	8.60844139239117\\
49	9.56010043567461\\
50	8.99900509730216\\
51	56.9981226867705\\
52	54.6346024226911\\
53	66.9474270630941\\
54	26.9010808436458\\
55	34.7615655835898\\
56	31.6887986411409\\
57	71.3866469743243\\
58	99.3688970529055\\
59	106.436790977691\\
60	73.5270933592442\\
61	117.71703058261\\
62	106.627347122476\\
63	112.646441201585\\
64	267.864106349447\\
65	25.0694408759163\\
66	99.3564319474633\\
67	33.3134222185612\\
68	25.4721041379239\\
69	364.354136184511\\
70	4047.2718828963\\
71	93.8366931491143\\
72	56.8500151474125\\
73	80.0360392952182\\
74	69.5583969164176\\
75	45.51932872627\\
76	79.6287558080528\\
77	106.263352972824\\
78	92.3842203330598\\
79	71.4388488447706\\
80	46.3499731952413\\
81	75.172396404123\\
82	76.0599834302751\\
83	38.7451813261172\\
84	34.1523492930666\\
85	18.7263463428295\\
86	44.659098877532\\
87	38.0597382291642\\
88	19.2437288900854\\
89	46.6496175907428\\
90	87.8670997140022\\
91	110.718717241532\\
92	80.2331708449372\\
93	38.6251157304354\\
94	36.4592944412691\\
95	60.6213376062026\\
96	49.1535096748234\\
97	37.5592759202941\\
98	9.24188775080396\\
99	5.77331103361365\\
100	1\\
};
\addlegendentry{$\ol{\eta}_s$ - \texttt{CUR-CR} $|$ (no oversamp.)}

\addplot [color=mycolor2,  line width=1.2pt, mark=o, mark options={solid, mycolor2}]
  table[row sep=crcr]{%
1	3.69276285917803\\
2	3.53833890253118\\
3	3.54889762765579\\
4	4.89254785677607\\
5	5.5677894075573\\
6	4.74824134873114\\
7	4.93225279761555\\
8	5.35891709588462\\
9	5.42548051082557\\
10	5.37565947954572\\
11	6.46992895118493\\
12	5.22896723367493\\
13	8.1258972464157\\
14	8.07570892147291\\
15	6.6153313180071\\
16	7.22322261904848\\
17	8.95868807612574\\
18	9.42748414324589\\
19	9.48582452756929\\
20	8.72803528438014\\
21	7.42176636388717\\
22	7.73209763341806\\
23	8.51491565487108\\
24	7.23406578318496\\
25	7.72760232737181\\
26	7.52368365193236\\
27	8.03685225228425\\
28	6.92216959304894\\
29	8.85784822124667\\
30	8.40598908916947\\
31	8.95864871002309\\
32	8.9321934832376\\
33	8.14156272173808\\
34	7.82876528540931\\
35	9.90950273585176\\
36	9.08406491349591\\
37	9.33883547887344\\
38	9.40244702380841\\
39	9.54733980195312\\
40	6.91397807944236\\
41	6.86987650323704\\
42	9.93653135126903\\
43	7.61288383068216\\
44	9.19394798073354\\
45	9.83087854727573\\
46	6.99800325198589\\
47	9.45341838478755\\
48	9.81422965221733\\
49	9.53376428549881\\
50	9.26592654588578\\
51	9.66296866311127\\
52	9.83056961024624\\
53	8.27304231775024\\
54	5.35157346490417\\
55	9.36639222085035\\
56	9.21844598878568\\
57	6.99727419312147\\
58	8.51159387240658\\
59	8.88153564346592\\
60	6.46828967409314\\
61	9.66102564872018\\
62	9.52090722017579\\
63	8.57377603935407\\
64	8.80662984628305\\
65	9.54362258261072\\
66	9.94562125238275\\
67	8.47760955275333\\
68	7.94590690526179\\
69	9.72439350017803\\
70	9.41204838503556\\
71	7.85528661292991\\
72	8.93650096105371\\
73	9.95086572331221\\
74	7.98969935068942\\
75	7.56326979864834\\
76	7.116214640521\\
77	9.79871626686473\\
78	9.62237377110663\\
79	8.37297964656143\\
80	7.28350703463973\\
81	9.13351680591494\\
82	8.73325011299853\\
83	8.32138066566465\\
84	8.95752378098393\\
85	9.45581522617159\\
86	3.62308377878904\\
87	8.13137979420162\\
88	6.69267648639974\\
89	6.29347187728055\\
90	5.5566380854894\\
91	7.62719984956429\\
92	7.81929934983657\\
93	8.85873390033983\\
94	9.71367667472451\\
95	1\\
96	2.89271246557596\\
97	1\\
98	4.46093424913814\\
99	1\\
100	1\\
};
\addlegendentry{$\ol{\eta}_p$ - \texttt{CUR-CR-OS}}

\addplot [color=mycolor1, line width=1.2pt, mark=o, mark options={solid, mycolor1}]
  table[row sep=crcr]{%
1	3.93824962202478\\
2	3.43293627900168\\
3	4.00120166011026\\
4	4.16562021179195\\
5	4.96780272664394\\
6	5.05817443488834\\
7	4.55167593224573\\
8	5.7113700946033\\
9	5.84142651683713\\
10	6.85149844347051\\
11	6.70280319828595\\
12	7.06931604062937\\
13	6.38332609927467\\
14	6.47924823836974\\
15	7.41124567726841\\
16	6.84283493305121\\
17	7.71499754065457\\
18	8.67887161231186\\
19	8.69277506702064\\
20	9.49033280245568\\
21	8.55929067572811\\
22	6.91646607845474\\
23	8.4699802299633\\
24	7.82830497758762\\
25	8.88941874277892\\
26	6.23868480592144\\
27	6.36606074448294\\
28	7.37627038535451\\
29	6.51590789653155\\
30	9.92466456500508\\
31	7.13883132886454\\
32	7.74642403667553\\
33	8.193144305457\\
34	8.30187169894542\\
35	9.56515479653703\\
36	9.06912141815402\\
37	9.01390655405411\\
38	9.20156140377768\\
39	9.59686242015802\\
40	9.2859511505521\\
41	9.00275876124246\\
42	6.96990799756096\\
43	8.74374478533552\\
44	8.65938746301282\\
45	8.95970159454677\\
46	9.77122438960001\\
47	8.3027749479553\\
48	8.60844139239117\\
49	9.56010043567461\\
50	8.99900509730216\\
51	9.2941740465908\\
52	8.1850435272066\\
53	5.9476201091964\\
54	8.68055369157493\\
55	9.19571570013975\\
56	6.13970338333632\\
57	9.03313553181083\\
58	9.67345438875446\\
59	6.76124214450997\\
60	8.05795782793695\\
61	8.76472510265867\\
62	7.82881112048427\\
63	6.62665253121894\\
64	9.28577847392152\\
65	9.52502606952686\\
66	8.85763663239168\\
67	9.76551329539195\\
68	9.5648828871579\\
69	5.68156430386341\\
70	9.99229022924109\\
71	9.38396164087602\\
72	9.72311649640518\\
73	8.49764456137816\\
74	9.86724176583895\\
75	8.27309675962873\\
76	7.34913162155015\\
77	7.19543636807384\\
78	8.35765390544566\\
79	9.46874445800739\\
80	9.95776574887451\\
81	7.81834876371496\\
82	8.8677021584863\\
83	7.80592847630227\\
84	9.00031493550549\\
85	9.4881528049428\\
86	7.97828938684104\\
87	5.6203004708624\\
88	9.27445323291218\\
89	8.89229633486839\\
90	6.98962104320029\\
91	7.40642617521955\\
92	7.1938749489333\\
93	8.63515754556799\\
94	9.90508053908695\\
95	6.08325394563719\\
96	8.92469627033582\\
97	8.31229359770124\\
98	9.24188775080396\\
99	5.77331103361365\\
100	1\\
};
\addlegendentry{$\ol{\eta}_s$ - \texttt{CUR-CR-OS}}
\end{axis}

\end{tikzpicture}%

%% file: figures/Amp_I_J_r_slow.tex
\definecolor{mycolor1}{rgb}{0.00000,0.44700,0.74100}%
\definecolor{mycolor2}{rgb}{0.85000,0.32500,0.09800}%
\definecolor{mycolor3}{rgb}{0.5,0.5,0.5}
\begin{tikzpicture}[scale=0.8]

\begin{axis}[%
width=2.2in,
height=2.0in,
at={(1.011in,0.642in)},
scale only axis,
xmin=0,
xmax=110,
xlabel style={font=\color{white!15!black}},
xlabel={$r$},
ymode=log,
ymin=1,
ymax=100000,
yminorticks=true,
ylabel style={font=\color{white!15!black}},
ylabel={Norm of Submatrices},
axis background/.style={fill=white},
xmajorgrids,
ymajorgrids,
%yminorgrids,
legend style={fill opacity=0.8, legend cell align=left,
font=\scriptsize,
  draw opacity=1,
  text opacity=1,
  at={(0.009,0.84)},
  anchor=west}
]

\addplot [color=mycolor2, dotted, line width=1.2pt, mark=+, mark options={solid, mycolor2}]
  table[row sep=crcr]{%
1	3.53456577962587\\
2	3.53272983487208\\
3	3.57903133795742\\
4	5.10630485091203\\
5	5.96761395574305\\
6	5.16352675029076\\
7	5.18561553488803\\
8	5.56080694687045\\
9	5.68605530889895\\
10	5.93546986490507\\
11	6.99061346396084\\
12	6.36574646391378\\
13	9.11022208265634\\
14	7.67213807461259\\
15	8.33946381370814\\
16	7.15209587507093\\
17	14.9518513993958\\
18	16.1384877429022\\
19	16.3781739544687\\
20	16.9915587157991\\
21	11.8134814814492\\
22	26.2854518639795\\
23	15.5679997540321\\
24	264.026989966732\\
25	12.4979531488088\\
26	800.311988393021\\
27	98.4328957679832\\
28	374.324191984357\\
29	25.4432914404503\\
30	17.1609694143225\\
31	25.0817783296692\\
32	19.2525634690911\\
33	21.9733072380675\\
34	13.06792469801\\
35	18.2268698245152\\
36	156.402740320244\\
37	39.3708800008736\\
38	21.296221519549\\
39	23.0888106142509\\
40	26.751116319159\\
41	41.8714820939611\\
42	48.3077244973391\\
43	66.7339483884878\\
44	25.3919404863174\\
45	15.9959133356408\\
46	24.7084935160585\\
47	27.4872117770039\\
48	17.6784505893205\\
49	14.1756016013744\\
50	31.1768343606316\\
51	44.0394588732413\\
52	117.704315765537\\
53	208.233394063312\\
54	47.5324082954206\\
55	53.0564976305507\\
56	28.2990353486796\\
57	27.9267003745456\\
58	24.112565881343\\
59	56.1947151956865\\
60	33.5099273844937\\
61	167.473127453882\\
62	59.6426839225109\\
63	38.0735115513019\\
64	133.151498305046\\
65	133.399839643503\\
66	54.6574118372449\\
67	67.9184596634221\\
68	34.4748374640653\\
69	25.2206319908344\\
70	38.6100899329658\\
71	351.257181034587\\
72	62.8785205262399\\
73	98.2617781443188\\
74	80.8545885450157\\
75	79.1180430735357\\
76	75.7078425801821\\
77	1350.0588577577\\
78	1940.38096937735\\
79	32.9608117476134\\
80	141.313724242329\\
81	100.550373428486\\
82	31.1738850376193\\
83	18.7000693676819\\
84	18.4004782455552\\
85	23.8698343584342\\
86	24.5890823446237\\
87	19.4676123858838\\
88	27.2385414596824\\
89	15336.9088115307\\
90	1574.24199143875\\
91	522.240396402209\\
92	1456.96947212862\\
93	88.3850254170754\\
94	95.4092959708452\\
95	27.5084896793668\\
96	14.8238182325177\\
97	58.6096053829792\\
98	15.2154599543324\\
99	5.20041749538385\\
100	1\\
};
\addlegendentry{$\ol{\eta}_p$ - \texttt{CUR-CR} $|$ (no oversamp.)}

\addplot [color=mycolor1, dotted, line width=1.2pt, mark=+, mark options={solid, mycolor1}]
  table[row sep=crcr]{%
1	3.98545416460458\\
2	3.38284287782404\\
3	3.97939912615686\\
4	4.19745307577015\\
5	5.57408130510737\\
6	5.16772210603135\\
7	5.02617031593698\\
8	5.79013655842614\\
9	5.5080431922631\\
10	6.57932459462923\\
11	7.66377140673162\\
12	6.65583732673056\\
13	7.90577389209253\\
14	6.42116798701518\\
15	9.27607370190683\\
16	8.0854452213\\
17	10.8149882992411\\
18	14.1603322717112\\
19	12.2678898512247\\
20	14.0833923017726\\
21	12.7581455331236\\
22	32.466890757446\\
23	16.4075923322566\\
24	186.03927961491\\
25	19.5842524807958\\
26	1046.05833055598\\
27	130.740055177494\\
28	375.236830170018\\
29	30.9369321523832\\
30	18.9760300447252\\
31	19.810953223125\\
32	19.7966635076136\\
33	25.2820147248946\\
34	13.5499936952959\\
35	18.3144910823426\\
36	162.876936637581\\
37	44.7456328679548\\
38	24.9986781303966\\
39	19.8930474129778\\
40	22.1823763466264\\
41	37.5701303985425\\
42	48.5074675371987\\
43	55.870810782782\\
44	18.1639057163274\\
45	15.6071833058249\\
46	26.4116301822446\\
47	31.5068961220475\\
48	19.0887280623321\\
49	15.0008566372388\\
50	28.5468208803667\\
51	39.9537708780279\\
52	90.3479872571629\\
53	216.328431446841\\
54	46.5251415064969\\
55	62.4707802921023\\
56	29.8169890059274\\
57	28.8187240531973\\
58	19.5765220557605\\
59	49.7667485926241\\
60	28.3503102844678\\
61	185.847666976715\\
62	44.2892677784801\\
63	35.8147389995636\\
64	104.361564801464\\
65	92.7149311081128\\
66	53.7782600754693\\
67	61.777132368695\\
68	27.2387876529052\\
69	22.4344764046151\\
70	39.9690683128324\\
71	388.519291367853\\
72	64.1985165533299\\
73	98.8291902820464\\
74	80.500400075454\\
75	89.4442613573424\\
76	87.0391413834679\\
77	928.036106706879\\
78	1537.74168719815\\
79	27.2684125753287\\
80	137.474750672154\\
81	111.289845305843\\
82	28.2204995460629\\
83	16.6032534266463\\
84	18.5439569259374\\
85	21.2653813770735\\
86	23.3653148004587\\
87	17.4232650605386\\
88	27.9764781837371\\
89	17294.8786009278\\
90	1784.25574359067\\
91	593.796024825618\\
92	1705.45008559048\\
93	90.0664436784451\\
94	101.436789600783\\
95	26.9057770558397\\
96	14.1189515627281\\
97	57.6348069729745\\
98	17.2968068355902\\
99	5.07390209734838\\
100	1\\
};
\addlegendentry{$\ol{\eta}_s$ - \texttt{CUR-CR} $|$ (no oversamp.)}

\addplot [color=mycolor2, line width=1.2pt, mark=o, mark options={solid, mycolor2}]
  table[row sep=crcr]{%
1	3.53456577962587\\
2	3.53272983487208\\
3	3.57903133795742\\
4	5.10630485091203\\
5	5.96761395574305\\
6	5.16352675029076\\
7	5.18561553488803\\
8	5.56080694687045\\
9	5.68605530889895\\
10	5.93546986490507\\
11	6.99061346396084\\
12	6.36574646391378\\
13	9.11022208265634\\
14	7.67213807461259\\
15	8.33946381370814\\
16	7.15209587507093\\
17	6.68906301994462\\
18	7.23246874161653\\
19	9.56370513600328\\
20	9.04527109367675\\
21	8.74347203885827\\
22	9.73022295290387\\
23	4.67865234262125\\
24	5.28981258711519\\
25	7.19336123048303\\
26	8.48478002642789\\
27	5.5758088599063\\
28	9.04331147088163\\
29	7.6354701637822\\
30	7.40489457460425\\
31	8.40691208922857\\
32	9.08074420841868\\
33	9.54771408024394\\
34	7.8099160633903\\
35	8.78393266685271\\
36	9.78108643743424\\
37	8.84473303915441\\
38	9.55953905208391\\
39	8.76308846427872\\
40	7.67586968232457\\
41	7.97861205314364\\
42	7.41606403417743\\
43	6.89664192622193\\
44	9.538430629982\\
45	9.55048420672919\\
46	9.89859212053971\\
47	7.42673356186299\\
48	7.54305478048709\\
49	9.06368525991419\\
50	8.6390065027433\\
51	9.10874445166314\\
52	8.69713675995252\\
53	9.22453412489386\\
54	6.95649596836826\\
55	8.64643200922266\\
56	9.96951249489921\\
57	9.99270016231322\\
58	8.66897153428372\\
59	9.11394099420779\\
60	6.94389615410669\\
61	6.56363768465611\\
62	5.67024596565343\\
63	8.39000201162442\\
64	7.57797971261731\\
65	7.29578877587616\\
66	9.23334396485778\\
67	9.44139271216136\\
68	9.98380066247266\\
69	6.82068599113201\\
70	9.59954338238093\\
71	7.5322613627097\\
72	5.30440230275672\\
73	9.70236412705304\\
74	7.97122268459983\\
75	6.04912503974008\\
76	9.64873081396142\\
77	8.3015350439036\\
78	6.23555829137287\\
79	8.93971674368628\\
80	9.29769899330843\\
81	7.54532445421813\\
82	9.06469945863363\\
83	9.09909684931804\\
84	7.63059701685664\\
85	7.22990493126227\\
86	9.01767186979421\\
87	7.56972285498337\\
88	8.07268854231828\\
89	6.20966307019671\\
90	9.26359438002029\\
91	9.7561782893714\\
92	5.59025301739735\\
93	5.95634355178483\\
94	8.79191222735646\\
95	6.7346878752702\\
96	8.47804626248313\\
97	5.77136418001656\\
98	1\\
99	5.20041749538385\\
100	1\\
};
\addlegendentry{$\ol{\eta}_p$ - \texttt{CUR-CR-OS}}

\addplot [color=mycolor1, line width=1.2pt, mark=o, mark options={solid, mycolor1}]
  table[row sep=crcr]{%
1	3.98545416460458\\
2	3.38284287782404\\
3	3.97939912615686\\
4	4.19745307577015\\
5	5.57408130510737\\
6	5.16772210603135\\
7	5.02617031593698\\
8	5.79013655842614\\
9	5.5080431922631\\
10	6.57932459462923\\
11	7.66377140673162\\
12	6.65583732673056\\
13	7.90577389209253\\
14	6.42116798701518\\
15	9.27607370190683\\
16	8.0854452213\\
17	8.1296983694795\\
18	8.64823591802751\\
19	9.65406333599301\\
20	7.93680846333341\\
21	8.55082242753825\\
22	7.33116108854232\\
23	6.04993830344466\\
24	6.96150047386599\\
25	6.63538832068593\\
26	8.41260912009345\\
27	8.54167751180103\\
28	9.74523714005642\\
29	8.33494318907065\\
30	7.76557113250981\\
31	8.25847665146372\\
32	7.50484515937132\\
33	7.16345315082824\\
34	6.99269384119891\\
35	8.56088986631181\\
36	8.97244038913234\\
37	9.34753463823979\\
38	9.47383744241415\\
39	6.84651336833781\\
40	9.85432036516462\\
41	9.61910023234825\\
42	7.57308043881831\\
43	9.20260868162932\\
44	6.59590236534179\\
45	8.65775430131711\\
46	9.82666123196691\\
47	6.60124136706531\\
48	8.85582753314232\\
49	9.08983877655439\\
50	9.76173789952972\\
51	7.25985168716568\\
52	9.70799050886814\\
53	8.18878589894268\\
54	7.11034311287122\\
55	7.68942653933445\\
56	9.49380893161032\\
57	9.70977118755657\\
58	9.56897692021933\\
59	8.76189146473348\\
60	9.57288105486387\\
61	9.33920114003838\\
62	6.64261154134567\\
63	9.86334829506263\\
64	8.25768545142639\\
65	9.2074705648869\\
66	6.04139130009072\\
67	7.73792345030571\\
68	9.36636980310854\\
69	7.57710903936235\\
70	7.43782834161071\\
71	9.6044914905633\\
72	8.65152445352239\\
73	9.12847302715515\\
74	8.38552652032718\\
75	8.25537878447889\\
76	6.94045759047569\\
77	7.4948587970807\\
78	9.88999271523068\\
79	8.33671385374168\\
80	7.36669227848906\\
81	6.2011836403687\\
82	6.24527568427248\\
83	8.84047461564903\\
84	8.16653453518879\\
85	9.08756463768366\\
86	9.40144632349058\\
87	9.53630869682972\\
88	9.88520171171887\\
89	8.1109905259229\\
90	5.79148317094228\\
91	4.78544139830749\\
92	7.89607096210637\\
93	5.51918809910307\\
94	7.16678593002389\\
95	8.28898287917187\\
96	7.28980521879005\\
97	4.92177920717198\\
98	4.93773050260402\\
99	5.07390209734838\\
100	1\\
};
\addlegendentry{$\ol{\eta}_s$ - \texttt{CUR-CR-OS}}
\end{axis}

\end{tikzpicture}%

%% file: figures/Toy_m_fast.tex
\definecolor{mycolor1}{rgb}{0.00000,0.44700,0.74100}%
\definecolor{mycolor2}{rgb}{0.85000,0.32500,0.09800}%
\definecolor{mycolor3}{rgb}{0.5,0.5,0.5}
\begin{tikzpicture}[scale=0.8]

\begin{axis}[%
width=2.2in,
height=2.0in,
at={(1.011in,0.642in)},
scale only axis,
xmin=0,
xmax=110,
xlabel style={font=\color{white!15!black}},
xlabel={$r$},
%ymode=log,
ymin=0,
ymax=10,
yminorticks=true,
ylabel style={font=\color{white!15!black}},
ylabel={Oversampling values},
axis background/.style={fill=white},
legend style={fill opacity=0.8, 
font=\scriptsize,
  draw opacity=1,
  text opacity=1,
  at={(0.04,0.9)},
  anchor=west}
  % draw=lightgray204}
]
\addplot [color=mycolor2, line width=1.2pt, mark=o, mark options={solid, mycolor2}]
  table[row sep=crcr]{%
1	0\\
2	0\\
3	0\\
4	0\\
5	0\\
6	0\\
7	0\\
8	0\\
9	0\\
10	0\\
11	0\\
12	0\\
13	0\\
14	0\\
15	0\\
16	0\\
17	0\\
18	0\\
19	0\\
20	0\\
21	0\\
22	0\\
23	0\\
24	1\\
25	0\\
26	0\\
27	0\\
28	1\\
29	0\\
30	0\\
31	0\\
32	0\\
33	1\\
34	1\\
35	0\\
36	1\\
37	1\\
38	1\\
39	0\\
40	1\\
41	1\\
42	0\\
43	2\\
44	1\\
45	1\\
46	2\\
47	1\\
48	1\\
49	1\\
50	1\\
51	2\\
52	2\\
53	2\\
54	8\\
55	4\\
56	5\\
57	5\\
58	2\\
59	9\\
60	7\\
61	4\\
62	5\\
63	2\\
64	2\\
65	4\\
66	3\\
67	7\\
68	7\\
69	6\\
70	6\\
71	6\\
72	7\\
73	5\\
74	4\\
75	4\\
76	4\\
77	4\\
78	3\\
79	6\\
80	5\\
81	3\\
82	5\\
83	5\\
84	4\\
85	4\\
86	6\\
87	3\\
88	4\\
89	5\\
90	6\\
91	3\\
92	3\\
93	4\\
94	4\\
95	5\\
96	3\\
97	3\\
98	1\\
99	1\\
100	0\\
};
\addlegendentry{ $m_r$ - \texttt{CUR-CR-OS} }

\addplot [color=mycolor1, line width=1.2pt, mark=o, mark options={solid, mycolor1}]
  table[row sep=crcr]{%
1	0\\
2	0\\
3	0\\
4	0\\
5	0\\
6	0\\
7	0\\
8	0\\
9	0\\
10	0\\
11	0\\
12	0\\
13	0\\
14	0\\
15	0\\
16	0\\
17	0\\
18	0\\
19	0\\
20	0\\
21	0\\
22	1\\
23	0\\
24	0\\
25	0\\
26	1\\
27	1\\
28	1\\
29	1\\
30	0\\
31	1\\
32	1\\
33	1\\
34	1\\
35	0\\
36	0\\
37	0\\
38	0\\
39	0\\
40	0\\
41	0\\
42	1\\
43	1\\
44	0\\
45	0\\
46	0\\
47	0\\
48	0\\
49	0\\
50	0\\
51	2\\
52	3\\
53	6\\
54	1\\
55	6\\
56	5\\
57	4\\
58	5\\
59	5\\
60	5\\
61	4\\
62	5\\
63	8\\
64	5\\
65	6\\
66	4\\
67	5\\
68	3\\
69	6\\
70	3\\
71	4\\
72	4\\
73	4\\
74	6\\
75	4\\
76	5\\
77	5\\
78	4\\
79	4\\
80	5\\
81	5\\
82	2\\
83	4\\
84	3\\
85	2\\
86	3\\
87	4\\
88	4\\
89	5\\
90	4\\
91	5\\
92	4\\
93	3\\
94	1\\
95	3\\
96	1\\
97	1\\
98	0\\
99	0\\
100	0\\
};

\addlegendentry{ $m_c$ - \texttt{CUR-CR-OS} }

\end{axis}

\end{tikzpicture}%

%% file: figures/Toy_m_slow.tex
\definecolor{mycolor1}{rgb}{0.00000,0.44700,0.74100}%
\definecolor{mycolor2}{rgb}{0.85000,0.32500,0.09800}%
\definecolor{mycolor3}{rgb}{0.5,0.5,0.5}
\begin{tikzpicture}[scale=0.8]

\begin{axis}[%
width=2.2in,
height=2.0in,
at={(1.011in,0.642in)},
scale only axis,
xmin=0,
xmax=110,
xlabel style={font=\color{white!15!black}},
xlabel={$r$},
%ymode=log,
ymin=0,
ymax=10,
yminorticks=true,
ylabel style={font=\color{white!15!black}},
ylabel={Oversampling values},
axis background/.style={fill=white},
legend style={fill opacity=0.8, 
font=\scriptsize,
  draw opacity=1,
  text opacity=1,
  at={(0.04,0.9)},
  anchor=west}
  % draw=lightgray204}
]
\addplot [color=mycolor2, line width=1.2pt, mark=o, mark options={solid, mycolor2}]
  table[row sep=crcr]{%
1	0\\
2	0\\
3	0\\
4	0\\
5	0\\
6	0\\
7	0\\
8	0\\
9	0\\
10	0\\
11	0\\
12	0\\
13	0\\
14	0\\
15	0\\
16	0\\
17	1\\
18	1\\
19	3\\
20	1\\
21	1\\
22	4\\
23	3\\
24	4\\
25	2\\
26	4\\
27	4\\
28	1\\
29	2\\
30	1\\
31	5\\
32	2\\
33	5\\
34	6\\
35	4\\
36	4\\
37	3\\
38	3\\
39	4\\
40	4\\
41	4\\
42	5\\
43	7\\
44	5\\
45	3\\
46	3\\
47	5\\
48	3\\
49	1\\
50	3\\
51	2\\
52	3\\
53	3\\
54	6\\
55	5\\
56	2\\
57	2\\
58	2\\
59	5\\
60	7\\
61	9\\
62	9\\
63	5\\
64	5\\
65	7\\
66	6\\
67	5\\
68	4\\
69	5\\
70	4\\
71	4\\
72	8\\
73	3\\
74	5\\
75	10\\
76	7\\
77	6\\
78	6\\
79	4\\
80	6\\
81	5\\
82	4\\
83	3\\
84	4\\
85	6\\
86	6\\
87	6\\
88	4\\
89	4\\
90	3\\
91	3\\
92	3\\
93	2\\
94	2\\
95	1\\
96	1\\
97	1\\
98	2\\
99	0\\
100	0\\
};
\addlegendentry{ $m_r$ - \texttt{CUR-CR-OS} }

\addplot [color=mycolor1, line width=1.2pt, mark=o, mark options={solid, mycolor1}]
  table[row sep=crcr]{%
1	0\\
2	0\\
3	0\\
4	0\\
5	0\\
6	0\\
7	0\\
8	0\\
9	0\\
10	0\\
11	0\\
12	0\\
13	0\\
14	0\\
15	0\\
16	0\\
17	1\\
18	2\\
19	1\\
20	2\\
21	1\\
22	2\\
23	3\\
24	3\\
25	1\\
26	1\\
27	1\\
28	5\\
29	1\\
30	1\\
31	2\\
32	2\\
33	2\\
34	2\\
35	2\\
36	2\\
37	1\\
38	3\\
39	7\\
40	4\\
41	2\\
42	4\\
43	6\\
44	7\\
45	5\\
46	7\\
47	4\\
48	3\\
49	3\\
50	4\\
51	8\\
52	4\\
53	8\\
54	7\\
55	7\\
56	7\\
57	7\\
58	4\\
59	5\\
60	4\\
61	5\\
62	6\\
63	8\\
64	6\\
65	7\\
66	5\\
67	5\\
68	4\\
69	6\\
70	5\\
71	6\\
72	5\\
73	6\\
74	5\\
75	6\\
76	7\\
77	6\\
78	6\\
79	6\\
80	5\\
81	5\\
82	4\\
83	4\\
84	5\\
85	7\\
86	5\\
87	4\\
88	3\\
89	3\\
90	4\\
91	4\\
92	3\\
93	3\\
94	2\\
95	1\\
96	1\\
97	1\\
98	1\\
99	0\\
100	0\\
};

\addlegendentry{ $m_c$ - \texttt{CUR-CR-OS} }

\end{axis}

\end{tikzpicture}%

%% file: figures/Burger_Sensitivity.tex
\definecolor{mycolor1}{rgb}{0.00000,0.44700,0.74100}%
\definecolor{mycolor2}{rgb}{0.85000,0.32500,0.09800}%
\definecolor{mycolor3}{rgb}{0.5,0.5,0.5}
\begin{tikzpicture}[scale=0.8]

\begin{axis}[%
width=2.4in,
height=2.2in,
at={(1.011in,0.642in)},
scale only axis,
xmin=0,
xmax=5,
xlabel style={font=\color{white!15!black}},
xlabel={Time},
ymode=log,
ymin=1e-18,
ymax=1e-06,
yminorticks=true,
ylabel style={font=\color{white!15!black}},
ylabel={$\mathcal{E}$},
axis background/.style={fill=white},
legend style={at={(0.35,0.07)}, font=\scriptsize, anchor=south west, legend cell align=left, align=left, draw=white!15!black}
]
\addplot [color=mycolor1, line width=1.2pt]
  table[row sep=crcr]{%
0	6.25305022111612e-19\\
0.025	4.84418336361229e-12\\
0.05	5.00338104968003e-12\\
0.075	4.83007823943687e-12\\
0.1	4.60129232524839e-12\\
0.125	4.72823933475612e-12\\
0.15	5.75890753956158e-12\\
0.175	6.76513103973087e-12\\
0.2	1.02135683110072e-11\\
0.225	1.8424529006758e-11\\
0.25	2.2358801987881e-11\\
0.275	2.50900365737414e-11\\
0.3	3.16774648645064e-11\\
0.325	3.72292316822725e-11\\
0.35	4.02515206973135e-11\\
0.375	4.5979517443859e-11\\
0.4	5.20127071515712e-11\\
0.425	5.9786926759319e-11\\
0.45	6.36554706152921e-11\\
0.475	7.0026760989017e-11\\
0.5	7.63326574292481e-11\\
0.525	8.37553176665008e-11\\
0.55	9.41410948525034e-11\\
0.575	1.044105720637e-10\\
0.6	1.17805458110726e-10\\
0.625	1.32710188103736e-10\\
0.65	1.54643744054597e-10\\
0.675	2.17345367479236e-10\\
0.7	2.56870226122177e-10\\
0.725	2.65148321858087e-10\\
0.75	2.82607662969992e-10\\
0.775	3.36889431734218e-10\\
0.8	3.77467382330632e-10\\
0.825	3.78992887220955e-10\\
0.85	3.76717191429742e-10\\
0.875	3.68643772601016e-10\\
0.9	3.75042967511847e-10\\
0.925	3.79658540032414e-10\\
0.95	3.83335809665147e-10\\
0.975	3.84174697103418e-10\\
1	3.80083459892834e-10\\
1.025	3.74876921732617e-10\\
1.05	3.69452506483339e-10\\
1.075	3.63367091790084e-10\\
1.1	3.57741163859561e-10\\
1.125	3.53351683397179e-10\\
1.15	3.48240922942776e-10\\
1.175	3.43281992487169e-10\\
1.2	3.38689745933788e-10\\
1.225	3.34599426316148e-10\\
1.25	3.31256763753513e-10\\
1.275	3.29745561521685e-10\\
1.3	3.28733044007194e-10\\
1.325	3.27222033088131e-10\\
1.35	3.27904327758105e-10\\
1.375	3.31495497181428e-10\\
1.4	3.30137712949311e-10\\
1.425	3.27495325230214e-10\\
1.45	3.24151967781056e-10\\
1.475	3.20896206216179e-10\\
1.5	3.18034153580086e-10\\
1.525	3.15362553089199e-10\\
1.55	3.12997884458775e-10\\
1.575	3.11848388380705e-10\\
1.6	3.12781493116931e-10\\
1.625	3.13370312404961e-10\\
1.65	3.1637941419085e-10\\
1.675	3.15247695906376e-10\\
1.7	3.21448941710022e-10\\
1.725	3.29825884888095e-10\\
1.75	3.36958393917949e-10\\
1.775	3.42964191003062e-10\\
1.8	3.55350578347948e-10\\
1.825	4.10178261908487e-10\\
1.85	4.11013530338429e-10\\
1.875	4.19508709860546e-10\\
1.9	4.12954081335394e-10\\
1.925	4.1047470820738e-10\\
1.95	4.07422533333232e-10\\
1.975	4.02941483866678e-10\\
2	3.9835048167519e-10\\
2.025	3.94288211336159e-10\\
2.05	3.90837001171595e-10\\
2.075	3.88005075413298e-10\\
2.1	3.86353594318163e-10\\
2.125	3.85651399708986e-10\\
2.15	3.83040697386446e-10\\
2.175	3.83532956828482e-10\\
2.2	3.87875782694779e-10\\
2.225	3.83527289678237e-10\\
2.25	3.82463771271744e-10\\
2.275	3.84575794143559e-10\\
2.3	3.89003974932395e-10\\
2.325	3.95553553062162e-10\\
2.35	4.06594779465629e-10\\
2.375	4.21114180397701e-10\\
2.4	4.36228677196089e-10\\
2.425	4.53097898598504e-10\\
2.45	4.6942188250978e-10\\
2.475	4.83309544418103e-10\\
2.5	4.94386787044707e-10\\
2.525	5.03165363742236e-10\\
2.55	5.10592368828991e-10\\
2.575	5.15469817940423e-10\\
2.6	5.20025190953739e-10\\
2.625	5.27930965517875e-10\\
2.65	5.34408035626271e-10\\
2.675	5.3684910504678e-10\\
2.7	5.46916149384499e-10\\
2.725	5.85462175906253e-10\\
2.75	6.17962057842261e-10\\
2.775	6.47491667269843e-10\\
2.8	6.43399481121976e-10\\
2.825	6.71738169553723e-10\\
2.85	6.75009351894978e-10\\
2.875	6.74933412394452e-10\\
2.9	7.02904762609204e-10\\
2.925	6.95283787182309e-10\\
2.95	6.90182006593936e-10\\
2.975	6.84044062959218e-10\\
3	6.77690104613794e-10\\
3.025	6.72555714671098e-10\\
3.05	6.83622688969217e-10\\
3.075	6.73508451388015e-10\\
3.1	6.65655874272813e-10\\
3.125	6.60297614405119e-10\\
3.15	6.55890642602329e-10\\
3.175	6.52778722275762e-10\\
3.2	6.50563904950563e-10\\
3.225	6.46722844564237e-10\\
3.25	6.43672294875356e-10\\
3.275	6.41698164184176e-10\\
3.3	6.41717676405914e-10\\
3.325	6.41012028249633e-10\\
3.35	6.35353287730107e-10\\
3.375	6.28393265140149e-10\\
3.4	6.32579698094653e-10\\
3.425	6.32274016928777e-10\\
3.45	6.32174748720736e-10\\
3.475	6.3103472201322e-10\\
3.5	6.28444087258469e-10\\
3.525	6.249129855e-10\\
3.55	6.20769458968355e-10\\
3.575	6.19391440743983e-10\\
3.6	6.21529260702606e-10\\
3.625	6.43574713278679e-10\\
3.65	6.45266294968902e-10\\
3.675	6.8270221700201e-10\\
3.7	6.77054292664953e-10\\
3.725	6.98189500828295e-10\\
3.75	8.05905074199481e-10\\
3.775	8.14345445821159e-10\\
3.8	8.34953882523123e-10\\
3.825	9.63218073414217e-10\\
3.85	1.30329586597543e-09\\
3.875	1.3504370439712e-09\\
3.9	1.27084706466765e-09\\
3.925	1.23154878114053e-09\\
3.95	1.21603051905559e-09\\
3.975	1.19700840648247e-09\\
4	1.16680310623403e-09\\
4.025	1.14521174883835e-09\\
4.05	1.1394320548118e-09\\
4.075	1.14024207926026e-09\\
4.1	1.14399189803414e-09\\
4.125	1.17513284173665e-09\\
4.15	1.18078559359629e-09\\
4.175	1.19754848174343e-09\\
4.2	1.19211626444248e-09\\
4.225	1.18826364494293e-09\\
4.25	1.1828286930461e-09\\
4.275	1.17793498517659e-09\\
4.3	1.16960115501923e-09\\
4.325	1.15568464875564e-09\\
4.35	1.14061355142726e-09\\
4.375	1.12771158349929e-09\\
4.4	1.10928928753872e-09\\
4.425	1.09252823338564e-09\\
4.45	1.07921029935789e-09\\
4.475	1.06827500717195e-09\\
4.5	1.05274947795355e-09\\
4.525	1.03584978227583e-09\\
4.55	1.01778620556297e-09\\
4.575	1.00031069798557e-09\\
4.6	9.87919203238332e-10\\
4.625	9.79073448622687e-10\\
4.65	9.74155096134002e-10\\
4.675	9.75164230662376e-10\\
4.7	9.74540160378675e-10\\
4.725	9.78315667036506e-10\\
4.75	9.85215390137776e-10\\
4.775	9.93981125865893e-10\\
4.8	1.00400308211502e-09\\
4.825	1.0134057375183e-09\\
4.85	1.02372804178683e-09\\
4.875	1.02984736558598e-09\\
4.9	1.03534510973013e-09\\
4.925	1.03876618926706e-09\\
4.95	1.04018719598998e-09\\
4.975	1.03965263166666e-09\\
5	1.03853458143577e-09\\
};
\addlegendentry{$\epsilon_{os} = 5$}

\addplot [color=mycolor2, line width=1.2pt]
  table[row sep=crcr]{%
0	6.25305022111612e-19\\
0.025	5.69096745972987e-12\\
0.05	5.87651583327269e-12\\
0.075	5.63681296952606e-12\\
0.1	5.49310292741789e-12\\
0.125	5.53478307758135e-12\\
0.15	6.85512018097151e-12\\
0.175	8.01439190702309e-12\\
0.2	1.11685792521109e-11\\
0.225	1.88834845387802e-11\\
0.25	2.28867796928877e-11\\
0.275	2.60187766906643e-11\\
0.3	3.74017488726199e-11\\
0.325	4.50760549450117e-11\\
0.35	4.89307210490719e-11\\
0.375	5.58412943731898e-11\\
0.4	6.10152762322836e-11\\
0.425	6.99769041217023e-11\\
0.45	7.19475892152266e-11\\
0.475	7.76684801178352e-11\\
0.5	8.32014024196618e-11\\
0.525	9.18625559259438e-11\\
0.55	1.02246527598012e-10\\
0.575	1.11766579112182e-10\\
0.6	1.25566972834276e-10\\
0.625	1.54751887977056e-10\\
0.65	1.83314384179735e-10\\
0.675	2.41057719484233e-10\\
0.7	2.79335951151002e-10\\
0.725	2.88262010608968e-10\\
0.75	3.07436707470317e-10\\
0.775	3.57438085461045e-10\\
0.8	4.30758655271148e-10\\
0.825	4.28197538005432e-10\\
0.85	4.25948844305566e-10\\
0.875	4.21275298290855e-10\\
0.9	4.25150398110267e-10\\
0.925	4.22908364954509e-10\\
0.95	4.21417633483144e-10\\
0.975	4.19690609195517e-10\\
1	4.16177846400022e-10\\
1.025	4.1116562796501e-10\\
1.05	4.05801471963483e-10\\
1.075	4.01242680278546e-10\\
1.1	3.96790785815289e-10\\
1.125	3.92001125161413e-10\\
1.15	3.87562864086871e-10\\
1.175	3.84991845684217e-10\\
1.2	3.8176423448352e-10\\
1.225	3.78804079924392e-10\\
1.25	3.76900721702475e-10\\
1.275	3.76690441466994e-10\\
1.3	3.76317664119081e-10\\
1.325	3.74399386287834e-10\\
1.35	3.90319464182517e-10\\
1.375	3.84318704902843e-10\\
1.4	3.81877959438464e-10\\
1.425	3.80467891202757e-10\\
1.45	3.78015576584159e-10\\
1.475	3.75167971163898e-10\\
1.5	3.72668911493393e-10\\
1.525	3.70625241618876e-10\\
1.55	3.68677686634641e-10\\
1.575	3.67318676574472e-10\\
1.6	3.67725769861902e-10\\
1.625	3.67658987391088e-10\\
1.65	3.69738672130698e-10\\
1.675	3.69749211345729e-10\\
1.7	3.73386039412466e-10\\
1.725	3.79939997932076e-10\\
1.75	3.83537015267285e-10\\
1.775	3.92552128820945e-10\\
1.8	4.03085251117664e-10\\
1.825	4.51429974243296e-10\\
1.85	4.54801593440906e-10\\
1.875	4.59410149747279e-10\\
1.9	4.50519007645105e-10\\
1.925	4.48221867313421e-10\\
1.95	4.46236427812924e-10\\
1.975	4.43899624280683e-10\\
2	4.51010141703764e-10\\
2.025	4.52472025561706e-10\\
2.05	4.52653698686914e-10\\
2.075	4.51076227695062e-10\\
2.1	4.49630867570448e-10\\
2.125	4.48677053811643e-10\\
2.15	4.45058246144408e-10\\
2.175	4.45098880642554e-10\\
2.2	4.50476111734346e-10\\
2.225	4.45776686807559e-10\\
2.25	4.43423422865668e-10\\
2.275	4.43260981526476e-10\\
2.3	4.45953096756605e-10\\
2.325	4.51034620751969e-10\\
2.35	4.6088660982208e-10\\
2.375	4.74786674391666e-10\\
2.4	4.90081936602579e-10\\
2.425	5.06851621297898e-10\\
2.45	5.23152385202036e-10\\
2.475	5.3736675681119e-10\\
2.5	5.48731563636845e-10\\
2.525	5.57225461096708e-10\\
2.55	5.63834780835845e-10\\
2.575	5.67511278429925e-10\\
2.6	5.7051720526354e-10\\
2.625	5.75475375268177e-10\\
2.65	5.80528945572205e-10\\
2.675	5.80021180914504e-10\\
2.7	5.85455402648858e-10\\
2.725	6.1707391679015e-10\\
2.75	6.45827234829841e-10\\
2.775	6.83381782478584e-10\\
2.8	6.60672659859868e-10\\
2.825	6.57632899645225e-10\\
2.85	6.49736680814883e-10\\
2.875	6.40682790602421e-10\\
2.9	6.43640596643247e-10\\
2.925	6.35773999375859e-10\\
2.95	6.32994145916786e-10\\
2.975	6.27410355571316e-10\\
3	1.73662917312528e-09\\
3.025	1.85420980344424e-09\\
3.05	1.83542708398123e-09\\
3.075	1.85436521348735e-09\\
3.1	1.90651553566966e-09\\
3.125	1.9452280374482e-09\\
3.15	1.99929384763656e-09\\
3.175	2.04278549435217e-09\\
3.2	2.08399461077725e-09\\
3.225	2.12183788671326e-09\\
3.25	2.1565405558714e-09\\
3.275	2.18871246860225e-09\\
3.3	2.2165213380057e-09\\
3.325	2.23868873057003e-09\\
3.35	2.254364169562e-09\\
3.375	2.26554212195198e-09\\
3.4	2.2737013050389e-09\\
3.425	2.28140490517693e-09\\
3.45	2.28356046028904e-09\\
3.475	2.28365639739195e-09\\
3.5	2.28210154728432e-09\\
3.525	2.28268747500704e-09\\
3.55	2.27265024237104e-09\\
3.575	2.26369432180691e-09\\
3.6	2.25455120035859e-09\\
3.625	2.24817614494271e-09\\
3.65	2.23789604529287e-09\\
3.675	2.23561454669387e-09\\
3.7	2.25604295151631e-09\\
3.725	2.80521088962888e-09\\
3.75	2.64748122647463e-09\\
3.775	2.69490269794283e-09\\
3.8	2.71671391702753e-09\\
3.825	2.75372426947146e-09\\
3.85	2.7588512018849e-09\\
3.875	2.74420404456133e-09\\
3.9	2.73159423331677e-09\\
3.925	2.72416437089105e-09\\
3.95	2.71028102687331e-09\\
3.975	2.68373617910093e-09\\
4	2.6660190398233e-09\\
4.025	2.66769096960186e-09\\
4.05	2.66096935505665e-09\\
4.075	2.65819087884145e-09\\
4.1	2.90880067768412e-09\\
4.125	2.83822425413668e-09\\
4.15	2.81596569878345e-09\\
4.175	2.81179270224427e-09\\
4.2	2.77033508128976e-09\\
4.225	2.77273934256692e-09\\
4.25	2.76794159243744e-09\\
4.275	2.75403725198205e-09\\
4.3	2.74553497106368e-09\\
4.325	2.73320931097992e-09\\
4.35	2.70829490693463e-09\\
4.375	2.68413564274159e-09\\
4.4	2.66186340094415e-09\\
4.425	2.64720772790371e-09\\
4.45	2.63435235205974e-09\\
4.475	2.61325494324552e-09\\
4.5	2.56707401836539e-09\\
4.525	2.47247679009104e-09\\
4.55	2.33636527434705e-09\\
4.575	2.19192105594553e-09\\
4.6	2.07725275619879e-09\\
4.625	2.00278375442376e-09\\
4.65	1.96077881302343e-09\\
4.675	1.9362350598924e-09\\
4.7	1.92057937787573e-09\\
4.725	1.90683497734424e-09\\
4.75	1.89695063380668e-09\\
4.775	1.88743921714353e-09\\
4.8	1.8801517269825e-09\\
4.825	1.87106045862147e-09\\
4.85	1.86564541608214e-09\\
4.875	1.8625508096685e-09\\
4.9	1.86169806546958e-09\\
4.925	1.85696918650319e-09\\
4.95	1.85202127601008e-09\\
4.975	1.84796843861989e-09\\
5	1.84602025547158e-09\\
};
\addlegendentry{$\epsilon_{os} = 10$}

\addplot [color=mycolor3, line width=1.2pt]
  table[row sep=crcr]{%
0	6.25305022111612e-19\\
0.025	5.60745301884619e-12\\
0.05	5.79831394136184e-12\\
0.075	5.55759158568115e-12\\
0.1	5.40917102948759e-12\\
0.125	5.45088158684646e-12\\
0.15	6.79253285499461e-12\\
0.175	8.1465234162711e-12\\
0.2	1.13779765755678e-11\\
0.225	6.3230261850176e-11\\
0.25	4.26558674495487e-11\\
0.275	3.99887491904058e-11\\
0.3	4.43678048092438e-11\\
0.325	5.12201112575699e-11\\
0.35	5.05266764648135e-11\\
0.375	5.8049638447657e-11\\
0.4	6.43727558893603e-11\\
0.425	6.91459510351302e-11\\
0.45	7.26102403873988e-11\\
0.475	7.91779072034606e-11\\
0.5	8.45770657190616e-11\\
0.525	9.21461180704877e-11\\
0.55	1.02047745598569e-10\\
0.575	1.1166582597781e-10\\
0.6	1.25329030014041e-10\\
0.625	1.43026913150272e-10\\
0.65	1.74327014854029e-10\\
0.675	2.32430675466079e-10\\
0.7	2.75549168712841e-10\\
0.725	1.51854173359346e-09\\
0.75	1.26395344247261e-09\\
0.775	1.20850500095643e-09\\
0.8	1.22250610790293e-09\\
0.825	1.23929784015131e-09\\
0.85	1.26907095679912e-09\\
0.875	1.26356495437405e-09\\
0.9	3.31296442502755e-09\\
0.925	3.01190179779599e-09\\
0.95	2.96072029209368e-09\\
0.975	2.8850463944749e-09\\
1	2.8018162432493e-09\\
1.025	2.74149908577275e-09\\
1.05	2.69272677587909e-09\\
1.075	2.66487376569192e-09\\
1.1	2.65285346163511e-09\\
1.125	2.65392545734138e-09\\
1.15	2.66873529283237e-09\\
1.175	2.69132217094031e-09\\
1.2	2.72168065107497e-09\\
1.225	2.75712339022326e-09\\
1.25	2.79243657528913e-09\\
1.275	2.82872969772917e-09\\
1.3	2.86410410766256e-09\\
1.325	2.89688632399329e-09\\
1.35	2.92555454487042e-09\\
1.375	2.94997125591995e-09\\
1.4	2.96953792825691e-09\\
1.425	2.98422258997384e-09\\
1.45	2.99559533094198e-09\\
1.475	3.0033881978988e-09\\
1.5	3.00715855932725e-09\\
1.525	3.00733523230055e-09\\
1.55	3.00432733743966e-09\\
1.575	2.99792672256449e-09\\
1.6	2.98871206559014e-09\\
1.625	2.97702099209601e-09\\
1.65	2.96302580212182e-09\\
1.675	2.94642156312389e-09\\
1.7	2.92704361123518e-09\\
1.725	2.9102987558636e-09\\
1.75	2.89651391813221e-09\\
1.775	2.87292529096582e-09\\
1.8	2.84975716848152e-09\\
1.825	2.8254132135839e-09\\
1.85	2.79929832083631e-09\\
1.875	3.26744257853039e-09\\
1.9	3.17305882299703e-09\\
1.925	3.13919075422812e-09\\
1.95	3.11405062757629e-09\\
1.975	3.09862529799631e-09\\
2	3.08814180364103e-09\\
2.025	3.94854644304991e-09\\
2.05	3.66602871964138e-09\\
2.075	3.60979930511216e-09\\
2.1	3.59398216233921e-09\\
2.125	3.60772206207674e-09\\
2.15	3.61622196909103e-09\\
2.175	3.63210806867839e-09\\
2.2	3.64838448697322e-09\\
2.225	3.65585837649612e-09\\
2.25	3.66293370913513e-09\\
2.275	3.67627142284468e-09\\
2.3	3.70725086583033e-09\\
2.325	3.76874027145243e-09\\
2.35	3.86838963293908e-09\\
2.375	4.00405867220515e-09\\
2.4	4.16032205839608e-09\\
2.425	4.31765080243762e-09\\
2.45	4.45629106637988e-09\\
2.475	4.56621753396737e-09\\
2.5	4.64594483328198e-09\\
2.525	4.69597431861441e-09\\
2.55	4.72199307250908e-09\\
2.575	4.72878336575914e-09\\
2.6	4.72057539406905e-09\\
2.625	4.70030442921901e-09\\
2.65	4.67496234246733e-09\\
2.675	4.64653249940786e-09\\
2.7	4.61030097996443e-09\\
2.725	4.5731033953322e-09\\
2.75	4.53546566718512e-09\\
2.775	4.49578134288474e-09\\
2.8	4.45726321870603e-09\\
2.825	4.41650679820055e-09\\
2.85	4.37521421143979e-09\\
2.875	4.34016960672101e-09\\
2.9	4.30276049728477e-09\\
2.925	4.26601707316556e-09\\
2.95	4.23080389256009e-09\\
2.975	4.19757251579923e-09\\
3	4.16449850644227e-09\\
3.025	4.13047263265251e-09\\
3.05	4.09519047687033e-09\\
3.075	4.05506142753741e-09\\
3.1	4.01789534633255e-09\\
3.125	3.98139062025219e-09\\
3.15	3.94571087574705e-09\\
3.175	3.90798613385303e-09\\
3.2	3.87311170181328e-09\\
3.225	3.83876777545355e-09\\
3.25	3.80478767522706e-09\\
3.275	3.77120593777195e-09\\
3.3	3.73844008861518e-09\\
3.325	3.70522232030824e-09\\
3.35	3.67206352028585e-09\\
3.375	3.64076673197286e-09\\
3.4	3.61224794462588e-09\\
3.425	3.58293627343037e-09\\
3.45	3.55463486470146e-09\\
3.475	3.52801071055727e-09\\
3.5	3.55182540037998e-09\\
3.525	3.51621781392813e-09\\
3.55	3.48616120447895e-09\\
3.575	3.45718258983729e-09\\
3.6	3.43004197434142e-09\\
3.625	3.4147897797002e-09\\
3.65	3.39372176988527e-09\\
3.675	3.37786972809334e-09\\
3.7	3.36904947844196e-09\\
3.725	3.35264444474662e-09\\
3.75	3.35132482658943e-09\\
3.775	3.33290209075834e-09\\
3.8	3.32855474193795e-09\\
3.825	3.35826400006395e-09\\
3.85	1.48525334424685e-08\\
3.875	1.36986748056409e-08\\
3.9	1.37726211994815e-08\\
3.925	1.4431066693815e-08\\
3.95	1.53291007237367e-08\\
3.975	1.61844211673215e-08\\
4	1.68933253847755e-08\\
4.025	1.74447684660837e-08\\
4.05	1.78505770746548e-08\\
4.075	1.81373848100014e-08\\
4.1	1.83257893598384e-08\\
4.125	1.84396772042212e-08\\
4.15	1.84896058064861e-08\\
4.175	1.84893535663544e-08\\
4.2	1.84491672521872e-08\\
4.225	1.83788317919773e-08\\
4.25	1.8283794362538e-08\\
4.275	1.81694876411753e-08\\
4.3	1.80363344600934e-08\\
4.325	1.78868913648975e-08\\
4.35	1.77252167711743e-08\\
4.375	1.75545995064955e-08\\
4.4	1.7379364677436e-08\\
4.425	1.7205197606311e-08\\
4.45	1.70355772115416e-08\\
4.475	1.68728513511942e-08\\
4.5	1.67165196904819e-08\\
4.525	1.65647323543874e-08\\
4.55	1.64161020987807e-08\\
4.575	1.62699914740737e-08\\
4.6	1.6126323710499e-08\\
4.625	1.59862373601245e-08\\
4.65	1.5848748911368e-08\\
4.675	1.57142812803715e-08\\
4.7	1.55816131581165e-08\\
4.725	1.5451395970101e-08\\
4.75	1.53236961602424e-08\\
4.775	1.51989716490999e-08\\
4.8	1.50763146449173e-08\\
4.825	1.49575402111903e-08\\
4.85	1.4839838115158e-08\\
4.875	1.47245620365486e-08\\
4.9	1.46113117318751e-08\\
4.925	1.45003268544061e-08\\
4.95	1.43915209050042e-08\\
4.975	1.42847903582663e-08\\
5	1.41799828436737e-08\\
};
\addlegendentry{$\epsilon_{os} = 20$}

\addplot [color=black, line width=1.2pt]
  table[row sep=crcr]{%
0	6.25305022111612e-19\\
0.025	5.60745301884619e-12\\
0.05	5.79831394136184e-12\\
0.075	5.55759158568115e-12\\
0.1	5.40917102948759e-12\\
0.125	5.45088158684646e-12\\
0.15	6.79253285499461e-12\\
0.175	8.1465234162711e-12\\
0.2	1.13779765755678e-11\\
0.225	6.3230261850176e-11\\
0.25	4.26558674495487e-11\\
0.275	3.99887491904058e-11\\
0.3	4.43678048092438e-11\\
0.325	5.12201112575699e-11\\
0.35	5.05266764648135e-11\\
0.375	5.8049638447657e-11\\
0.4	6.43727558893603e-11\\
0.425	6.91459510351302e-11\\
0.45	7.26102403873988e-11\\
0.475	7.91779072034606e-11\\
0.5	8.45770657190616e-11\\
0.525	9.21461180704877e-11\\
0.55	1.02047745598569e-10\\
0.575	1.1166582597781e-10\\
0.6	1.25329030014041e-10\\
0.625	1.43026913150272e-10\\
0.65	1.74327014854029e-10\\
0.675	2.32430675466079e-10\\
0.7	2.75549168712841e-10\\
0.725	1.51855602857986e-09\\
0.75	1.26398124926609e-09\\
0.775	1.20853406269709e-09\\
0.8	1.22253456319653e-09\\
0.825	1.23932563240718e-09\\
0.85	1.26909738862874e-09\\
0.875	1.2599670769349e-09\\
0.9	1.26386242094807e-09\\
0.925	1.26851791587079e-09\\
0.95	1.28431126218711e-09\\
0.975	1.29799863911958e-09\\
1	1.30935777246853e-09\\
1.025	1.32906933458431e-09\\
1.05	1.36244167915202e-09\\
1.075	1.47503912184211e-09\\
1.1	1.5444460664169e-09\\
1.125	1.5950921693642e-09\\
1.15	1.6233027276611e-09\\
1.175	1.63973434236658e-09\\
1.2	1.64910695294197e-09\\
1.225	1.65388234494248e-09\\
1.25	1.65587061180512e-09\\
1.275	1.65690886079296e-09\\
1.3	1.65457070382494e-09\\
1.325	1.65116037031506e-09\\
1.35	1.64635858891793e-09\\
1.375	1.64113598604197e-09\\
1.4	1.63568317632981e-09\\
1.425	1.63012383341414e-09\\
1.45	1.62374099703657e-09\\
1.475	1.61558541760905e-09\\
1.5	1.60700170855361e-09\\
1.525	1.59673584993072e-09\\
1.55	1.58791505555934e-09\\
1.575	1.57959286960087e-09\\
1.6	1.57057677599162e-09\\
1.625	1.56118748531783e-09\\
1.65	1.55192265602256e-09\\
1.675	1.5420190910756e-09\\
1.7	1.53452434864072e-09\\
1.725	1.52594622058728e-09\\
1.75	1.52252218983673e-09\\
1.775	1.50867975652756e-09\\
1.8	1.49765477763216e-09\\
1.825	1.49729655563626e-09\\
1.85	1.49778221204316e-09\\
1.875	1.4887376592192e-09\\
1.9	1.49635746065285e-09\\
1.925	1.49643032890445e-09\\
1.95	1.49025614889327e-09\\
1.975	1.47429673427215e-09\\
2	1.45292803887335e-09\\
2.025	1.43468723949956e-09\\
2.05	1.41772934959249e-09\\
2.075	1.4023474184318e-09\\
2.1	1.38817900028731e-09\\
2.125	1.37539359211911e-09\\
2.15	1.3654092670341e-09\\
2.175	1.35195985252429e-09\\
2.2	1.33944396480586e-09\\
2.225	1.32560008257612e-09\\
2.25	1.31388053979899e-09\\
2.275	1.31232901949913e-09\\
2.3	1.31139898632427e-09\\
2.325	1.32258765576757e-09\\
2.35	1.35241648814077e-09\\
2.375	1.4019303764193e-09\\
2.4	1.46426189563251e-09\\
2.425	1.53099651701497e-09\\
2.45	1.59347422829541e-09\\
2.475	1.64563127388131e-09\\
2.5	1.68443446600209e-09\\
2.525	1.71030492794703e-09\\
2.55	1.72612349177152e-09\\
2.575	1.73449784648474e-09\\
2.6	1.73746415728826e-09\\
2.625	1.73428956408675e-09\\
2.65	1.72819295575068e-09\\
2.675	1.72065455410055e-09\\
2.7	1.71207129441408e-09\\
2.725	1.71450546558317e-09\\
2.75	1.71403784198e-09\\
2.775	1.71139694478006e-09\\
2.8	1.69315156470591e-09\\
2.825	1.68570969757853e-09\\
2.85	1.66799286532411e-09\\
2.875	1.63736049076925e-09\\
2.9	1.62520243153595e-09\\
2.925	1.59833266062607e-09\\
2.95	1.56434489552127e-09\\
2.975	1.53383683614402e-09\\
3	1.50684946817779e-09\\
3.025	1.48551935516781e-09\\
3.05	1.46834689678288e-09\\
3.075	1.44258391284984e-09\\
3.1	1.42172203770708e-09\\
3.125	1.40362718789089e-09\\
3.15	1.38722555351103e-09\\
3.175	1.37242875638311e-09\\
3.2	1.35838522479966e-09\\
3.225	1.34398468921807e-09\\
3.25	1.33051262924166e-09\\
3.275	1.31813482712263e-09\\
3.3	1.30803560876717e-09\\
3.325	1.29990898242979e-09\\
3.35	1.28766538340914e-09\\
3.375	1.27654613917913e-09\\
3.4	1.27682035576638e-09\\
3.425	1.27068769016714e-09\\
3.45	1.26408646640328e-09\\
3.475	1.25824534536126e-09\\
3.5	1.24981052564561e-09\\
3.525	1.24141814440791e-09\\
3.55	1.23313208904829e-09\\
3.575	1.22561976184656e-09\\
3.6	1.22200984542248e-09\\
3.625	1.22916332765394e-09\\
3.65	1.22906178074373e-09\\
3.675	1.24578469161028e-09\\
3.7	1.23726335394396e-09\\
3.725	1.23559822886072e-09\\
3.75	1.27236111596271e-09\\
3.775	1.26843905895319e-09\\
3.8	1.29463903649923e-09\\
3.825	1.40604279503885e-09\\
3.85	1.68781220758251e-09\\
3.875	1.73917499422142e-09\\
3.9	1.66734727546659e-09\\
3.925	1.64522468301011e-09\\
3.95	1.63343549188108e-09\\
3.975	1.59696392890241e-09\\
4	1.56936510804212e-09\\
4.025	1.55600451344239e-09\\
4.05	1.53182822126786e-09\\
4.075	1.52199432493817e-09\\
4.1	1.50204011877815e-09\\
4.125	1.49847167636161e-09\\
4.15	1.48272599361841e-09\\
4.175	1.46922514934742e-09\\
4.2	1.46443621784971e-09\\
4.225	1.45761714747531e-09\\
4.25	1.43885292721039e-09\\
4.275	1.4133083470242e-09\\
4.3	1.37768958511355e-09\\
4.325	1.3282776447652e-09\\
4.35	1.26420656564626e-09\\
4.375	1.18583322219351e-09\\
4.4	1.10322284910503e-09\\
4.425	1.02550404287679e-09\\
4.45	9.60483717589209e-10\\
4.475	9.11307033897855e-10\\
4.5	8.75271574234791e-10\\
4.525	8.51346203395927e-10\\
4.55	8.34465452973682e-10\\
4.575	8.23380672476127e-10\\
4.6	8.17155293596423e-10\\
4.625	8.12553355240454e-10\\
4.65	8.14296255752783e-10\\
4.675	8.28834654964315e-10\\
4.7	8.4271275047264e-10\\
4.725	8.59704846684266e-10\\
4.75	8.79820768116937e-10\\
4.775	9.03229533751924e-10\\
4.8	9.2936777895428e-10\\
4.825	9.61933954871513e-10\\
4.85	9.80574548977144e-10\\
4.875	9.99657266661885e-10\\
4.9	1.01555265691074e-09\\
4.925	1.02692045228406e-09\\
4.95	1.03397076477243e-09\\
4.975	1.04198399852528e-09\\
5	1.04249288137018e-09\\
};
\addlegendentry{$\epsilon_{os} = 30$}

\end{axis}
\end{tikzpicture}%

%% file: figures/Error_Burgers.tex
\definecolor{mycolor1}{rgb}{0.00000,0.44700,0.74100}%
\definecolor{mycolor2}{rgb}{0.85000,0.32500,0.09800}%
\definecolor{mycolor3}{rgb}{0.5,0.5,0.5}
\begin{tikzpicture}[scale=0.8]

\begin{axis}[%
title={Burgers},
width=1.5in,
height=1.5in,
scale only axis,
xmin=0,
xmax=5,
xlabel style={font=\color{white!15!black}},
xlabel={Time},
ymode=log,
log basis y=10,
ymin=1e-15,
ymax=1e-05,
yminorticks=true,
ylabel style={font=\color{white!15!black}},
ylabel={$\mathcal{E}$},
axis background/.style={fill=white},
legend style={fill opacity=0.8, 
font=\tiny,
  draw opacity=1,
  text opacity=1,
  at={(0.04,0.22)},
  anchor=west}
  % draw=lightgray204}
]
\addplot [color=mycolor1, line width=1.2pt,y filter/.expression={log10(y)}]
  table[row sep=crcr]{%
0	6.25305022111612e-19\\
0.0125	4.74558749574291e-09\\
0.025	7.1404817448325e-09\\
0.0375	8.41832358688583e-09\\
0.05	1.18526221504213e-08\\
0.0625	1.57889334997273e-08\\
0.075	1.90894212245822e-08\\
0.0875	2.12776729105025e-08\\
0.1	2.24765472242708e-08\\
0.1125	2.33219776966992e-08\\
0.125	2.34998292364818e-08\\
0.1375	2.31200903333664e-08\\
0.15	2.25524667366513e-08\\
0.1625	2.20996021835405e-08\\
0.175	2.21133231127903e-08\\
0.1875	2.28009659292086e-08\\
0.2	2.44249343835678e-08\\
0.2125	2.75869858287993e-08\\
0.225	3.17331393812747e-08\\
0.2375	3.65864404356045e-08\\
0.25	4.27333759632771e-08\\
0.2625	5.10329307500994e-08\\
0.275	6.11315443975405e-08\\
0.2875	6.53302549922977e-08\\
0.3	6.60392310832461e-08\\
0.3125	6.63411521341859e-08\\
0.325	6.57891319649309e-08\\
0.3375	6.47350116700401e-08\\
0.35	6.40614564091693e-08\\
0.3625	6.329944816595e-08\\
0.375	6.23468490985603e-08\\
0.3875	6.12121742318163e-08\\
0.4	6.00911997309666e-08\\
0.4125	5.88272694810917e-08\\
0.425	5.74592572176885e-08\\
0.4375	5.62106401560078e-08\\
0.45	5.52283201129799e-08\\
0.4625	5.43653978802957e-08\\
0.475	5.32094259667982e-08\\
0.4875	5.22923366150077e-08\\
0.5	5.15684531418479e-08\\
0.5125	5.10761620190862e-08\\
0.525	4.99769370565053e-08\\
0.5375	4.87596245868138e-08\\
0.55	4.7747037587343e-08\\
0.5625	4.68334995946479e-08\\
0.575	4.88887371500047e-08\\
0.5875	5.60041347063583e-08\\
0.6	6.20729409726389e-08\\
0.6125	6.61130258760891e-08\\
0.625	6.89670607045497e-08\\
0.6375	7.10462729873247e-08\\
0.65	7.54135980411398e-08\\
0.6625	8.25150508454171e-08\\
0.675	8.9285259464561e-08\\
0.6875	9.50906142678833e-08\\
0.7	9.99021235608862e-08\\
0.7125	1.01364529077259e-07\\
0.725	1.02946757861101e-07\\
0.7375	1.04692961818009e-07\\
0.75	1.07052724279704e-07\\
0.7625	1.09967393411362e-07\\
0.775	1.15195191449742e-07\\
0.7875	1.21782514449105e-07\\
0.8	1.26308820809739e-07\\
0.8125	1.3153505623443e-07\\
0.825	1.37685235657236e-07\\
0.8375	1.44848409253245e-07\\
0.85	1.52369500083918e-07\\
0.8625	1.58420686502307e-07\\
0.875	1.64661448149989e-07\\
0.8875	1.72406062452603e-07\\
0.9	1.79460563872439e-07\\
0.9125	1.86725191633087e-07\\
0.925	1.93910574383413e-07\\
0.9375	2.00815566569262e-07\\
0.95	2.06888365303201e-07\\
0.9625	2.12158518380452e-07\\
0.975	2.16543314388177e-07\\
0.9875	2.20111142542406e-07\\
1	2.22235115274322e-07\\
1.0125	2.24571057440713e-07\\
1.025	2.26676282113186e-07\\
1.0375	2.27781812815746e-07\\
1.05	2.28866734470823e-07\\
1.0625	2.30091821378013e-07\\
1.075	2.31587604507675e-07\\
1.0875	2.3306595207053e-07\\
1.1	2.33247269650232e-07\\
1.1125	2.32891372831345e-07\\
1.125	2.32672728129724e-07\\
1.1375	2.32246796116514e-07\\
1.15	2.31711509861566e-07\\
1.1625	2.30846622297101e-07\\
1.175	2.29595018627383e-07\\
1.1875	2.27821498174314e-07\\
1.2	2.25874883304909e-07\\
1.2125	8.09376793044725e-07\\
1.225	6.58774011563893e-07\\
1.2375	5.71362095343447e-07\\
1.25	5.09735403085149e-07\\
1.2625	4.64075852309988e-07\\
1.275	4.29073460775176e-07\\
1.2875	4.00920241859788e-07\\
1.3	3.77616591600158e-07\\
1.3125	3.57940260937862e-07\\
1.325	3.41274659752358e-07\\
1.3375	3.26611178480304e-07\\
1.35	3.14101489238306e-07\\
1.3625	3.03241059439066e-07\\
1.375	2.93869724864658e-07\\
1.3875	2.85980919920606e-07\\
1.4	2.79151041435763e-07\\
1.4125	2.7329390744588e-07\\
1.425	2.69264717008497e-07\\
1.4375	2.66097071364189e-07\\
1.45	2.63288405742778e-07\\
1.4625	2.61634526645894e-07\\
1.475	2.6061598848134e-07\\
1.4875	2.60069610745846e-07\\
1.5	2.822891807427e-07\\
1.5125	2.77721603067227e-07\\
1.525	2.74127067879451e-07\\
1.5375	2.71212904164683e-07\\
1.55	2.68774035155565e-07\\
1.5625	2.67397275893719e-07\\
1.575	2.6536386558993e-07\\
1.5875	2.64119234013792e-07\\
1.6	2.63586094202621e-07\\
1.6125	2.63686338875899e-07\\
1.625	2.64546181865037e-07\\
1.6375	2.6622103778709e-07\\
1.65	2.68605479684959e-07\\
1.6625	2.71688694711687e-07\\
1.675	2.75397324930839e-07\\
1.6875	2.79682386248152e-07\\
1.7	2.84633506930605e-07\\
1.7125	2.90346017292994e-07\\
1.725	2.96753142458067e-07\\
1.7375	3.03899866680518e-07\\
1.75	3.11710591164935e-07\\
1.7625	3.20123618926722e-07\\
1.775	3.29792167737562e-07\\
1.7875	3.37946576483465e-07\\
1.8	3.46361756163784e-07\\
1.8125	3.54906940388776e-07\\
1.825	3.62719502781876e-07\\
1.8375	3.70448362520523e-07\\
1.85	3.77567927524694e-07\\
1.8625	3.84457637347008e-07\\
1.875	3.90895399705701e-07\\
1.8875	3.96959397443103e-07\\
1.9	4.02490622964663e-07\\
1.9125	4.07528882747138e-07\\
1.925	4.12035107168964e-07\\
1.9375	4.16021062375671e-07\\
1.95	4.19438977185547e-07\\
1.9625	4.2234068708096e-07\\
1.975	4.24926716261325e-07\\
1.9875	4.27072793611749e-07\\
2	4.2906644774016e-07\\
2.0125	4.30492331203121e-07\\
2.025	4.31474893780579e-07\\
2.0375	4.32246500379674e-07\\
2.05	4.32791361975202e-07\\
2.0625	4.33115674242671e-07\\
2.075	4.33214131777605e-07\\
2.0875	4.33155081222598e-07\\
2.1	4.32947366002036e-07\\
2.1125	4.32623303635115e-07\\
2.125	4.32231854003879e-07\\
2.1375	4.31590255952096e-07\\
2.15	4.30944854654344e-07\\
2.1625	4.30229259962832e-07\\
2.175	4.29473128412319e-07\\
2.1875	4.28683165468609e-07\\
2.2	4.27780362381799e-07\\
2.2125	4.2684960066506e-07\\
2.225	4.25915017867394e-07\\
2.2375	4.25070147432411e-07\\
2.25	4.24402259523679e-07\\
2.2625	4.23993380825911e-07\\
2.275	4.23942300574245e-07\\
2.2875	4.24419262573016e-07\\
2.3	4.25670198055287e-07\\
2.3125	4.27362185274487e-07\\
2.325	4.29887623314741e-07\\
2.3375	4.3327456404986e-07\\
2.35	4.37513699290399e-07\\
2.3625	4.42470593177339e-07\\
2.375	4.48135519516744e-07\\
2.3875	4.54324445250156e-07\\
2.4	4.60700541714798e-07\\
2.4125	4.67193021537198e-07\\
2.425	4.73575499022861e-07\\
2.4375	4.79573366402205e-07\\
2.45	4.85230699197387e-07\\
2.4625	4.90309827635448e-07\\
2.475	4.94830934863025e-07\\
2.4875	4.98660531691111e-07\\
2.5	5.01976856894092e-07\\
2.5125	5.04971119730733e-07\\
2.525	8.78234231562658e-07\\
2.5375	7.70714777369339e-07\\
2.55	6.98059362322544e-07\\
2.5625	6.48823027459876e-07\\
2.575	6.1710428280848e-07\\
2.5875	5.98557696103295e-07\\
2.6	5.84124364397108e-07\\
2.6125	5.73674066592562e-07\\
2.625	5.68505835859718e-07\\
2.6375	5.66886236039758e-07\\
2.65	5.67519155889569e-07\\
2.6625	5.70215088010948e-07\\
2.675	5.74980700197888e-07\\
2.6875	5.81855290968675e-07\\
2.7	5.90170620712168e-07\\
2.7125	5.99946587373183e-07\\
2.725	6.1212964502144e-07\\
2.7375	6.25673101865931e-07\\
2.75	6.39815454811545e-07\\
2.7625	6.53999088623055e-07\\
2.775	6.6777046905875e-07\\
2.7875	6.83706677556635e-07\\
2.8	7.00820782106852e-07\\
2.8125	7.19213530892352e-07\\
2.825	7.38639630582167e-07\\
2.8375	7.58363628882649e-07\\
2.85	7.78050723134366e-07\\
2.8625	7.97984623271425e-07\\
2.875	8.17917137570754e-07\\
2.8875	8.35162592555972e-07\\
2.9	8.50499981655053e-07\\
2.9125	8.65203986341821e-07\\
2.925	8.79015331569784e-07\\
2.9375	8.91335929595233e-07\\
2.95	9.01518439766003e-07\\
2.9625	9.10511531013867e-07\\
2.975	9.18613109300798e-07\\
2.9875	9.25854886727775e-07\\
3	9.32374078138853e-07\\
3.0125	9.39140691325436e-07\\
3.025	9.42610145988515e-07\\
3.0375	9.46502356124381e-07\\
3.05	9.50020484821914e-07\\
3.0625	9.527911776751e-07\\
3.075	9.60545706737554e-07\\
3.0875	9.66958309802767e-07\\
3.1	9.72192157616836e-07\\
3.1125	9.7594026513748e-07\\
3.125	9.7920242518262e-07\\
3.1375	9.81905546952716e-07\\
3.15	9.79023333896625e-07\\
3.1625	9.69449697469462e-07\\
3.175	9.63456986213922e-07\\
3.1875	9.59151204610991e-07\\
3.2	9.55694343955365e-07\\
3.2125	9.52704335783769e-07\\
3.225	9.49877713607636e-07\\
3.2375	9.47252211824749e-07\\
3.25	9.44788314642862e-07\\
3.2625	9.42236535007873e-07\\
3.275	9.39542289491932e-07\\
3.2875	9.36912954880272e-07\\
3.3	9.33801089604815e-07\\
3.3125	9.32990623278906e-07\\
3.325	9.33483880869426e-07\\
3.3375	9.30720495018022e-07\\
3.35	9.28336684613923e-07\\
3.3625	9.2602598742079e-07\\
3.375	9.24029757957857e-07\\
3.3875	9.22239596050348e-07\\
3.4	9.20612921921742e-07\\
3.4125	9.18277365748488e-07\\
3.425	9.16113423830708e-07\\
3.4375	9.18916521564571e-07\\
3.45	9.18110309371392e-07\\
3.4625	9.17794649687578e-07\\
3.475	9.16950711274734e-07\\
3.4875	9.16224358962299e-07\\
3.5	9.15533110025519e-07\\
3.5125	9.14632228082065e-07\\
3.525	9.13928527683594e-07\\
3.5375	9.13149320219973e-07\\
3.55	9.12709957430668e-07\\
3.5625	9.12475378299e-07\\
3.575	9.13876589865886e-07\\
3.5875	9.14618286844015e-07\\
3.6	9.16849679526867e-07\\
3.6125	9.15386192128464e-07\\
3.625	9.15619527316139e-07\\
3.6375	9.18633923005438e-07\\
3.65	9.24574291186677e-07\\
3.6625	9.32784252487697e-07\\
3.675	9.42331083499433e-07\\
3.6875	9.54298144462823e-07\\
3.7	9.67475842949e-07\\
3.7125	9.81758802915935e-07\\
3.725	9.97523319603998e-07\\
3.7375	1.01463016168501e-06\\
3.75	1.03341177466801e-06\\
3.7625	1.05455053052992e-06\\
3.775	1.07747929021437e-06\\
3.7875	1.1026150278517e-06\\
3.8	1.12953885158492e-06\\
3.8125	1.1582324507009e-06\\
3.825	1.18827870298694e-06\\
3.8375	1.21588728510724e-06\\
3.85	1.23927008986903e-06\\
3.8625	1.26234677408408e-06\\
3.875	1.28367666992715e-06\\
3.8875	1.30328249157864e-06\\
3.9	1.3215018345625e-06\\
3.9125	1.338314351944e-06\\
3.925	1.35367972569191e-06\\
3.9375	1.36794928002551e-06\\
3.95	1.38139481937758e-06\\
3.9625	1.39343047025277e-06\\
3.975	1.40407613323507e-06\\
3.9875	1.41342769743336e-06\\
4	1.4215684203541e-06\\
4.0125	1.4286002079552e-06\\
4.025	1.43461294674321e-06\\
4.0375	1.4396956960684e-06\\
4.05	1.4439340997591e-06\\
4.0625	1.44739786706095e-06\\
4.075	1.45018476989002e-06\\
4.0875	1.45189756506095e-06\\
4.1	1.45290018599359e-06\\
4.1125	1.45353076891355e-06\\
4.125	1.45365138277074e-06\\
4.1375	1.45333022597958e-06\\
4.15	1.45266039618728e-06\\
4.1625	1.45168564908512e-06\\
4.175	1.45040597573394e-06\\
4.1875	1.44902269923191e-06\\
4.2	1.44694720444125e-06\\
4.2125	1.44628632144476e-06\\
4.225	1.44502436141595e-06\\
4.2375	1.44423738226133e-06\\
4.25	1.44331380174065e-06\\
4.2625	1.44278618275948e-06\\
4.275	1.44217117108046e-06\\
4.2875	1.44164842718684e-06\\
4.3	1.43959525767976e-06\\
4.3125	1.43252558100386e-06\\
4.325	1.42595128450285e-06\\
4.3375	1.42073747047293e-06\\
4.35	1.41427710815592e-06\\
4.3625	1.40831307359611e-06\\
4.375	1.40340277867814e-06\\
4.3875	1.39762971778891e-06\\
4.4	1.39326063716526e-06\\
4.4125	1.38180569039731e-06\\
4.425	1.37632329555195e-06\\
4.4375	1.36747332029995e-06\\
4.45	1.35797436880403e-06\\
4.4625	1.34697365859694e-06\\
4.475	1.33498614266004e-06\\
4.4875	1.62365803932431e-06\\
4.5	1.56190209622255e-06\\
4.5125	1.51799919397516e-06\\
4.525	1.46073827998577e-06\\
4.5375	1.41563728113715e-06\\
4.55	1.36919690624627e-06\\
4.5625	1.31354329360397e-06\\
4.575	1.26717458388515e-06\\
4.5875	1.23207222490884e-06\\
4.6	1.20731994974707e-06\\
4.6125	1.19041583757321e-06\\
4.625	1.18075614448847e-06\\
4.6375	1.17671972316637e-06\\
4.65	1.17668832649811e-06\\
4.6625	1.17970526048777e-06\\
4.675	1.18537865494155e-06\\
4.6875	1.19248001647355e-06\\
4.7	1.20085606407702e-06\\
4.7125	1.21005595744974e-06\\
4.725	1.21993579276601e-06\\
4.7375	1.23097185082075e-06\\
4.75	1.24308464237565e-06\\
4.7625	1.25596220620917e-06\\
4.775	1.26983128340074e-06\\
4.7875	1.28417490177393e-06\\
4.8	1.29766525731731e-06\\
4.8125	1.3110399292084e-06\\
4.825	1.32422422548447e-06\\
4.8375	1.33696357967183e-06\\
4.85	1.34878333330876e-06\\
4.8625	1.36084362052017e-06\\
4.875	1.37158109225633e-06\\
4.8875	1.3816077361179e-06\\
4.9	1.39101965293266e-06\\
4.9125	1.40013181957914e-06\\
4.925	1.40868985225342e-06\\
4.9375	1.41633460287212e-06\\
4.95	1.42271051834148e-06\\
4.9625	1.42834749121244e-06\\
4.975	1.43302122756703e-06\\
4.9875	1.43656112257392e-06\\
5	1.43912657019324e-06\\
};
\addlegendentry{\texttt{TDB-CUR-CR-OS} $|$ $\epsilon_u = 10^{-8}$ }

\addplot [color=mycolor2, line width=1.2pt,y filter/.expression={log10(y)}]
  table[row sep=crcr]{%
0	6.25305022111612e-19\\
0.0125	4.57396484974063e-09\\
0.025	5.40957073080839e-09\\
0.0375	7.0152195366919e-09\\
0.05	8.98680047527615e-09\\
0.0625	1.07657161429857e-08\\
0.075	1.20492038806194e-08\\
0.0875	1.28086885099942e-08\\
0.1	1.31177288263468e-08\\
0.1125	1.30698967722278e-08\\
0.125	1.27602384314494e-08\\
0.1375	1.22767174726221e-08\\
0.15	1.16988290274225e-08\\
0.1625	1.11131063366971e-08\\
0.175	1.06265497682019e-08\\
0.1875	1.03451822777607e-08\\
0.2	1.03241649397481e-08\\
0.2125	1.05474430302258e-08\\
0.225	1.09675274061253e-08\\
0.2375	1.15555991787762e-08\\
0.25	1.23086117111691e-08\\
0.2625	1.32250163988259e-08\\
0.275	1.42902159928152e-08\\
0.2875	1.421351415492e-08\\
0.3	1.35374647898614e-08\\
0.3125	1.30396423300771e-08\\
0.325	1.26446703298485e-08\\
0.3375	1.23105913611871e-08\\
0.35	1.20194513554798e-08\\
0.3625	1.17685880886914e-08\\
0.375	1.1561278170078e-08\\
0.3875	1.14012410058994e-08\\
0.4	1.12930506762087e-08\\
0.4125	1.12450580194504e-08\\
0.425	1.12700984911966e-08\\
0.4375	1.13830942442702e-08\\
0.45	1.15982705280768e-08\\
0.4625	1.19280554649971e-08\\
0.475	1.2383449382253e-08\\
0.4875	1.29756828101742e-08\\
0.5	1.37193503846784e-08\\
0.5125	1.46344074954115e-08\\
0.525	1.57422998509876e-08\\
0.5375	1.70569996003278e-08\\
0.55	1.85802228215856e-08\\
0.5625	2.03066630760877e-08\\
0.575	2.57532074060301e-08\\
0.5875	2.93199473027604e-08\\
0.6	3.36244780396673e-08\\
0.6125	3.86102752577227e-08\\
0.625	4.4164802337937e-08\\
0.6375	5.01430355793386e-08\\
0.65	5.63930518321573e-08\\
0.6625	6.27689913263082e-08\\
0.675	6.91364140723322e-08\\
0.6875	7.53798164750261e-08\\
0.7	8.14126750633997e-08\\
0.7125	8.73309751584386e-08\\
0.725	9.3231669568397e-08\\
0.7375	9.89904719967285e-08\\
0.75	1.04564549936953e-07\\
0.7625	1.09993829391041e-07\\
0.775	1.15384375195463e-07\\
0.7875	1.20846863182847e-07\\
0.8	1.26439159526797e-07\\
0.8125	1.32377384969003e-07\\
0.825	1.38890962159882e-07\\
0.8375	1.46239828664546e-07\\
0.85	1.55035331914782e-07\\
0.8625	1.64239974549854e-07\\
0.875	1.73561777606197e-07\\
0.8875	1.82176978759543e-07\\
0.9	1.90587057217649e-07\\
0.9125	1.98507162180798e-07\\
0.925	2.05524033566255e-07\\
0.9375	2.11586104643542e-07\\
0.95	2.16742426396213e-07\\
0.9625	2.21076020053735e-07\\
0.975	2.2468202950049e-07\\
0.9875	2.27811925216181e-07\\
1	2.30380048934613e-07\\
1.0125	2.32712322600105e-07\\
1.025	2.3461035360636e-07\\
1.0375	2.36085940380319e-07\\
1.05	2.37158584194355e-07\\
1.0625	2.37858294898673e-07\\
1.075	2.38222241607946e-07\\
1.0875	2.38295934661385e-07\\
1.1	2.38215536948931e-07\\
1.1125	2.37980049548629e-07\\
1.125	2.37588369518424e-07\\
1.1375	2.3705392527422e-07\\
1.15	2.36393798802372e-07\\
1.1625	2.35625550759101e-07\\
1.175	2.34765848186978e-07\\
1.1875	2.3370223579115e-07\\
1.2	2.32597987786089e-07\\
1.2125	2.31465171041977e-07\\
1.225	2.30316969628012e-07\\
1.2375	2.29151635254794e-07\\
1.25	2.27957350618231e-07\\
1.2625	2.26748537828842e-07\\
1.275	2.25535236629674e-07\\
1.2875	2.24288412236702e-07\\
1.3	2.23040244937955e-07\\
1.3125	2.21792526600911e-07\\
1.325	2.20547394113382e-07\\
1.3375	2.19768186080982e-07\\
1.35	2.18704088428588e-07\\
1.3625	2.17656790849915e-07\\
1.375	2.16618008313827e-07\\
1.3875	2.15586147222143e-07\\
1.4	2.1456111649578e-07\\
1.4125	2.13543282451248e-07\\
1.425	2.12448970404476e-07\\
1.4375	2.11393631108235e-07\\
1.45	2.10372947070221e-07\\
1.4625	2.0939516464278e-07\\
1.475	2.08471026473202e-07\\
1.4875	2.07614592895094e-07\\
1.5	2.06821689740204e-07\\
1.5125	2.06053787940224e-07\\
1.525	2.05366599524546e-07\\
1.5375	2.04792760648891e-07\\
1.55	2.04370380477799e-07\\
1.5625	2.04396615884888e-07\\
1.575	2.04694848578163e-07\\
1.5875	2.05418378229329e-07\\
1.6	2.06639284944203e-07\\
1.6125	2.08390578157559e-07\\
1.625	2.10659372879775e-07\\
1.6375	2.13389397472597e-07\\
1.65	2.16498463028419e-07\\
1.6625	2.19900724984273e-07\\
1.675	2.23522854558193e-07\\
1.6875	2.27310789042448e-07\\
1.7	2.31229565584003e-07\\
1.7125	2.35259494896801e-07\\
1.725	2.39391121492374e-07\\
1.7375	2.43620492897492e-07\\
1.75	2.47944079089479e-07\\
1.7625	2.52351312548963e-07\\
1.775	2.56815760308841e-07\\
1.7875	2.61253339688353e-07\\
1.8	2.65616768517748e-07\\
1.8125	2.69846551597564e-07\\
1.825	2.7386921978907e-07\\
1.8375	2.77606926768358e-07\\
1.85	2.8085629154997e-07\\
1.8625	2.83578508989894e-07\\
1.875	2.85769878137186e-07\\
1.8875	2.87447024751089e-07\\
1.9	2.88643826450379e-07\\
1.9125	2.89402983936246e-07\\
1.925	2.89770093447344e-07\\
1.9375	2.89791039477466e-07\\
1.95	2.8951128255459e-07\\
1.9625	2.88979242778339e-07\\
1.975	2.88255351303203e-07\\
1.9875	2.87510651180549e-07\\
2	2.86649876281594e-07\\
2.0125	2.86267696065675e-07\\
2.025	2.85406794929514e-07\\
2.0375	2.84542659269817e-07\\
2.05	2.83673168772301e-07\\
2.0625	2.82800022707605e-07\\
2.075	2.81925518957049e-07\\
2.0875	2.81051609425753e-07\\
2.1	2.80180148453086e-07\\
2.1125	2.79313869408319e-07\\
2.125	2.78396518948794e-07\\
2.1375	2.77598411171128e-07\\
2.15	2.76826850366516e-07\\
2.1625	2.76100358483057e-07\\
2.175	2.7544620342862e-07\\
2.1875	2.7490213883878e-07\\
2.2	2.74518774190852e-07\\
2.2125	2.74362584104112e-07\\
2.225	2.74518960112328e-07\\
2.2375	2.75094105959373e-07\\
2.25	2.76214316402635e-07\\
2.2625	2.78021184313731e-07\\
2.275	2.80661549969336e-07\\
2.2875	2.84271737046174e-07\\
2.3	2.89325773624224e-07\\
2.3125	2.95154956202098e-07\\
2.325	3.02093622832707e-07\\
2.3375	3.10020622240115e-07\\
2.35	3.18741061367328e-07\\
2.3625	3.28003642719442e-07\\
2.375	3.37526274574812e-07\\
2.3875	3.47024980819274e-07\\
2.4	3.5624129017802e-07\\
2.4125	3.6496280530941e-07\\
2.425	3.73031910304667e-07\\
2.4375	3.80343090429872e-07\\
2.45	3.8683699847877e-07\\
2.4625	3.92498741144003e-07\\
2.475	3.9735055063017e-07\\
2.4875	4.01444139088597e-07\\
2.5	4.04849400934072e-07\\
2.5125	4.07603461499746e-07\\
2.525	4.09604602255582e-07\\
2.5375	4.11195288251504e-07\\
2.55	4.12446696329277e-07\\
2.5625	4.13429507853473e-07\\
2.575	4.14212832787125e-07\\
2.5875	4.14861844467532e-07\\
2.6	4.15432928536919e-07\\
2.6125	4.15967780757328e-07\\
2.625	4.16487407243834e-07\\
2.6375	4.17015942852523e-07\\
2.65	4.17561106098832e-07\\
2.6625	4.18080567975907e-07\\
2.675	4.18588189461113e-07\\
2.6875	4.19089455189606e-07\\
2.7	4.21291475048208e-07\\
2.7125	4.22222192947918e-07\\
2.725	4.2321217072221e-07\\
2.7375	4.24243402776284e-07\\
2.75	4.25323311588895e-07\\
2.7625	4.2647805634686e-07\\
2.775	4.27739812096568e-07\\
2.7875	4.29120620879148e-07\\
2.8	4.30602532569218e-07\\
2.8125	4.32145491607327e-07\\
2.825	4.33692152081563e-07\\
2.8375	4.35438574846361e-07\\
2.85	4.37014424857201e-07\\
2.8625	4.38545301949231e-07\\
2.875	4.40001097336476e-07\\
2.8875	4.41343064454894e-07\\
2.9	4.42540444667099e-07\\
2.9125	4.43519916943176e-07\\
2.925	4.44264843342844e-07\\
2.9375	4.44760771767054e-07\\
2.95	4.45031382413467e-07\\
2.9625	4.45137579601591e-07\\
2.975	4.45069665784831e-07\\
2.9875	4.44834727399203e-07\\
3	4.44448516492538e-07\\
3.0125	4.44386577370026e-07\\
3.025	4.4757849297593e-07\\
3.0375	4.47050698849695e-07\\
3.05	4.46481332072198e-07\\
3.0625	4.45855022877605e-07\\
3.075	4.55016506188333e-07\\
3.0875	4.54312318167345e-07\\
3.1	4.53552782645853e-07\\
3.1125	4.52738940919038e-07\\
3.125	4.51875101033597e-07\\
3.1375	4.5096772960684e-07\\
3.15	4.48910691386142e-07\\
3.1625	4.45061170105331e-07\\
3.175	4.42184748492155e-07\\
3.1875	4.39848020184664e-07\\
3.2	4.37850269467572e-07\\
3.2125	4.36069433235462e-07\\
3.225	4.34438563594205e-07\\
3.2375	4.32917012318576e-07\\
3.25	4.31471921403827e-07\\
3.2625	4.3009276610371e-07\\
3.275	4.2874688130546e-07\\
3.2875	4.2743831955141e-07\\
3.3	4.26160364396537e-07\\
3.3125	4.24904335927303e-07\\
3.325	4.23669563907031e-07\\
3.3375	4.22414216567701e-07\\
3.35	4.21259026916321e-07\\
3.3625	4.200770953226e-07\\
3.375	4.18935210907737e-07\\
3.3875	4.17798563209535e-07\\
3.4	4.18304968932132e-07\\
3.4125	4.1723323072291e-07\\
3.425	4.16180325356157e-07\\
3.4375	4.24155235412264e-07\\
3.45	4.23191192491723e-07\\
3.4625	4.22244526400837e-07\\
3.475	4.21316029161141e-07\\
3.4875	4.20409001868371e-07\\
3.5	4.19529437918874e-07\\
3.5125	4.18686253070908e-07\\
3.525	4.17891876055392e-07\\
3.5375	4.17162815079899e-07\\
3.55	4.16519561368528e-07\\
3.5625	4.15985359059335e-07\\
3.575	4.15583452686123e-07\\
3.5875	4.15332537394098e-07\\
3.6	4.15241073953187e-07\\
3.6125	4.13609458831526e-07\\
3.625	4.12374655528444e-07\\
3.6375	4.11639358092408e-07\\
3.65	4.11281899285743e-07\\
3.6625	4.11219823289554e-07\\
3.675	4.11453670191018e-07\\
3.6875	4.119590543628e-07\\
3.7	4.12689609724182e-07\\
3.7125	4.13647279940776e-07\\
3.725	4.14818500923685e-07\\
3.7375	4.16169637824912e-07\\
3.75	4.17675095869512e-07\\
3.7625	4.19560900052114e-07\\
3.775	4.21365268257777e-07\\
3.7875	4.23321550694144e-07\\
3.8	4.25411811431394e-07\\
3.8125	4.2760635402435e-07\\
3.825	4.29866554780149e-07\\
3.8375	4.32144767647771e-07\\
3.85	4.34383668064204e-07\\
3.8625	4.38206690809208e-07\\
3.875	4.40017401614012e-07\\
3.8875	4.41623600277906e-07\\
3.9	4.42993779396893e-07\\
3.9125	4.44076057354658e-07\\
3.925	4.44849656834475e-07\\
3.9375	4.45877669192413e-07\\
3.95	4.46130718466344e-07\\
3.9625	4.46141455595527e-07\\
3.975	4.45966081159358e-07\\
3.9875	4.45669740928711e-07\\
4	4.45316061301642e-07\\
4.0125	4.44952174205091e-07\\
4.025	4.44599626664984e-07\\
4.0375	4.44258485391122e-07\\
4.05	4.43920366788871e-07\\
4.0625	4.43579587183777e-07\\
4.075	4.43236494443046e-07\\
4.0875	4.42744928591048e-07\\
4.1	4.42262879668705e-07\\
4.1125	4.41829755669622e-07\\
4.125	4.414554214367e-07\\
4.1375	4.41121786383974e-07\\
4.15	4.40832463828702e-07\\
4.1625	4.40582567510977e-07\\
4.175	4.40995194677049e-07\\
4.1875	4.40736626630304e-07\\
4.2	4.40637185467e-07\\
4.2125	4.47064399706302e-07\\
4.225	4.46269156119848e-07\\
4.2375	4.45226936313378e-07\\
4.25	4.43849708099327e-07\\
4.2625	4.42040554118588e-07\\
4.275	4.39699613729594e-07\\
4.2875	4.36733786484464e-07\\
4.3	4.32529163767926e-07\\
4.3125	4.27052172084025e-07\\
4.325	4.21287882517193e-07\\
4.3375	4.15097582560444e-07\\
4.35	4.08442406697829e-07\\
4.3625	4.0197700280544e-07\\
4.375	3.94774364310322e-07\\
4.3875	3.88461217518489e-07\\
4.4	3.81443379037176e-07\\
4.4125	3.52321391007274e-07\\
4.425	3.4217370308792e-07\\
4.4375	3.32033320545024e-07\\
4.45	3.21999108803966e-07\\
4.4625	3.12169766888386e-07\\
4.475	3.01669755674372e-07\\
4.4875	1.15213501286718e-06\\
4.5	1.06152489717725e-06\\
4.5125	9.93478701007781e-07\\
4.525	9.0946249061018e-07\\
4.5375	8.41963634219937e-07\\
4.55	7.70404555638985e-07\\
4.5625	6.95754317092175e-07\\
4.575	6.33379887513555e-07\\
4.5875	5.84678293662208e-07\\
4.6	5.48663597293232e-07\\
4.6125	5.23325253582455e-07\\
4.625	5.06201705166643e-07\\
4.6375	4.94935615418636e-07\\
4.65	4.8759046721506e-07\\
4.6625	4.82763954675726e-07\\
4.675	4.79510310665424e-07\\
4.6875	4.77229914647568e-07\\
4.7	4.75559336879248e-07\\
4.7125	4.74297734444275e-07\\
4.725	4.73417562272847e-07\\
4.7375	4.7281857032571e-07\\
4.75	4.72428870853415e-07\\
4.7625	4.72224832282258e-07\\
4.775	4.72211704137533e-07\\
4.7875	4.72401861333309e-07\\
4.8	4.72799261244036e-07\\
4.8125	4.73392253634387e-07\\
4.825	4.74152172104501e-07\\
4.8375	4.74976428754087e-07\\
4.85	4.75810674191837e-07\\
4.8625	4.76605041034665e-07\\
4.875	4.77297992773442e-07\\
4.8875	4.77834791042005e-07\\
4.9	4.7816799603727e-07\\
4.9125	4.78258653812022e-07\\
4.925	4.78079046284053e-07\\
4.9375	4.77614744649522e-07\\
4.95	4.76864323136911e-07\\
4.9625	4.75839030511556e-07\\
4.975	4.74564721823815e-07\\
4.9875	4.73082206474227e-07\\
5	4.71439265481534e-07\\
};
\addlegendentry{\texttt{TDB-SVD} $|$ $\epsilon_u = 10^{-8}$}

\addplot [color=mycolor1, dotted, line width=1.2pt,y filter/.expression={log10(y)}]
  table[row sep=crcr]{%
0	6.25305022111612e-19\\
0.0125	5.18747011939031e-12\\
0.025	5.70068007147524e-12\\
0.0375	5.87556444215407e-12\\
0.05	5.86652495580232e-12\\
0.0625	5.71581258500695e-12\\
0.075	5.72744526776692e-12\\
0.0875	5.63255809535912e-12\\
0.1	5.65800267209923e-12\\
0.1125	5.97987822696888e-12\\
0.125	7.23276594947053e-12\\
0.1375	8.60195910418864e-12\\
0.15	1.07628154359891e-11\\
0.1625	1.47561130333986e-11\\
0.175	2.04799251995241e-11\\
0.1875	2.50607485560757e-11\\
0.2	2.85942822252962e-11\\
0.2125	3.2325043192954e-11\\
0.225	3.44653840630988e-11\\
0.2375	3.46415716719e-11\\
0.25	3.55160803572537e-11\\
0.2625	3.85336425694492e-11\\
0.275	4.47544507321231e-11\\
0.2875	4.63528537351341e-11\\
0.3	4.87595620796389e-11\\
0.3125	4.97697204103345e-11\\
0.325	5.33714475097115e-11\\
0.3375	5.68130368541525e-11\\
0.35	6.3270518524848e-11\\
0.3625	6.81892993912627e-11\\
0.375	7.08603359120889e-11\\
0.3875	7.53565728073004e-11\\
0.4	7.74968493557484e-11\\
0.4125	7.95807815056109e-11\\
0.425	8.52025853811024e-11\\
0.4375	9.34037779937754e-11\\
0.45	1.04302257327811e-10\\
0.4625	1.17550116676037e-10\\
0.475	1.21898405821122e-10\\
0.4875	1.35564554781884e-10\\
0.5	1.31080015627624e-10\\
0.5125	1.31471472261558e-10\\
0.525	1.368883541705e-10\\
0.5375	1.40133681644619e-10\\
0.55	1.44919628061274e-10\\
0.5625	1.48631324606984e-10\\
0.575	1.52898695311794e-10\\
0.5875	1.60234567809485e-10\\
0.6	1.74238986941928e-10\\
0.6125	1.7637425646831e-10\\
0.625	1.80863056121284e-10\\
0.6375	1.84469793682843e-10\\
0.65	1.85957801750338e-10\\
0.6625	2.01396286178787e-10\\
0.675	2.10906362337127e-10\\
0.6875	2.09613530077846e-10\\
0.7	2.09932287954833e-10\\
0.7125	2.08559079756448e-10\\
0.725	2.09278827118527e-10\\
0.7375	2.09625574922357e-10\\
0.75	2.10576301926434e-10\\
0.7625	2.11982465421739e-10\\
0.775	2.13880682817342e-10\\
0.7875	2.16184546420881e-10\\
0.8	2.17639574313429e-10\\
0.8125	2.18752192984288e-10\\
0.825	2.20057057033505e-10\\
0.8375	2.21713035920322e-10\\
0.85	2.24230721697219e-10\\
0.8625	2.27908164282168e-10\\
0.875	2.29375039963071e-10\\
0.8875	2.31663752706657e-10\\
0.9	2.34267440097961e-10\\
0.9125	2.36631652379683e-10\\
0.925	2.38690890378786e-10\\
0.9375	2.40930787045474e-10\\
0.95	2.42809255874723e-10\\
0.9625	2.44623468774208e-10\\
0.975	2.45652656192112e-10\\
0.9875	2.46949697677364e-10\\
1	2.48493067219326e-10\\
1.0125	2.49855285668473e-10\\
1.025	2.51385971305823e-10\\
1.0375	2.52947626760915e-10\\
1.05	2.54337983359838e-10\\
1.0625	2.55828225322944e-10\\
1.075	2.57244887264415e-10\\
1.0875	2.58450981522302e-10\\
1.1	2.59486699254454e-10\\
1.1125	2.60650447608704e-10\\
1.125	2.61777518024278e-10\\
1.1375	2.62982721891622e-10\\
1.15	2.64110275874936e-10\\
1.1625	2.65117937323821e-10\\
1.175	2.6607773520256e-10\\
1.1875	2.66898328084365e-10\\
1.2	2.6763911622453e-10\\
1.2125	2.68312291937346e-10\\
1.225	2.68907090865635e-10\\
1.2375	2.69448304906909e-10\\
1.25	2.69937577439056e-10\\
1.2625	2.70350890612467e-10\\
1.275	2.70709827615938e-10\\
1.2875	2.71019547232599e-10\\
1.3	2.71307410367583e-10\\
1.3125	2.71551972941526e-10\\
1.325	2.7175710858247e-10\\
1.3375	2.71918556379381e-10\\
1.35	2.7203222146735e-10\\
1.3625	2.72107304283966e-10\\
1.375	2.72148651137472e-10\\
1.3875	2.72171432416827e-10\\
1.4	2.72166782438681e-10\\
1.4125	2.72061668388389e-10\\
1.425	2.71956317371674e-10\\
1.4375	2.71787051354702e-10\\
1.45	2.71605802448338e-10\\
1.4625	2.71415808875704e-10\\
1.475	2.71149338946305e-10\\
1.4875	2.70836183504167e-10\\
1.5	2.70482644723755e-10\\
1.5125	2.70102497247977e-10\\
1.525	2.69688124786775e-10\\
1.5375	2.6924982538682e-10\\
1.55	2.68786592258622e-10\\
1.5625	2.68027259958223e-10\\
1.575	2.67540111498072e-10\\
1.5875	2.67058717487846e-10\\
1.6	2.66489143083613e-10\\
1.6125	2.65933185218162e-10\\
1.625	2.65368074847672e-10\\
1.6375	2.64855863580705e-10\\
1.65	2.64394564659065e-10\\
1.6625	2.63775131643229e-10\\
1.675	2.63055547746794e-10\\
1.6875	2.62336469546309e-10\\
1.7	2.61655474193573e-10\\
1.7125	2.61036079676564e-10\\
1.725	2.6038302755477e-10\\
1.7375	2.59798696006608e-10\\
1.75	2.59162076037906e-10\\
1.7625	2.583745137048e-10\\
1.775	2.57553966138246e-10\\
1.7875	2.56781329610181e-10\\
1.8	2.55973002487706e-10\\
1.8125	2.54907853809952e-10\\
1.825	2.53917589313153e-10\\
1.8375	2.5300145673291e-10\\
1.85	2.51928263122116e-10\\
1.8625	2.50747799443927e-10\\
1.875	2.49600313402427e-10\\
1.8875	2.484270227505e-10\\
1.9	2.47172312058068e-10\\
1.9125	2.45850823607741e-10\\
1.925	2.4445470956557e-10\\
1.9375	2.43037597280182e-10\\
1.95	2.41619174730326e-10\\
1.9625	2.40160593666702e-10\\
1.975	2.38725679437693e-10\\
1.9875	2.37247672522351e-10\\
2	2.35722217908703e-10\\
2.0125	2.34254647779162e-10\\
2.025	2.32762106034533e-10\\
2.0375	2.31046376344705e-10\\
2.05	2.29295883566355e-10\\
2.0625	2.27677653529377e-10\\
2.075	2.26040726069235e-10\\
2.0875	2.24412350854692e-10\\
2.1	2.22905349790869e-10\\
2.1125	2.21334532878987e-10\\
2.125	2.19814709997665e-10\\
2.1375	2.18250754066781e-10\\
2.15	2.16765452259585e-10\\
2.1625	2.15217796864406e-10\\
2.175	2.13825423961976e-10\\
2.1875	2.1245668621035e-10\\
2.2	2.11177743775403e-10\\
2.2125	2.10063780308509e-10\\
2.225	2.09104291870387e-10\\
2.2375	2.08379756167683e-10\\
2.25	2.07973799980608e-10\\
2.2625	2.07946383726093e-10\\
2.275	2.08440341360477e-10\\
2.2875	2.09566544895454e-10\\
2.3	2.11392407879029e-10\\
2.3125	2.13998743770847e-10\\
2.325	2.17417545354142e-10\\
2.3375	2.216504854374e-10\\
2.35	2.26635652435283e-10\\
2.3625	2.3416864214972e-10\\
2.375	2.39679307635293e-10\\
2.3875	2.45634832831759e-10\\
2.4	2.51644383669402e-10\\
2.4125	2.57501363710677e-10\\
2.425	2.63084439117607e-10\\
2.4375	2.68279302558075e-10\\
2.45	2.72917828346021e-10\\
2.4625	2.76983935090123e-10\\
2.475	2.80514195800514e-10\\
2.4875	2.83514842953023e-10\\
2.5	2.86016084284009e-10\\
2.5125	2.88049486328902e-10\\
2.525	2.89645973726358e-10\\
2.5375	2.90846791461586e-10\\
2.55	2.91739754668666e-10\\
2.5625	2.92441581522438e-10\\
2.575	2.92799456177061e-10\\
2.5875	2.92991135460134e-10\\
2.6	2.9303842334659e-10\\
2.6125	2.92918551572642e-10\\
2.625	2.92691482353715e-10\\
2.6375	2.92442689740295e-10\\
2.65	2.92257674256104e-10\\
2.6625	2.92016904917734e-10\\
2.675	2.91771731448756e-10\\
2.6875	2.91521260952926e-10\\
2.7	2.91354717474079e-10\\
2.7125	2.91101640110612e-10\\
2.725	2.90986909604512e-10\\
2.7375	2.90998547716747e-10\\
2.75	2.91018509364583e-10\\
2.7625	2.91074881700534e-10\\
2.775	2.91175675480982e-10\\
2.7875	2.91335123011023e-10\\
2.8	2.91558457574942e-10\\
2.8125	2.9190616375709e-10\\
2.825	2.91960895031002e-10\\
2.8375	2.92039761681516e-10\\
2.85	2.92057253042438e-10\\
2.8625	2.92175630488038e-10\\
2.875	2.92311266053828e-10\\
2.8875	2.92285853295425e-10\\
2.9	2.92217435384331e-10\\
2.9125	2.92030873817259e-10\\
2.925	2.91818699719105e-10\\
2.9375	2.91539088932939e-10\\
2.95	2.91206898024746e-10\\
2.9625	2.90779965905101e-10\\
2.975	2.9029380040318e-10\\
2.9875	2.89800950279614e-10\\
3	2.89423802821507e-10\\
3.0125	2.88815461048243e-10\\
3.025	2.88217804849579e-10\\
3.0375	2.8777724350888e-10\\
3.05	2.87179120468539e-10\\
3.0625	2.86484375069879e-10\\
3.075	2.85738330987562e-10\\
3.0875	2.84952982862084e-10\\
3.1	2.84251419597824e-10\\
3.1125	2.835394753939e-10\\
3.125	2.82757884646394e-10\\
3.1375	2.81987635094272e-10\\
3.15	2.81203191624648e-10\\
3.1625	2.80401018575522e-10\\
3.175	2.79574802144186e-10\\
3.1875	2.78727433754862e-10\\
3.2	2.77872236611714e-10\\
3.2125	2.77019464212848e-10\\
3.225	2.76151825176867e-10\\
3.2375	2.75283547025301e-10\\
3.25	2.74414493811533e-10\\
3.2625	2.7354649854263e-10\\
3.275	2.72699585296093e-10\\
3.2875	2.71856771469734e-10\\
3.3	2.71025618770081e-10\\
3.3125	2.70209375849673e-10\\
3.325	2.69393479323635e-10\\
3.3375	2.68603472441278e-10\\
3.35	2.6792438857889e-10\\
3.3625	2.67130779264406e-10\\
3.375	2.66311924389543e-10\\
3.3875	2.65506632471619e-10\\
3.4	2.6467934146381e-10\\
3.4125	2.63861986711479e-10\\
3.425	2.63051690055224e-10\\
3.4375	2.62261687119036e-10\\
3.45	2.61502563350264e-10\\
3.4625	2.60758603410562e-10\\
3.475	2.59992466038974e-10\\
3.4875	2.59298784473107e-10\\
3.5	2.58630660590339e-10\\
3.5125	2.57954134085492e-10\\
3.525	2.57258802360719e-10\\
3.5375	2.56712996737478e-10\\
3.55	2.56039173886982e-10\\
3.5625	2.55415744339021e-10\\
3.575	2.54851060999934e-10\\
3.5875	2.54344888657659e-10\\
3.6	2.53861559037224e-10\\
3.6125	2.53458790663202e-10\\
3.625	2.53080616031612e-10\\
3.6375	2.52795184734046e-10\\
3.65	2.52640274331848e-10\\
3.6625	2.52618868624358e-10\\
3.675	2.52676254407087e-10\\
3.6875	2.52669182802472e-10\\
3.7	2.52816993296183e-10\\
3.7125	2.53062957811836e-10\\
3.725	2.5339435048446e-10\\
3.7375	2.53840250698261e-10\\
3.75	2.54505572601616e-10\\
3.7625	2.55384235837122e-10\\
3.775	2.56082017622366e-10\\
3.7875	2.56528683493881e-10\\
3.8	2.57357909913733e-10\\
3.8125	2.58245234034301e-10\\
3.825	2.58915245763382e-10\\
3.8375	2.59670217171487e-10\\
3.85	2.60394093022469e-10\\
3.8625	2.61117615059641e-10\\
3.875	2.61825314097008e-10\\
3.8875	2.62482718847196e-10\\
3.9	2.63078668989832e-10\\
3.9125	2.63622485874282e-10\\
3.925	2.64093476778774e-10\\
3.9375	2.64506433594757e-10\\
3.95	2.64879497580266e-10\\
3.9625	2.65206024100886e-10\\
3.975	2.65470201645479e-10\\
3.9875	2.65680788799635e-10\\
4	2.65882526731805e-10\\
4.0125	2.65975364487138e-10\\
4.025	2.66033032221146e-10\\
4.0375	2.66070029227393e-10\\
4.05	2.660905239291e-10\\
4.0625	2.66082411725656e-10\\
4.075	2.66078480835771e-10\\
4.0875	2.66087875263285e-10\\
4.1	2.66118708767457e-10\\
4.1125	2.66137682522884e-10\\
4.125	2.66141443575829e-10\\
4.1375	2.66146023706483e-10\\
4.15	2.66153059253673e-10\\
4.1625	2.66164587478058e-10\\
4.175	2.66150606971473e-10\\
4.1875	2.66070076829711e-10\\
4.2	2.6594277705322e-10\\
4.2125	2.65727284038337e-10\\
4.225	2.65343092285146e-10\\
4.2375	2.64737691461847e-10\\
4.25	2.63871230145591e-10\\
4.2625	2.62639893602557e-10\\
4.275	2.6099202409162e-10\\
4.2875	2.58808233804332e-10\\
4.3	2.56013429592102e-10\\
4.3125	2.5258998761535e-10\\
4.325	2.48485267533793e-10\\
4.3375	2.43704530891549e-10\\
4.35	2.38314640355968e-10\\
4.3625	2.32417892612797e-10\\
4.375	2.26120504294364e-10\\
4.3875	2.19595384287288e-10\\
4.4	2.12993385991781e-10\\
4.4125	2.06487649305119e-10\\
4.425	2.00192365717215e-10\\
4.4375	1.94166517488953e-10\\
4.45	1.8843676507609e-10\\
4.4625	1.8296100182286e-10\\
4.475	1.77636377191144e-10\\
4.4875	1.72335264946283e-10\\
4.5	1.66941603556902e-10\\
4.5125	1.61382773021664e-10\\
4.525	1.55677660258285e-10\\
4.5375	1.49934821353843e-10\\
4.55	1.44364982161862e-10\\
4.5625	1.39230110139795e-10\\
4.575	1.3476144873757e-10\\
4.5875	1.3109722049029e-10\\
4.6	1.28300908702202e-10\\
4.6125	1.26339508822012e-10\\
4.625	1.25089136782123e-10\\
4.6375	1.24461998582062e-10\\
4.65	1.24326327172931e-10\\
4.6625	1.24592397390146e-10\\
4.675	1.25199240825366e-10\\
4.6875	1.26086851061892e-10\\
4.7	1.27218876560834e-10\\
4.7125	1.28561133652766e-10\\
4.725	1.30100264059673e-10\\
4.7375	1.31800181427303e-10\\
4.75	1.33640088037867e-10\\
4.7625	1.35595185610762e-10\\
4.775	1.37639566041648e-10\\
4.7875	1.39752007165206e-10\\
4.8	1.41901640827436e-10\\
4.8125	1.44059958378202e-10\\
4.825	1.46200186262485e-10\\
4.8375	1.48304828153097e-10\\
4.85	1.50346660631576e-10\\
4.8625	1.5230395208525e-10\\
4.875	1.54143920535228e-10\\
4.8875	1.55873741439661e-10\\
4.9	1.5747357995632e-10\\
4.9125	1.5894472268081e-10\\
4.925	1.60282148884713e-10\\
4.9375	1.61485076869744e-10\\
4.95	1.62553146550832e-10\\
4.9625	1.63487574802597e-10\\
4.975	1.64294901858055e-10\\
4.9875	1.64981387277337e-10\\
5	1.65553864374066e-10\\
};
\addlegendentry{\texttt{TDB-CUR-CR-OS} $|$ $\epsilon_u = 10^{-11}$}

\addplot [color=mycolor2, dotted, line width=1.2pt,y filter/.expression={log10(y)}]
  table[row sep=crcr]{%
0	6.25305022111612e-19\\
0.0125	2.57821786147377e-14\\
0.025	4.9768591903554e-14\\
0.0375	1.46155800819145e-13\\
0.05	2.06909838438981e-13\\
0.0625	3.11407814489249e-13\\
0.075	5.19036042075774e-13\\
0.0875	6.46673875066252e-13\\
0.1	8.14925234540622e-13\\
0.1125	1.11465845998137e-12\\
0.125	1.58794280708104e-12\\
0.1375	2.24715092452125e-12\\
0.15	3.04176492656747e-12\\
0.1625	3.93173847593347e-12\\
0.175	5.11883504895187e-12\\
0.1875	6.89561333936423e-12\\
0.2	9.65655317901314e-12\\
0.2125	1.34721155397483e-11\\
0.225	1.64939614812536e-11\\
0.2375	1.80360890602349e-11\\
0.25	1.89347452006755e-11\\
0.2625	1.9940719378409e-11\\
0.275	2.13091022874618e-11\\
0.2875	2.2608475273703e-11\\
0.3	2.37744267694863e-11\\
0.3125	2.53358731906743e-11\\
0.325	2.80421024543617e-11\\
0.3375	3.2256803143307e-11\\
0.35	3.65126714736412e-11\\
0.3625	3.90911954589358e-11\\
0.375	4.04786044077894e-11\\
0.3875	4.15232380712102e-11\\
0.4	4.26599249236053e-11\\
0.4125	4.42757685215506e-11\\
0.425	4.68719291710486e-11\\
0.4375	5.11256873757258e-11\\
0.45	5.76976408099581e-11\\
0.4625	6.70177389969487e-11\\
0.475	7.98410685282691e-11\\
0.4875	9.78272994399686e-11\\
0.5	1.10547389132086e-10\\
0.5125	1.16208731271789e-10\\
0.525	1.23530567349189e-10\\
0.5375	1.33882770808135e-10\\
0.55	1.4918917936293e-10\\
0.5625	1.70496304056769e-10\\
0.575	1.97845321332289e-10\\
0.5875	2.31431557429502e-10\\
0.6	2.71041587953965e-10\\
0.6125	3.13019452407016e-10\\
0.625	3.57295745119992e-10\\
0.6375	4.03470333805182e-10\\
0.65	4.51190506197088e-10\\
0.6625	5.01405061829858e-10\\
0.675	5.54670701999263e-10\\
0.6875	6.08515010496503e-10\\
0.7	6.61920026731604e-10\\
0.7125	7.16080262978518e-10\\
0.725	7.712160708598e-10\\
0.7375	8.26853990366644e-10\\
0.75	8.83021966838653e-10\\
0.7625	9.40255154082214e-10\\
0.775	9.99338942479845e-10\\
0.7875	1.06066326869238e-09\\
0.8	1.12425195098155e-09\\
0.8125	1.18984490876633e-09\\
0.825	1.25804423203108e-09\\
0.8375	1.33034462169224e-09\\
0.85	1.40755418121071e-09\\
0.8625	1.48851774788902e-09\\
0.875	1.57054691184513e-09\\
0.8875	1.65046204567327e-09\\
0.9	1.72587167451559e-09\\
0.9125	1.7955053323273e-09\\
0.925	1.85882829399606e-09\\
0.9375	1.91571446462037e-09\\
0.95	1.96627472725165e-09\\
0.9625	2.01076243628149e-09\\
0.975	2.04976177388453e-09\\
0.9875	2.08362986605243e-09\\
1	2.11277199592984e-09\\
1.0125	2.13758658684918e-09\\
1.025	2.15845603080373e-09\\
1.0375	2.17575266686472e-09\\
1.05	2.1898348544215e-09\\
1.0625	2.20103053820503e-09\\
1.075	2.20963173306207e-09\\
1.0875	2.21590960967327e-09\\
1.1	2.22011390170358e-09\\
1.1125	2.22244774548385e-09\\
1.125	2.22311775626731e-09\\
1.1375	2.22231478191946e-09\\
1.15	2.22020716135688e-09\\
1.1625	2.21693965812107e-09\\
1.175	2.21263181542093e-09\\
1.1875	2.20740500093483e-09\\
1.2	2.20137877014665e-09\\
1.2125	2.19465954494001e-09\\
1.225	2.18734041409542e-09\\
1.2375	2.17950164867086e-09\\
1.25	2.17121232673668e-09\\
1.2625	2.16253403805361e-09\\
1.275	2.15352313979732e-09\\
1.2875	2.14423202205433e-09\\
1.3	2.13471080271581e-09\\
1.3125	2.12500749235588e-09\\
1.325	2.11516326872489e-09\\
1.3375	2.10520534761562e-09\\
1.35	2.09515544525535e-09\\
1.3625	2.08503708061417e-09\\
1.375	2.07487344490627e-09\\
1.3875	2.06468501895912e-09\\
1.4	2.05449049463381e-09\\
1.4125	2.0443082868925e-09\\
1.425	2.03415745795257e-09\\
1.4375	2.02405689337261e-09\\
1.45	2.01402671381273e-09\\
1.4625	2.0040907565763e-09\\
1.475	1.99427867858443e-09\\
1.4875	1.9846228115737e-09\\
1.5	1.97515774574348e-09\\
1.5125	1.96593135561278e-09\\
1.525	1.95700320857938e-09\\
1.5375	1.9484533693648e-09\\
1.55	1.94040584302906e-09\\
1.5625	1.93301664529254e-09\\
1.575	1.92649204802735e-09\\
1.5875	1.92105581197183e-09\\
1.6	1.91691847149298e-09\\
1.6125	1.91427888641771e-09\\
1.625	1.91322838797914e-09\\
1.6375	1.91368309250708e-09\\
1.65	1.91575342001979e-09\\
1.6625	1.91962458722631e-09\\
1.675	1.92534903945495e-09\\
1.6875	1.93294242669213e-09\\
1.7	1.9423722316129e-09\\
1.7125	1.95353992828331e-09\\
1.725	1.96628122576847e-09\\
1.7375	1.98036908348464e-09\\
1.75	1.99584406897292e-09\\
1.7625	2.01260000603698e-09\\
1.775	2.0304371557752e-09\\
1.7875	2.04919405081901e-09\\
1.8	2.06870731186843e-09\\
1.8125	2.08858315097674e-09\\
1.825	2.10838911673507e-09\\
1.8375	2.12782731870687e-09\\
1.85	2.14670416323052e-09\\
1.8625	2.16487156407523e-09\\
1.875	2.18214425185511e-09\\
1.8875	2.19835533766007e-09\\
1.9	2.21338029558071e-09\\
1.9125	2.227126204481e-09\\
1.925	2.23952204173577e-09\\
1.9375	2.25051439611533e-09\\
1.95	2.26006418733738e-09\\
1.9625	2.26814463238045e-09\\
1.975	2.27478396397418e-09\\
1.9875	2.28015220174638e-09\\
2	2.28438471555699e-09\\
2.0125	2.28755321362381e-09\\
2.025	2.28974333281004e-09\\
2.0375	2.2910473618708e-09\\
2.05	2.29155338485432e-09\\
2.0625	2.29134386001665e-09\\
2.075	2.29050335838125e-09\\
2.0875	2.28910942267264e-09\\
2.1	2.28724012713204e-09\\
2.1125	2.28499982000905e-09\\
2.125	2.28248355376531e-09\\
2.1375	2.2798349733356e-09\\
2.15	2.27724318845888e-09\\
2.1625	2.27494476928652e-09\\
2.175	2.27325898988496e-09\\
2.1875	2.27262391739894e-09\\
2.2	2.27361949363623e-09\\
2.2125	2.27700330057596e-09\\
2.225	2.2837451332769e-09\\
2.2375	2.2950416002103e-09\\
2.25	2.31230163384783e-09\\
2.2625	2.33708653394655e-09\\
2.275	2.37098727815608e-09\\
2.2875	2.41543431993813e-09\\
2.3	2.47146093647598e-09\\
2.3125	2.53946789267135e-09\\
2.325	2.61903233769136e-09\\
2.3375	2.70883477307606e-09\\
2.35	2.8067347164715e-09\\
2.3625	2.90998741134881e-09\\
2.375	3.01553803317054e-09\\
2.3875	3.12034931117035e-09\\
2.4	3.22168322177029e-09\\
2.4125	3.31729613367035e-09\\
2.425	3.40553730784458e-09\\
2.4375	3.48535870276913e-09\\
2.45	3.55625907345534e-09\\
2.4625	3.61818937451151e-09\\
2.475	3.67144297301909e-09\\
2.4875	3.71655155599779e-09\\
2.5	3.75419366849024e-09\\
2.5125	3.78512162502358e-09\\
2.525	3.81011364722116e-09\\
2.5375	3.82993709682067e-09\\
2.55	3.84532374308041e-09\\
2.5625	3.85695786514458e-09\\
2.575	3.86549319490856e-09\\
2.5875	3.87152334886914e-09\\
2.6	3.87558154995564e-09\\
2.6125	3.8781426335962e-09\\
2.625	3.8796196547714e-09\\
2.6375	3.88033901830202e-09\\
2.65	3.88058212578051e-09\\
2.6625	3.88058737674534e-09\\
2.675	3.88054939387344e-09\\
2.6875	3.88063696332225e-09\\
2.7	3.88099811296971e-09\\
2.7125	3.88173745220879e-09\\
2.725	3.88293014976311e-09\\
2.7375	3.88463347149835e-09\\
2.75	3.88686859107564e-09\\
2.7625	3.88959977051785e-09\\
2.775	3.89276560028585e-09\\
2.7875	3.89626972612671e-09\\
2.8	3.89999757940189e-09\\
2.8125	3.90384170823331e-09\\
2.825	3.90769278157687e-09\\
2.8375	3.91144598813016e-09\\
2.85	3.91501138745011e-09\\
2.8625	3.91827700444594e-09\\
2.875	3.92110049784116e-09\\
2.8875	3.92335759811341e-09\\
2.9	3.924953775611e-09\\
2.9125	3.92581772881001e-09\\
2.925	3.925897119403e-09\\
2.9375	3.9251562160074e-09\\
2.95	3.92357204914569e-09\\
2.9625	3.92112721439841e-09\\
2.975	3.91785503249374e-09\\
2.9875	3.9138656870374e-09\\
3	3.90919147993497e-09\\
3.0125	3.90387251332499e-09\\
3.025	3.89796294437729e-09\\
3.0375	3.89152125283385e-09\\
3.05	3.88460432050699e-09\\
3.0625	3.87726514730265e-09\\
3.075	3.86955100485457e-09\\
3.0875	3.86150272979573e-09\\
3.1	3.85315800654607e-09\\
3.1125	3.8445489791975e-09\\
3.125	3.83569986817846e-09\\
3.1375	3.82663878949433e-09\\
3.15	3.8173898849005e-09\\
3.1625	3.80797244004574e-09\\
3.175	3.79840934945257e-09\\
3.1875	3.78871871213108e-09\\
3.2	3.77891728383409e-09\\
3.2125	3.76902084355872e-09\\
3.225	3.75904384023344e-09\\
3.2375	3.74900031575363e-09\\
3.25	3.73890345671277e-09\\
3.2625	3.7287653488998e-09\\
3.275	3.71859700972488e-09\\
3.2875	3.70840797879043e-09\\
3.3	3.69820703434562e-09\\
3.3125	3.68800259858765e-09\\
3.325	3.67780274745731e-09\\
3.3375	3.66761499428556e-09\\
3.35	3.65744587067915e-09\\
3.3625	3.64729933281431e-09\\
3.375	3.63718155499987e-09\\
3.3875	3.62711814222209e-09\\
3.4	3.61710325978597e-09\\
3.4125	3.60713820790757e-09\\
3.425	3.59722779847e-09\\
3.4375	3.58737792154291e-09\\
3.45	3.57759565045373e-09\\
3.4625	3.56788997134362e-09\\
3.475	3.55827215596484e-09\\
3.4875	3.54875240744618e-09\\
3.5	3.53934719828583e-09\\
3.5125	3.53008022922215e-09\\
3.525	3.52098030949947e-09\\
3.5375	3.512088851519e-09\\
3.55	3.50346392418865e-09\\
3.5625	3.49518626330178e-09\\
3.575	3.48734257246249e-09\\
3.5875	3.48003590824241e-09\\
3.6	3.47338476743607e-09\\
3.6125	3.46750317236506e-09\\
3.625	3.46248891459394e-09\\
3.6375	3.45838998496968e-09\\
3.65	3.45521484317885e-09\\
3.6625	3.45298562052695e-09\\
3.675	3.45180554032805e-09\\
3.6875	3.45174298250069e-09\\
3.7	3.4528108552994e-09\\
3.7125	3.45500953290992e-09\\
3.725	3.4582927042369e-09\\
3.7375	3.46257196569269e-09\\
3.75	3.46773566942895e-09\\
3.7625	3.47370027656568e-09\\
3.775	3.4804447920049e-09\\
3.7875	3.48782796998263e-09\\
3.8	3.49571265805462e-09\\
3.8125	3.50394617448183e-09\\
3.825	3.51228731413009e-09\\
3.8375	3.52055845939783e-09\\
3.85	3.52859774451118e-09\\
3.8625	3.53625889168512e-09\\
3.875	3.54342306506143e-09\\
3.8875	3.55002131097349e-09\\
3.9	3.55600163085894e-09\\
3.9125	3.56131747352725e-09\\
3.925	3.5659388714783e-09\\
3.9375	3.56985081893992e-09\\
3.95	3.57305309162354e-09\\
3.9625	3.57556152146069e-09\\
3.975	3.57741379574642e-09\\
3.9875	3.57868007377993e-09\\
4	3.57944239342768e-09\\
4.0125	3.57978401633959e-09\\
4.025	3.57979688706604e-09\\
4.0375	3.57957733520053e-09\\
4.05	3.5792218400886e-09\\
4.0625	3.57882413537711e-09\\
4.075	3.57847355604752e-09\\
4.0875	3.57825220953236e-09\\
4.1	3.57823210812885e-09\\
4.1125	3.57846906195268e-09\\
4.125	3.57899830973997e-09\\
4.1375	3.5798291769115e-09\\
4.15	3.58093238011296e-09\\
4.1625	3.58222512060716e-09\\
4.175	3.58355485811391e-09\\
4.1875	3.5846804110684e-09\\
4.2	3.58525141594485e-09\\
4.2125	3.58478837147457e-09\\
4.225	3.58266679889873e-09\\
4.2375	3.57810877879123e-09\\
4.25	3.57018877020028e-09\\
4.2625	3.55785977888373e-09\\
4.275	3.54000596216124e-09\\
4.2875	3.51552541345742e-09\\
4.3	3.48344192998645e-09\\
4.3125	3.44303724923592e-09\\
4.325	3.39398598198948e-09\\
4.3375	3.3364666713042e-09\\
4.35	3.27121710978248e-09\\
4.3625	3.1995040943333e-09\\
4.375	3.12299025458885e-09\\
4.3875	3.04350509324968e-09\\
4.4	2.9627527303316e-09\\
4.4125	2.88201001083356e-09\\
4.425	2.80188781371072e-09\\
4.4375	2.72221269391724e-09\\
4.45	2.64206234192312e-09\\
4.4625	2.55995841400729e-09\\
4.475	2.47419910496774e-09\\
4.4875	2.38330506024523e-09\\
4.5	2.28653949221084e-09\\
4.5125	2.18441542467705e-09\\
4.525	2.07903257150953e-09\\
4.5375	1.97403573517632e-09\\
4.55	1.87406025078492e-09\\
4.5625	1.78375560680707e-09\\
4.575	1.70673150807737e-09\\
4.5875	1.64485655549055e-09\\
4.6	1.5981619127351e-09\\
4.6125	1.56527459009479e-09\\
4.625	1.54407589008725e-09\\
4.6375	1.53229318820378e-09\\
4.65	1.52788923537109e-09\\
4.6625	1.52922982814121e-09\\
4.675	1.53509963338066e-09\\
4.6875	1.54463466593203e-09\\
4.7	1.55722609557323e-09\\
4.7125	1.5724304395614e-09\\
4.725	1.58990261559539e-09\\
4.7375	1.60934457090631e-09\\
4.75	1.63047052231014e-09\\
4.7625	1.65298906611527e-09\\
4.775	1.67659630136186e-09\\
4.7875	1.70097863762288e-09\\
4.8	1.7258176949348e-09\\
4.8125	1.75079734138123e-09\\
4.825	1.77561211244684e-09\\
4.8375	1.79997261112288e-09\\
4.85	1.82361678812973e-09\\
4.8625	1.84631267010701e-09\\
4.875	1.86786079546231e-09\\
4.8875	1.88809496065256e-09\\
4.9	1.90688242581884e-09\\
4.9125	1.92413093709537e-09\\
4.925	1.93978437040679e-09\\
4.9375	1.95381733179759e-09\\
4.95	1.96623273830047e-09\\
4.9625	1.97706906812262e-09\\
4.975	1.98638838716605e-09\\
4.9875	1.99426708648481e-09\\
5	2.00079667922965e-09\\
};
\addlegendentry{\texttt{TDB-SVD} $|$ $\epsilon_u = 10^{-11}$}

\end{axis}
\end{tikzpicture}%

%% file: figures/Error_AC.tex
\definecolor{mycolor1}{rgb}{0.00000,0.44700,0.74100}%
\definecolor{mycolor2}{rgb}{0.85000,0.32500,0.09800}%
\definecolor{mycolor3}{rgb}{0.5,0.5,0.5}
\begin{tikzpicture}[scale=0.8]

\begin{axis}[%
title={Allen-Cahn},
width=1.5in,
height=1.5in,
scale only axis,
xmin=0,
xmax=200,
xlabel style={font=\color{white!15!black}},
xlabel={Time},
ymin=5.91173206668798e-15, ymax=3.28323957774326e-05,
ymode=log,
yminorticks=true,
ylabel style={font=\color{white!15!black}},
axis background/.style={fill=white},
legend style={fill opacity=0.8, 
font=\tiny,
  draw opacity=1,
  text opacity=1,
  at={(0.04,0.22)},
  anchor=west}
  % draw=lightgray204}
]
\addplot [line width=1.2pt,  mycolor1]
table {%
0 2.45692690877759e-07
1 2.31492830389881e-06
2 1.60013520054813e-06
3 1.95632965617386e-06
4 1.76580187513535e-06
5 1.81394451004179e-06
6 1.90469608197097e-06
7 1.94395632978442e-06
8 1.95519964529717e-06
9 1.83986676809932e-06
10 1.82495665907816e-06
11 1.88516407675114e-06
12 1.84423489002756e-06
13 1.87534611670293e-06
14 1.87000490269911e-06
15 1.87867382368284e-06
16 1.86660463478596e-06
17 1.87237343154608e-06
18 1.87528476971749e-06
19 1.88357720202079e-06
20 1.90617600154736e-06
21 1.90944321225898e-06
22 1.90367107011149e-06
23 1.90118292402996e-06
24 1.90429376587688e-06
25 1.9099572166601e-06
26 1.9168257987714e-06
27 1.92148967438605e-06
28 1.92314181627005e-06
29 1.92799310994158e-06
30 1.93340960377825e-06
31 1.9384792902849e-06
32 1.94563536186549e-06
33 1.95174523703661e-06
34 1.9576301253963e-06
35 1.96335271579285e-06
36 1.96863357965562e-06
37 1.97294293036079e-06
38 1.97784875979678e-06
39 1.9860421305778e-06
40 1.9917440649471e-06
41 2.00767448425999e-06
42 2.00883989534068e-06
43 2.0208393189785e-06
44 2.03110933471337e-06
45 2.04842252635996e-06
46 2.06676697126539e-06
47 2.06186889329657e-06
48 2.08499432946442e-06
49 2.09266773115163e-06
50 2.11264746326696e-06
51 2.13136403987683e-06
52 2.14117588037773e-06
53 2.17476561113623e-06
54 2.17558459811477e-06
55 2.19823400635348e-06
56 2.23251624180234e-06
57 2.26527692620129e-06
58 2.29323448829829e-06
59 2.33466355804783e-06
60 2.35310080884023e-06
61 2.38157868104812e-06
62 2.42785799860717e-06
63 2.48652190001763e-06
64 2.55659935276099e-06
65 2.64113148430892e-06
66 2.76013769740897e-06
67 2.90505480214301e-06
68 3.12998371930904e-06
69 3.53715288191181e-06
70 4.25735062576737e-06
71 5.96470716708511e-06
72 8.43627972995395e-06
73 7.74991495295134e-06
74 1.01103919373136e-05
75 6.89224780489232e-06
76 7.38625938632393e-06
77 1.18384531402632e-05
78 1.50393917206993e-05
79 1.03504750168038e-05
80 1.1509041396335e-05
81 1.17210962390972e-05
82 1.26743726072088e-05
83 1.04671313012719e-05
84 1.36561996615961e-05
85 2.04127159871033e-05
86 2.29101101108319e-05
87 1.87479892209314e-05
88 1.66890542688994e-05
89 2.00038732198428e-05
90 2.10732493929208e-05
91 2.33737755099174e-05
92 2.29715441650544e-05
93 1.79946783471127e-05
94 2.12383060127869e-05
95 2.37250902333471e-05
96 2.60057120046332e-05
97 2.47724989979681e-05
98 2.07929450253628e-05
99 1.71309732319405e-05
100 1.67924040126071e-05
101 1.73440219688274e-05
102 1.71928329709084e-05
103 1.92777951568176e-05
104 2.01522874782918e-05
105 2.29071502066295e-05
106 2.0955051616437e-05
107 1.93509559656269e-05
108 1.59846139971382e-05
109 1.63781719118537e-05
110 1.84343751207363e-05
111 1.38780247920476e-05
112 1.7357216754528e-05
113 1.93021738985999e-05
114 2.04878379477482e-05
115 1.84151582234665e-05
116 1.52632117845208e-05
117 1.37574628751232e-05
118 1.39931113470523e-05
119 1.10757133522181e-05
120 1.0237436662577e-05
121 1.26114680769938e-05
122 1.37802488765248e-05
123 1.44333524599091e-05
124 1.42449283908134e-05
125 1.28202246300293e-05
126 1.3652994136451e-05
127 1.15750809210106e-05
128 1.13501339480018e-05
129 1.24574552213138e-05
130 1.07198049301808e-05
131 6.98608509251919e-06
132 6.06334746505588e-06
133 6.87395947243776e-06
134 9.26426038947838e-06
135 1.11713015144013e-05
136 6.42946155819278e-06
137 6.15503000020453e-06
138 7.27768449147432e-06
139 8.52561185104977e-06
140 1.20992495561818e-05
141 1.13297292340361e-05
142 1.4168090758142e-05
143 2.12544107413314e-05
144 1.71943754078273e-05
145 9.8503153456132e-06
146 9.89732933565548e-06
147 9.70984321675819e-06
148 6.75450463157172e-06
149 9.88338916661575e-06
150 6.9805017240674e-06
151 3.09827748793244e-06
152 2.66853718145731e-06
153 2.91309737000688e-06
154 3.51435032735708e-06
155 4.99695954125211e-06
156 6.20386034608748e-06
157 3.96658172976882e-06
158 4.80380703948384e-06
159 7.68629082599287e-06
160 8.85664364571583e-06
161 4.18020517993044e-06
162 4.82258312243e-06
163 6.86034534302664e-06
164 4.17694157713955e-06
165 2.60536562788688e-06
166 2.48278103196182e-06
167 2.49574271319872e-06
168 2.51294379801763e-06
169 2.53117413526118e-06
170 2.55155579600029e-06
171 2.57267046339756e-06
172 2.59196045313681e-06
173 2.62394668534617e-06
174 2.63498940229489e-06
175 2.65935319403959e-06
176 2.68402523640171e-06
177 2.71363363922286e-06
178 2.74685056829829e-06
179 2.78766784152086e-06
180 2.8203704712154e-06
181 2.86305213094336e-06
182 2.89966666166539e-06
183 2.96690177564672e-06
184 3.02756180375538e-06
185 3.08450894491334e-06
186 3.17750304188138e-06
187 3.2675464508897e-06
188 3.35558766973549e-06
189 3.58713251944853e-06
190 3.67107348457806e-06
191 3.92390992669065e-06
192 4.38651818374826e-06
193 4.94501854835843e-06
194 6.1198360993347e-06
195 8.76599639399967e-06
196 1.43111711981322e-05
197 1.53482524374279e-05
198 7.11816092082601e-06
199 7.29948823446378e-06
};

\addplot [line width=1.2pt,  mycolor2, dotted]
table {%
0 2.45692690877759e-07
1 3.07214811697117e-07
2 2.76212960330526e-07
3 3.48616850392662e-07
4 4.03862326491645e-07
5 4.21260644039921e-07
6 4.25619957165308e-07
7 4.26721068372356e-07
8 4.27031375610754e-07
9 4.27142391062207e-07
10 4.27201412203719e-07
11 4.27248864032735e-07
12 4.27291952478064e-07
13 4.27336838546079e-07
14 4.27381088353726e-07
15 4.27426251195897e-07
16 4.27471887531287e-07
17 4.27518213094374e-07
18 4.2756469789371e-07
19 4.27613760180709e-07
20 4.27661919369334e-07
21 4.27712015955611e-07
22 4.27763829947623e-07
23 4.27815006518817e-07
24 4.27868132817493e-07
25 4.27922510583836e-07
26 4.27978117048744e-07
27 4.28034904231643e-07
28 4.28093284343812e-07
29 4.28153628972408e-07
30 4.28215279646196e-07
31 4.28278638567383e-07
32 4.28343846518116e-07
33 4.28410873195163e-07
34 4.28479706757821e-07
35 4.28550199036806e-07
36 4.28622449718175e-07
37 4.28697314479529e-07
38 4.28775423813986e-07
39 4.28855363389027e-07
40 4.28938539246213e-07
41 4.29025400365056e-07
42 4.29114559838772e-07
43 4.29207832318492e-07
44 4.29305818839735e-07
45 4.29407314423514e-07
46 4.29513488237437e-07
47 4.29624665124123e-07
48 4.29741438267607e-07
49 4.29864386500186e-07
50 4.299941375992e-07
51 4.30131419636359e-07
52 4.30277105411823e-07
53 4.30432085350226e-07
54 4.30597771377347e-07
55 4.30774904750049e-07
56 4.30966465817885e-07
57 4.31174024036358e-07
58 4.31399318865429e-07
59 4.31646725765954e-07
60 4.31919326341051e-07
61 4.32223765229284e-07
62 4.32567157727946e-07
63 4.32960080736217e-07
64 4.33418081759821e-07
65 4.3396582056813e-07
66 4.34642566386208e-07
67 4.3551920595284e-07
68 4.36736791780775e-07
69 4.38620385725427e-07
70 4.42080928195035e-07
71 4.50141928319952e-07
72 4.66328977615934e-07
73 4.9590372928248e-07
74 5.5988854795259e-07
75 4.87703093570683e-07
76 4.96448127278546e-07
77 5.91162468218918e-07
78 5.81030613846887e-07
79 4.95347004026226e-07
80 5.20759002799715e-07
81 5.45933123095681e-07
82 5.5162262983119e-07
83 5.16407013200795e-07
84 6.03742009796836e-07
85 8.2525393266322e-07
86 8.00891726026059e-07
87 5.62351542687905e-07
88 6.14285149497008e-07
89 7.23927582033275e-07
90 8.40119831760348e-07
91 9.6058045874697e-07
92 7.07180456275756e-07
93 5.92897939065438e-07
94 6.9400628201373e-07
95 8.60885533834589e-07
96 1.09235417486363e-06
97 1.19316237309864e-06
98 1.02065779552945e-06
99 1.00644484465694e-06
100 1.20560801354431e-06
101 8.73107025654227e-07
102 9.60564911053833e-07
103 1.10696135615232e-06
104 1.34180198241136e-06
105 1.25455862323298e-06
106 1.38732384828004e-06
107 1.41248478665568e-06
108 1.1599581434943e-06
109 9.77369658369268e-07
110 1.04965557355242e-06
111 1.11691311765537e-06
112 1.4092952908868e-06
113 1.48728210388127e-06
114 1.61222087932111e-06
115 1.27572278264852e-06
116 1.1194790533709e-06
117 9.68866938178286e-07
118 1.19909696536969e-06
119 9.84450747595668e-07
120 9.12368722968346e-07
121 1.00357827884792e-06
122 1.2464012195845e-06
123 1.28324989591906e-06
124 1.12665035836706e-06
125 9.83117115047415e-07
126 1.09300281889776e-06
127 8.66198713691264e-07
128 9.64384774580002e-07
129 1.03550870564306e-06
130 7.9017077636512e-07
131 6.70887354356776e-07
132 7.67565200892894e-07
133 8.11048082692192e-07
134 1.15584137004389e-06
135 1.31085887440205e-06
136 7.20315350529021e-07
137 6.02087088181223e-07
138 6.49779859647184e-07
139 6.4596172117946e-07
140 7.11347034526967e-07
141 7.79623511718911e-07
142 1.01317510880625e-06
143 1.63311110468963e-06
144 2.04567388967209e-06
145 1.09095466703355e-06
146 9.08684609701765e-07
147 9.75852951263612e-07
148 6.57565702518176e-07
149 7.34837835184995e-07
150 6.12320342582652e-07
151 4.42603018188568e-07
152 4.34089700668402e-07
153 4.34781281706722e-07
154 4.36049112130562e-07
155 4.38211719066536e-07
156 4.4156661788139e-07
157 4.4468532929866e-07
158 4.53617645395165e-07
159 4.75329659629327e-07
160 4.98472152600397e-07
161 5.5150574766193e-07
162 7.34812221617771e-07
163 9.66853024268555e-07
164 7.26737005094283e-07
165 5.01210780910255e-07
166 4.33950290838891e-07
167 4.32221692047435e-07
168 4.32309353390576e-07
169 4.32431562251238e-07
170 4.32559600538543e-07
171 4.32693720421257e-07
172 4.32834244421631e-07
173 4.32981681791414e-07
174 4.33136855535795e-07
175 4.33300655537214e-07
176 4.33474859295386e-07
177 4.33659201179524e-07
178 4.33855981600326e-07
179 4.34066805476345e-07
180 4.34293503274084e-07
181 4.34538461019986e-07
182 4.34804717127373e-07
183 4.35096135022763e-07
184 4.35417662120977e-07
185 4.3577588843072e-07
186 4.36178850247678e-07
187 4.36638396591803e-07
188 4.37172217223044e-07
189 4.37806720551768e-07
190 4.3857997537623e-07
191 4.39559033686424e-07
192 4.40865122957314e-07
193 4.42766559796968e-07
194 4.46005638107091e-07
195 4.53213372017856e-07
196 4.70112600080035e-07
197 4.66641335544134e-07
198 4.69357308478595e-07
199 4.94605438263497e-07
};

\addplot [line width=1.2pt,  mycolor1, dotted]
table {%
0 1.97629294434362e-18
1 5.29922398286269e-10
2 6.82634670584001e-10
3 1.18962623282018e-09
4 9.53589989441056e-10
5 1.86082876457661e-09
6 1.00690565276711e-09
7 1.03983023842192e-09
8 9.75234453268053e-10
9 9.73967566494279e-10
10 1.1313770322985e-09
11 1.0761813269453e-09
12 1.11961816467441e-09
13 1.12980525715614e-09
14 1.20156804130447e-09
15 1.12120023958051e-09
16 1.36983936550348e-09
17 1.38239992986596e-09
18 1.21730381558946e-09
19 1.21936337093582e-09
20 1.16021065331653e-09
21 1.15433837357836e-09
22 1.15605647967188e-09
23 1.17777432976375e-09
24 1.18633690833037e-09
25 1.55431194921437e-09
26 1.39130313500349e-09
27 1.47439840095426e-09
28 1.61741899554693e-09
29 1.8872869081976e-09
30 1.92907830768305e-09
31 1.8181353134766e-09
32 1.91756737039289e-09
33 1.83815921079305e-09
34 1.79404063984803e-09
35 1.85106671171827e-09
36 1.79391872586359e-09
37 1.76221508478448e-09
38 1.78263705742186e-09
39 1.78160160833939e-09
40 1.78493166157452e-09
41 1.81562097635785e-09
42 1.829256465179e-09
43 1.8917204864197e-09
44 1.8688099973637e-09
45 1.87745243864262e-09
46 1.89119510546251e-09
47 1.89500136367307e-09
48 1.90638770518664e-09
49 1.92200922839778e-09
50 1.94901267367662e-09
51 1.96662794670531e-09
52 1.99512036648879e-09
53 2.00756825745658e-09
54 2.0259944847235e-09
55 2.04531651477845e-09
56 2.05873287521151e-09
57 2.07744320160712e-09
58 2.10658275014285e-09
59 2.14002764766618e-09
60 2.20016975683699e-09
61 2.20425002176042e-09
62 2.24274268022618e-09
63 2.25220984410012e-09
64 2.28606551946124e-09
65 2.3249584476279e-09
66 2.38325138463418e-09
67 2.42278836337012e-09
68 2.58295568581719e-09
69 2.80908805820179e-09
70 3.26457092676837e-09
71 3.40299883108109e-09
72 3.45859891667815e-09
73 4.08500874913883e-09
74 5.34880852875442e-09
75 4.13376153381202e-09
76 4.48740054556066e-09
77 6.42887994105991e-09
78 6.6784458242028e-09
79 6.44064060443941e-09
80 8.04980460636247e-09
81 8.03768560428644e-09
82 8.16901799136048e-09
83 5.79230266115352e-09
84 6.45684628188177e-09
85 9.49082849441232e-09
86 1.01122272842743e-08
87 7.95483557773064e-09
88 7.55693905451011e-09
89 9.35981260059094e-09
90 1.17883784230693e-08
91 1.21439725941547e-08
92 1.27932425514142e-08
93 1.25688206417915e-08
94 1.62458246618053e-08
95 1.651207693391e-08
96 1.75210761799779e-08
97 1.55672470699844e-08
98 1.57219439295109e-08
99 1.72877094026648e-08
100 1.7802299392485e-08
101 1.66399664989995e-08
102 1.70995968515852e-08
103 1.83859865940096e-08
104 1.9081640749168e-08
105 2.31226741659203e-08
106 2.6406135142496e-08
107 2.50235423655134e-08
108 2.47900040601998e-08
109 2.6985947274039e-08
110 2.78987451743008e-08
111 3.13984907826505e-08
112 4.29451153431085e-08
113 4.71322952490029e-08
114 5.11901609720644e-08
115 4.97390320215183e-08
116 5.17275097289779e-08
117 5.09782466880552e-08
118 5.0202851968867e-08
119 3.89620927233008e-08
120 3.42655108152015e-08
121 4.70785650480366e-08
122 6.00596602401288e-08
123 5.3643667117845e-08
124 4.54216942658873e-08
125 4.31318527756312e-08
126 3.28971635762762e-08
127 3.34706947610602e-08
128 4.2026887468523e-08
129 4.45589059955136e-08
130 3.54185501032351e-08
131 2.44108610589117e-08
132 2.20021632010677e-08
133 2.36814006511709e-08
134 2.57339171021687e-08
135 2.65418764715355e-08
136 1.42382880504674e-08
137 1.5643845907745e-08
138 1.24915718105493e-08
139 1.17853832649851e-08
140 1.28383182628951e-08
141 1.22893779339868e-08
142 1.60789308307867e-08
143 2.58382567776247e-08
144 2.5351325824017e-08
145 1.76453326476933e-08
146 1.56460170789008e-08
147 1.46233385284235e-08
148 1.14697943867855e-08
149 1.65059613085398e-08
150 1.96349813337511e-08
151 1.01043226047267e-08
152 6.86392872955355e-09
153 7.65972098598358e-09
154 9.77634314232718e-09
155 1.49754392697222e-08
156 2.04969377339437e-08
157 1.01626282142359e-08
158 8.06413758724704e-09
159 1.07675935633392e-08
160 1.11698454171813e-08
161 1.05558242734166e-08
162 1.60231021601342e-08
163 2.29886772879533e-08
164 1.87796715625927e-08
165 1.01499556618461e-08
166 5.95668879958483e-09
167 5.82178292164473e-09
168 5.83843819298895e-09
169 5.85760494782414e-09
170 5.87387165093448e-09
171 5.88551041557946e-09
172 5.9406013952171e-09
173 5.97093359287431e-09
174 5.95107088744406e-09
175 5.97701905842603e-09
176 5.98788053401674e-09
177 6.02678645405219e-09
178 6.02670352329914e-09
179 6.0521823905534e-09
180 6.08101138635492e-09
181 6.09375578732023e-09
182 6.10631435643226e-09
183 6.12433673802562e-09
184 6.14700258949188e-09
185 6.16817595967288e-09
186 6.1934342137099e-09
187 6.21764267912855e-09
188 6.23805123734276e-09
189 6.26609464916441e-09
190 6.29659003744072e-09
191 6.33676619296396e-09
192 6.37773012238342e-09
193 6.42590474013333e-09
194 6.50455408371754e-09
195 6.67090707733345e-09
196 7.06373679642864e-09
197 6.8889870627373e-09
198 6.76992550252463e-09
199 6.98785854080359e-09
};

\addplot [line width=1.2pt,  mycolor2,dotted]
table {%
0 1.97629294434362e-18
1 7.70291461582861e-12
2 2.0340560860084e-11
3 3.71496167283403e-11
4 3.46305572661578e-11
5 3.66355680857342e-11
6 3.25617981567896e-11
7 3.24679355788506e-11
8 3.24828226479156e-11
9 3.2487271315178e-11
10 3.2502235437067e-11
11 3.25387738935225e-11
12 3.25699739679594e-11
13 3.26011584136411e-11
14 3.26420490045898e-11
15 3.26702168410749e-11
16 3.27131234057992e-11
17 3.2756728617608e-11
18 3.27803662680924e-11
19 3.28212027957312e-11
20 3.2857604452965e-11
21 3.28958864553375e-11
22 3.29275055884222e-11
23 3.29580679556487e-11
24 3.3001368656125e-11
25 3.37233573571725e-11
26 3.33784549876121e-11
27 3.32709712680377e-11
28 3.32909908391187e-11
29 3.33455499822057e-11
30 3.33955131611262e-11
31 3.34449805708521e-11
32 3.35167168197156e-11
33 3.35694507284568e-11
34 3.36399263894715e-11
35 3.36957948906303e-11
36 3.37709795608255e-11
37 3.38339856730802e-11
38 3.38110013372888e-11
39 3.38490903265943e-11
40 3.3929657212367e-11
41 3.39900383687188e-11
42 3.40641128535243e-11
43 3.41210938347512e-11
44 3.41240147746103e-11
45 3.41982697256884e-11
46 3.42787148857026e-11
47 3.43388185260589e-11
48 3.44138250409074e-11
49 3.44951745202728e-11
50 3.45879671263468e-11
51 3.46815533837454e-11
52 3.47856140940405e-11
53 3.48789544347105e-11
54 3.49738370560443e-11
55 3.50973635213409e-11
56 3.52192152509648e-11
57 3.53524461769981e-11
58 3.5511716680514e-11
59 3.56732541252011e-11
60 3.58545867906258e-11
61 3.60446822301907e-11
62 3.62798609934526e-11
63 3.65145508607864e-11
64 3.6807993904818e-11
65 3.71654646519309e-11
66 3.76235118667051e-11
67 3.82059714157682e-11
68 3.8996901111566e-11
69 4.03063646373893e-11
70 4.27675439651941e-11
71 4.78953981019536e-11
72 5.57487410690011e-11
73 5.99594253972334e-11
74 7.43455277104872e-11
75 6.46201033245092e-11
76 8.16740624058536e-11
77 1.28108912618808e-10
78 1.55394803622757e-10
79 1.04560269768934e-10
80 1.1690889392768e-10
81 1.59146706425484e-10
82 1.33523935963281e-10
83 6.15913409175711e-11
84 7.03896196195267e-11
85 9.43085439940273e-11
86 1.10086757791938e-10
87 7.63152869959569e-11
88 8.64102728356826e-11
89 1.06660653176732e-10
90 1.00071117088522e-10
91 1.19400949568987e-10
92 1.20858732323741e-10
93 1.25521327427783e-10
94 1.40554058062963e-10
95 1.51939010250181e-10
96 1.69759656223706e-10
97 1.7406381430233e-10
98 1.69274150476432e-10
99 1.23812322890183e-10
100 1.49301382089848e-10
101 2.00215683482059e-10
102 1.99038364916202e-10
103 1.47309361496131e-10
104 2.0823647481462e-10
105 2.95558742666829e-10
106 2.72220301662261e-10
107 2.05363307498519e-10
108 1.98279128923627e-10
109 1.69116530279666e-10
110 1.99885261052323e-10
111 1.96761985173649e-10
112 2.12610017517568e-10
113 2.07787997566115e-10
114 2.05580515112946e-10
115 2.00817575917919e-10
116 1.50432280540037e-10
117 1.31948637835824e-10
118 1.77913490024511e-10
119 1.85340478008849e-10
120 1.88474722181709e-10
121 1.86277620309189e-10
122 1.79436353900853e-10
123 2.41225188185799e-10
124 3.24882221470726e-10
125 2.14241131799089e-10
126 1.48768640865401e-10
127 1.55318379596218e-10
128 1.74647998997368e-10
129 1.93324765830317e-10
130 1.70825254685611e-10
131 1.6601659824794e-10
132 1.86517693581444e-10
133 2.02816630152224e-10
134 2.34705037986134e-10
135 2.72135951970838e-10
136 1.84866089457814e-10
137 1.73080305333142e-10
138 1.40874235153631e-10
139 1.25283251070479e-10
140 1.1075479595353e-10
141 9.01496749587576e-11
142 9.85327725191392e-11
143 1.33506095446058e-10
144 1.67084385126569e-10
145 1.9187164843738e-10
146 2.32952771222168e-10
147 1.66197105988974e-10
148 1.25502037259152e-10
149 1.89937163378475e-10
150 1.97686771472182e-10
151 1.01178947791468e-10
152 7.86168914948589e-11
153 9.73057309457773e-11
154 1.36600318561017e-10
155 2.22276520882376e-10
156 3.11813591978085e-10
157 1.84659217683225e-10
158 2.31867268893973e-10
159 4.00411110031455e-10
160 3.83235466094246e-10
161 1.37947374637368e-10
162 1.65480665271672e-10
163 2.71778131190179e-10
164 3.9606024341839e-10
165 2.04972197929057e-10
166 5.0871240529196e-11
167 3.97103374403822e-11
168 3.96446763268645e-11
169 3.97908031424799e-11
170 3.99469381982067e-11
171 4.01105172368e-11
172 4.0282292284809e-11
173 4.04628393547056e-11
174 4.0653281626476e-11
175 4.08549809009276e-11
176 4.10692124749976e-11
177 4.1297311282004e-11
178 4.15411956524114e-11
179 4.18029936962988e-11
180 4.20851714501769e-11
181 4.23907659387177e-11
182 4.27234413566223e-11
183 4.30879340553008e-11
184 4.348964417269e-11
185 4.39356921608579e-11
186 4.44346952269458e-11
187 4.49985494998259e-11
188 4.56421031249105e-11
189 4.63856863503223e-11
190 4.72570626668164e-11
191 4.8295361025671e-11
192 4.95574618302473e-11
193 5.11292744459703e-11
194 5.31462228699573e-11
195 5.58346209365568e-11
196 5.95931995500561e-11
197 6.51247941402299e-11
198 7.39889189406865e-11
199 8.96598200981604e-11
};
25.3 1.33947237426163e-09

\end{axis}
\end{tikzpicture}%

%% file: figures/rank_Burgers.tex
\definecolor{mycolor1}{rgb}{0.00000,0.44700,0.74100}%
\definecolor{mycolor2}{rgb}{0.85000,0.32500,0.09800}%
\definecolor{mycolor3}{rgb}{0.5,0.5,0.5}
\begin{tikzpicture}[scale=0.8]

\begin{axis}[%
title={Burgers},
width=1.5in,
height=1.5in,
at={(1.011in,0.642in)},
scale only axis,
xmin=0,
xmax=5,
xlabel style={font=\color{white!15!black}},
xlabel={Time},
%ymode=log,
ymin=0,
ymax=28,
yminorticks=true,
ylabel style={font=\color{white!15!black}},
ylabel={$r$},
axis background/.style={fill=white},
legend style={fill opacity=0.8, 
font=\tiny,
  draw opacity=1,
  text opacity=1,
  at={(0.3,0.13)},
  anchor=west}
  % draw=lightgray204}
]

\addplot [color=mycolor1, solid, line width=1.2pt]
  table[row sep=crcr]{%
0	17\\
0.0125	5\\
0.025	5\\
0.0375	5\\
0.05	5\\
0.0625	5\\
0.075	5\\
0.0875	5\\
0.1	5\\
0.1125	5\\
0.125	5\\
0.1375	5\\
0.15	5\\
0.1625	5\\
0.175	5\\
0.1875	5\\
0.2	5\\
0.2125	5\\
0.225	5\\
0.2375	5\\
0.25	5\\
0.2625	5\\
0.275	5\\
0.2875	7\\
0.3	7\\
0.3125	7\\
0.325	7\\
0.3375	7\\
0.35	7\\
0.3625	7\\
0.375	7\\
0.3875	7\\
0.4	7\\
0.4125	7\\
0.425	7\\
0.4375	7\\
0.45	7\\
0.4625	7\\
0.475	7\\
0.4875	7\\
0.5	7\\
0.5125	7\\
0.525	7\\
0.5375	7\\
0.55	7\\
0.5625	7\\
0.575	6\\
0.5875	6\\
0.6	6\\
0.6125	6\\
0.625	6\\
0.6375	6\\
0.65	6\\
0.6625	6\\
0.675	6\\
0.6875	6\\
0.7	6\\
0.7125	7\\
0.725	7\\
0.7375	7\\
0.75	7\\
0.7625	7\\
0.775	7\\
0.7875	8\\
0.8	8\\
0.8125	8\\
0.825	8\\
0.8375	8\\
0.85	7\\
0.8625	7\\
0.875	7\\
0.8875	8\\
0.9	8\\
0.9125	8\\
0.925	8\\
0.9375	8\\
0.95	8\\
0.9625	8\\
0.975	8\\
0.9875	6\\
1	6\\
1.0125	6\\
1.025	6\\
1.0375	6\\
1.05	6\\
1.0625	6\\
1.075	6\\
1.0875	8\\
1.1	8\\
1.1125	8\\
1.125	8\\
1.1375	8\\
1.15	8\\
1.1625	8\\
1.175	7\\
1.1875	7\\
1.2	7\\
1.2125	8\\
1.225	8\\
1.2375	8\\
1.25	8\\
1.2625	8\\
1.275	7\\
1.2875	7\\
1.3	7\\
1.3125	7\\
1.325	7\\
1.3375	7\\
1.35	7\\
1.3625	7\\
1.375	7\\
1.3875	7\\
1.4	7\\
1.4125	7\\
1.425	6\\
1.4375	6\\
1.45	6\\
1.4625	6\\
1.475	6\\
1.4875	6\\
1.5	7\\
1.5125	7\\
1.525	7\\
1.5375	7\\
1.55	7\\
1.5625	6\\
1.575	7\\
1.5875	7\\
1.6	7\\
1.6125	7\\
1.625	7\\
1.6375	7\\
1.65	7\\
1.6625	7\\
1.675	7\\
1.6875	7\\
1.7	7\\
1.7125	7\\
1.725	7\\
1.7375	7\\
1.75	7\\
1.7625	7\\
1.775	7\\
1.7875	8\\
1.8	8\\
1.8125	8\\
1.825	8\\
1.8375	9\\
1.85	9\\
1.8625	9\\
1.875	9\\
1.8875	9\\
1.9	9\\
1.9125	9\\
1.925	9\\
1.9375	9\\
1.95	9\\
1.9625	9\\
1.975	9\\
1.9875	8\\
2	8\\
2.0125	8\\
2.025	8\\
2.0375	8\\
2.05	8\\
2.0625	8\\
2.075	8\\
2.0875	8\\
2.1	8\\
2.1125	8\\
2.125	7\\
2.1375	7\\
2.15	7\\
2.1625	7\\
2.175	7\\
2.1875	7\\
2.2	7\\
2.2125	7\\
2.225	7\\
2.2375	7\\
2.25	7\\
2.2625	7\\
2.275	7\\
2.2875	7\\
2.3	6\\
2.3125	6\\
2.325	6\\
2.3375	6\\
2.35	6\\
2.3625	6\\
2.375	6\\
2.3875	6\\
2.4	6\\
2.4125	6\\
2.425	6\\
2.4375	6\\
2.45	6\\
2.4625	6\\
2.475	6\\
2.4875	6\\
2.5	6\\
2.5125	5\\
2.525	5\\
2.5375	5\\
2.55	5\\
2.5625	5\\
2.575	5\\
2.5875	5\\
2.6	5\\
2.6125	5\\
2.625	6\\
2.6375	6\\
2.65	5\\
2.6625	5\\
2.675	5\\
2.6875	5\\
2.7	6\\
2.7125	6\\
2.725	6\\
2.7375	6\\
2.75	6\\
2.7625	6\\
2.775	6\\
2.7875	7\\
2.8	7\\
2.8125	7\\
2.825	7\\
2.8375	6\\
2.85	6\\
2.8625	6\\
2.875	6\\
2.8875	6\\
2.9	6\\
2.9125	6\\
2.925	6\\
2.9375	6\\
2.95	7\\
2.9625	7\\
2.975	7\\
2.9875	7\\
3	7\\
3.0125	5\\
3.025	4\\
3.0375	4\\
3.05	4\\
3.0625	4\\
3.075	3\\
3.0875	3\\
3.1	3\\
3.1125	3\\
3.125	3\\
3.1375	3\\
3.15	6\\
3.1625	5\\
3.175	6\\
3.1875	5\\
3.2	6\\
3.2125	6\\
3.225	5\\
3.2375	6\\
3.25	5\\
3.2625	6\\
3.275	6\\
3.2875	5\\
3.3	5\\
3.3125	5\\
3.325	5\\
3.3375	5\\
3.35	6\\
3.3625	6\\
3.375	6\\
3.3875	6\\
3.4	4\\
3.4125	4\\
3.425	4\\
3.4375	3\\
3.45	3\\
3.4625	3\\
3.475	3\\
3.4875	3\\
3.5	3\\
3.5125	3\\
3.525	3\\
3.5375	3\\
3.55	3\\
3.5625	3\\
3.575	3\\
3.5875	3\\
3.6	3\\
3.6125	5\\
3.625	5\\
3.6375	5\\
3.65	5\\
3.6625	7\\
3.675	6\\
3.6875	6\\
3.7	7\\
3.7125	7\\
3.725	7\\
3.7375	7\\
3.75	7\\
3.7625	6\\
3.775	6\\
3.7875	6\\
3.8	6\\
3.8125	6\\
3.825	6\\
3.8375	6\\
3.85	6\\
3.8625	5\\
3.875	6\\
3.8875	6\\
3.9	6\\
3.9125	6\\
3.925	6\\
3.9375	5\\
3.95	5\\
3.9625	5\\
3.975	5\\
3.9875	5\\
4	5\\
4.0125	5\\
4.025	5\\
4.0375	5\\
4.05	5\\
4.0625	5\\
4.075	5\\
4.0875	6\\
4.1	7\\
4.1125	7\\
4.125	6\\
4.1375	6\\
4.15	6\\
4.1625	6\\
4.175	5\\
4.1875	6\\
4.2	6\\
4.2125	4\\
4.225	4\\
4.2375	4\\
4.25	4\\
4.2625	4\\
4.275	4\\
4.2875	4\\
4.3	6\\
4.3125	6\\
4.325	6\\
4.3375	6\\
4.35	5\\
4.3625	5\\
4.375	5\\
4.3875	6\\
4.4	5\\
4.4125	6\\
4.425	6\\
4.4375	6\\
4.45	6\\
4.4625	4\\
4.475	4\\
4.4875	4\\
4.5	4\\
4.5125	4\\
4.525	4\\
4.5375	3\\
4.55	3\\
4.5625	5\\
4.575	4\\
4.5875	4\\
4.6	4\\
4.6125	5\\
4.625	5\\
4.6375	5\\
4.65	4\\
4.6625	5\\
4.675	4\\
4.6875	5\\
4.7	4\\
4.7125	5\\
4.725	5\\
4.7375	5\\
4.75	5\\
4.7625	5\\
4.775	5\\
4.7875	5\\
4.8	5\\
4.8125	5\\
4.825	5\\
4.8375	6\\
4.85	6\\
4.8625	6\\
4.875	6\\
4.8875	6\\
4.9	6\\
4.9125	6\\
4.925	6\\
4.9375	6\\
4.95	6\\
4.9625	6\\
4.975	6\\
4.9875	6\\
};
\addlegendentry{$\epsilon_u = 10^{-8}$}

\addplot [color=mycolor2, solid, line width=1.2pt]
  table[row sep=crcr]{%
0	18\\
0.0125	17\\
0.025	17\\
0.0375	17\\
0.05	17\\
0.0625	17\\
0.075	17\\
0.0875	17\\
0.1	17\\
0.1125	17\\
0.125	17\\
0.1375	17\\
0.15	17\\
0.1625	17\\
0.175	17\\
0.1875	17\\
0.2	17\\
0.2125	17\\
0.225	17\\
0.2375	17\\
0.25	17\\
0.2625	17\\
0.275	17\\
0.2875	17\\
0.3	17\\
0.3125	17\\
0.325	17\\
0.3375	17\\
0.35	17\\
0.3625	17\\
0.375	17\\
0.3875	17\\
0.4	17\\
0.4125	17\\
0.425	17\\
0.4375	17\\
0.45	17\\
0.4625	17\\
0.475	17\\
0.4875	17\\
0.5	18\\
0.5125	18\\
0.525	18\\
0.5375	18\\
0.55	18\\
0.5625	18\\
0.575	18\\
0.5875	18\\
0.6	19\\
0.6125	19\\
0.625	19\\
0.6375	19\\
0.65	19\\
0.6625	19\\
0.675	20\\
0.6875	20\\
0.7	21\\
0.7125	21\\
0.725	21\\
0.7375	21\\
0.75	21\\
0.7625	21\\
0.775	21\\
0.7875	21\\
0.8	21\\
0.8125	21\\
0.825	21\\
0.8375	21\\
0.85	21\\
0.8625	21\\
0.875	21\\
0.8875	21\\
0.9	21\\
0.9125	21\\
0.925	21\\
0.9375	21\\
0.95	21\\
0.9625	22\\
0.975	22\\
0.9875	22\\
1	22\\
1.0125	22\\
1.025	22\\
1.0375	22\\
1.05	22\\
1.0625	22\\
1.075	22\\
1.0875	22\\
1.1	22\\
1.1125	22\\
1.125	22\\
1.1375	22\\
1.15	22\\
1.1625	22\\
1.175	22\\
1.1875	22\\
1.2	22\\
1.2125	22\\
1.225	22\\
1.2375	22\\
1.25	22\\
1.2625	22\\
1.275	22\\
1.2875	22\\
1.3	22\\
1.3125	22\\
1.325	22\\
1.3375	22\\
1.35	22\\
1.3625	22\\
1.375	22\\
1.3875	22\\
1.4	22\\
1.4125	22\\
1.425	22\\
1.4375	22\\
1.45	22\\
1.4625	22\\
1.475	22\\
1.4875	22\\
1.5	22\\
1.5125	22\\
1.525	22\\
1.5375	22\\
1.55	22\\
1.5625	22\\
1.575	22\\
1.5875	22\\
1.6	22\\
1.6125	22\\
1.625	22\\
1.6375	22\\
1.65	22\\
1.6625	22\\
1.675	22\\
1.6875	22\\
1.7	22\\
1.7125	22\\
1.725	22\\
1.7375	23\\
1.75	23\\
1.7625	23\\
1.775	23\\
1.7875	23\\
1.8	23\\
1.8125	23\\
1.825	23\\
1.8375	23\\
1.85	23\\
1.8625	23\\
1.875	23\\
1.8875	23\\
1.9	23\\
1.9125	23\\
1.925	23\\
1.9375	23\\
1.95	23\\
1.9625	23\\
1.975	23\\
1.9875	23\\
2	23\\
2.0125	23\\
2.025	23\\
2.0375	23\\
2.05	23\\
2.0625	23\\
2.075	23\\
2.0875	23\\
2.1	23\\
2.1125	23\\
2.125	23\\
2.1375	23\\
2.15	23\\
2.1625	23\\
2.175	23\\
2.1875	23\\
2.2	23\\
2.2125	23\\
2.225	23\\
2.2375	23\\
2.25	23\\
2.2625	23\\
2.275	23\\
2.2875	23\\
2.3	23\\
2.3125	23\\
2.325	23\\
2.3375	23\\
2.35	23\\
2.3625	24\\
2.375	24\\
2.3875	24\\
2.4	24\\
2.4125	24\\
2.425	24\\
2.4375	24\\
2.45	24\\
2.4625	24\\
2.475	24\\
2.4875	24\\
2.5	24\\
2.5125	24\\
2.525	24\\
2.5375	24\\
2.55	24\\
2.5625	24\\
2.575	24\\
2.5875	24\\
2.6	24\\
2.6125	24\\
2.625	24\\
2.6375	24\\
2.65	24\\
2.6625	24\\
2.675	24\\
2.6875	24\\
2.7	24\\
2.7125	24\\
2.725	24\\
2.7375	24\\
2.75	24\\
2.7625	24\\
2.775	24\\
2.7875	24\\
2.8	24\\
2.8125	24\\
2.825	24\\
2.8375	24\\
2.85	24\\
2.8625	24\\
2.875	24\\
2.8875	24\\
2.9	24\\
2.9125	24\\
2.925	24\\
2.9375	24\\
2.95	24\\
2.9625	24\\
2.975	24\\
2.9875	24\\
3	24\\
3.0125	24\\
3.025	24\\
3.0375	24\\
3.05	24\\
3.0625	24\\
3.075	24\\
3.0875	24\\
3.1	24\\
3.1125	24\\
3.125	24\\
3.1375	24\\
3.15	24\\
3.1625	24\\
3.175	24\\
3.1875	24\\
3.2	24\\
3.2125	24\\
3.225	24\\
3.2375	24\\
3.25	24\\
3.2625	24\\
3.275	24\\
3.2875	24\\
3.3	24\\
3.3125	24\\
3.325	24\\
3.3375	24\\
3.35	24\\
3.3625	24\\
3.375	24\\
3.3875	24\\
3.4	24\\
3.4125	24\\
3.425	24\\
3.4375	24\\
3.45	24\\
3.4625	24\\
3.475	24\\
3.4875	24\\
3.5	24\\
3.5125	24\\
3.525	24\\
3.5375	24\\
3.55	24\\
3.5625	24\\
3.575	25\\
3.5875	25\\
3.6	25\\
3.6125	25\\
3.625	25\\
3.6375	25\\
3.65	25\\
3.6625	25\\
3.675	25\\
3.6875	25\\
3.7	25\\
3.7125	25\\
3.725	25\\
3.7375	25\\
3.75	25\\
3.7625	25\\
3.775	26\\
3.7875	26\\
3.8	26\\
3.8125	27\\
3.825	27\\
3.8375	27\\
3.85	27\\
3.8625	27\\
3.875	27\\
3.8875	27\\
3.9	27\\
3.9125	27\\
3.925	27\\
3.9375	27\\
3.95	27\\
3.9625	27\\
3.975	27\\
3.9875	27\\
4	27\\
4.0125	27\\
4.025	27\\
4.0375	27\\
4.05	27\\
4.0625	27\\
4.075	27\\
4.0875	27\\
4.1	27\\
4.1125	27\\
4.125	27\\
4.1375	27\\
4.15	27\\
4.1625	27\\
4.175	27\\
4.1875	27\\
4.2	27\\
4.2125	27\\
4.225	27\\
4.2375	27\\
4.25	27\\
4.2625	27\\
4.275	27\\
4.2875	27\\
4.3	27\\
4.3125	27\\
4.325	27\\
4.3375	27\\
4.35	27\\
4.3625	27\\
4.375	27\\
4.3875	27\\
4.4	27\\
4.4125	27\\
4.425	27\\
4.4375	27\\
4.45	27\\
4.4625	27\\
4.475	27\\
4.4875	27\\
4.5	27\\
4.5125	27\\
4.525	27\\
4.5375	27\\
4.55	27\\
4.5625	27\\
4.575	27\\
4.5875	27\\
4.6	27\\
4.6125	27\\
4.625	27\\
4.6375	27\\
4.65	27\\
4.6625	27\\
4.675	27\\
4.6875	27\\
4.7	27\\
4.7125	27\\
4.725	27\\
4.7375	27\\
4.75	27\\
4.7625	27\\
4.775	27\\
4.7875	27\\
4.8	27\\
4.8125	27\\
4.825	27\\
4.8375	26\\
4.85	26\\
4.8625	26\\
4.875	26\\
4.8875	26\\
4.9	26\\
4.9125	26\\
4.925	26\\
4.9375	25\\
4.95	25\\
4.9625	25\\
4.975	25\\
4.9875	25\\
}; 
\addlegendentry{$\epsilon_u = 10^{-11}$}  
  
\end{axis}

\end{tikzpicture}%

%% file: figures/rank_AC.tex
% This file was created with tikzplotlib v0.10.1.
\definecolor{mycolor1}{rgb}{0.00000,0.44700,0.74100}%
\definecolor{mycolor2}{rgb}{0.85000,0.32500,0.09800}%
\definecolor{mycolor3}{rgb}{0.5,0.5,0.5}
\begin{tikzpicture}[scale=0.8]

\begin{axis}[%
title={Allen-Cahn},
width=1.5in,
height=1.5in,
at={(1.011in,0.642in)},
scale only axis,
xmin=0,
xmax=200,
xlabel style={font=\color{white!15!black}},
xlabel={Time},
%ymode=log,
ymin=0,
ymax=50,
yminorticks=true,
ylabel style={font=\color{white!15!black}},
axis background/.style={fill=white},
legend style={fill opacity=0.8, 
font=\tiny,
  draw opacity=1,
  text opacity=1,
  at={(0.3,0.13)},
  anchor=west}
  % draw=lightgray204}
]

\addplot [line width=1.2pt,  mycolor1]
table {%
0 4
0.5 8
1 8
1.5 9
2 9
2.5 9
3 10
3.5 10
4 11
4.5 11
5 11
5.5 11
6 11
6.5 11
7 12
7.5 12
8 11
8.5 13
9 13
9.5 13
10 13
10.5 13
11 12
11.5 13
12 13
12.5 12
13 12
13.5 12
14 11
14.5 11
15 12
15.5 12
16 12
16.5 12
17 11
17.5 12
18 12
18.5 12
19 12
19.5 12
20 12
20.5 12
21 12
21.5 11
22 12
22.5 12
23 12
23.5 12
24 13
24.5 13
25 13
25.5 13
26 12
26.5 12
27 13
27.5 13
28 13
28.5 13
29 13
29.5 13
30 13
30.5 14
31 13
31.5 15
32 13
32.5 14
33 12
33.5 13
34 13
34.5 13
35 13
35.5 13
36 13
36.5 13
37 13
37.5 14
38 13
38.5 13
39 14
39.5 13
40 13
40.5 13
41 13
41.5 14
42 14
42.5 12
43 13
43.5 14
44 13
44.5 12
45 12
45.5 12
46 13
46.5 13
47 13
47.5 13
48 13
48.5 13
49 13
49.5 13
50 13
50.5 13
51 13
51.5 13
52 13
52.5 13
53 12
53.5 13
54 13
54.5 12
55 13
55.5 13
56 13
56.5 13
57 12
57.5 13
58 13
58.5 13
59 13
59.5 13
60 13
60.5 14
61 14
61.5 14
62 13
62.5 13
63 13
63.5 14
64 13
64.5 13
65 13
65.5 14
66 14
66.5 14
67 14
67.5 14
68 13
68.5 14
69 14
69.5 14
70 14
70.5 14
71 14
71.5 14
72 14
72.5 14
73 14
73.5 13
74 14
74.5 14
75 15
75.5 16
76 18
76.5 18
77 16
77.5 17
78 18
78.5 18
79 18
79.5 18
80 16
80.5 17
81 17
81.5 18
82 19
82.5 18
83 19
83.5 18
84 19
84.5 19
85 19
85.5 19
86 18
86.5 19
87 20
87.5 20
88 20
88.5 20
89 18
89.5 19
90 18
90.5 18
91 18
91.5 19
92 19
92.5 20
93 19
93.5 20
94 20
94.5 20
95 20
95.5 20
96 20
96.5 19
97 20
97.5 20
98 19
98.5 21
99 20
99.5 22
100 22
100.5 23
101 22
101.5 21
102 21
102.5 21
103 21
103.5 20
104 21
104.5 22
105 22
105.5 21
106 21
106.5 20
107 20
107.5 20
108 20
108.5 19
109 18
109.5 20
110 20
110.5 20
111 20
111.5 19
112 20
112.5 20
113 19
113.5 20
114 20
114.5 20
115 21
115.5 21
116 22
116.5 22
117 22
117.5 22
118 22
118.5 22
119 21
119.5 22
120 21
120.5 21
121 21
121.5 22
122 22
122.5 22
123 22
123.5 22
124 22
124.5 21
125 22
125.5 23
126 22
126.5 22
127 23
127.5 23
128 22
128.5 22
129 23
129.5 23
130 21
130.5 21
131 22
131.5 22
132 21
132.5 21
133 21
133.5 21
134 21
134.5 21
135 21
135.5 21
136 21
136.5 20
137 21
137.5 21
138 21
138.5 21
139 21
139.5 21
140 22
140.5 21
141 21
141.5 22
142 21
142.5 21
143 23
143.5 23
144 23
144.5 22
145 23
145.5 22
146 21
146.5 21
147 21
147.5 21
148 21
148.5 21
149 21
149.5 21
150 20
150.5 21
151 21
151.5 20
152 20
152.5 19
153 19
153.5 19
154 19
154.5 19
155 18
155.5 19
156 19
156.5 19
157 19
157.5 19
158 19
158.5 18
159 18
159.5 18
160 18
160.5 18
161 17
161.5 16
162 16
162.5 18
163 18
163.5 17
164 18
164.5 18
165 18
165.5 18
166 19
166.5 19
167 19
167.5 19
168 18
168.5 19
169 17
169.5 18
170 18
170.5 18
171 16
171.5 16
172 15
172.5 16
173 16
173.5 16
174 16
174.5 16
175 16
175.5 15
176 14
176.5 14
177 14
177.5 14
178 14
178.5 14
179 14
179.5 14
180 14
180.5 14
181 14
181.5 14
182 14
182.5 14
183 14
183.5 14
184 14
184.5 14
185 15
185.5 15
186 14
186.5 15
187 15
187.5 14
188 14
188.5 14
189 15
189.5 16
190 15
190.5 16
191 16
191.5 16
192 15
192.5 15
193 15
193.5 16
194 15
194.5 15
195 15
195.5 16
196 16
196.5 16
197 15
197.5 16
198 16
198.5 16
199 15
199.5 16
};

%\addlegendentry{$\epsilon_u = 10^{-8}$}
\addplot [line width=1.2pt,  mycolor2]
table {%
0 12
0.5 14
1 16
1.5 19
2 19
2.5 19
3 22
3.5 22
4 22
4.5 22
5 23
5.5 25
6 25
6.5 25
7 25
7.5 25
8 25
8.5 26
9 26
9.5 25
10 27
10.5 27
11 26
11.5 26
12 25
12.5 26
13 26
13.5 25
14 24
14.5 25
15 25
15.5 25
16 24
16.5 25
17 24
17.5 25
18 25
18.5 25
19 25
19.5 24
20 25
20.5 25
21 25
21.5 26
22 25
22.5 25
23 25
23.5 24
24 24
24.5 23
25 22
25.5 23
26 23
26.5 24
27 23
27.5 24
28 24
28.5 24
29 24
29.5 24
30 24
30.5 24
31 23
31.5 24
32 24
32.5 23
33 23
33.5 24
34 24
34.5 23
35 23
35.5 23
36 24
36.5 24
37 23
37.5 24
38 25
38.5 25
39 25
39.5 25
40 25
40.5 24
41 25
41.5 24
42 25
42.5 24
43 26
43.5 27
44 27
44.5 27
45 27
45.5 25
46 25
46.5 26
47 27
47.5 27
48 27
48.5 27
49 27
49.5 27
50 27
50.5 27
51 26
51.5 26
52 26
52.5 27
53 28
53.5 27
54 29
54.5 28
55 28
55.5 28
56 28
56.5 28
57 27
57.5 26
58 27
58.5 27
59 27
59.5 26
60 27
60.5 26
61 26
61.5 27
62 27
62.5 26
63 27
63.5 27
64 28
64.5 28
65 28
65.5 28
66 27
66.5 27
67 28
67.5 28
68 29
68.5 29
69 28
69.5 28
70 28
70.5 28
71 30
71.5 30
72 30
72.5 29
73 30
73.5 30
74 30
74.5 31
75 31
75.5 32
76 32
76.5 34
77 34
77.5 34
78 34
78.5 35
79 35
79.5 36
80 35
80.5 36
81 35
81.5 36
82 36
82.5 36
83 35
83.5 36
84 36
84.5 36
85 37
85.5 38
86 38
86.5 38
87 38
87.5 38
88 37
88.5 38
89 39
89.5 39
90 38
90.5 39
91 39
91.5 39
92 39
92.5 39
93 39
93.5 39
94 39
94.5 39
95 39
95.5 40
96 40
96.5 40
97 41
97.5 42
98 42
98.5 44
99 44
99.5 44
100 45
100.5 45
101 45
101.5 45
102 45
102.5 45
103 45
103.5 45
104 45
104.5 44
105 44
105.5 44
106 45
106.5 45
107 45
107.5 45
108 45
108.5 46
109 45
109.5 45
110 44
110.5 44
111 45
111.5 44
112 44
112.5 44
113 43
113.5 45
114 45
114.5 46
115 46
115.5 46
116 46
116.5 45
117 45
117.5 44
118 43
118.5 44
119 43
119.5 45
120 47
120.5 46
121 47
121.5 46
122 47
122.5 47
123 46
123.5 45
124 45
124.5 45
125 44
125.5 46
126 45
126.5 46
127 45
127.5 45
128 44
128.5 45
129 45
129.5 45
130 45
130.5 44
131 45
131.5 45
132 45
132.5 44
133 44
133.5 43
134 43
134.5 43
135 44
135.5 44
136 45
136.5 45
137 44
137.5 44
138 44
138.5 44
139 43
139.5 43
140 43
140.5 44
141 43
141.5 42
142 43
142.5 43
143 43
143.5 42
144 41
144.5 41
145 41
145.5 42
146 42
146.5 44
147 44
147.5 40
148 41
148.5 39
149 40
149.5 39
150 39
150.5 39
151 38
151.5 37
152 37
152.5 36
153 36
153.5 34
154 34
154.5 34
155 32
155.5 31
156 32
156.5 32
157 31
157.5 31
158 32
158.5 31
159 32
159.5 32
160 31
160.5 32
161 32
161.5 31
162 30
162.5 31
163 31
163.5 31
164 31
164.5 30
165 31
165.5 31
166 31
166.5 30
167 30
167.5 30
168 30
168.5 31
169 30
169.5 31
170 31
170.5 30
171 30
171.5 30
172 30
172.5 29
173 29
173.5 29
174 28
174.5 29
175 30
175.5 30
176 30
176.5 30
177 30
177.5 30
178 29
178.5 30
179 29
179.5 29
180 29
180.5 29
181 29
181.5 29
182 30
182.5 29
183 30
183.5 30
184 31
184.5 30
185 28
185.5 30
186 30
186.5 30
187 30
187.5 30
188 30
188.5 31
189 30
189.5 30
190 30
190.5 29
191 30
191.5 30
192 29
192.5 30
193 30
193.5 31
194 32
194.5 30
195 31
195.5 33
196 32
196.5 31
197 31
197.5 32
198 31
198.5 32
199 32
199.5 33
};
%\addlegendentry{$\epsilon_u = 10^{-11}$}
\end{axis}

\end{tikzpicture}

%% file: figures/mr_Burgers.tex
\definecolor{mycolor1}{rgb}{0.00000,0.44700,0.74100}%
\definecolor{mycolor2}{rgb}{0.85000,0.32500,0.09800}%
\definecolor{mycolor3}{rgb}{0.5,0.5,0.5}
\begin{tikzpicture}[scale=0.8]

\begin{axis}[%
width=1.5in,
height=1.5in,
at={(1.011in,0.642in)},
scale only axis,
xmin=0,
xmax=5,
xlabel style={font=\color{white!15!black}},
xlabel={Time},
%ymode=log,
ymin=0,
ymax=40,
yminorticks=true,
ylabel style={font=\color{white!15!black}},
ylabel={$m_r$},
axis background/.style={fill=white},
legend style={fill opacity=0.8, 
font=\tiny,
  draw opacity=1,
  text opacity=1,
  at={(0.02,0.87)},
  anchor=west}
  % draw=lightgray204}
]

\addplot [color=mycolor1, line width=1.2pt]
  table[row sep=crcr]{%
0	1\\
0.0025	1\\
0.005	1\\
0.0075	1\\
0.01	1\\
0.0125	1\\
0.015	1\\
0.0175	0\\
0.02	0\\
0.0225	0\\
0.025	0\\
0.0275	0\\
0.03	0\\
0.0325	0\\
0.035	0\\
0.0375	0\\
0.04	0\\
0.0425	0\\
0.045	0\\
0.0475	0\\
0.05	0\\
0.0525	0\\
0.055	0\\
0.0575	0\\
0.06	0\\
0.0625	0\\
0.065	0\\
0.0675	0\\
0.07	0\\
0.0725	0\\
0.075	0\\
0.0775	0\\
0.08	0\\
0.0825	0\\
0.085	0\\
0.0875	0\\
0.09	0\\
0.0925	0\\
0.095	0\\
0.0975	0\\
0.1	0\\
0.1025	0\\
0.105	0\\
0.1075	0\\
0.11	0\\
0.1125	0\\
0.115	0\\
0.1175	1\\
0.12	0\\
0.1225	0\\
0.125	0\\
0.1275	0\\
0.13	0\\
0.1325	0\\
0.135	0\\
0.1375	0\\
0.14	0\\
0.1425	0\\
0.145	0\\
0.1475	0\\
0.15	0\\
0.1525	0\\
0.155	0\\
0.1575	0\\
0.16	0\\
0.1625	0\\
0.165	0\\
0.1675	0\\
0.17	0\\
0.1725	0\\
0.175	0\\
0.1775	0\\
0.18	0\\
0.1825	0\\
0.185	0\\
0.1875	0\\
0.19	0\\
0.1925	0\\
0.195	0\\
0.1975	0\\
0.2	0\\
0.2025	0\\
0.205	0\\
0.2075	0\\
0.21	0\\
0.2125	0\\
0.215	0\\
0.2175	0\\
0.22	0\\
0.2225	0\\
0.225	0\\
0.2275	0\\
0.23	0\\
0.2325	0\\
0.235	0\\
0.2375	0\\
0.24	0\\
0.2425	0\\
0.245	0\\
0.2475	0\\
0.25	0\\
0.2525	0\\
0.255	0\\
0.2575	0\\
0.26	0\\
0.2625	0\\
0.265	0\\
0.2675	0\\
0.27	0\\
0.2725	0\\
0.275	0\\
0.2775	0\\
0.28	3\\
0.2825	1\\
0.285	1\\
0.2875	1\\
0.29	1\\
0.2925	1\\
0.295	1\\
0.2975	1\\
0.3	1\\
0.3025	1\\
0.305	1\\
0.3075	1\\
0.31	1\\
0.3125	1\\
0.315	1\\
0.3175	1\\
0.32	1\\
0.3225	1\\
0.325	1\\
0.3275	1\\
0.33	1\\
0.3325	1\\
0.335	1\\
0.3375	1\\
0.34	1\\
0.3425	1\\
0.345	1\\
0.3475	1\\
0.35	1\\
0.3525	1\\
0.355	1\\
0.3575	1\\
0.36	1\\
0.3625	1\\
0.365	1\\
0.3675	1\\
0.37	1\\
0.3725	1\\
0.375	1\\
0.3775	1\\
0.38	1\\
0.3825	1\\
0.385	1\\
0.3875	1\\
0.39	1\\
0.3925	1\\
0.395	1\\
0.3975	1\\
0.4	1\\
0.4025	1\\
0.405	1\\
0.4075	1\\
0.41	1\\
0.4125	0\\
0.415	0\\
0.4175	0\\
0.42	0\\
0.4225	0\\
0.425	0\\
0.4275	0\\
0.43	0\\
0.4325	0\\
0.435	0\\
0.4375	0\\
0.44	0\\
0.4425	0\\
0.445	0\\
0.4475	0\\
0.45	0\\
0.4525	0\\
0.455	0\\
0.4575	0\\
0.46	0\\
0.4625	1\\
0.465	1\\
0.4675	0\\
0.47	1\\
0.4725	1\\
0.475	1\\
0.4775	1\\
0.48	0\\
0.4825	0\\
0.485	0\\
0.4875	0\\
0.49	0\\
0.4925	0\\
0.495	0\\
0.4975	0\\
0.5	0\\
0.5025	0\\
0.505	0\\
0.5075	0\\
0.51	0\\
0.5125	0\\
0.515	1\\
0.5175	1\\
0.52	1\\
0.5225	1\\
0.525	1\\
0.5275	1\\
0.53	1\\
0.5325	1\\
0.535	1\\
0.5375	1\\
0.54	1\\
0.5425	1\\
0.545	1\\
0.5475	1\\
0.55	1\\
0.5525	1\\
0.555	1\\
0.5575	1\\
0.56	1\\
0.5625	1\\
0.565	1\\
0.5675	1\\
0.57	1\\
0.5725	1\\
0.575	2\\
0.5775	1\\
0.58	0\\
0.5825	0\\
0.585	0\\
0.5875	0\\
0.59	0\\
0.5925	0\\
0.595	0\\
0.5975	0\\
0.6	0\\
0.6025	0\\
0.605	0\\
0.6075	0\\
0.61	0\\
0.6125	0\\
0.615	0\\
0.6175	0\\
0.62	0\\
0.6225	0\\
0.625	0\\
0.6275	0\\
0.63	0\\
0.6325	0\\
0.635	0\\
0.6375	0\\
0.64	0\\
0.6425	0\\
0.645	0\\
0.6475	0\\
0.65	0\\
0.6525	0\\
0.655	0\\
0.6575	0\\
0.66	0\\
0.6625	0\\
0.665	0\\
0.6675	0\\
0.67	0\\
0.6725	0\\
0.675	0\\
0.6775	0\\
0.68	0\\
0.6825	0\\
0.685	0\\
0.6875	0\\
0.69	0\\
0.6925	0\\
0.695	0\\
0.6975	0\\
0.7	0\\
0.7025	0\\
0.705	0\\
0.7075	0\\
0.71	0\\
0.7125	0\\
0.715	0\\
0.7175	0\\
0.72	0\\
0.7225	0\\
0.725	0\\
0.7275	0\\
0.73	0\\
0.7325	0\\
0.735	0\\
0.7375	0\\
0.74	0\\
0.7425	0\\
0.745	0\\
0.7475	0\\
0.75	0\\
0.7525	0\\
0.755	0\\
0.7575	0\\
0.76	0\\
0.7625	0\\
0.765	0\\
0.7675	0\\
0.77	0\\
0.7725	0\\
0.775	0\\
0.7775	0\\
0.78	0\\
0.7825	0\\
0.785	0\\
0.7875	0\\
0.79	0\\
0.7925	0\\
0.795	0\\
0.7975	0\\
0.8	0\\
0.8025	0\\
0.805	0\\
0.8075	0\\
0.81	0\\
0.8125	0\\
0.815	0\\
0.8175	0\\
0.82	0\\
0.8225	0\\
0.825	0\\
0.8275	1\\
0.83	1\\
0.8325	1\\
0.835	1\\
0.8375	1\\
0.84	1\\
0.8425	1\\
0.845	1\\
0.8475	1\\
0.85	0\\
0.8525	0\\
0.855	0\\
0.8575	0\\
0.86	0\\
0.8625	0\\
0.865	0\\
0.8675	0\\
0.87	0\\
0.8725	0\\
0.875	0\\
0.8775	0\\
0.88	0\\
0.8825	0\\
0.885	6\\
0.8875	0\\
0.89	0\\
0.8925	0\\
0.895	0\\
0.8975	0\\
0.9	0\\
0.9025	0\\
0.905	0\\
0.9075	0\\
0.91	0\\
0.9125	0\\
0.915	0\\
0.9175	0\\
0.92	0\\
0.9225	0\\
0.925	0\\
0.9275	0\\
0.93	0\\
0.9325	0\\
0.935	0\\
0.9375	0\\
0.94	0\\
0.9425	0\\
0.945	0\\
0.9475	0\\
0.95	0\\
0.9525	0\\
0.955	0\\
0.9575	0\\
0.96	0\\
0.9625	0\\
0.965	0\\
0.9675	1\\
0.97	1\\
0.9725	1\\
0.975	1\\
0.9775	1\\
0.98	1\\
0.9825	1\\
0.985	1\\
0.9875	1\\
0.99	1\\
0.9925	1\\
0.995	1\\
0.9975	1\\
1	1\\
1.0025	1\\
1.005	1\\
1.0075	1\\
1.01	1\\
1.0125	1\\
1.015	1\\
1.0175	1\\
1.02	1\\
1.0225	1\\
1.025	1\\
1.0275	1\\
1.03	1\\
1.0325	0\\
1.035	0\\
1.0375	0\\
1.04	0\\
1.0425	0\\
1.045	0\\
1.0475	0\\
1.05	0\\
1.0525	0\\
1.055	0\\
1.0575	0\\
1.06	0\\
1.0625	0\\
1.065	0\\
1.0675	0\\
1.07	0\\
1.0725	0\\
1.075	0\\
1.0775	0\\
1.08	0\\
1.0825	0\\
1.085	0\\
1.0875	0\\
1.09	0\\
1.0925	0\\
1.095	0\\
1.0975	0\\
1.1	0\\
1.1025	0\\
1.105	0\\
1.1075	0\\
1.11	1\\
1.1125	1\\
1.115	1\\
1.1175	1\\
1.12	1\\
1.1225	1\\
1.125	1\\
1.1275	1\\
1.13	1\\
1.1325	1\\
1.135	1\\
1.1375	1\\
1.14	1\\
1.1425	1\\
1.145	1\\
1.1475	1\\
1.15	1\\
1.1525	1\\
1.155	1\\
1.1575	1\\
1.16	1\\
1.1625	1\\
1.165	1\\
1.1675	1\\
1.17	1\\
1.1725	1\\
1.175	1\\
1.1775	1\\
1.18	0\\
1.1825	0\\
1.185	0\\
1.1875	0\\
1.19	0\\
1.1925	0\\
1.195	0\\
1.1975	0\\
1.2	0\\
1.2025	0\\
1.205	0\\
1.2075	0\\
1.21	8\\
1.2125	0\\
1.215	0\\
1.2175	0\\
1.22	0\\
1.2225	0\\
1.225	0\\
1.2275	0\\
1.23	0\\
1.2325	0\\
1.235	0\\
1.2375	0\\
1.24	0\\
1.2425	0\\
1.245	0\\
1.2475	0\\
1.25	0\\
1.2525	0\\
1.255	0\\
1.2575	0\\
1.26	0\\
1.2625	0\\
1.265	0\\
1.2675	0\\
1.27	0\\
1.2725	0\\
1.275	0\\
1.2775	0\\
1.28	0\\
1.2825	0\\
1.285	0\\
1.2875	0\\
1.29	0\\
1.2925	0\\
1.295	0\\
1.2975	0\\
1.3	0\\
1.3025	0\\
1.305	0\\
1.3075	0\\
1.31	0\\
1.3125	0\\
1.315	0\\
1.3175	0\\
1.32	0\\
1.3225	0\\
1.325	0\\
1.3275	0\\
1.33	0\\
1.3325	0\\
1.335	0\\
1.3375	1\\
1.34	0\\
1.3425	0\\
1.345	0\\
1.3475	0\\
1.35	0\\
1.3525	0\\
1.355	0\\
1.3575	0\\
1.36	0\\
1.3625	0\\
1.365	0\\
1.3675	0\\
1.37	0\\
1.3725	0\\
1.375	0\\
1.3775	0\\
1.38	0\\
1.3825	0\\
1.385	0\\
1.3875	0\\
1.39	0\\
1.3925	0\\
1.395	0\\
1.3975	0\\
1.4	0\\
1.4025	0\\
1.405	0\\
1.4075	0\\
1.41	0\\
1.4125	0\\
1.415	0\\
1.4175	0\\
1.42	0\\
1.4225	0\\
1.425	0\\
1.4275	0\\
1.43	0\\
1.4325	0\\
1.435	0\\
1.4375	0\\
1.44	0\\
1.4425	0\\
1.445	0\\
1.4475	0\\
1.45	0\\
1.4525	0\\
1.455	0\\
1.4575	0\\
1.46	0\\
1.4625	0\\
1.465	0\\
1.4675	0\\
1.47	0\\
1.4725	0\\
1.475	0\\
1.4775	0\\
1.48	0\\
1.4825	0\\
1.485	0\\
1.4875	0\\
1.49	0\\
1.4925	0\\
1.495	0\\
1.4975	0\\
1.5	0\\
1.5025	0\\
1.505	0\\
1.5075	0\\
1.51	0\\
1.5125	0\\
1.515	0\\
1.5175	0\\
1.52	0\\
1.5225	0\\
1.525	0\\
1.5275	0\\
1.53	0\\
1.5325	0\\
1.535	0\\
1.5375	0\\
1.54	0\\
1.5425	0\\
1.545	0\\
1.5475	0\\
1.55	0\\
1.5525	0\\
1.555	1\\
1.5575	1\\
1.56	1\\
1.5625	1\\
1.565	0\\
1.5675	0\\
1.57	0\\
1.5725	0\\
1.575	0\\
1.5775	0\\
1.58	0\\
1.5825	0\\
1.585	0\\
1.5875	0\\
1.59	0\\
1.5925	0\\
1.595	0\\
1.5975	0\\
1.6	0\\
1.6025	0\\
1.605	0\\
1.6075	0\\
1.61	0\\
1.6125	0\\
1.615	0\\
1.6175	0\\
1.62	0\\
1.6225	0\\
1.625	0\\
1.6275	0\\
1.63	0\\
1.6325	0\\
1.635	0\\
1.6375	0\\
1.64	0\\
1.6425	0\\
1.645	0\\
1.6475	0\\
1.65	0\\
1.6525	0\\
1.655	0\\
1.6575	0\\
1.66	0\\
1.6625	0\\
1.665	0\\
1.6675	0\\
1.67	0\\
1.6725	0\\
1.675	0\\
1.6775	0\\
1.68	0\\
1.6825	0\\
1.685	0\\
1.6875	0\\
1.69	0\\
1.6925	0\\
1.695	0\\
1.6975	0\\
1.7	0\\
1.7025	0\\
1.705	0\\
1.7075	0\\
1.71	0\\
1.7125	0\\
1.715	0\\
1.7175	0\\
1.72	0\\
1.7225	0\\
1.725	0\\
1.7275	0\\
1.73	0\\
1.7325	0\\
1.735	0\\
1.7375	0\\
1.74	0\\
1.7425	0\\
1.745	0\\
1.7475	0\\
1.75	0\\
1.7525	0\\
1.755	0\\
1.7575	0\\
1.76	0\\
1.7625	0\\
1.765	0\\
1.7675	0\\
1.77	0\\
1.7725	0\\
1.775	1\\
1.7775	2\\
1.78	1\\
1.7825	1\\
1.785	1\\
1.7875	1\\
1.79	1\\
1.7925	1\\
1.795	1\\
1.7975	1\\
1.8	1\\
1.8025	1\\
1.805	1\\
1.8075	1\\
1.81	1\\
1.8125	1\\
1.815	1\\
1.8175	1\\
1.82	0\\
1.8225	0\\
1.825	0\\
1.8275	0\\
1.83	0\\
1.8325	0\\
1.835	0\\
1.8375	0\\
1.84	0\\
1.8425	0\\
1.845	0\\
1.8475	0\\
1.85	0\\
1.8525	0\\
1.855	0\\
1.8575	0\\
1.86	0\\
1.8625	0\\
1.865	0\\
1.8675	0\\
1.87	0\\
1.8725	0\\
1.875	0\\
1.8775	0\\
1.88	0\\
1.8825	0\\
1.885	0\\
1.8875	0\\
1.89	0\\
1.8925	0\\
1.895	1\\
1.8975	1\\
1.9	1\\
1.9025	1\\
1.905	1\\
1.9075	1\\
1.91	1\\
1.9125	1\\
1.915	1\\
1.9175	1\\
1.92	1\\
1.9225	1\\
1.925	1\\
1.9275	1\\
1.93	1\\
1.9325	1\\
1.935	1\\
1.9375	1\\
1.94	1\\
1.9425	1\\
1.945	1\\
1.9475	1\\
1.95	1\\
1.9525	1\\
1.955	1\\
1.9575	0\\
1.96	0\\
1.9625	0\\
1.965	0\\
1.9675	1\\
1.97	0\\
1.9725	0\\
1.975	0\\
1.9775	0\\
1.98	0\\
1.9825	0\\
1.985	0\\
1.9875	0\\
1.99	0\\
1.9925	0\\
1.995	0\\
1.9975	0\\
2	0\\
2.0025	1\\
2.005	1\\
2.0075	1\\
2.01	1\\
2.0125	1\\
2.015	1\\
2.0175	1\\
2.02	1\\
2.0225	1\\
2.025	1\\
2.0275	1\\
2.03	1\\
2.0325	1\\
2.035	1\\
2.0375	1\\
2.04	1\\
2.0425	1\\
2.045	1\\
2.0475	1\\
2.05	1\\
2.0525	1\\
2.055	1\\
2.0575	1\\
2.06	1\\
2.0625	1\\
2.065	1\\
2.0675	0\\
2.07	1\\
2.0725	0\\
2.075	1\\
2.0775	0\\
2.08	1\\
2.0825	0\\
2.085	0\\
2.0875	0\\
2.09	0\\
2.0925	0\\
2.095	0\\
2.0975	0\\
2.1	0\\
2.1025	0\\
2.105	0\\
2.1075	0\\
2.11	0\\
2.1125	0\\
2.115	0\\
2.1175	0\\
2.12	0\\
2.1225	0\\
2.125	0\\
2.1275	0\\
2.13	0\\
2.1325	0\\
2.135	0\\
2.1375	0\\
2.14	0\\
2.1425	0\\
2.145	0\\
2.1475	0\\
2.15	0\\
2.1525	0\\
2.155	0\\
2.1575	0\\
2.16	0\\
2.1625	0\\
2.165	0\\
2.1675	0\\
2.17	0\\
2.1725	0\\
2.175	0\\
2.1775	0\\
2.18	0\\
2.1825	0\\
2.185	0\\
2.1875	0\\
2.19	0\\
2.1925	0\\
2.195	0\\
2.1975	0\\
2.2	0\\
2.2025	0\\
2.205	0\\
2.2075	0\\
2.21	0\\
2.2125	0\\
2.215	0\\
2.2175	0\\
2.22	0\\
2.2225	0\\
2.225	0\\
2.2275	0\\
2.23	0\\
2.2325	0\\
2.235	0\\
2.2375	0\\
2.24	0\\
2.2425	0\\
2.245	0\\
2.2475	0\\
2.25	0\\
2.2525	0\\
2.255	0\\
2.2575	0\\
2.26	0\\
2.2625	0\\
2.265	0\\
2.2675	0\\
2.27	0\\
2.2725	0\\
2.275	0\\
2.2775	0\\
2.28	0\\
2.2825	0\\
2.285	0\\
2.2875	0\\
2.29	0\\
2.2925	0\\
2.295	0\\
2.2975	0\\
2.3	0\\
2.3025	0\\
2.305	0\\
2.3075	0\\
2.31	0\\
2.3125	0\\
2.315	0\\
2.3175	0\\
2.32	0\\
2.3225	0\\
2.325	0\\
2.3275	0\\
2.33	0\\
2.3325	0\\
2.335	0\\
2.3375	0\\
2.34	0\\
2.3425	0\\
2.345	0\\
2.3475	0\\
2.35	0\\
2.3525	0\\
2.355	0\\
2.3575	0\\
2.36	0\\
2.3625	0\\
2.365	0\\
2.3675	0\\
2.37	0\\
2.3725	0\\
2.375	0\\
2.3775	0\\
2.38	0\\
2.3825	0\\
2.385	0\\
2.3875	0\\
2.39	0\\
2.3925	0\\
2.395	0\\
2.3975	0\\
2.4	0\\
2.4025	0\\
2.405	0\\
2.4075	0\\
2.41	0\\
2.4125	0\\
2.415	0\\
2.4175	0\\
2.42	0\\
2.4225	0\\
2.425	0\\
2.4275	0\\
2.43	0\\
2.4325	0\\
2.435	0\\
2.4375	0\\
2.44	0\\
2.4425	21\\
2.445	0\\
2.4475	0\\
2.45	0\\
2.4525	0\\
2.455	0\\
2.4575	0\\
2.46	0\\
2.4625	0\\
2.465	0\\
2.4675	0\\
2.47	0\\
2.4725	0\\
2.475	0\\
2.4775	0\\
2.48	0\\
2.4825	0\\
2.485	0\\
2.4875	0\\
2.49	0\\
2.4925	0\\
2.495	0\\
2.4975	0\\
2.5	0\\
2.5025	0\\
2.505	0\\
2.5075	0\\
2.51	0\\
2.5125	0\\
2.515	0\\
2.5175	0\\
2.52	0\\
2.5225	0\\
2.525	0\\
2.5275	0\\
2.53	0\\
2.5325	0\\
2.535	0\\
2.5375	0\\
2.54	0\\
2.5425	0\\
2.545	0\\
2.5475	0\\
2.55	0\\
2.5525	0\\
2.555	0\\
2.5575	0\\
2.56	0\\
2.5625	0\\
2.565	0\\
2.5675	0\\
2.57	0\\
2.5725	0\\
2.575	0\\
2.5775	0\\
2.58	0\\
2.5825	0\\
2.585	0\\
2.5875	0\\
2.59	0\\
2.5925	0\\
2.595	0\\
2.5975	0\\
2.6	0\\
2.6025	0\\
2.605	0\\
2.6075	0\\
2.61	0\\
2.6125	0\\
2.615	0\\
2.6175	0\\
2.62	0\\
2.6225	0\\
2.625	0\\
2.6275	0\\
2.63	1\\
2.6325	0\\
2.635	0\\
2.6375	0\\
2.64	0\\
2.6425	0\\
2.645	0\\
2.6475	0\\
2.65	0\\
2.6525	0\\
2.655	0\\
2.6575	0\\
2.66	0\\
2.6625	0\\
2.665	0\\
2.6675	0\\
2.67	0\\
2.6725	0\\
2.675	0\\
2.6775	0\\
2.68	0\\
2.6825	0\\
2.685	0\\
2.6875	0\\
2.69	0\\
2.6925	0\\
2.695	0\\
2.6975	0\\
2.7	0\\
2.7025	0\\
2.705	0\\
2.7075	0\\
2.71	0\\
2.7125	0\\
2.715	1\\
2.7175	0\\
2.72	1\\
2.7225	0\\
2.725	1\\
2.7275	0\\
2.73	1\\
2.7325	0\\
2.735	0\\
2.7375	0\\
2.74	0\\
2.7425	0\\
2.745	0\\
2.7475	0\\
2.75	0\\
2.7525	0\\
2.755	0\\
2.7575	0\\
2.76	0\\
2.7625	0\\
2.765	0\\
2.7675	0\\
2.77	0\\
2.7725	0\\
2.775	0\\
2.7775	0\\
2.78	0\\
2.7825	0\\
2.785	0\\
2.7875	0\\
2.79	0\\
2.7925	0\\
2.795	0\\
2.7975	0\\
2.8	0\\
2.8025	0\\
2.805	0\\
2.8075	0\\
2.81	0\\
2.8125	0\\
2.815	0\\
2.8175	0\\
2.82	0\\
2.8225	0\\
2.825	0\\
2.8275	0\\
2.83	0\\
2.8325	0\\
2.835	0\\
2.8375	0\\
2.84	1\\
2.8425	1\\
2.845	1\\
2.8475	1\\
2.85	1\\
2.8525	1\\
2.855	1\\
2.8575	1\\
2.86	1\\
2.8625	1\\
2.865	1\\
2.8675	1\\
2.87	1\\
2.8725	1\\
2.875	1\\
2.8775	1\\
2.88	1\\
2.8825	1\\
2.885	1\\
2.8875	1\\
2.89	1\\
2.8925	1\\
2.895	1\\
2.8975	1\\
2.9	1\\
2.9025	1\\
2.905	1\\
2.9075	1\\
2.91	1\\
2.9125	1\\
2.915	1\\
2.9175	1\\
2.92	1\\
2.9225	1\\
2.925	0\\
2.9275	1\\
2.93	0\\
2.9325	0\\
2.935	0\\
2.9375	0\\
2.94	0\\
2.9425	0\\
2.945	1\\
2.9475	1\\
2.95	1\\
2.9525	1\\
2.955	1\\
2.9575	1\\
2.96	1\\
2.9625	1\\
2.965	1\\
2.9675	1\\
2.97	1\\
2.9725	1\\
2.975	1\\
2.9775	1\\
2.98	1\\
2.9825	1\\
2.985	1\\
2.9875	1\\
2.99	1\\
2.9925	1\\
2.995	1\\
2.9975	1\\
3	1\\
3.0025	1\\
3.005	0\\
3.0075	0\\
3.01	0\\
3.0125	0\\
3.015	0\\
3.0175	0\\
3.02	0\\
3.0225	0\\
3.025	0\\
3.0275	0\\
3.03	0\\
3.0325	0\\
3.035	0\\
3.0375	0\\
3.04	0\\
3.0425	0\\
3.045	0\\
3.0475	0\\
3.05	0\\
3.0525	0\\
3.055	0\\
3.0575	0\\
3.06	0\\
3.0625	0\\
3.065	0\\
3.0675	0\\
3.07	0\\
3.0725	0\\
3.075	0\\
3.0775	0\\
3.08	0\\
3.0825	0\\
3.085	0\\
3.0875	0\\
3.09	0\\
3.0925	0\\
3.095	0\\
3.0975	0\\
3.1	0\\
3.1025	0\\
3.105	0\\
3.1075	0\\
3.11	0\\
3.1125	0\\
3.115	0\\
3.1175	0\\
3.12	0\\
3.1225	0\\
3.125	0\\
3.1275	0\\
3.13	0\\
3.1325	0\\
3.135	0\\
3.1375	0\\
3.14	0\\
3.1425	0\\
3.145	0\\
3.1475	0\\
3.15	0\\
3.1525	0\\
3.155	0\\
3.1575	0\\
3.16	0\\
3.1625	0\\
3.165	1\\
3.1675	0\\
3.17	0\\
3.1725	0\\
3.175	0\\
3.1775	1\\
3.18	0\\
3.1825	0\\
3.185	0\\
3.1875	1\\
3.19	0\\
3.1925	0\\
3.195	0\\
3.1975	0\\
3.2	0\\
3.2025	0\\
3.205	0\\
3.2075	1\\
3.21	0\\
3.2125	0\\
3.215	0\\
3.2175	0\\
3.22	1\\
3.2225	0\\
3.225	0\\
3.2275	0\\
3.23	0\\
3.2325	3\\
3.235	0\\
3.2375	0\\
3.24	0\\
3.2425	0\\
3.245	0\\
3.2475	0\\
3.25	0\\
3.2525	0\\
3.255	0\\
3.2575	0\\
3.26	0\\
3.2625	0\\
3.265	0\\
3.2675	0\\
3.27	0\\
3.2725	0\\
3.275	0\\
3.2775	0\\
3.28	1\\
3.2825	0\\
3.285	1\\
3.2875	0\\
3.29	2\\
3.2925	0\\
3.295	0\\
3.2975	0\\
3.3	0\\
3.3025	0\\
3.305	0\\
3.3075	0\\
3.31	0\\
3.3125	0\\
3.315	0\\
3.3175	0\\
3.32	0\\
3.3225	0\\
3.325	0\\
3.3275	0\\
3.33	0\\
3.3325	0\\
3.335	0\\
3.3375	1\\
3.34	0\\
3.3425	0\\
3.345	0\\
3.3475	1\\
3.35	1\\
3.3525	0\\
3.355	0\\
3.3575	1\\
3.36	1\\
3.3625	0\\
3.365	0\\
3.3675	1\\
3.37	0\\
3.3725	1\\
3.375	1\\
3.3775	0\\
3.38	1\\
3.3825	0\\
3.385	1\\
3.3875	1\\
3.39	1\\
3.3925	1\\
3.395	1\\
3.3975	1\\
3.4	1\\
3.4025	1\\
3.405	1\\
3.4075	1\\
3.41	1\\
3.4125	1\\
3.415	1\\
3.4175	1\\
3.42	1\\
3.4225	1\\
3.425	1\\
3.4275	1\\
3.43	1\\
3.4325	1\\
3.435	1\\
3.4375	1\\
3.44	1\\
3.4425	1\\
3.445	1\\
3.4475	1\\
3.45	1\\
3.4525	1\\
3.455	1\\
3.4575	1\\
3.46	1\\
3.4625	1\\
3.465	1\\
3.4675	1\\
3.47	1\\
3.4725	1\\
3.475	1\\
3.4775	1\\
3.48	1\\
3.4825	1\\
3.485	1\\
3.4875	1\\
3.49	1\\
3.4925	1\\
3.495	1\\
3.4975	1\\
3.5	1\\
3.5025	1\\
3.505	1\\
3.5075	1\\
3.51	1\\
3.5125	1\\
3.515	1\\
3.5175	1\\
3.52	1\\
3.5225	1\\
3.525	1\\
3.5275	1\\
3.53	1\\
3.5325	1\\
3.535	1\\
3.5375	1\\
3.54	1\\
3.5425	1\\
3.545	1\\
3.5475	1\\
3.55	1\\
3.5525	1\\
3.555	1\\
3.5575	1\\
3.56	1\\
3.5625	1\\
3.565	1\\
3.5675	1\\
3.57	1\\
3.5725	1\\
3.575	1\\
3.5775	1\\
3.58	1\\
3.5825	1\\
3.585	1\\
3.5875	1\\
3.59	1\\
3.5925	1\\
3.595	1\\
3.5975	1\\
3.6	1\\
3.6025	3\\
3.605	1\\
3.6075	1\\
3.61	1\\
3.6125	1\\
3.615	1\\
3.6175	0\\
3.62	1\\
3.6225	1\\
3.625	1\\
3.6275	0\\
3.63	0\\
3.6325	1\\
3.635	5\\
3.6375	1\\
3.64	0\\
3.6425	0\\
3.645	0\\
3.6475	0\\
3.65	0\\
3.6525	0\\
3.655	0\\
3.6575	0\\
3.66	0\\
3.6625	0\\
3.665	1\\
3.6675	1\\
3.67	0\\
3.6725	0\\
3.675	0\\
3.6775	0\\
3.68	0\\
3.6825	0\\
3.685	0\\
3.6875	0\\
3.69	0\\
3.6925	1\\
3.695	0\\
3.6975	1\\
3.7	1\\
3.7025	1\\
3.705	1\\
3.7075	1\\
3.71	1\\
3.7125	1\\
3.715	1\\
3.7175	1\\
3.72	1\\
3.7225	1\\
3.725	1\\
3.7275	1\\
3.73	1\\
3.7325	1\\
3.735	1\\
3.7375	1\\
3.74	1\\
3.7425	1\\
3.745	1\\
3.7475	1\\
3.75	1\\
3.7525	1\\
3.755	1\\
3.7575	1\\
3.76	1\\
3.7625	1\\
3.765	1\\
3.7675	1\\
3.77	1\\
3.7725	1\\
3.775	1\\
3.7775	1\\
3.78	1\\
3.7825	1\\
3.785	1\\
3.7875	1\\
3.79	1\\
3.7925	1\\
3.795	1\\
3.7975	1\\
3.8	1\\
3.8025	1\\
3.805	1\\
3.8075	1\\
3.81	1\\
3.8125	1\\
3.815	1\\
3.8175	1\\
3.82	1\\
3.8225	1\\
3.825	1\\
3.8275	1\\
3.83	1\\
3.8325	1\\
3.835	1\\
3.8375	1\\
3.84	1\\
3.8425	1\\
3.845	1\\
3.8475	1\\
3.85	1\\
3.8525	1\\
3.855	1\\
3.8575	1\\
3.86	1\\
3.8625	0\\
3.865	0\\
3.8675	0\\
3.87	0\\
3.8725	0\\
3.875	0\\
3.8775	0\\
3.88	0\\
3.8825	0\\
3.885	0\\
3.8875	0\\
3.89	0\\
3.8925	0\\
3.895	0\\
3.8975	0\\
3.9	0\\
3.9025	0\\
3.905	0\\
3.9075	0\\
3.91	0\\
3.9125	0\\
3.915	0\\
3.9175	0\\
3.92	0\\
3.9225	0\\
3.925	0\\
3.9275	0\\
3.93	0\\
3.9325	0\\
3.935	1\\
3.9375	1\\
3.94	1\\
3.9425	1\\
3.945	1\\
3.9475	1\\
3.95	1\\
3.9525	1\\
3.955	1\\
3.9575	1\\
3.96	1\\
3.9625	1\\
3.965	1\\
3.9675	1\\
3.97	1\\
3.9725	1\\
3.975	1\\
3.9775	1\\
3.98	1\\
3.9825	1\\
3.985	1\\
3.9875	1\\
3.99	1\\
3.9925	1\\
3.995	1\\
3.9975	1\\
4	1\\
4.0025	1\\
4.005	1\\
4.0075	1\\
4.01	1\\
4.0125	1\\
4.015	1\\
4.0175	1\\
4.02	1\\
4.0225	1\\
4.025	1\\
4.0275	1\\
4.03	1\\
4.0325	1\\
4.035	1\\
4.0375	1\\
4.04	1\\
4.0425	1\\
4.045	1\\
4.0475	1\\
4.05	1\\
4.0525	1\\
4.055	1\\
4.0575	1\\
4.06	1\\
4.0625	1\\
4.065	1\\
4.0675	1\\
4.07	1\\
4.0725	1\\
4.075	1\\
4.0775	1\\
4.08	1\\
4.0825	1\\
4.085	1\\
4.0875	1\\
4.09	1\\
4.0925	1\\
4.095	1\\
4.0975	1\\
4.1	1\\
4.1025	1\\
4.105	1\\
4.1075	1\\
4.11	1\\
4.1125	1\\
4.115	1\\
4.1175	1\\
4.12	1\\
4.1225	1\\
4.125	1\\
4.1275	1\\
4.13	1\\
4.1325	1\\
4.135	1\\
4.1375	1\\
4.14	1\\
4.1425	1\\
4.145	1\\
4.1475	1\\
4.15	1\\
4.1525	1\\
4.155	1\\
4.1575	1\\
4.16	1\\
4.1625	1\\
4.165	1\\
4.1675	1\\
4.17	0\\
4.1725	0\\
4.175	1\\
4.1775	2\\
4.18	1\\
4.1825	1\\
4.185	0\\
4.1875	0\\
4.19	0\\
4.1925	0\\
4.195	1\\
4.1975	0\\
4.2	1\\
4.2025	0\\
4.205	0\\
4.2075	0\\
4.21	0\\
4.2125	0\\
4.215	0\\
4.2175	0\\
4.22	0\\
4.2225	0\\
4.225	0\\
4.2275	0\\
4.23	0\\
4.2325	0\\
4.235	0\\
4.2375	0\\
4.24	0\\
4.2425	0\\
4.245	0\\
4.2475	0\\
4.25	0\\
4.2525	0\\
4.255	0\\
4.2575	0\\
4.26	0\\
4.2625	0\\
4.265	0\\
4.2675	0\\
4.27	0\\
4.2725	0\\
4.275	0\\
4.2775	0\\
4.28	0\\
4.2825	0\\
4.285	0\\
4.2875	0\\
4.29	0\\
4.2925	0\\
4.295	0\\
4.2975	1\\
4.3	0\\
4.3025	0\\
4.305	0\\
4.3075	0\\
4.31	0\\
4.3125	0\\
4.315	0\\
4.3175	0\\
4.32	0\\
4.3225	0\\
4.325	0\\
4.3275	0\\
4.33	0\\
4.3325	0\\
4.335	0\\
4.3375	0\\
4.34	0\\
4.3425	1\\
4.345	1\\
4.3475	1\\
4.35	1\\
4.3525	1\\
4.355	1\\
4.3575	1\\
4.36	1\\
4.3625	1\\
4.365	1\\
4.3675	1\\
4.37	1\\
4.3725	1\\
4.375	1\\
4.3775	27\\
4.38	1\\
4.3825	2\\
4.385	1\\
4.3875	1\\
4.39	1\\
4.3925	1\\
4.395	1\\
4.3975	1\\
4.4	1\\
4.4025	1\\
4.405	1\\
4.4075	1\\
4.41	1\\
4.4125	1\\
4.415	1\\
4.4175	1\\
4.42	1\\
4.4225	1\\
4.425	1\\
4.4275	1\\
4.43	1\\
4.4325	1\\
4.435	2\\
4.4375	1\\
4.44	2\\
4.4425	2\\
4.445	2\\
4.4475	2\\
4.45	2\\
4.4525	1\\
4.455	1\\
4.4575	1\\
4.46	1\\
4.4625	1\\
4.465	1\\
4.4675	1\\
4.47	1\\
4.4725	1\\
4.475	1\\
4.4775	1\\
4.48	1\\
4.4825	1\\
4.485	1\\
4.4875	1\\
4.49	1\\
4.4925	1\\
4.495	1\\
4.4975	1\\
4.5	1\\
4.5025	1\\
4.505	1\\
4.5075	1\\
4.51	1\\
4.5125	1\\
4.515	8\\
4.5175	1\\
4.52	1\\
4.5225	1\\
4.525	1\\
4.5275	1\\
4.53	1\\
4.5325	1\\
4.535	16\\
4.5375	10\\
4.54	1\\
4.5425	1\\
4.545	1\\
4.5475	35\\
4.55	1\\
4.5525	1\\
4.555	1\\
4.5575	1\\
4.56	1\\
4.5625	1\\
4.565	1\\
4.5675	1\\
4.57	1\\
4.5725	1\\
4.575	1\\
4.5775	1\\
4.58	1\\
4.5825	1\\
4.585	2\\
4.5875	3\\
4.59	1\\
4.5925	1\\
4.595	1\\
4.5975	1\\
4.6	1\\
4.6025	1\\
4.605	0\\
4.6075	0\\
4.61	0\\
4.6125	0\\
4.615	0\\
4.6175	0\\
4.62	0\\
4.6225	1\\
4.625	1\\
4.6275	1\\
4.63	1\\
4.6325	1\\
4.635	1\\
4.6375	0\\
4.64	2\\
4.6425	1\\
4.645	0\\
4.6475	3\\
4.65	1\\
4.6525	0\\
4.655	3\\
4.6575	1\\
4.66	0\\
4.6625	2\\
4.665	1\\
4.6675	1\\
4.67	1\\
4.6725	3\\
4.675	1\\
4.6775	0\\
4.68	7\\
4.6825	1\\
4.685	1\\
4.6875	5\\
4.69	1\\
4.6925	1\\
4.695	1\\
4.6975	4\\
4.7	1\\
4.7025	1\\
4.705	8\\
4.7075	1\\
4.71	1\\
4.7125	1\\
4.715	1\\
4.7175	1\\
4.72	1\\
4.7225	1\\
4.725	1\\
4.7275	1\\
4.73	1\\
4.7325	1\\
4.735	1\\
4.7375	1\\
4.74	1\\
4.7425	1\\
4.745	1\\
4.7475	1\\
4.75	1\\
4.7525	1\\
4.755	2\\
4.7575	1\\
4.76	2\\
4.7625	1\\
4.765	1\\
4.7675	2\\
4.77	1\\
4.7725	2\\
4.775	1\\
4.7775	2\\
4.78	1\\
4.7825	2\\
4.785	1\\
4.7875	2\\
4.79	1\\
4.7925	1\\
4.795	1\\
4.7975	1\\
4.8	1\\
4.8025	1\\
4.805	1\\
4.8075	1\\
4.81	1\\
4.8125	1\\
4.815	1\\
4.8175	1\\
4.82	1\\
4.8225	1\\
4.825	1\\
4.8275	1\\
4.83	1\\
4.8325	1\\
4.835	1\\
4.8375	1\\
4.84	1\\
4.8425	1\\
4.845	1\\
4.8475	1\\
4.85	1\\
4.8525	1\\
4.855	1\\
4.8575	1\\
4.86	1\\
4.8625	1\\
4.865	1\\
4.8675	1\\
4.87	1\\
4.8725	1\\
4.875	1\\
4.8775	1\\
4.88	1\\
4.8825	1\\
4.885	1\\
4.8875	1\\
4.89	1\\
4.8925	1\\
4.895	1\\
4.8975	1\\
4.9	1\\
4.9025	1\\
4.905	1\\
4.9075	1\\
4.91	1\\
4.9125	1\\
4.915	1\\
4.9175	1\\
4.92	1\\
4.9225	1\\
4.925	1\\
4.9275	1\\
4.93	1\\
4.9325	1\\
4.935	1\\
4.9375	1\\
4.94	1\\
4.9425	1\\
4.945	1\\
4.9475	1\\
4.95	1\\
4.9525	1\\
4.955	1\\
4.9575	1\\
4.96	1\\
4.9625	1\\
4.965	1\\
4.9675	1\\
4.97	1\\
4.9725	1\\
4.975	1\\
4.9775	1\\
4.98	1\\
4.9825	1\\
4.985	1\\
4.9875	1\\
4.99	1\\
4.9925	1\\
4.995	1\\
4.9975	1\\
};
\addlegendentry{$\epsilon_u = 10^{-8}$}

\addplot [color=mycolor2, line width=1.2pt,opacity=0.5]
table[row sep=crcr]{%
0	1\\
0.0025	0\\
0.005	0\\
0.0075	0\\
0.01	0\\
0.0125	0\\
0.015	0\\
0.0175	0\\
0.02	0\\
0.0225	0\\
0.025	0\\
0.0275	0\\
0.03	0\\
0.0325	0\\
0.035	0\\
0.0375	0\\
0.04	0\\
0.0425	0\\
0.045	0\\
0.0475	0\\
0.05	0\\
0.0525	0\\
0.055	0\\
0.0575	0\\
0.06	0\\
0.0625	0\\
0.065	0\\
0.0675	0\\
0.07	0\\
0.0725	0\\
0.075	0\\
0.0775	0\\
0.08	0\\
0.0825	1\\
0.085	1\\
0.0875	1\\
0.09	0\\
0.0925	0\\
0.095	0\\
0.0975	0\\
0.1	1\\
0.1025	1\\
0.105	1\\
0.1075	1\\
0.11	1\\
0.1125	1\\
0.115	1\\
0.1175	1\\
0.12	1\\
0.1225	1\\
0.125	1\\
0.1275	1\\
0.13	1\\
0.1325	1\\
0.135	1\\
0.1375	1\\
0.14	1\\
0.1425	1\\
0.145	1\\
0.1475	1\\
0.15	1\\
0.1525	1\\
0.155	1\\
0.1575	1\\
0.16	1\\
0.1625	1\\
0.165	1\\
0.1675	1\\
0.17	1\\
0.1725	1\\
0.175	1\\
0.1775	1\\
0.18	1\\
0.1825	0\\
0.185	0\\
0.1875	1\\
0.19	1\\
0.1925	1\\
0.195	1\\
0.1975	1\\
0.2	1\\
0.2025	1\\
0.205	1\\
0.2075	1\\
0.21	1\\
0.2125	1\\
0.215	1\\
0.2175	1\\
0.22	1\\
0.2225	1\\
0.225	1\\
0.2275	1\\
0.23	1\\
0.2325	1\\
0.235	1\\
0.2375	1\\
0.24	1\\
0.2425	1\\
0.245	1\\
0.2475	1\\
0.25	1\\
0.2525	1\\
0.255	1\\
0.2575	1\\
0.26	1\\
0.2625	1\\
0.265	1\\
0.2675	1\\
0.27	1\\
0.2725	1\\
0.275	1\\
0.2775	1\\
0.28	1\\
0.2825	1\\
0.285	1\\
0.2875	1\\
0.29	1\\
0.2925	1\\
0.295	1\\
0.2975	1\\
0.3	1\\
0.3025	1\\
0.305	1\\
0.3075	1\\
0.31	1\\
0.3125	1\\
0.315	1\\
0.3175	1\\
0.32	1\\
0.3225	1\\
0.325	1\\
0.3275	1\\
0.33	1\\
0.3325	1\\
0.335	1\\
0.3375	1\\
0.34	1\\
0.3425	1\\
0.345	1\\
0.3475	1\\
0.35	1\\
0.3525	1\\
0.355	1\\
0.3575	1\\
0.36	1\\
0.3625	1\\
0.365	1\\
0.3675	1\\
0.37	1\\
0.3725	1\\
0.375	1\\
0.3775	1\\
0.38	1\\
0.3825	1\\
0.385	1\\
0.3875	1\\
0.39	1\\
0.3925	1\\
0.395	1\\
0.3975	1\\
0.4	1\\
0.4025	1\\
0.405	1\\
0.4075	1\\
0.41	1\\
0.4125	1\\
0.415	1\\
0.4175	1\\
0.42	1\\
0.4225	1\\
0.425	1\\
0.4275	1\\
0.43	1\\
0.4325	1\\
0.435	1\\
0.4375	1\\
0.44	1\\
0.4425	1\\
0.445	1\\
0.4475	1\\
0.45	1\\
0.4525	1\\
0.455	1\\
0.4575	1\\
0.46	1\\
0.4625	1\\
0.465	1\\
0.4675	1\\
0.47	1\\
0.4725	1\\
0.475	1\\
0.4775	1\\
0.48	1\\
0.4825	1\\
0.485	1\\
0.4875	1\\
0.49	1\\
0.4925	1\\
0.495	1\\
0.4975	1\\
0.5	1\\
0.5025	1\\
0.505	1\\
0.5075	1\\
0.51	1\\
0.5125	1\\
0.515	1\\
0.5175	1\\
0.52	1\\
0.5225	1\\
0.525	1\\
0.5275	1\\
0.53	1\\
0.5325	1\\
0.535	1\\
0.5375	1\\
0.54	1\\
0.5425	1\\
0.545	1\\
0.5475	1\\
0.55	1\\
0.5525	1\\
0.555	1\\
0.5575	1\\
0.56	1\\
0.5625	1\\
0.565	1\\
0.5675	1\\
0.57	1\\
0.5725	1\\
0.575	1\\
0.5775	1\\
0.58	1\\
0.5825	1\\
0.585	1\\
0.5875	1\\
0.59	1\\
0.5925	1\\
0.595	1\\
0.5975	3\\
0.6	1\\
0.6025	1\\
0.605	1\\
0.6075	1\\
0.61	1\\
0.6125	1\\
0.615	1\\
0.6175	1\\
0.62	1\\
0.6225	1\\
0.625	1\\
0.6275	1\\
0.63	1\\
0.6325	1\\
0.635	1\\
0.6375	1\\
0.64	1\\
0.6425	1\\
0.645	1\\
0.6475	1\\
0.65	1\\
0.6525	1\\
0.655	1\\
0.6575	1\\
0.66	1\\
0.6625	1\\
0.665	1\\
0.6675	1\\
0.67	1\\
0.6725	1\\
0.675	1\\
0.6775	1\\
0.68	1\\
0.6825	1\\
0.685	1\\
0.6875	1\\
0.69	1\\
0.6925	1\\
0.695	1\\
0.6975	1\\
0.7	1\\
0.7025	1\\
0.705	1\\
0.7075	1\\
0.71	1\\
0.7125	1\\
0.715	1\\
0.7175	1\\
0.72	1\\
0.7225	1\\
0.725	1\\
0.7275	1\\
0.73	1\\
0.7325	1\\
0.735	1\\
0.7375	1\\
0.74	1\\
0.7425	1\\
0.745	1\\
0.7475	1\\
0.75	1\\
0.7525	1\\
0.755	1\\
0.7575	1\\
0.76	1\\
0.7625	1\\
0.765	1\\
0.7675	1\\
0.77	1\\
0.7725	1\\
0.775	0\\
0.7775	0\\
0.78	0\\
0.7825	0\\
0.785	0\\
0.7875	0\\
0.79	0\\
0.7925	0\\
0.795	0\\
0.7975	0\\
0.8	0\\
0.8025	0\\
0.805	0\\
0.8075	0\\
0.81	0\\
0.8125	0\\
0.815	0\\
0.8175	0\\
0.82	0\\
0.8225	0\\
0.825	0\\
0.8275	0\\
0.83	0\\
0.8325	0\\
0.835	0\\
0.8375	0\\
0.84	0\\
0.8425	0\\
0.845	0\\
0.8475	0\\
0.85	0\\
0.8525	0\\
0.855	0\\
0.8575	0\\
0.86	0\\
0.8625	0\\
0.865	0\\
0.8675	0\\
0.87	0\\
0.8725	1\\
0.875	1\\
0.8775	1\\
0.88	1\\
0.8825	1\\
0.885	1\\
0.8875	1\\
0.89	1\\
0.8925	1\\
0.895	1\\
0.8975	1\\
0.9	1\\
0.9025	1\\
0.905	1\\
0.9075	1\\
0.91	1\\
0.9125	1\\
0.915	1\\
0.9175	1\\
0.92	1\\
0.9225	1\\
0.925	1\\
0.9275	1\\
0.93	1\\
0.9325	1\\
0.935	1\\
0.9375	1\\
0.94	1\\
0.9425	1\\
0.945	1\\
0.9475	0\\
0.95	0\\
0.9525	0\\
0.955	0\\
0.9575	0\\
0.96	1\\
0.9625	1\\
0.965	1\\
0.9675	1\\
0.97	1\\
0.9725	0\\
0.975	0\\
0.9775	0\\
0.98	0\\
0.9825	0\\
0.985	0\\
0.9875	0\\
0.99	0\\
0.9925	0\\
0.995	0\\
0.9975	0\\
1	0\\
1.0025	1\\
1.005	0\\
1.0075	0\\
1.01	0\\
1.0125	0\\
1.015	0\\
1.0175	0\\
1.02	0\\
1.0225	0\\
1.025	0\\
1.0275	0\\
1.03	0\\
1.0325	0\\
1.035	0\\
1.0375	0\\
1.04	0\\
1.0425	0\\
1.045	1\\
1.0475	1\\
1.05	1\\
1.0525	1\\
1.055	1\\
1.0575	1\\
1.06	1\\
1.0625	1\\
1.065	1\\
1.0675	1\\
1.07	1\\
1.0725	1\\
1.075	0\\
1.0775	0\\
1.08	0\\
1.0825	0\\
1.085	0\\
1.0875	0\\
1.09	0\\
1.0925	0\\
1.095	0\\
1.0975	1\\
1.1	0\\
1.1025	0\\
1.105	0\\
1.1075	0\\
1.11	0\\
1.1125	0\\
1.115	0\\
1.1175	0\\
1.12	0\\
1.1225	0\\
1.125	0\\
1.1275	0\\
1.13	0\\
1.1325	0\\
1.135	0\\
1.1375	0\\
1.14	0\\
1.1425	0\\
1.145	0\\
1.1475	0\\
1.15	0\\
1.1525	0\\
1.155	0\\
1.1575	0\\
1.16	0\\
1.1625	0\\
1.165	0\\
1.1675	0\\
1.17	0\\
1.1725	0\\
1.175	0\\
1.1775	0\\
1.18	0\\
1.1825	0\\
1.185	0\\
1.1875	0\\
1.19	0\\
1.1925	0\\
1.195	0\\
1.1975	0\\
1.2	0\\
1.2025	0\\
1.205	0\\
1.2075	0\\
1.21	0\\
1.2125	0\\
1.215	0\\
1.2175	0\\
1.22	0\\
1.2225	1\\
1.225	0\\
1.2275	1\\
1.23	1\\
1.2325	1\\
1.235	1\\
1.2375	0\\
1.24	0\\
1.2425	1\\
1.245	1\\
1.2475	0\\
1.25	0\\
1.2525	0\\
1.255	0\\
1.2575	0\\
1.26	0\\
1.2625	0\\
1.265	0\\
1.2675	0\\
1.27	0\\
1.2725	0\\
1.275	0\\
1.2775	0\\
1.28	0\\
1.2825	0\\
1.285	0\\
1.2875	0\\
1.29	1\\
1.2925	1\\
1.295	1\\
1.2975	0\\
1.3	0\\
1.3025	0\\
1.305	0\\
1.3075	0\\
1.31	0\\
1.3125	0\\
1.315	0\\
1.3175	0\\
1.32	0\\
1.3225	0\\
1.325	0\\
1.3275	0\\
1.33	0\\
1.3325	0\\
1.335	0\\
1.3375	0\\
1.34	0\\
1.3425	0\\
1.345	0\\
1.3475	0\\
1.35	0\\
1.3525	0\\
1.355	0\\
1.3575	0\\
1.36	0\\
1.3625	0\\
1.365	0\\
1.3675	0\\
1.37	0\\
1.3725	0\\
1.375	0\\
1.3775	0\\
1.38	0\\
1.3825	0\\
1.385	0\\
1.3875	0\\
1.39	0\\
1.3925	0\\
1.395	0\\
1.3975	0\\
1.4	0\\
1.4025	0\\
1.405	0\\
1.4075	0\\
1.41	0\\
1.4125	0\\
1.415	0\\
1.4175	0\\
1.42	0\\
1.4225	0\\
1.425	0\\
1.4275	0\\
1.43	0\\
1.4325	0\\
1.435	0\\
1.4375	0\\
1.44	0\\
1.4425	0\\
1.445	0\\
1.4475	0\\
1.45	0\\
1.4525	0\\
1.455	0\\
1.4575	0\\
1.46	0\\
1.4625	0\\
1.465	0\\
1.4675	0\\
1.47	0\\
1.4725	0\\
1.475	0\\
1.4775	0\\
1.48	0\\
1.4825	0\\
1.485	0\\
1.4875	0\\
1.49	0\\
1.4925	0\\
1.495	0\\
1.4975	0\\
1.5	0\\
1.5025	0\\
1.505	0\\
1.5075	0\\
1.51	0\\
1.5125	0\\
1.515	0\\
1.5175	0\\
1.52	0\\
1.5225	0\\
1.525	0\\
1.5275	0\\
1.53	0\\
1.5325	0\\
1.535	0\\
1.5375	0\\
1.54	0\\
1.5425	0\\
1.545	0\\
1.5475	0\\
1.55	0\\
1.5525	0\\
1.555	0\\
1.5575	0\\
1.56	1\\
1.5625	1\\
1.565	0\\
1.5675	0\\
1.57	0\\
1.5725	0\\
1.575	0\\
1.5775	0\\
1.58	0\\
1.5825	0\\
1.585	0\\
1.5875	0\\
1.59	0\\
1.5925	0\\
1.595	0\\
1.5975	0\\
1.6	0\\
1.6025	0\\
1.605	0\\
1.6075	0\\
1.61	0\\
1.6125	0\\
1.615	0\\
1.6175	0\\
1.62	0\\
1.6225	0\\
1.625	0\\
1.6275	0\\
1.63	0\\
1.6325	0\\
1.635	0\\
1.6375	0\\
1.64	0\\
1.6425	0\\
1.645	0\\
1.6475	0\\
1.65	0\\
1.6525	0\\
1.655	0\\
1.6575	0\\
1.66	0\\
1.6625	0\\
1.665	0\\
1.6675	0\\
1.67	0\\
1.6725	0\\
1.675	0\\
1.6775	0\\
1.68	0\\
1.6825	0\\
1.685	0\\
1.6875	0\\
1.69	0\\
1.6925	0\\
1.695	0\\
1.6975	0\\
1.7	0\\
1.7025	1\\
1.705	0\\
1.7075	0\\
1.71	0\\
1.7125	0\\
1.715	0\\
1.7175	0\\
1.72	1\\
1.7225	0\\
1.725	0\\
1.7275	0\\
1.73	0\\
1.7325	0\\
1.735	9\\
1.7375	0\\
1.74	1\\
1.7425	0\\
1.745	0\\
1.7475	0\\
1.75	0\\
1.7525	0\\
1.755	0\\
1.7575	0\\
1.76	0\\
1.7625	0\\
1.765	0\\
1.7675	0\\
1.77	0\\
1.7725	0\\
1.775	0\\
1.7775	0\\
1.78	0\\
1.7825	0\\
1.785	0\\
1.7875	0\\
1.79	0\\
1.7925	0\\
1.795	0\\
1.7975	1\\
1.8	0\\
1.8025	0\\
1.805	0\\
1.8075	0\\
1.81	0\\
1.8125	0\\
1.815	0\\
1.8175	0\\
1.82	0\\
1.8225	1\\
1.825	1\\
1.8275	1\\
1.83	1\\
1.8325	1\\
1.835	0\\
1.8375	0\\
1.84	0\\
1.8425	0\\
1.845	0\\
1.8475	0\\
1.85	0\\
1.8525	0\\
1.855	0\\
1.8575	0\\
1.86	0\\
1.8625	0\\
1.865	0\\
1.8675	0\\
1.87	0\\
1.8725	0\\
1.875	0\\
1.8775	0\\
1.88	0\\
1.8825	0\\
1.885	0\\
1.8875	0\\
1.89	0\\
1.8925	0\\
1.895	0\\
1.8975	0\\
1.9	0\\
1.9025	0\\
1.905	0\\
1.9075	0\\
1.91	0\\
1.9125	0\\
1.915	0\\
1.9175	0\\
1.92	0\\
1.9225	0\\
1.925	0\\
1.9275	0\\
1.93	0\\
1.9325	0\\
1.935	0\\
1.9375	0\\
1.94	0\\
1.9425	0\\
1.945	0\\
1.9475	0\\
1.95	0\\
1.9525	0\\
1.955	0\\
1.9575	0\\
1.96	0\\
1.9625	0\\
1.965	0\\
1.9675	0\\
1.97	0\\
1.9725	0\\
1.975	0\\
1.9775	0\\
1.98	0\\
1.9825	0\\
1.985	0\\
1.9875	0\\
1.99	0\\
1.9925	0\\
1.995	0\\
1.9975	0\\
2	0\\
2.0025	0\\
2.005	0\\
2.0075	0\\
2.01	0\\
2.0125	0\\
2.015	0\\
2.0175	0\\
2.02	0\\
2.0225	0\\
2.025	0\\
2.0275	0\\
2.03	0\\
2.0325	0\\
2.035	1\\
2.0375	1\\
2.04	1\\
2.0425	0\\
2.045	0\\
2.0475	0\\
2.05	0\\
2.0525	0\\
2.055	0\\
2.0575	0\\
2.06	0\\
2.0625	0\\
2.065	0\\
2.0675	0\\
2.07	0\\
2.0725	0\\
2.075	0\\
2.0775	0\\
2.08	0\\
2.0825	0\\
2.085	0\\
2.0875	0\\
2.09	0\\
2.0925	0\\
2.095	0\\
2.0975	0\\
2.1	0\\
2.1025	0\\
2.105	0\\
2.1075	0\\
2.11	0\\
2.1125	0\\
2.115	0\\
2.1175	0\\
2.12	0\\
2.1225	0\\
2.125	0\\
2.1275	0\\
2.13	0\\
2.1325	0\\
2.135	0\\
2.1375	0\\
2.14	0\\
2.1425	0\\
2.145	0\\
2.1475	0\\
2.15	0\\
2.1525	0\\
2.155	0\\
2.1575	0\\
2.16	0\\
2.1625	0\\
2.165	0\\
2.1675	0\\
2.17	0\\
2.1725	1\\
2.175	1\\
2.1775	1\\
2.18	1\\
2.1825	1\\
2.185	1\\
2.1875	1\\
2.19	1\\
2.1925	1\\
2.195	1\\
2.1975	0\\
2.2	0\\
2.2025	0\\
2.205	0\\
2.2075	0\\
2.21	0\\
2.2125	0\\
2.215	0\\
2.2175	0\\
2.22	0\\
2.2225	0\\
2.225	0\\
2.2275	0\\
2.23	0\\
2.2325	0\\
2.235	0\\
2.2375	0\\
2.24	0\\
2.2425	0\\
2.245	0\\
2.2475	0\\
2.25	0\\
2.2525	0\\
2.255	0\\
2.2575	0\\
2.26	0\\
2.2625	0\\
2.265	0\\
2.2675	0\\
2.27	0\\
2.2725	0\\
2.275	0\\
2.2775	0\\
2.28	0\\
2.2825	0\\
2.285	0\\
2.2875	0\\
2.29	0\\
2.2925	0\\
2.295	0\\
2.2975	0\\
2.3	0\\
2.3025	0\\
2.305	0\\
2.3075	0\\
2.31	0\\
2.3125	0\\
2.315	0\\
2.3175	0\\
2.32	0\\
2.3225	0\\
2.325	0\\
2.3275	0\\
2.33	0\\
2.3325	0\\
2.335	1\\
2.3375	1\\
2.34	1\\
2.3425	0\\
2.345	0\\
2.3475	0\\
2.35	0\\
2.3525	1\\
2.355	0\\
2.3575	0\\
2.36	0\\
2.3625	0\\
2.365	1\\
2.3675	0\\
2.37	1\\
2.3725	1\\
2.375	1\\
2.3775	1\\
2.38	1\\
2.3825	1\\
2.385	0\\
2.3875	0\\
2.39	0\\
2.3925	0\\
2.395	0\\
2.3975	0\\
2.4	0\\
2.4025	1\\
2.405	0\\
2.4075	1\\
2.41	1\\
2.4125	1\\
2.415	1\\
2.4175	1\\
2.42	1\\
2.4225	1\\
2.425	1\\
2.4275	1\\
2.43	1\\
2.4325	1\\
2.435	0\\
2.4375	1\\
2.44	0\\
2.4425	0\\
2.445	0\\
2.4475	0\\
2.45	0\\
2.4525	0\\
2.455	0\\
2.4575	0\\
2.46	0\\
2.4625	0\\
2.465	0\\
2.4675	0\\
2.47	0\\
2.4725	0\\
2.475	0\\
2.4775	0\\
2.48	0\\
2.4825	0\\
2.485	0\\
2.4875	0\\
2.49	0\\
2.4925	0\\
2.495	0\\
2.4975	0\\
2.5	0\\
2.5025	0\\
2.505	0\\
2.5075	0\\
2.51	0\\
2.5125	0\\
2.515	0\\
2.5175	0\\
2.52	0\\
2.5225	0\\
2.525	1\\
2.5275	0\\
2.53	0\\
2.5325	0\\
2.535	0\\
2.5375	0\\
2.54	0\\
2.5425	0\\
2.545	0\\
2.5475	0\\
2.55	0\\
2.5525	0\\
2.555	0\\
2.5575	0\\
2.56	0\\
2.5625	1\\
2.565	1\\
2.5675	1\\
2.57	1\\
2.5725	1\\
2.575	1\\
2.5775	0\\
2.58	0\\
2.5825	0\\
2.585	0\\
2.5875	0\\
2.59	0\\
2.5925	0\\
2.595	0\\
2.5975	1\\
2.6	0\\
2.6025	1\\
2.605	1\\
2.6075	1\\
2.61	1\\
2.6125	1\\
2.615	1\\
2.6175	1\\
2.62	0\\
2.6225	0\\
2.625	0\\
2.6275	0\\
2.63	0\\
2.6325	0\\
2.635	0\\
2.6375	0\\
2.64	0\\
2.6425	0\\
2.645	0\\
2.6475	0\\
2.65	0\\
2.6525	0\\
2.655	0\\
2.6575	0\\
2.66	0\\
2.6625	0\\
2.665	0\\
2.6675	0\\
2.67	0\\
2.6725	0\\
2.675	0\\
2.6775	0\\
2.68	0\\
2.6825	0\\
2.685	0\\
2.6875	0\\
2.69	0\\
2.6925	0\\
2.695	0\\
2.6975	0\\
2.7	0\\
2.7025	1\\
2.705	1\\
2.7075	1\\
2.71	1\\
2.7125	1\\
2.715	1\\
2.7175	1\\
2.72	1\\
2.7225	1\\
2.725	1\\
2.7275	1\\
2.73	0\\
2.7325	1\\
2.735	1\\
2.7375	1\\
2.74	1\\
2.7425	1\\
2.745	1\\
2.7475	1\\
2.75	1\\
2.7525	1\\
2.755	1\\
2.7575	1\\
2.76	1\\
2.7625	1\\
2.765	1\\
2.7675	1\\
2.77	1\\
2.7725	1\\
2.775	1\\
2.7775	1\\
2.78	1\\
2.7825	1\\
2.785	1\\
2.7875	1\\
2.79	1\\
2.7925	1\\
2.795	1\\
2.7975	0\\
2.8	2\\
2.8025	0\\
2.805	0\\
2.8075	0\\
2.81	1\\
2.8125	0\\
2.815	0\\
2.8175	1\\
2.82	0\\
2.8225	0\\
2.825	1\\
2.8275	1\\
2.83	0\\
2.8325	0\\
2.835	1\\
2.8375	0\\
2.84	0\\
2.8425	0\\
2.845	0\\
2.8475	0\\
2.85	0\\
2.8525	0\\
2.855	0\\
2.8575	0\\
2.86	0\\
2.8625	0\\
2.865	0\\
2.8675	0\\
2.87	0\\
2.8725	0\\
2.875	0\\
2.8775	0\\
2.88	0\\
2.8825	0\\
2.885	0\\
2.8875	0\\
2.89	0\\
2.8925	0\\
2.895	0\\
2.8975	0\\
2.9	0\\
2.9025	0\\
2.905	0\\
2.9075	0\\
2.91	0\\
2.9125	0\\
2.915	0\\
2.9175	0\\
2.92	0\\
2.9225	0\\
2.925	0\\
2.9275	0\\
2.93	0\\
2.9325	0\\
2.935	0\\
2.9375	0\\
2.94	0\\
2.9425	0\\
2.945	0\\
2.9475	0\\
2.95	0\\
2.9525	0\\
2.955	0\\
2.9575	0\\
2.96	0\\
2.9625	0\\
2.965	0\\
2.9675	0\\
2.97	0\\
2.9725	0\\
2.975	0\\
2.9775	0\\
2.98	0\\
2.9825	0\\
2.985	0\\
2.9875	0\\
2.99	0\\
2.9925	0\\
2.995	0\\
2.9975	0\\
3	0\\
3.0025	0\\
3.005	0\\
3.0075	0\\
3.01	0\\
3.0125	0\\
3.015	0\\
3.0175	0\\
3.02	0\\
3.0225	0\\
3.025	0\\
3.0275	0\\
3.03	0\\
3.0325	1\\
3.035	1\\
3.0375	0\\
3.04	0\\
3.0425	0\\
3.045	0\\
3.0475	0\\
3.05	0\\
3.0525	0\\
3.055	0\\
3.0575	0\\
3.06	0\\
3.0625	0\\
3.065	0\\
3.0675	0\\
3.07	1\\
3.0725	1\\
3.075	0\\
3.0775	1\\
3.08	1\\
3.0825	1\\
3.085	1\\
3.0875	1\\
3.09	0\\
3.0925	0\\
3.095	0\\
3.0975	0\\
3.1	0\\
3.1025	0\\
3.105	0\\
3.1075	0\\
3.11	0\\
3.1125	0\\
3.115	0\\
3.1175	0\\
3.12	0\\
3.1225	0\\
3.125	0\\
3.1275	0\\
3.13	0\\
3.1325	0\\
3.135	0\\
3.1375	0\\
3.14	0\\
3.1425	0\\
3.145	0\\
3.1475	0\\
3.15	0\\
3.1525	0\\
3.155	0\\
3.1575	0\\
3.16	0\\
3.1625	0\\
3.165	0\\
3.1675	0\\
3.17	0\\
3.1725	0\\
3.175	0\\
3.1775	0\\
3.18	0\\
3.1825	0\\
3.185	0\\
3.1875	0\\
3.19	0\\
3.1925	0\\
3.195	0\\
3.1975	0\\
3.2	0\\
3.2025	0\\
3.205	0\\
3.2075	0\\
3.21	0\\
3.2125	0\\
3.215	0\\
3.2175	0\\
3.22	0\\
3.2225	0\\
3.225	0\\
3.2275	0\\
3.23	0\\
3.2325	0\\
3.235	0\\
3.2375	0\\
3.24	0\\
3.2425	0\\
3.245	0\\
3.2475	0\\
3.25	0\\
3.2525	0\\
3.255	0\\
3.2575	0\\
3.26	0\\
3.2625	0\\
3.265	0\\
3.2675	0\\
3.27	0\\
3.2725	0\\
3.275	0\\
3.2775	0\\
3.28	0\\
3.2825	0\\
3.285	0\\
3.2875	0\\
3.29	0\\
3.2925	0\\
3.295	0\\
3.2975	0\\
3.3	0\\
3.3025	0\\
3.305	0\\
3.3075	0\\
3.31	0\\
3.3125	0\\
3.315	0\\
3.3175	0\\
3.32	0\\
3.3225	0\\
3.325	0\\
3.3275	0\\
3.33	0\\
3.3325	0\\
3.335	0\\
3.3375	0\\
3.34	0\\
3.3425	0\\
3.345	1\\
3.3475	1\\
3.35	1\\
3.3525	1\\
3.355	1\\
3.3575	1\\
3.36	0\\
3.3625	0\\
3.365	0\\
3.3675	0\\
3.37	0\\
3.3725	0\\
3.375	0\\
3.3775	0\\
3.38	0\\
3.3825	0\\
3.385	0\\
3.3875	0\\
3.39	0\\
3.3925	0\\
3.395	0\\
3.3975	0\\
3.4	0\\
3.4025	0\\
3.405	0\\
3.4075	0\\
3.41	0\\
3.4125	0\\
3.415	0\\
3.4175	0\\
3.42	0\\
3.4225	0\\
3.425	0\\
3.4275	0\\
3.43	0\\
3.4325	0\\
3.435	0\\
3.4375	0\\
3.44	0\\
3.4425	0\\
3.445	0\\
3.4475	0\\
3.45	0\\
3.4525	0\\
3.455	0\\
3.4575	0\\
3.46	0\\
3.4625	0\\
3.465	0\\
3.4675	0\\
3.47	0\\
3.4725	0\\
3.475	0\\
3.4775	0\\
3.48	0\\
3.4825	0\\
3.485	0\\
3.4875	0\\
3.49	0\\
3.4925	0\\
3.495	0\\
3.4975	0\\
3.5	0\\
3.5025	0\\
3.505	0\\
3.5075	0\\
3.51	0\\
3.5125	0\\
3.515	0\\
3.5175	0\\
3.52	0\\
3.5225	0\\
3.525	0\\
3.5275	0\\
3.53	0\\
3.5325	0\\
3.535	0\\
3.5375	0\\
3.54	0\\
3.5425	0\\
3.545	0\\
3.5475	0\\
3.55	0\\
3.5525	0\\
3.555	0\\
3.5575	0\\
3.56	0\\
3.5625	0\\
3.565	0\\
3.5675	0\\
3.57	0\\
3.5725	0\\
3.575	0\\
3.5775	0\\
3.58	0\\
3.5825	0\\
3.585	0\\
3.5875	0\\
3.59	0\\
3.5925	0\\
3.595	0\\
3.5975	0\\
3.6	0\\
3.6025	0\\
3.605	0\\
3.6075	0\\
3.61	0\\
3.6125	0\\
3.615	0\\
3.6175	0\\
3.62	0\\
3.6225	0\\
3.625	0\\
3.6275	0\\
3.63	0\\
3.6325	0\\
3.635	0\\
3.6375	0\\
3.64	0\\
3.6425	0\\
3.645	0\\
3.6475	0\\
3.65	0\\
3.6525	0\\
3.655	0\\
3.6575	0\\
3.66	0\\
3.6625	0\\
3.665	0\\
3.6675	0\\
3.67	0\\
3.6725	0\\
3.675	0\\
3.6775	0\\
3.68	0\\
3.6825	0\\
3.685	0\\
3.6875	0\\
3.69	0\\
3.6925	0\\
3.695	0\\
3.6975	0\\
3.7	0\\
3.7025	0\\
3.705	0\\
3.7075	0\\
3.71	0\\
3.7125	0\\
3.715	0\\
3.7175	0\\
3.72	0\\
3.7225	0\\
3.725	0\\
3.7275	0\\
3.73	0\\
3.7325	0\\
3.735	0\\
3.7375	1\\
3.74	1\\
3.7425	0\\
3.745	0\\
3.7475	0\\
3.75	0\\
3.7525	0\\
3.755	0\\
3.7575	0\\
3.76	0\\
3.7625	0\\
3.765	0\\
3.7675	0\\
3.77	0\\
3.7725	0\\
3.775	0\\
3.7775	0\\
3.78	0\\
3.7825	0\\
3.785	0\\
3.7875	0\\
3.79	0\\
3.7925	0\\
3.795	0\\
3.7975	0\\
3.8	0\\
3.8025	0\\
3.805	0\\
3.8075	0\\
3.81	0\\
3.8125	0\\
3.815	0\\
3.8175	0\\
3.82	0\\
3.8225	0\\
3.825	0\\
3.8275	0\\
3.83	0\\
3.8325	0\\
3.835	0\\
3.8375	0\\
3.84	0\\
3.8425	0\\
3.845	0\\
3.8475	0\\
3.85	0\\
3.8525	0\\
3.855	0\\
3.8575	0\\
3.86	0\\
3.8625	0\\
3.865	0\\
3.8675	0\\
3.87	0\\
3.8725	0\\
3.875	0\\
3.8775	0\\
3.88	0\\
3.8825	0\\
3.885	0\\
3.8875	0\\
3.89	0\\
3.8925	0\\
3.895	0\\
3.8975	0\\
3.9	0\\
3.9025	0\\
3.905	0\\
3.9075	0\\
3.91	0\\
3.9125	0\\
3.915	0\\
3.9175	0\\
3.92	0\\
3.9225	0\\
3.925	0\\
3.9275	0\\
3.93	0\\
3.9325	0\\
3.935	0\\
3.9375	0\\
3.94	0\\
3.9425	0\\
3.945	0\\
3.9475	0\\
3.95	0\\
3.9525	0\\
3.955	0\\
3.9575	0\\
3.96	0\\
3.9625	0\\
3.965	0\\
3.9675	0\\
3.97	0\\
3.9725	0\\
3.975	0\\
3.9775	0\\
3.98	0\\
3.9825	0\\
3.985	0\\
3.9875	0\\
3.99	0\\
3.9925	0\\
3.995	0\\
3.9975	0\\
4	0\\
4.0025	0\\
4.005	0\\
4.0075	0\\
4.01	0\\
4.0125	0\\
4.015	0\\
4.0175	0\\
4.02	0\\
4.0225	0\\
4.025	0\\
4.0275	0\\
4.03	0\\
4.0325	0\\
4.035	0\\
4.0375	0\\
4.04	0\\
4.0425	0\\
4.045	0\\
4.0475	0\\
4.05	0\\
4.0525	0\\
4.055	0\\
4.0575	0\\
4.06	0\\
4.0625	0\\
4.065	0\\
4.0675	0\\
4.07	0\\
4.0725	0\\
4.075	0\\
4.0775	0\\
4.08	0\\
4.0825	0\\
4.085	0\\
4.0875	0\\
4.09	0\\
4.0925	1\\
4.095	0\\
4.0975	1\\
4.1	0\\
4.1025	0\\
4.105	0\\
4.1075	0\\
4.11	0\\
4.1125	0\\
4.115	0\\
4.1175	1\\
4.12	0\\
4.1225	0\\
4.125	0\\
4.1275	0\\
4.13	0\\
4.1325	0\\
4.135	0\\
4.1375	0\\
4.14	0\\
4.1425	0\\
4.145	0\\
4.1475	0\\
4.15	0\\
4.1525	0\\
4.155	0\\
4.1575	0\\
4.16	0\\
4.1625	0\\
4.165	0\\
4.1675	0\\
4.17	0\\
4.1725	0\\
4.175	0\\
4.1775	0\\
4.18	0\\
4.1825	0\\
4.185	0\\
4.1875	0\\
4.19	0\\
4.1925	0\\
4.195	0\\
4.1975	0\\
4.2	0\\
4.2025	0\\
4.205	0\\
4.2075	0\\
4.21	0\\
4.2125	0\\
4.215	1\\
4.2175	0\\
4.22	0\\
4.2225	0\\
4.225	0\\
4.2275	0\\
4.23	0\\
4.2325	0\\
4.235	0\\
4.2375	0\\
4.24	0\\
4.2425	0\\
4.245	0\\
4.2475	0\\
4.25	0\\
4.2525	0\\
4.255	0\\
4.2575	0\\
4.26	0\\
4.2625	0\\
4.265	0\\
4.2675	0\\
4.27	0\\
4.2725	0\\
4.275	0\\
4.2775	0\\
4.28	0\\
4.2825	0\\
4.285	0\\
4.2875	0\\
4.29	0\\
4.2925	0\\
4.295	0\\
4.2975	0\\
4.3	0\\
4.3025	0\\
4.305	0\\
4.3075	0\\
4.31	0\\
4.3125	0\\
4.315	0\\
4.3175	0\\
4.32	0\\
4.3225	0\\
4.325	0\\
4.3275	0\\
4.33	0\\
4.3325	0\\
4.335	0\\
4.3375	0\\
4.34	0\\
4.3425	0\\
4.345	0\\
4.3475	0\\
4.35	0\\
4.3525	0\\
4.355	0\\
4.3575	0\\
4.36	0\\
4.3625	0\\
4.365	0\\
4.3675	0\\
4.37	0\\
4.3725	0\\
4.375	0\\
4.3775	0\\
4.38	0\\
4.3825	0\\
4.385	0\\
4.3875	0\\
4.39	0\\
4.3925	0\\
4.395	0\\
4.3975	0\\
4.4	0\\
4.4025	0\\
4.405	0\\
4.4075	0\\
4.41	0\\
4.4125	0\\
4.415	0\\
4.4175	0\\
4.42	0\\
4.4225	0\\
4.425	0\\
4.4275	0\\
4.43	0\\
4.4325	0\\
4.435	0\\
4.4375	0\\
4.44	0\\
4.4425	0\\
4.445	0\\
4.4475	0\\
4.45	0\\
4.4525	0\\
4.455	0\\
4.4575	0\\
4.46	0\\
4.4625	0\\
4.465	0\\
4.4675	0\\
4.47	0\\
4.4725	0\\
4.475	0\\
4.4775	0\\
4.48	0\\
4.4825	0\\
4.485	0\\
4.4875	0\\
4.49	0\\
4.4925	0\\
4.495	0\\
4.4975	0\\
4.5	0\\
4.5025	0\\
4.505	0\\
4.5075	0\\
4.51	0\\
4.5125	0\\
4.515	0\\
4.5175	0\\
4.52	0\\
4.5225	0\\
4.525	0\\
4.5275	0\\
4.53	0\\
4.5325	0\\
4.535	0\\
4.5375	0\\
4.54	0\\
4.5425	0\\
4.545	1\\
4.5475	0\\
4.55	0\\
4.5525	0\\
4.555	0\\
4.5575	0\\
4.56	0\\
4.5625	0\\
4.565	0\\
4.5675	0\\
4.57	0\\
4.5725	0\\
4.575	0\\
4.5775	0\\
4.58	0\\
4.5825	0\\
4.585	0\\
4.5875	0\\
4.59	0\\
4.5925	0\\
4.595	0\\
4.5975	0\\
4.6	0\\
4.6025	0\\
4.605	0\\
4.6075	0\\
4.61	0\\
4.6125	0\\
4.615	0\\
4.6175	0\\
4.62	0\\
4.6225	0\\
4.625	0\\
4.6275	0\\
4.63	0\\
4.6325	0\\
4.635	0\\
4.6375	0\\
4.64	0\\
4.6425	0\\
4.645	0\\
4.6475	0\\
4.65	0\\
4.6525	0\\
4.655	0\\
4.6575	0\\
4.66	0\\
4.6625	0\\
4.665	0\\
4.6675	0\\
4.67	1\\
4.6725	0\\
4.675	1\\
4.6775	1\\
4.68	1\\
4.6825	0\\
4.685	1\\
4.6875	0\\
4.69	1\\
4.6925	0\\
4.695	0\\
4.6975	0\\
4.7	0\\
4.7025	0\\
4.705	0\\
4.7075	0\\
4.71	0\\
4.7125	0\\
4.715	0\\
4.7175	0\\
4.72	0\\
4.7225	0\\
4.725	0\\
4.7275	0\\
4.73	0\\
4.7325	0\\
4.735	0\\
4.7375	0\\
4.74	0\\
4.7425	0\\
4.745	0\\
4.7475	0\\
4.75	0\\
4.7525	0\\
4.755	0\\
4.7575	0\\
4.76	1\\
4.7625	0\\
4.765	0\\
4.7675	0\\
4.77	0\\
4.7725	0\\
4.775	0\\
4.7775	0\\
4.78	0\\
4.7825	0\\
4.785	0\\
4.7875	0\\
4.79	0\\
4.7925	0\\
4.795	0\\
4.7975	0\\
4.8	0\\
4.8025	0\\
4.805	0\\
4.8075	0\\
4.81	0\\
4.8125	1\\
4.815	0\\
4.8175	1\\
4.82	0\\
4.8225	0\\
4.825	0\\
4.8275	0\\
4.83	1\\
4.8325	1\\
4.835	1\\
4.8375	1\\
4.84	1\\
4.8425	0\\
4.845	0\\
4.8475	0\\
4.85	0\\
4.8525	0\\
4.855	0\\
4.8575	0\\
4.86	0\\
4.8625	0\\
4.865	0\\
4.8675	0\\
4.87	0\\
4.8725	0\\
4.875	0\\
4.8775	0\\
4.88	1\\
4.8825	1\\
4.885	1\\
4.8875	1\\
4.89	1\\
4.8925	1\\
4.895	1\\
4.8975	1\\
4.9	0\\
4.9025	0\\
4.905	0\\
4.9075	0\\
4.91	0\\
4.9125	0\\
4.915	0\\
4.9175	0\\
4.92	0\\
4.9225	0\\
4.925	0\\
4.9275	0\\
4.93	0\\
4.9325	0\\
4.935	0\\
4.9375	0\\
4.94	0\\
4.9425	0\\
4.945	0\\
4.9475	0\\
4.95	0\\
4.9525	0\\
4.955	0\\
4.9575	0\\
4.96	0\\
4.9625	0\\
4.965	0\\
4.9675	0\\
4.97	0\\
4.9725	0\\
4.975	0\\
4.9775	0\\
4.98	0\\
4.9825	0\\
4.985	0\\
4.9875	0\\
4.99	0\\
4.9925	0\\
4.995	0\\
4.9975	0\\
};
\addlegendentry{$\epsilon_u = 10^{-11}$}

\end{axis}
\end{tikzpicture}%

%% file: figures/mr_AC.tex
% This file was created with tikzplotlib v0.10.1.
\definecolor{mycolor1}{rgb}{0.00000,0.44700,0.74100}%
\definecolor{mycolor2}{rgb}{0.85000,0.32500,0.09800}%
\definecolor{mycolor3}{rgb}{0.5,0.5,0.5}
\begin{tikzpicture}[scale=0.8]

\begin{axis}[%
width=1.5in,
height=1.5in,
at={(1.011in,0.642in)},
scale only axis,
xmin=0,
xmax=200,
xlabel style={font=\color{white!15!black}},
xlabel={Time},
%ymode=log,
ymin=0,
ymax=60,
yminorticks=true,
ylabel style={font=\color{white!15!black}},
axis background/.style={fill=white},
legend style={fill opacity=0.8, 
font=\tiny,
  draw opacity=1,
  text opacity=1,
  at={(0.02,0.87)}, %(0.3,0.13)
  anchor=west}
  % draw=lightgray204}
]

\addplot [line width=1.2pt,  mycolor1]
table {%
0 5
0.5 0
1 0
1.5 1
2 1
2.5 1
3 1
3.5 1
4 1
4.5 1
5 1
5.5 1
6 1
6.5 1
7 1
7.5 11
8 2
8.5 1
9 1
9.5 1
10 2
10.5 1
11 0
11.5 1
12 1
12.5 1
13 1
13.5 1
14 2
14.5 2
15 5
15.5 1
16 2
16.5 3
17 2
17.5 1
18 1
18.5 6
19 2
19.5 1
20 1
20.5 7
21 1
21.5 2
22 2
22.5 1
23 2
23.5 1
24 7
24.5 48
25 4
25.5 2
26 2
26.5 2
27 2
27.5 2
28 1
28.5 1
29 1
29.5 1
30 2
30.5 9
31 2
31.5 7
32 2
32.5 18
33 2
33.5 1
34 6
34.5 1
35 2
35.5 5
36 1
36.5 1
37 2
37.5 12
38 19
38.5 6
39 1
39.5 2
40 2
40.5 2
41 6
41.5 15
42 10
42.5 2
43 2
43.5 0
44 1
44.5 2
45 2
45.5 2
46 7
46.5 1
47 2
47.5 1
48 4
48.5 6
49 9
49.5 1
50 4
50.5 16
51 2
51.5 5
52 18
52.5 9
53 2
53.5 1
54 14
54.5 2
55 1
55.5 3
56 1
56.5 3
57 2
57.5 1
58 1
58.5 4
59 1
59.5 1
60 1
60.5 1
61 4
61.5 17
62 2
62.5 2
63 2
63.5 12
64 2
64.5 2
65 2
65.5 3
66 5
66.5 2
67 4
67.5 1
68 2
68.5 3
69 1
69.5 1
70 1
70.5 1
71 1
71.5 1
72 1
72.5 1
73 1
73.5 2
74 1
74.5 1
75 1
75.5 1
76 5
76.5 1
77 2
77.5 2
78 1
78.5 3
79 3
79.5 1
80 2
80.5 1
81 1
81.5 2
82 1
82.5 2
83 1
83.5 2
84 1
84.5 1
85 1
85.5 1
86 2
86.5 1
87 16
87.5 10
88 1
88.5 1
89 2
89.5 1
90 1
90.5 9
91 1
91.5 1
92 6
92.5 5
93 2
93.5 1
94 2
94.5 1
95 1
95.5 1
96 1
96.5 2
97 1
97.5 1
98 2
98.5 1
99 2
99.5 2
100 1
100.5 1
101 1
101.5 2
102 2
102.5 23
103 1
103.5 1
104 2
104.5 2
105 4
105.5 11
106 4
106.5 2
107 1
107.5 5
108 1
108.5 1
109 2
109.5 1
110 9
110.5 1
111 1
111.5 2
112 2
112.5 1
113 2
113.5 1
114 1
114.5 1
115 1
115.5 6
116 5
116.5 3
117 1
117.5 1
118 1
118.5 1
119 3
119.5 4
120 1
120.5 0
121 11
121.5 1
122 1
122.5 0
123 2
123.5 1
124 7
124.5 2
125 2
125.5 1
126 1
126.5 2
127 7
127.5 1
128 2
128.5 2
129 5
129.5 6
130 2
130.5 2
131 1
131.5 16
132 8
132.5 1
133 5
133.5 1
134 1
134.5 1
135 8
135.5 1
136 1
136.5 2
137 6
137.5 1
138 1
138.5 1
139 1
139.5 0
140 1
140.5 2
141 2
141.5 3
142 2
142.5 2
143 7
143.5 9
144 12
144.5 2
145 1
145.5 2
146 4
146.5 3
147 1
147.5 1
148 1
148.5 4
149 6
149.5 1
150 1
150.5 1
151 1
151.5 5
152 3
152.5 1
153 2
153.5 3
154 18
154.5 1
155 2
155.5 5
156 7
156.5 7
157 1
157.5 5
158 1
158.5 5
159 1
159.5 1
160 3
160.5 6
161 1
161.5 2
162 2
162.5 4
163 1
163.5 2
164 7
164.5 1
165 1
165.5 1
166 8
166.5 11
167 9
167.5 4
168 2
168.5 1
169 2
169.5 1
170 1
170.5 5
171 2
171.5 2
172 2
172.5 2
173 1
173.5 1
174 4
174.5 3
175 2
175.5 7
176 2
176.5 2
177 2
177.5 2
178 2
178.5 2
179 2
179.5 2
180 2
180.5 2
181 2
181.5 2
182 2
182.5 2
183 2
183.5 2
184 2
184.5 2
185 4
185.5 3
186 2
186.5 3
187 1
187.5 2
188 2
188.5 2
189 5
189.5 10
190 2
190.5 3
191 11
191.5 2
192 2
192.5 2
193 2
193.5 2
194 2
194.5 2
195 2
195.5 3
196 1
196.5 4
197 2
197.5 1
198 1
198.5 1
199 2
199.5 4
};
%\addlegendentry{$\epsilon_u = 10^{-8}$}

\addplot [line width=1.2pt,  mycolor2,opacity=0.5]
table {%
0 10
0.5 0
1 0
1.5 0
2 0
2.5 0
3 0
3.5 1
4 1
4.5 1
5 1
5.5 1
6 1
6.5 1
7 1
7.5 1
8 1
8.5 1
9 1
9.5 2
10 1
10.5 1
11 7
11.5 1
12 2
12.5 5
13 1
13.5 4
14 2
14.5 2
15 2
15.5 1
16 2
16.5 1
17 2
17.5 1
18 1
18.5 1
19 1
19.5 2
20 1
20.5 6
21 1
21.5 8
22 1
22.5 1
23 1
23.5 2
24 1
24.5 5
25 2
25.5 8
26 1
26.5 7
27 2
27.5 1
28 1
28.5 1
29 1
29.5 1
30 1
30.5 1
31 2
31.5 15
32 1
32.5 2
33 2
33.5 1
34 2
34.5 2
35 2
35.5 2
36 10
36.5 1
37 2
37.5 1
38 1
38.5 1
39 1
39.5 1
40 1
40.5 2
41 1
41.5 2
42 5
42.5 2
43 1
43.5 1
44 1
44.5 1
45 2
45.5 2
46 2
46.5 2
47 1
47.5 2
48 1
48.5 1
49 1
49.5 3
50 1
50.5 3
51 2
51.5 4
52 1
52.5 1
53 3
53.5 2
54 1
54.5 6
55 4
55.5 2
56 2
56.5 1
57 2
57.5 2
58 1
58.5 10
59 1
59.5 2
60 7
60.5 2
61 2
61.5 1
62 1
62.5 2
63 1
63.5 2
64 1
64.5 6
65 1
65.5 1
66 2
66.5 3
67 1
67.5 1
68 6
68.5 15
69 1
69.5 4
70 3
70.5 1
71 1
71.5 1
72 3
72.5 2
73 1
73.5 1
74 2
74.5 8
75 1
75.5 1
76 1
76.5 1
77 5
77.5 10
78 1
78.5 7
79 3
79.5 7
80 2
80.5 1
81 2
81.5 1
82 4
82.5 1
83 2
83.5 4
84 1
84.5 2
85 1
85.5 1
86 1
86.5 7
87 1
87.5 1
88 2
88.5 2
89 1
89.5 1
90 2
90.5 1
91 1
91.5 1
92 1
92.5 1
93 1
93.5 1
94 8
94.5 1
95 1
95.5 1
96 1
96.5 2
97 1
97.5 1
98 2
98.5 1
99 0
99.5 9
100 1
100.5 3
101 1
101.5 1
102 1
102.5 1
103 1
103.5 9
104 6
104.5 1
105 1
105.5 1
106 1
106.5 5
107 1
107.5 1
108 1
108.5 3
109 6
109.5 0
110 2
110.5 1
111 1
111.5 1
112 1
112.5 1
113 2
113.5 1
114 4
114.5 1
115 1
115.5 5
116 3
116.5 5
117 1
117.5 1
118 2
118.5 11
119 2
119.5 2
120 1
120.5 1
121 0
121.5 1
122 1
122.5 1
123 2
123.5 2
124 1
124.5 1
125 2
125.5 1
126 1
126.5 1
127 1
127.5 2
128 2
128.5 0
129 1
129.5 1
130 1
130.5 2
131 1
131.5 2
132 1
132.5 6
133 1
133.5 1
134 2
134.5 1
135 1
135.5 1
136 7
136.5 3
137 5
137.5 1
138 1
138.5 1
139 2
139.5 1
140 2
140.5 1
141 1
141.5 2
142 2
142.5 1
143 6
143.5 3
144 6
144.5 4
145 2
145.5 1
146 2
146.5 0
147 8
147.5 1
148 4
148.5 1
149 9
149.5 2
150 7
150.5 1
151 2
151.5 8
152 8
152.5 1
153 1
153.5 8
154 4
154.5 2
155 2
155.5 2
156 1
156.5 1
157 2
157.5 2
158 1
158.5 2
159 10
159.5 8
160 2
160.5 10
161 7
161.5 3
162 2
162.5 1
163 1
163.5 5
164 9
164.5 2
165 1
165.5 1
166 1
166.5 2
167 2
167.5 2
168 2
168.5 1
169 2
169.5 0
170 7
170.5 2
171 1
171.5 1
172 1
172.5 1
173 1
173.5 3
174 2
174.5 2
175 1
175.5 2
176 1
176.5 11
177 1
177.5 4
178 2
178.5 1
179 1
179.5 2
180 4
180.5 9
181 1
181.5 1
182 8
182.5 2
183 2
183.5 1
184 1
184.5 11
185 2
185.5 1
186 2
186.5 2
187 3
187.5 2
188 1
188.5 1
189 1
189.5 2
190 1
190.5 2
191 1
191.5 2
192 2
192.5 8
193 5
193.5 2
194 17
194.5 2
195 2
195.5 2
196 2
196.5 3
197 2
197.5 1
198 2
198.5 1
199 1
199.5 1
};
%\addlegendentry{$\epsilon_u = 10^{-11}$}
\end{axis}

\end{tikzpicture}

%% file: figures/mr_KdV.tex
% This file was created with tikzplotlib v0.10.1.
\definecolor{mycolor1}{rgb}{0.00000,0.44700,0.74100}%
\definecolor{mycolor2}{rgb}{0.85000,0.32500,0.09800}%
\definecolor{mycolor3}{rgb}{0.5,0.5,0.5}
\begin{tikzpicture}[scale=0.8]

\begin{axis}[%
width=1.5in,
height=1.5in,
at={(1.011in,0.642in)},
scale only axis,
xmin=0,
xmax=1,
xlabel style={font=\color{white!15!black}},
xlabel={Time},
%ymode=log,
ymin=0,
ymax=601,
yminorticks=true,
ylabel style={font=\color{white!15!black}},
axis background/.style={fill=white},
legend style={fill opacity=0.8, 
font=\tiny,
  draw opacity=1,
  text opacity=1,
  at={(0.02,0.87)},
  anchor=west}
  % draw=lightgray204}
]

\addplot [line width=1.2pt,  mycolor1]
table {%
0 43
0.0005 7
0.001 8
0.0015 8
0.002 8
0.0025 8
0.003 8
0.0035 8
0.004 8
0.0045 8
0.005 8
0.0055 8
0.006 8
0.0065 8
0.007 8
0.0075 8
0.008 8
0.0085 8
0.009 8
0.0095 8
0.01 8
0.0105 8
0.011 8
0.0115 8
0.012 8
0.0125 8
0.013 8
0.0135 8
0.014 8
0.0145 8
0.015 8
0.0155 8
0.016 8
0.0165 8
0.017 8
0.0175 8
0.018 8
0.0185 8
0.019 7
0.0195 7
0.02 7
0.0205 7
0.021 7
0.0215 7
0.022 7
0.0225 7
0.023 7
0.0235 7
0.024 7
0.0245 7
0.025 7
0.0255 7
0.026 7
0.0265 7
0.027 7
0.0275 7
0.028 7
0.0285 7
0.029 7
0.0295 7
0.03 7
0.0305 7
0.031 7
0.0315 7
0.032 7
0.0325 7
0.033 7
0.0335 7
0.034 7
0.0345 8
0.035 7
0.0355 7
0.036 7
0.0365 7
0.037 7
0.0375 7
0.038 7
0.0385 7
0.039 7
0.0395 7
0.04 7
0.0405 7
0.041 7
0.0415 7
0.042 7
0.0425 7
0.043 7
0.0435 7
0.044 7
0.0445 7
0.045 7
0.0455 7
0.046 7
0.0465 7
0.047 7
0.0475 7
0.048 7
0.0485 7
0.049 7
0.0495 7
0.05 7
0.0505 7
0.051 7
0.0515 7
0.052 7
0.0525 7
0.053 7
0.0535 7
0.054 7
0.0545 7
0.055 7
0.0555 7
0.056 7
0.0565 7
0.057 7
0.0575 7
0.058 7
0.0585 7
0.059 7
0.0595 7
0.06 303
0.0605 7
0.061 7
0.0615 7
0.062 7
0.0625 7
0.063 7
0.0635 7
0.064 7
0.0645 7
0.065 7
0.0655 7
0.066 7
0.0665 7
0.067 7
0.0675 7
0.068 7
0.0685 7
0.069 7
0.0695 7
0.07 7
0.0705 7
0.071 7
0.0715 7
0.072 7
0.0725 7
0.073 7
0.0735 7
0.074 7
0.0745 7
0.075 7
0.0755 7
0.076 7
0.0765 7
0.077 7
0.0775 7
0.078 7
0.0785 7
0.079 7
0.0795 7
0.08 7
0.0805 7
0.081 7
0.0815 7
0.082 7
0.0825 7
0.083 7
0.0835 7
0.084 7
0.0845 7
0.085 7
0.0855 7
0.086 7
0.0865 7
0.087 7
0.0875 7
0.088 7
0.0885 7
0.089 7
0.0895 7
0.09 7
0.0905 7
0.091 7
0.0915 7
0.092 7
0.0925 7
0.093 7
0.0935 7
0.094 7
0.0945 7
0.095 7
0.0955 7
0.096 7
0.0965 7
0.097 7
0.0975 7
0.098 7
0.0985 7
0.099 7
0.0995 7
0.1 7
0.1005 7
0.101 7
0.1015 7
0.102 7
0.1025 7
0.103 7
0.1035 7
0.104 7
0.1045 7
0.105 7
0.1055 7
0.106 7
0.1065 7
0.107 7
0.1075 7
0.108 7
0.1085 7
0.109 7
0.1095 7
0.11 7
0.1105 7
0.111 7
0.1115 7
0.112 7
0.1125 7
0.113 7
0.1135 7
0.114 7
0.1145 7
0.115 7
0.1155 7
0.116 7
0.1165 7
0.117 7
0.1175 7
0.118 7
0.1185 7
0.119 7
0.1195 7
0.12 7
0.1205 7
0.121 7
0.1215 7
0.122 7
0.1225 7
0.123 7
0.1235 7
0.124 7
0.1245 7
0.125 7
0.1255 7
0.126 7
0.1265 7
0.127 7
0.1275 7
0.128 7
0.1285 7
0.129 7
0.1295 7
0.13 7
0.1305 7
0.131 7
0.1315 7
0.132 7
0.1325 7
0.133 7
0.1335 7
0.134 7
0.1345 7
0.135 7
0.1355 7
0.136 6
0.1365 6
0.137 6
0.1375 6
0.138 6
0.1385 6
0.139 6
0.1395 6
0.14 6
0.1405 6
0.141 6
0.1415 6
0.142 6
0.1425 6
0.143 6
0.1435 6
0.144 6
0.1445 6
0.145 6
0.1455 6
0.146 6
0.1465 6
0.147 6
0.1475 6
0.148 6
0.1485 6
0.149 6
0.1495 6
0.15 6
0.1505 6
0.151 6
0.1515 6
0.152 6
0.1525 6
0.153 6
0.1535 6
0.154 6
0.1545 6
0.155 6
0.1555 6
0.156 6
0.1565 6
0.157 6
0.1575 6
0.158 6
0.1585 6
0.159 6
0.1595 6
0.16 6
0.1605 6
0.161 6
0.1615 6
0.162 6
0.1625 6
0.163 6
0.1635 6
0.164 6
0.1645 6
0.165 6
0.1655 6
0.166 6
0.1665 6
0.167 6
0.1675 6
0.168 6
0.1685 6
0.169 6
0.1695 6
0.17 6
0.1705 6
0.171 6
0.1715 6
0.172 6
0.1725 6
0.173 6
0.1735 6
0.174 6
0.1745 6
0.175 6
0.1755 6
0.176 6
0.1765 6
0.177 6
0.1775 6
0.178 6
0.1785 6
0.179 6
0.1795 6
0.18 6
0.1805 6
0.181 6
0.1815 6
0.182 6
0.1825 6
0.183 6
0.1835 6
0.184 6
0.1845 6
0.185 6
0.1855 6
0.186 6
0.1865 6
0.187 6
0.1875 6
0.188 6
0.1885 6
0.189 6
0.1895 6
0.19 6
0.1905 6
0.191 6
0.1915 6
0.192 6
0.1925 6
0.193 6
0.1935 6
0.194 6
0.1945 6
0.195 6
0.1955 6
0.196 6
0.1965 6
0.197 6
0.1975 6
0.198 6
0.1985 6
0.199 6
0.1995 6
0.2 6
0.2005 6
0.201 6
0.2015 6
0.202 6
0.2025 6
0.203 6
0.2035 6
0.204 6
0.2045 6
0.205 6
0.2055 6
0.206 6
0.2065 6
0.207 6
0.2075 6
0.208 6
0.2085 6
0.209 6
0.2095 6
0.21 6
0.2105 6
0.211 6
0.2115 6
0.212 6
0.2125 6
0.213 6
0.2135 6
0.214 6
0.2145 6
0.215 6
0.2155 6
0.216 6
0.2165 6
0.217 6
0.2175 6
0.218 6
0.2185 5
0.219 5
0.2195 5
0.22 5
0.2205 5
0.221 5
0.2215 5
0.222 5
0.2225 5
0.223 5
0.2235 5
0.224 5
0.2245 5
0.225 5
0.2255 5
0.226 5
0.2265 5
0.227 5
0.2275 5
0.228 5
0.2285 5
0.229 5
0.2295 5
0.23 5
0.2305 5
0.231 5
0.2315 5
0.232 5
0.2325 5
0.233 5
0.2335 5
0.234 5
0.2345 5
0.235 5
0.2355 5
0.236 5
0.2365 5
0.237 5
0.2375 5
0.238 5
0.2385 5
0.239 5
0.2395 5
0.24 5
0.2405 5
0.241 5
0.2415 5
0.242 5
0.2425 5
0.243 5
0.2435 5
0.244 5
0.2445 5
0.245 5
0.2455 5
0.246 5
0.2465 5
0.247 5
0.2475 5
0.248 5
0.2485 5
0.249 5
0.2495 5
0.25 5
0.2505 5
0.251 5
0.2515 5
0.252 5
0.2525 5
0.253 5
0.2535 5
0.254 5
0.2545 5
0.255 5
0.2555 5
0.256 5
0.2565 5
0.257 5
0.2575 5
0.258 5
0.2585 5
0.259 5
0.2595 5
0.26 5
0.2605 5
0.261 5
0.2615 5
0.262 5
0.2625 5
0.263 5
0.2635 5
0.264 5
0.2645 5
0.265 5
0.2655 5
0.266 5
0.2665 5
0.267 5
0.2675 5
0.268 5
0.2685 5
0.269 5
0.2695 5
0.27 5
0.2705 5
0.271 5
0.2715 5
0.272 5
0.2725 5
0.273 5
0.2735 5
0.274 5
0.2745 5
0.275 5
0.2755 5
0.276 5
0.2765 5
0.277 5
0.2775 5
0.278 5
0.2785 5
0.279 5
0.2795 5
0.28 5
0.2805 5
0.281 5
0.2815 5
0.282 5
0.2825 5
0.283 5
0.2835 5
0.284 5
0.2845 5
0.285 5
0.2855 5
0.286 5
0.2865 5
0.287 5
0.2875 5
0.288 5
0.2885 5
0.289 5
0.2895 5
0.29 5
0.2905 5
0.291 5
0.2915 5
0.292 5
0.2925 5
0.293 5
0.2935 5
0.294 5
0.2945 5
0.295 5
0.2955 5
0.296 5
0.2965 5
0.297 5
0.2975 5
0.298 5
0.2985 5
0.299 5
0.2995 5
0.3 5
0.3005 5
0.301 5
0.3015 5
0.302 5
0.3025 5
0.303 5
0.3035 5
0.304 5
0.3045 5
0.305 5
0.3055 5
0.306 5
0.3065 5
0.307 4
0.3075 4
0.308 4
0.3085 4
0.309 4
0.3095 4
0.31 4
0.3105 4
0.311 4
0.3115 4
0.312 4
0.3125 4
0.313 4
0.3135 4
0.314 4
0.3145 4
0.315 4
0.3155 4
0.316 4
0.3165 4
0.317 4
0.3175 4
0.318 4
0.3185 4
0.319 4
0.3195 4
0.32 4
0.3205 4
0.321 4
0.3215 4
0.322 4
0.3225 4
0.323 4
0.3235 4
0.324 4
0.3245 4
0.325 4
0.3255 4
0.326 4
0.3265 4
0.327 4
0.3275 4
0.328 4
0.3285 4
0.329 4
0.3295 4
0.33 4
0.3305 4
0.331 4
0.3315 4
0.332 4
0.3325 4
0.333 4
0.3335 4
0.334 4
0.3345 4
0.335 4
0.3355 4
0.336 4
0.3365 4
0.337 4
0.3375 4
0.338 4
0.3385 4
0.339 4
0.3395 4
0.34 4
0.3405 4
0.341 4
0.3415 4
0.342 4
0.3425 4
0.343 4
0.3435 4
0.344 4
0.3445 4
0.345 4
0.3455 4
0.346 4
0.3465 4
0.347 4
0.3475 4
0.348 4
0.3485 4
0.349 4
0.3495 4
0.35 4
0.3505 4
0.351 4
0.3515 4
0.352 4
0.3525 4
0.353 4
0.3535 4
0.354 4
0.3545 4
0.355 4
0.3555 4
0.356 4
0.3565 4
0.357 4
0.3575 4
0.358 4
0.3585 4
0.359 4
0.3595 4
0.36 4
0.3605 4
0.361 4
0.3615 4
0.362 4
0.3625 4
0.363 4
0.3635 4
0.364 4
0.3645 4
0.365 4
0.3655 4
0.366 4
0.3665 4
0.367 4
0.3675 4
0.368 4
0.3685 4
0.369 4
0.3695 4
0.37 4
0.3705 4
0.371 4
0.3715 4
0.372 4
0.3725 4
0.373 4
0.3735 4
0.374 4
0.3745 4
0.375 4
0.3755 4
0.376 4
0.3765 4
0.377 4
0.3775 4
0.378 4
0.3785 4
0.379 4
0.3795 4
0.38 4
0.3805 4
0.381 4
0.3815 4
0.382 4
0.3825 3
0.383 3
0.3835 3
0.384 3
0.3845 3
0.385 3
0.3855 3
0.386 3
0.3865 3
0.387 3
0.3875 3
0.388 3
0.3885 3
0.389 3
0.3895 3
0.39 3
0.3905 3
0.391 3
0.3915 3
0.392 3
0.3925 3
0.393 3
0.3935 3
0.394 3
0.3945 3
0.395 3
0.3955 3
0.396 3
0.3965 3
0.397 3
0.3975 3
0.398 3
0.3985 3
0.399 3
0.3995 3
0.4 3
0.4005 3
0.401 3
0.4015 3
0.402 3
0.4025 3
0.403 3
0.4035 3
0.404 3
0.4045 3
0.405 3
0.4055 3
0.406 3
0.4065 3
0.407 3
0.4075 3
0.408 3
0.4085 3
0.409 3
0.4095 3
0.41 3
0.4105 3
0.411 3
0.4115 3
0.412 3
0.4125 3
0.413 3
0.4135 3
0.414 3
0.4145 3
0.415 3
0.4155 3
0.416 3
0.4165 3
0.417 3
0.4175 3
0.418 3
0.4185 3
0.419 3
0.4195 3
0.42 3
0.4205 3
0.421 3
0.4215 3
0.422 3
0.4225 3
0.423 3
0.4235 3
0.424 3
0.4245 3
0.425 3
0.4255 3
0.426 3
0.4265 3
0.427 3
0.4275 3
0.428 3
0.4285 3
0.429 3
0.4295 3
0.43 3
0.4305 3
0.431 3
0.4315 3
0.432 3
0.4325 3
0.433 3
0.4335 3
0.434 3
0.4345 3
0.435 3
0.4355 3
0.436 3
0.4365 3
0.437 3
0.4375 3
0.438 3
0.4385 3
0.439 3
0.4395 3
0.44 3
0.4405 3
0.441 3
0.4415 3
0.442 3
0.4425 3
0.443 3
0.4435 3
0.444 3
0.4445 3
0.445 3
0.4455 3
0.446 3
0.4465 3
0.447 3
0.4475 3
0.448 3
0.4485 3
0.449 3
0.4495 3
0.45 3
0.4505 3
0.451 3
0.4515 3
0.452 3
0.4525 3
0.453 3
0.4535 3
0.454 3
0.4545 3
0.455 3
0.4555 3
0.456 3
0.4565 3
0.457 3
0.4575 3
0.458 3
0.4585 3
0.459 3
0.4595 3
0.46 3
0.4605 3
0.461 3
0.4615 3
0.462 3
0.4625 3
0.463 2
0.4635 2
0.464 2
0.4645 2
0.465 2
0.4655 2
0.466 2
0.4665 2
0.467 2
0.4675 2
0.468 2
0.4685 2
0.469 2
0.4695 2
0.47 2
0.4705 2
0.471 2
0.4715 2
0.472 3
0.4725 3
0.473 3
0.4735 3
0.474 3
0.4745 3
0.475 3
0.4755 3
0.476 3
0.4765 3
0.477 3
0.4775 3
0.478 3
0.4785 3
0.479 3
0.4795 3
0.48 3
0.4805 3
0.481 3
0.4815 3
0.482 3
0.4825 3
0.483 3
0.4835 3
0.484 3
0.4845 3
0.485 3
0.4855 3
0.486 3
0.4865 3
0.487 3
0.4875 3
0.488 3
0.4885 3
0.489 3
0.4895 3
0.49 3
0.4905 3
0.491 3
0.4915 3
0.492 3
0.4925 3
0.493 3
0.4935 3
0.494 3
0.4945 3
0.495 3
0.4955 3
0.496 3
0.4965 3
0.497 3
0.4975 3
0.498 3
0.4985 3
0.499 3
0.4995 3
0.5 3
0.5005 3
0.501 3
0.5015 3
0.502 3
0.5025 3
0.503 3
0.5035 3
0.504 3
0.5045 3
0.505 3
0.5055 3
0.506 3
0.5065 3
0.507 3
0.5075 3
0.508 3
0.5085 3
0.509 3
0.5095 3
0.51 3
0.5105 3
0.511 3
0.5115 3
0.512 3
0.5125 3
0.513 3
0.5135 3
0.514 3
0.5145 3
0.515 3
0.5155 3
0.516 3
0.5165 3
0.517 3
0.5175 3
0.518 3
0.5185 3
0.519 3
0.5195 3
0.52 3
0.5205 3
0.521 3
0.5215 3
0.522 3
0.5225 3
0.523 3
0.5235 3
0.524 3
0.5245 3
0.525 3
0.5255 3
0.526 3
0.5265 3
0.527 3
0.5275 3
0.528 3
0.5285 3
0.529 3
0.5295 3
0.53 2
0.5305 3
0.531 3
0.5315 3
0.532 2
0.5325 2
0.533 3
0.5335 2
0.534 3
0.5345 2
0.535 2
0.5355 2
0.536 2
0.5365 2
0.537 2
0.5375 2
0.538 3
0.5385 2
0.539 2
0.5395 2
0.54 3
0.5405 3
0.541 3
0.5415 3
0.542 3
0.5425 3
0.543 2
0.5435 3
0.544 2
0.5445 3
0.545 3
0.5455 3
0.546 3
0.5465 3
0.547 2
0.5475 2
0.548 2
0.5485 2
0.549 2
0.5495 2
0.55 2
0.5505 2
0.551 2
0.5515 2
0.552 2
0.5525 2
0.553 2
0.5535 2
0.554 2
0.5545 2
0.555 2
0.5555 2
0.556 2
0.5565 2
0.557 2
0.5575 2
0.558 2
0.5585 2
0.559 2
0.5595 2
0.56 2
0.5605 2
0.561 2
0.5615 2
0.562 2
0.5625 2
0.563 2
0.5635 2
0.564 2
0.5645 2
0.565 2
0.5655 2
0.566 2
0.5665 2
0.567 2
0.5675 2
0.568 2
0.5685 2
0.569 2
0.5695 2
0.57 2
0.5705 2
0.571 2
0.5715 2
0.572 2
0.5725 2
0.573 2
0.5735 2
0.574 2
0.5745 2
0.575 2
0.5755 2
0.576 2
0.5765 2
0.577 2
0.5775 2
0.578 2
0.5785 2
0.579 2
0.5795 2
0.58 2
0.5805 2
0.581 2
0.5815 2
0.582 2
0.5825 2
0.583 2
0.5835 2
0.584 2
0.5845 2
0.585 2
0.5855 2
0.586 2
0.5865 2
0.587 2
0.5875 2
0.588 2
0.5885 2
0.589 2
0.5895 2
0.59 2
0.5905 2
0.591 2
0.5915 2
0.592 2
0.5925 2
0.593 2
0.5935 2
0.594 2
0.5945 2
0.595 2
0.5955 2
0.596 2
0.5965 2
0.597 2
0.5975 2
0.598 2
0.5985 2
0.599 2
0.5995 2
0.6 2
0.6005 1
0.601 2
0.6015 2
0.602 2
0.6025 2
0.603 2
0.6035 2
0.604 1
0.6045 1
0.605 1
0.6055 1
0.606 1
0.6065 1
0.607 2
0.6075 1
0.608 2
0.6085 2
0.609 1
0.6095 0
0.61 1
0.6105 1
0.611 0
0.6115 0
0.612 0
0.6125 0
0.613 0
0.6135 0
0.614 1
0.6145 0
0.615 0
0.6155 0
0.616 0
0.6165 0
0.617 0
0.6175 0
0.618 1
0.6185 0
0.619 0
0.6195 0
0.62 0
0.6205 0
0.621 0
0.6215 0
0.622 0
0.6225 0
0.623 0
0.6235 0
0.624 0
0.6245 0
0.625 1
0.6255 1
0.626 1
0.6265 1
0.627 1
0.6275 1
0.628 1
0.6285 1
0.629 2
0.6295 2
0.63 0
0.6305 0
0.631 0
0.6315 0
0.632 0
0.6325 0
0.633 0
0.6335 0
0.634 0
0.6345 0
0.635 0
0.6355 0
0.636 0
0.6365 0
0.637 0
0.6375 0
0.638 0
0.6385 1
0.639 1
0.6395 1
0.64 1
0.6405 0
0.641 0
0.6415 0
0.642 1
0.6425 1
0.643 1
0.6435 1
0.644 1
0.6445 1
0.645 1
0.6455 1
0.646 1
0.6465 1
0.647 1
0.6475 1
0.648 1
0.6485 1
0.649 1
0.6495 1
0.65 0
% 0.6505 0
% 0.651 0
% 0.6515 0
% 0.652 0
% 0.6525 0
% 0.653 0
% 0.6535 0
% 0.654 0
0.6545 0
0.655 0
0.6555 0
0.656 0
0.6565 1
0.657 1
0.6575 1
0.658 1
0.6585 0
0.659 0
0.6595 1
0.66 1
0.6605 1
0.661 1
0.6615 1
0.662 1
% 0.6625 1
% 0.663 1
% 0.6635 1
% 0.664 0
% 0.6645 1
% 0.665 1
% 0.6655 1
% 0.666 1
% 0.6665 1
% 0.667 0
% 0.6675 1
% 0.668 1
% 0.6685 1
% 0.669 1
% 0.6695 1
% 0.67 1
% 0.6705 1
% 0.671 1
% 0.6715 1
% 0.672 1
% 0.6725 1
% 0.673 1
% 0.6735 1
% 0.674 1
% 0.6745 1
% 0.675 1
% 0.6755 1
% 0.676 1
% 0.6765 1
% 0.677 1
% 0.6775 1
% 0.678 1
0.6785 1
0.679 1
0.6795 1
0.68 1
0.6805 464
0.681 0
0.6815 1
0.682 1
0.6825 1
0.683 0
0.6835 0
% 0.684 0
% 0.6845 0
% 0.685 0
% 0.6855 0
% 0.686 0
% 0.6865 0
% 0.687 0
% 0.6875 0
% 0.688 0
0.6885 0
0.689 0
0.6895 0
0.69 0
0.6905 1
0.691 1
0.6915 1
0.692 0
0.6925 0
0.693 0
0.6935 0
0.694 0
0.6945 0
0.695 0
0.6955 0
0.696 0
0.6965 0
0.697 0
0.6975 0
0.698 0
0.6985 1
0.699 1
0.6995 1
0.7 1
0.7005 1
0.701 0
0.7015 0
0.702 1
0.7025 1
0.703 1
0.7035 1
0.704 1
0.7045 1
0.705 1
0.7055 1
0.706 0
0.7065 1
0.707 1
0.7075 1
0.708 1
0.7085 1
0.709 1
0.7095 1
0.71 1
0.7105 1
0.711 1
0.7115 1
0.712 1
0.7125 1
0.713 0
% 0.7135 0
% 0.714 0
% 0.7145 0
% 0.715 0
% 0.7155 0
% 0.716 0
% 0.7165 0
% 0.717 0
% 0.7175 0
% 0.718 0
% 0.7185 0
% 0.719 0
% 0.7195 0
% 0.72 0
% 0.7205 0
% 0.721 0
% 0.7215 0
% 0.722 0
% 0.7225 0
% 0.723 0
% 0.7235 0
% 0.724 0
% 0.7245 0
% 0.725 0
% 0.7255 0
% 0.726 0
0.7265 0
0.727 1
 0.7275 0
% 0.728 0
% 0.7285 0
% 0.729 0
% 0.7295 0
% 0.73 0
% 0.7305 0
% 0.731 0
% 0.7315 0
% 0.732 0
% 0.7325 0
% 0.733 0
% 0.7335 0
% 0.734 0
% 0.7345 0
% 0.735 0
% 0.7355 0
% 0.736 0
% 0.7365 0
% 0.737 0
% 0.7375 0
% 0.738 0
% 0.7385 0
% 0.739 0
% 0.7395 0
% 0.74 0
% 0.7405 0
% 0.741 0
% 0.7415 0
% 0.742 0
% 0.7425 0
% 0.743 0
% 0.7435 0
% 0.744 0
% 0.7445 0
% 0.745 0
% 0.7455 0
% 0.746 0
% 0.7465 0
% 0.747 0
% 0.7475 0
% 0.748 0
% 0.7485 0
% 0.749 0
% 0.7495 0
% 0.75 0
% 0.7505 0
% 0.751 0
% 0.7515 0
% 0.752 0
% 0.7525 0
% 0.753 0
% 0.7535 0
% 0.754 0
% 0.7545 0
% 0.755 0
% 0.7555 0
% 0.756 0
% 0.7565 0
% 0.757 0
% 0.7575 0
% 0.758 0
% 0.7585 0
% 0.759 0
% 0.7595 0
% 0.76 0
% 0.7605 0
% 0.761 0
% 0.7615 0
% 0.762 0
% 0.7625 0
% 0.763 0
% 0.7635 0
% 0.764 0
% 0.7645 0
% 0.765 0
% 0.7655 0
% 0.766 0
% 0.7665 0
% 0.767 0
% 0.7675 0
% 0.768 0
% 0.7685 0
% 0.769 0
% 0.7695 0
% 0.77 0
% 0.7705 0
% 0.771 0
% 0.7715 0
% 0.772 0
% 0.7725 0
% 0.773 0
% 0.7735 0
% 0.774 0
% 0.7745 0
% 0.775 0
% 0.7755 0
% 0.776 0
% 0.7765 0
% 0.777 0
% 0.7775 0
% 0.778 0
% 0.7785 0
% 0.779 0
% 0.7795 0
% 0.78 0
% 0.7805 0
% 0.781 0
% 0.7815 0
% 0.782 0
% 0.7825 0
% 0.783 0
% 0.7835 0
% 0.784 0
% 0.7845 0
% 0.785 0
% 0.7855 0
% 0.786 0
% 0.7865 0
% 0.787 0
% 0.7875 0
% 0.788 0
% 0.7885 0
% 0.789 0
% 0.7895 0
% 0.79 0
% 0.7905 0
% 0.791 0
% 0.7915 0
% 0.792 0
% 0.7925 0
0.793 0
0.7935 249
0.794 0
0.7945 0
% 0.795 0
% 0.7955 0
% 0.796 0
% 0.7965 0
% 0.797 0
% 0.7975 0
% 0.798 0
% 0.7985 0
0.799 0
0.7995 0
0.936 0
0.9365 0
% 0.937 0
% 0.9375 0
% 0.938 0
% 0.9385 0
% 0.939 0
% 0.9395 0
% 0.94 0
% 0.9405 0
% 0.941 0
% 0.9415 0
% 0.942 0
% 0.9425 0
% 0.943 0
% 0.9435 0
% 0.944 0
% 0.9445 0
% 0.945 0
% 0.9455 0
% 0.946 0
% 0.9465 0
% 0.947 0
% 0.9475 0
% 0.948 0
% 0.9485 0
% 0.949 0
% 0.9495 0
% 0.95 0
% 0.9505 0
% 0.951 0
% 0.9515 0
% 0.952 0
% 0.9525 0
% 0.953 0
% 0.9535 0
% 0.954 0
% 0.9545 0
% 0.955 0
% 0.9555 0
% 0.956 0
% 0.9565 0
% 0.957 0
% 0.9575 0
% 0.958 0
% 0.9585 0
% 0.959 0
% 0.9595 0
% 0.96 0
% 0.9605 0
% 0.961 0
% 0.9615 0
% 0.962 0
% 0.9625 0
% 0.963 0
% 0.9635 0
% 0.964 0
% 0.9645 0
% 0.965 0
% 0.9655 0
% 0.966 0
% 0.9665 0
% 0.967 0
% 0.9675 0
% 0.968 0
% 0.9685 0
% 0.969 0
% 0.9695 0
% 0.97 0
% 0.9705 0
% 0.971 0
% 0.9715 0
% 0.972 0
% 0.9725 0
% 0.973 0
% 0.9735 0
% 0.974 0
% 0.9745 0
% 0.975 0
% 0.9755 0
% 0.976 0
% 0.9765 0
% 0.977 0
% 0.9775 0
% 0.978 0
% 0.9785 0
% 0.979 0
0.9795 0
0.98 0
0.9805 1
0.981 0
% 0.9815 0
% 0.982 0
% 0.9825 0
% 0.983 0
% 0.9835 0
% 0.984 0
% 0.9845 0
% 0.985 0
% 0.9855 0
% 0.986 0
% 0.9865 0
% 0.987 0
% 0.9875 0
% 0.988 0
% 0.9885 0
% 0.989 0
% 0.9895 0
% 0.99 0
% 0.9905 0
% 0.991 0
% 0.9915 0
0.992 0
0.9925 0
0.993 1
0.9935 1
0.994 0
0.9945 0
% 0.995 0
% 0.9955 0
% 0.996 0
% 0.9965 0
0.997 0
0.9975 1
0.998 0
0.9985 0
0.999 0
0.9995 0
};
%\addlegendentry{$\epsilon_u = 10^{-8}$}

\addplot [line width=1.2pt,  mycolor2,opacity=0.5]
table {%
0 43
0.0005 8
0.001 7
0.0015 7
0.002 7
0.0025 7
0.003 7
0.0035 7
0.004 7
0.0045 7
0.005 12
0.0055 7
0.006 7
0.0065 7
0.007 7
0.0075 7
0.008 7
0.0085 7
0.009 7
0.0095 6
0.01 6
0.0105 7
0.011 7
0.0115 7
0.012 7
0.0125 7
0.013 7
0.0135 7
0.014 7
0.0145 8
0.015 7
0.0155 7
0.016 7
0.0165 7
0.017 7
0.0175 7
0.018 7
0.0185 8
0.019 7
0.0195 7
0.02 7
0.0205 7
0.021 7
0.0215 7
0.022 8
0.0225 8
0.023 8
0.0235 8
0.024 8
0.0245 8
0.025 8
0.0255 8
0.026 8
0.0265 8
0.027 7
0.0275 8
0.028 7
0.0285 8
0.029 7
0.0295 7
0.03 7
0.0305 8
0.031 7
0.0315 7
0.032 7
0.0325 7
0.033 7
0.0335 8
0.034 8
0.0345 8
0.035 8
0.0355 8
0.036 8
0.0365 8
0.037 8
0.0375 8
0.038 8
0.0385 7
0.039 7
0.0395 7
0.04 7
0.0405 6
0.041 7
0.0415 7
0.042 8
0.0425 7
0.043 7
0.0435 7
0.044 7
0.0445 7
0.045 7
0.0455 6
0.046 8
0.0465 8
0.047 8
0.0475 7
0.048 8
0.0485 8
0.049 8
0.0495 7
0.05 7
0.0505 6
0.051 6
0.0515 7
0.052 7
0.0525 6
0.053 6
0.0535 7
0.054 7
0.0545 7
0.055 7
0.0555 7
0.056 6
0.0565 7
0.057 7
0.0575 6
0.058 8
0.0585 7
0.059 7
0.0595 7
0.06 7
0.0605 7
0.061 6
0.0615 6
0.062 6
0.0625 6
0.063 6
0.0635 6
0.064 6
0.0645 6
0.065 6
0.0655 6
0.066 6
0.0665 6
0.067 6
0.0675 6
0.068 6
0.0685 6
0.069 6
0.0695 6
0.07 6
0.0705 6
0.071 6
0.0715 6
0.072 6
0.0725 6
0.073 6
0.0735 6
0.074 6
0.0745 6
0.075 6
0.0755 6
0.076 6
0.0765 6
0.077 6
0.0775 6
0.078 6
0.0785 6
0.079 6
0.0795 6
0.08 6
0.0805 6
0.081 6
0.0815 6
0.082 7
0.0825 6
0.083 6
0.0835 6
0.084 6
0.0845 5
0.085 6
0.0855 6
0.086 6
0.0865 6
0.087 6
0.0875 6
0.088 6
0.0885 6
0.089 6
0.0895 6
0.09 6
0.0905 6
0.091 5
0.0915 6
0.092 6
0.0925 6
0.093 6
0.0935 7
0.094 5
0.0945 5
0.095 6
0.0955 6
0.096 6
0.0965 6
0.097 6
0.0975 6
0.098 6
0.0985 5
0.099 6
0.0995 5
0.1 6
0.1005 6
0.101 6
0.1015 6
0.102 6
0.1025 5
0.103 6
0.1035 5
0.104 6
0.1045 6
0.105 6
0.1055 5
0.106 6
0.1065 6
0.107 6
0.1075 5
0.108 6
0.1085 6
0.109 6
0.1095 6
0.11 5
0.1105 5
0.111 5
0.1115 5
0.112 5
0.1125 4
0.113 5
0.1135 5
0.114 5
0.1145 6
0.115 5
0.1155 5
0.116 5
0.1165 6
0.117 6
0.1175 5
0.118 7
0.1185 5
0.119 5
0.1195 5
0.12 6
0.1205 6
0.121 6
0.1215 6
0.122 6
0.1225 5
0.123 5
0.1235 5
0.124 5
0.1245 6
0.125 5
0.1255 5
0.126 6
0.1265 6
0.127 7
0.1275 6
0.128 7
0.1285 5
0.129 5
0.1295 6
0.13 6
0.1305 5
0.131 5
0.1315 6
0.132 6
0.1325 6
0.133 6
0.1335 6
0.134 6
0.1345 5
0.135 6
0.1355 6
0.136 6
0.1365 7
0.137 6
0.1375 6
0.138 5
0.1385 7
0.139 7
0.1395 6
0.14 7
0.1405 6
0.141 6
0.1415 6
0.142 6
0.1425 6
0.143 6
0.1435 6
0.144 6
0.1445 6
0.145 6
0.1455 6
0.146 6
0.1465 6
0.147 6
0.1475 6
0.148 6
0.1485 7
0.149 7
0.1495 6
0.15 7
0.1505 7
0.151 7
0.1515 6
0.152 6
0.1525 6
0.153 6
0.1535 6
0.154 6
0.1545 6
0.155 6
0.1555 6
0.156 7
0.1565 7
0.157 6
0.1575 6
0.158 6
0.1585 6
0.159 7
0.1595 7
0.16 5
0.1605 6
0.161 5
0.1615 6
0.162 7
0.1625 6
0.163 5
0.1635 5
0.164 7
0.1645 6
0.165 6
0.1655 6
0.166 6
0.1665 6
0.167 6
0.1675 7
0.168 7
0.1685 7
0.169 6
0.1695 6
0.17 6
0.1705 5
0.171 7
0.1715 5
0.172 7
0.1725 7
0.173 5
0.1735 6
0.174 7
0.1745 7
0.175 6
0.1755 6
0.176 6
0.1765 6
0.177 6
0.1775 6
0.178 7
0.1785 7
0.179 6
0.1795 6
0.18 7
0.1805 7
0.181 6
0.1815 6
0.182 6
0.1825 7
0.183 7
0.1835 5
0.184 7
0.1845 6
0.185 7
0.1855 7
0.186 6
0.1865 7
0.187 6
0.1875 6
0.188 6
0.1885 6
0.189 7
0.1895 6
0.19 5
0.1905 7
0.191 7
0.1915 6
0.192 5
0.1925 7
0.193 6
0.1935 7
0.194 7
0.1945 6
0.195 5
0.1955 6
0.196 5
0.1965 5
0.197 6
0.1975 5
0.198 6
0.1985 7
0.199 6
0.1995 7
0.2 7
0.2005 7
0.201 6
0.2015 5
0.202 8
0.2025 7
0.203 7
0.2035 6
0.204 6
0.2045 7
0.205 6
0.2055 6
0.206 7
0.2065 5
0.207 7
0.2075 7
0.208 8
0.2085 6
0.209 6
0.2095 7
0.21 5
0.2105 6
0.211 5
0.2115 5
0.212 5
0.2125 6
0.213 6
0.2135 7
0.214 7
0.2145 6
0.215 6
0.2155 6
0.216 6
0.2165 8
0.217 8
0.2175 6
0.218 7
0.2185 7
0.219 7
0.2195 6
0.22 6
0.2205 7
0.221 6
0.2215 6
0.222 6
0.2225 7
0.223 7
0.2235 7
0.224 6
0.2245 6
0.225 7
0.2255 6
0.226 8
0.2265 7
0.227 8
0.2275 7
0.228 7
0.2285 7
0.229 7
0.2295 7
0.23 7
0.2305 6
0.231 7
0.2315 6
0.232 6
0.2325 7
0.233 7
0.2335 6
0.234 6
0.2345 7
0.235 8
0.2355 6
0.236 6
0.2365 8
0.237 6
0.2375 6
0.238 7
0.2385 7
0.239 7
0.2395 6
0.24 6
0.2405 6
0.241 7
0.2415 8
0.242 6
0.2425 7
0.243 7
0.2435 6
0.244 6
0.2445 8
0.245 7
0.2455 8
0.246 8
0.2465 6
0.247 7
0.2475 7
0.248 6
0.2485 6
0.249 5
0.2495 7
0.25 7
0.2505 6
0.251 6
0.2515 6
0.252 6
0.2525 6
0.253 7
0.2535 7
0.254 6
0.2545 8
0.255 6
0.2555 7
0.256 6
0.2565 7
0.257 7
0.2575 7
0.258 6
0.2585 6
0.259 6
0.2595 8
0.26 6
0.2605 6
0.261 6
0.2615 7
0.262 7
0.2625 8
0.263 6
0.2635 7
0.264 7
0.2645 7
0.265 7
0.2655 7
0.266 6
0.2665 7
0.267 6
0.2675 7
0.268 6
0.2685 6
0.269 6
0.2695 7
0.27 6
0.2705 6
0.271 7
0.2715 6
0.272 7
0.2725 6
0.273 8
0.2735 6
0.274 6
0.2745 6
0.275 7
0.2755 6
0.276 6
0.2765 6
0.277 7
0.2775 8
0.278 7
0.2785 7
0.279 8
0.2795 9
0.28 6
0.2805 7
0.281 6
0.2815 7
0.282 7
0.2825 7
0.283 7
0.2835 7
0.284 7
0.2845 8
0.285 7
0.2855 7
0.286 7
0.2865 6
0.287 6
0.2875 6
0.288 7
0.2885 7
0.289 7
0.2895 6
0.29 7
0.2905 7
0.291 7
0.2915 6
0.292 7
0.2925 7
0.293 7
0.2935 7
0.294 6
0.2945 6
0.295 8
0.2955 6
0.296 5
0.2965 5
0.297 6
0.2975 7
0.298 7
0.2985 6
0.299 8
0.2995 7
0.3 7
0.3005 7
0.301 7
0.3015 7
0.302 6
0.3025 6
0.303 6
0.3035 7
0.304 6
0.3045 5
0.305 6
0.3055 7
0.306 7
0.3065 7
0.307 6
0.3075 7
0.308 6
0.3085 6
0.309 6
0.3095 9
0.31 7
0.3105 5
0.311 7
0.3115 8
0.312 7
0.3125 7
0.313 6
0.3135 7
0.314 7
0.3145 7
0.315 6
0.3155 8
0.316 7
0.3165 7
0.317 6
0.3175 8
0.318 6
0.3185 6
0.319 7
0.3195 7
0.32 6
0.3205 6
0.321 6
0.3215 6
0.322 7
0.3225 6
0.323 7
0.3235 6
0.324 6
0.3245 7
0.325 7
0.3255 7
0.326 8
0.3265 7
0.327 6
0.3275 6
0.328 6
0.3285 8
0.329 6
0.3295 7
0.33 7
0.3305 6
0.331 6
0.3315 7
0.332 7
0.3325 7
0.333 6
0.3335 6
0.334 8
0.3345 7
0.335 6
0.3355 6
0.336 7
0.3365 6
0.337 7
0.3375 7
0.338 7
0.3385 6
0.339 6
0.3395 6
0.34 5
0.3405 6
0.341 7
0.3415 6
0.342 7
0.3425 7
0.343 6
0.3435 6
0.344 6
0.3445 7
0.345 6
0.3455 6
0.346 6
0.3465 6
0.347 6
0.3475 6
0.348 7
0.3485 8
0.349 5
0.3495 6
0.35 6
0.3505 6
0.351 6
0.3515 6
0.352 6
0.3525 6
0.353 6
0.3535 7
0.354 7
0.3545 6
0.355 6
0.3555 6
0.356 6
0.3565 6
0.357 5
0.3575 7
0.358 6
0.3585 6
0.359 6
0.3595 7
0.36 8
0.3605 7
0.361 7
0.3615 8
0.362 7
0.3625 6
0.363 6
0.3635 6
0.364 6
0.3645 7
0.365 6
0.3655 6
0.366 7
0.3665 7
0.367 7
0.3675 8
0.368 6
0.3685 7
0.369 7
0.3695 7
0.37 6
0.3705 7
0.371 7
0.3715 6
0.372 5
0.3725 7
0.373 8
0.3735 7
0.374 6
0.3745 6
0.375 7
0.3755 7
0.376 7
0.3765 6
0.377 7
0.3775 9
0.378 7
0.3785 7
0.379 8
0.3795 6
0.38 7
0.3805 7
0.381 7
0.3815 7
0.382 7
0.3825 8
0.383 7
0.3835 7
0.384 5
0.3845 6
0.385 8
0.3855 8
0.386 7
0.3865 7
0.387 6
0.3875 8
0.388 6
0.3885 6
0.389 9
0.3895 7
0.39 8
0.3905 7
0.391 6
0.3915 5
0.392 7
0.3925 7
0.393 7
0.3935 7
0.394 7
0.3945 7
0.395 7
0.3955 6
0.396 7
0.3965 8
0.397 7
0.3975 8
0.398 7
0.3985 7
0.399 7
0.3995 7
0.4 7
0.4005 6
0.401 7
0.4015 8
0.402 7
0.4025 6
0.403 7
0.4035 7
0.404 7
0.4045 8
0.405 7
0.4055 7
0.406 7
0.4065 8
0.407 5
0.4075 7
0.408 8
0.4085 6
0.409 7
0.4095 7
0.41 7
0.4105 8
0.411 6
0.4115 7
0.412 8
0.4125 7
0.413 7
0.4135 6
0.414 6
0.4145 8
0.415 7
0.4155 7
0.416 8
0.4165 8
0.417 7
0.4175 7
0.418 7
0.4185 6
0.419 6
0.4195 8
0.42 6
0.4205 5
0.421 5
0.4215 7
0.422 6
0.4225 4
0.423 5
0.4235 7
0.424 5
0.4245 6
0.425 6
0.4255 6
0.426 6
0.4265 6
0.427 6
0.4275 7
0.428 8
0.4285 7
0.429 6
0.4295 6
0.43 6
0.4305 6
0.431 8
0.4315 6
0.432 6
0.4325 7
0.433 6
0.4335 7
0.434 6
0.4345 6
0.435 7
0.4355 8
0.436 7
0.4365 7
0.437 7
0.4375 6
0.438 7
0.4385 6
0.439 7
0.4395 6
0.44 6
0.4405 7
0.441 7
0.4415 6
0.442 7
0.4425 7
0.443 5
0.4435 5
0.444 6
0.4445 6
0.445 6
0.4455 6
0.446 5
0.4465 6
0.447 5
0.4475 6
0.448 7
0.4485 6
0.449 7
0.4495 6
0.45 6
0.4505 4
0.451 5
0.4515 7
0.452 6
0.4525 6
0.453 6
0.4535 8
0.454 6
0.4545 7
0.455 6
0.4555 7
0.456 7
0.4565 5
0.457 5
0.4575 6
0.458 7
0.4585 7
0.459 4
0.4595 5
0.46 6
0.4605 5
0.461 5
0.4615 7
0.462 6
0.4625 6
0.463 4
0.4635 5
0.464 4
0.4645 6
0.465 4
0.4655 6
0.466 6
0.4665 5
0.467 4
0.4675 5
0.468 4
0.4685 4
0.469 7
0.4695 7
0.47 6
0.4705 6
0.471 4
0.4715 4
0.472 5
0.4725 5
0.473 5
0.4735 5
0.474 4
0.4745 4
0.475 6
0.4755 5
0.476 5
0.4765 3
0.477 4
0.4775 6
0.478 4
0.4785 5
0.479 6
0.4795 7
0.48 6
0.4805 5
0.481 4
0.4815 6
0.482 7
0.4825 4
0.483 5
0.4835 4
0.484 4
0.4845 4
0.485 3
0.4855 4
0.486 4
0.4865 4
0.487 4
0.4875 4
0.488 4
0.4885 5
0.489 4
0.4895 3
0.49 3
0.4905 4
0.491 3
0.4915 4
0.492 3
0.4925 3
0.493 3
0.4935 4
0.494 4
0.4945 2
0.495 3
0.4955 4
0.496 3
0.4965 3
0.497 3
0.4975 3
0.498 6
0.4985 4
0.499 3
0.4995 3
0.5 3
0.5005 4
0.501 3
0.5015 3
0.502 3
0.5025 2
0.503 3
0.5035 6
0.504 3
0.5045 3
0.505 3
0.5055 2
0.506 3
0.5065 3
0.507 3
0.5075 3
0.508 5
0.5085 2
0.509 3
0.5095 2
0.51 3
0.5105 2
0.511 2
0.5115 2
0.512 2
0.5125 3
0.513 2
0.5135 3
0.514 2
0.5145 3
0.515 2
0.5155 2
0.516 3
0.5165 3
0.517 3
0.5175 3
0.518 2
0.5185 3
0.519 3
0.5195 3
0.52 3
0.5205 3
0.521 3
0.5215 3
0.522 3
0.5225 2
0.523 2
0.5235 4
0.524 3
0.5245 2
0.525 2
0.5255 3
0.526 2
0.5265 2
0.527 2
0.5275 3
0.528 2
0.5285 4
0.529 2
0.5295 2
0.53 2
0.5305 2
0.531 2
0.5315 2
0.532 2
0.5325 3
0.533 3
0.5335 3
0.534 2
0.5345 2
0.535 3
0.5355 3
0.536 3
0.5365 5
0.537 4
0.5375 3
0.538 3
0.5385 2
0.539 3
0.5395 3
0.54 2
0.5405 3
0.541 2
0.5415 3
0.542 3
0.5425 3
0.543 3
0.5435 3
0.544 2
0.5445 3
0.545 2
0.5455 4
0.546 3
0.5465 3
0.547 2
0.5475 2
0.548 2
0.5485 4
0.549 3
0.5495 4
0.55 3
0.5505 2
0.551 4
0.5515 3
0.552 4
0.5525 3
0.553 3
0.5535 3
0.554 4
0.5545 3
0.555 2
0.5555 3
0.556 3
0.5565 3
0.557 3
0.5575 4
0.558 4
0.5585 4
0.559 4
0.5595 4
0.56 3
0.5605 3
0.561 3
0.5615 4
0.562 4
0.5625 4
0.563 2
0.5635 4
0.564 2
0.5645 3
0.565 4
0.5655 4
0.566 3
0.5665 3
0.567 3
0.5675 2
0.568 2
0.5685 2
0.569 4
0.5695 4
0.57 3
0.5705 3
0.571 4
0.5715 3
0.572 3
0.5725 3
0.573 3
0.5735 2
0.574 3
0.5745 3
0.575 3
0.5755 3
0.576 2
0.5765 3
0.577 4
0.5775 3
0.578 3
0.5785 1
0.579 1
0.5795 3
0.58 2
0.5805 4
0.581 4
0.5815 2
0.582 2
0.5825 2
0.583 2
0.5835 1
0.584 1
0.5845 3
0.585 3
0.5855 2
0.586 2
0.5865 1
0.587 3
0.5875 3
0.588 3
0.5885 4
0.589 3
0.5895 3
0.59 2
0.5905 2
0.591 3
0.5915 2
0.592 4
0.5925 3
0.593 4
0.5935 3
0.594 3
0.5945 2
0.595 3
0.5955 3
0.596 2
0.5965 3
0.597 3
0.5975 4
0.598 2
0.5985 2
0.599 3
0.5995 3
0.6 3
0.6005 4
0.601 3
0.6015 2
0.602 2
0.6025 2
0.603 2
0.6035 3
0.604 3
0.6045 3
0.605 3
0.6055 3
0.606 3
0.6065 1
0.607 2
0.6075 2
0.608 2
0.6085 2
0.609 3
0.6095 4
0.61 1
0.6105 1
0.611 2
0.6115 2
0.612 2
0.6125 3
0.613 1
0.6135 1
0.614 2
0.6145 1
0.615 1
0.6155 1
0.616 1
0.6165 2
0.617 2
0.6175 2
0.618 2
0.6185 1
0.619 1
0.6195 1
0.62 1
0.6205 2
0.621 2
0.6215 1
0.622 2
0.6225 1
0.623 2
0.6235 2
0.624 3
0.6245 2
0.625 2
0.6255 1
0.626 1
0.6265 1
0.627 1
0.6275 3
0.628 4
0.6285 2
0.629 1
0.6295 3
0.63 2
0.6305 1
0.631 2
0.6315 2
0.632 2
0.6325 3
0.633 2
0.6335 2
0.634 3
0.6345 4
0.635 1
0.6355 3
0.636 2
0.6365 1
0.637 2
0.6375 2
0.638 2
0.6385 2
0.639 2
0.6395 1
0.64 1
0.6405 1
0.641 2
0.6415 1
0.642 2
0.6425 2
0.643 2
0.6435 2
0.644 1
0.6445 1
0.645 1
0.6455 2
0.646 1
0.6465 2
0.647 2
0.6475 2
0.648 1
0.6485 1
0.649 1
0.6495 2
0.65 2
0.6505 2
0.651 2
0.6515 2
0.652 1
0.6525 2
0.653 2
0.6535 2
0.654 2
0.6545 3
0.655 1
0.6555 1
0.656 1
0.6565 2
0.657 1
0.6575 1
0.658 2
0.6585 2
0.659 1
0.6595 2
0.66 2
0.6605 2
0.661 2
0.6615 2
0.662 2
0.6625 2
0.663 2
0.6635 3
0.664 3
0.6645 4
0.665 4
0.6655 4
0.666 3
0.6665 4
0.667 2
0.6675 3
0.668 3
0.6685 2
0.669 3
0.6695 2
0.67 3
0.6705 3
0.671 3
0.6715 2
0.672 3
0.6725 3
0.673 1
0.6735 1
0.674 2
0.6745 2
0.675 2
0.6755 2
0.676 2
0.6765 2
0.677 2
0.6775 3
0.678 3
0.6785 3
0.679 3
0.6795 2
0.68 2
0.6805 1
0.681 3
0.6815 3
0.682 2
0.6825 2
0.683 3
0.6835 2
0.684 2
0.6845 2
0.685 2
0.6855 2
0.686 3
0.6865 3
0.687 2
0.6875 2
0.688 3
0.6885 2
0.689 4
0.6895 2
0.69 2
0.6905 4
0.691 3
0.6915 4
0.692 3
0.6925 4
0.693 3
0.6935 2
0.694 1
0.6945 3
0.695 3
0.6955 3
0.696 3
0.6965 2
0.697 2
0.6975 3
0.698 2
0.6985 3
0.699 2
0.6995 4
0.7 2
0.7005 3
0.701 2
0.7015 2
0.702 3
0.7025 3
0.703 3
0.7035 2
0.704 3
0.7045 3
0.705 2
0.7055 2
0.706 4
0.7065 2
0.707 2
0.7075 4
0.708 2
0.7085 2
0.709 2
0.7095 2
0.71 2
0.7105 2
0.711 2
0.7115 3
0.712 2
0.7125 3
0.713 2
0.7135 2
0.714 3
0.7145 4
0.715 4
0.7155 4
0.716 4
0.7165 2
0.717 3
0.7175 3
0.718 3
0.7185 3
0.719 3
0.7195 3
0.72 3
0.7205 2
0.721 2
0.7215 3
0.722 3
0.7225 4
0.723 4
0.7235 2
0.724 3
0.7245 3
0.725 3
0.7255 5
0.726 3
0.7265 2
0.727 3
0.7275 2
0.728 3
0.7285 3
0.729 3
0.7295 2
0.73 4
0.7305 4
0.731 3
0.7315 3
0.732 2
0.7325 3
0.733 3
0.7335 3
0.734 4
0.7345 4
0.735 4
0.7355 3
0.736 4
0.7365 3
0.737 3
0.7375 3
0.738 2
0.7385 2
0.739 3
0.7395 4
0.74 4
0.7405 4
0.741 4
0.7415 4
0.742 4
0.7425 4
0.743 4
0.7435 5
0.744 3
0.7445 5
0.745 4
0.7455 3
0.746 3
0.7465 3
0.747 2
0.7475 3
0.748 4
0.7485 3
0.749 4
0.7495 3
0.75 3
0.7505 4
0.751 4
0.7515 4
0.752 4
0.7525 2
0.753 2
0.7535 2
0.754 3
0.7545 3
0.755 3
0.7555 3
0.756 4
0.7565 3
0.757 3
0.7575 4
0.758 4
0.7585 3
0.759 3
0.7595 4
0.76 3
0.7605 2
0.761 2
0.7615 4
0.762 3
0.7625 3
0.763 3
0.7635 3
0.764 3
0.7645 4
0.765 4
0.7655 4
0.766 3
0.7665 3
0.767 2
0.7675 3
0.768 3
0.7685 4
0.769 5
0.7695 3
0.77 4
0.7705 4
0.771 3
0.7715 4
0.772 3
0.7725 4
0.773 5
0.7735 3
0.774 4
0.7745 3
0.775 2
0.7755 2
0.776 5
0.7765 4
0.777 4
0.7775 4
0.778 4
0.7785 4
0.779 3
0.7795 3
0.78 3
0.7805 4
0.781 4
0.7815 4
0.782 4
0.7825 4
0.783 5
0.7835 3
0.784 3
0.7845 5
0.785 4
0.7855 2
0.786 4
0.7865 5
0.787 2
0.7875 3
0.788 4
0.7885 4
0.789 5
0.7895 5
0.79 3
0.7905 3
0.791 3
0.7915 3
0.792 4
0.7925 4
0.793 4
0.7935 3
0.794 3
0.7945 3
0.795 4
0.7955 4
0.796 4
0.7965 3
0.797 4
0.7975 3
0.798 5
0.7985 4
0.799 5
0.7995 2
0.8 1
0.8005 2
0.801 3
0.8015 3
0.802 3
0.8025 5
0.803 4
0.8035 3
0.804 3
0.8045 3
0.805 2
0.8055 2
0.806 4
0.8065 3
0.807 5
0.8075 4
0.808 3
0.8085 4
0.809 3
0.8095 6
0.81 4
0.8105 4
0.811 4
0.8115 4
0.812 5
0.8125 5
0.813 4
0.8135 4
0.814 5
0.8145 5
0.815 5
0.8155 4
0.816 4
0.8165 3
0.817 3
0.8175 4
0.818 3
0.8185 3
0.819 4
0.8195 2
0.82 3
0.8205 3
0.821 4
0.8215 5
0.822 6
0.8225 4
0.823 5
0.8235 4
0.824 4
0.8245 3
0.825 4
0.8255 5
0.826 3
0.8265 2
0.827 4
0.8275 3
0.828 3
0.8285 5
0.829 4
0.8295 4
0.83 4
0.8305 4
0.831 5
0.8315 4
0.832 5
0.8325 4
0.833 5
0.8335 5
0.834 5
0.8345 3
0.835 4
0.8355 5
0.836 5
0.8365 4
0.837 4
0.8375 4
0.838 5
0.8385 5
0.839 6
0.8395 3
0.84 3
0.8405 2
0.841 4
0.8415 4
0.842 4
0.8425 4
0.843 4
0.8435 4
0.844 3
0.8445 4
0.845 3
0.8455 4
0.846 3
0.8465 3
0.847 4
0.8475 5
0.848 4
0.8485 4
0.849 3
0.8495 4
0.85 4
0.8505 4
0.851 5
0.8515 5
0.852 5
0.8525 4
0.853 4
0.8535 4
0.854 3
0.8545 4
0.855 4
0.8555 4
0.856 4
0.8565 5
0.857 6
0.8575 4
0.858 4
0.8585 6
0.859 4
0.8595 3
0.86 4
0.8605 4
0.861 4
0.8615 6
0.862 4
0.8625 4
0.863 3
0.8635 3
0.864 3
0.8645 4
0.865 5
0.8655 3
0.866 3
0.8665 4
0.867 4
0.8675 3
0.868 4
0.8685 4
0.869 3
0.8695 4
0.87 3
0.8705 4
0.871 2
0.8715 3
0.872 3
0.8725 4
0.873 4
0.8735 2
0.874 3
0.8745 4
0.875 3
0.8755 4
0.876 4
0.8765 4
0.877 4
0.8775 4
0.878 5
0.8785 2
0.879 2
0.8795 6
0.88 5
0.8805 3
0.881 4
0.8815 4
0.882 3
0.8825 3
0.883 4
0.8835 4
0.884 4
0.8845 3
0.885 4
0.8855 5
0.886 5
0.8865 2
0.887 2
0.8875 2
0.888 1
0.8885 3
0.889 1
0.8895 4
0.89 4
0.8905 4
0.891 3
0.8915 3
0.892 5
0.8925 4
0.893 3
0.8935 2
0.894 2
0.8945 2
0.895 3
0.8955 4
0.896 2
0.8965 2
0.897 3
0.8975 3
0.898 2
0.8985 4
0.899 4
0.8995 3
0.9 3
0.9005 5
0.901 5
0.9015 3
0.902 3
0.9025 2
0.903 4
0.9035 4
0.904 3
0.9045 3
0.905 2
0.9055 3
0.906 4
0.9065 2
0.907 3
0.9075 3
0.908 4
0.9085 4
0.909 3
0.9095 2
0.91 3
0.9105 3
0.911 2
0.9115 4
0.912 3
0.9125 4
0.913 5
0.9135 4
0.914 7
0.9145 4
0.915 4
0.9155 4
0.916 3
0.9165 3
0.917 4
0.9175 2
0.918 3
0.9185 4
0.919 3
0.9195 4
0.92 4
0.9205 4
0.921 1
0.9215 4
0.922 4
0.9225 5
0.923 5
0.9235 4
0.924 4
0.9245 2
0.925 3
0.9255 4
0.926 3
0.9265 2
0.927 7
0.9275 3
0.928 3
0.9285 3
0.929 4
0.9295 4
0.93 5
0.9305 3
0.931 4
0.9315 3
0.932 3
0.9325 3
0.933 2
0.9335 4
0.934 5
0.9345 2
0.935 3
0.9355 3
0.936 6
0.9365 2
0.937 4
0.9375 5
0.938 3
0.9385 2
0.939 4
0.9395 2
0.94 4
0.9405 5
0.941 3
0.9415 4
0.942 4
0.9425 2
0.943 3
0.9435 2
0.944 2
0.9445 2
0.945 3
0.9455 5
0.946 4
0.9465 2
0.947 3
0.9475 3
0.948 3
0.9485 2
0.949 3
0.9495 4
0.95 4
0.9505 4
0.951 4
0.9515 2
0.952 2
0.9525 1
0.953 2
0.9535 3
0.954 3
0.9545 4
0.955 2
0.9555 2
0.956 18
0.9565 4
0.957 5
0.9575 3
0.958 3
0.9585 3
0.959 2
0.9595 4
0.96 4
0.9605 4
0.961 4
0.9615 3
0.962 5
0.9625 6
0.963 6
0.9635 2
0.964 4
0.9645 4
0.965 3
0.9655 3
0.966 3
0.9665 3
0.967 3
0.9675 3
0.968 3
0.9685 3
0.969 3
0.9695 4
0.97 5
0.9705 2
0.971 3
0.9715 3
0.972 4
0.9725 4
0.973 3
0.9735 3
0.974 4
0.9745 5
0.975 3
0.9755 4
0.976 3
0.9765 5
0.977 2
0.9775 4
0.978 3
0.9785 5
0.979 5
0.9795 4
0.98 3
0.9805 21
0.981 5
0.9815 4
0.982 4
0.9825 4
0.983 3
0.9835 2
0.984 3
0.9845 3
0.985 4
0.9855 5
0.986 3
0.9865 2
0.987 4
0.9875 3
0.988 4
0.9885 3
0.989 3
0.9895 2
0.99 2
0.9905 3
0.991 2
0.9915 2
0.992 3
0.9925 3
0.993 2
0.9935 3
0.994 4
0.9945 5
0.995 3
0.9955 5
0.996 3
0.9965 4
0.997 4
0.9975 3
0.998 3
0.9985 2
0.999 3
0.9995 2
};
%\addlegendentry{$\epsilon_u = 10^{-11}$}
\end{axis}

\end{tikzpicture}

%% file: figures/mc_Burgers.tex
\definecolor{mycolor1}{rgb}{0.00000,0.44700,0.74100}%
\definecolor{mycolor2}{rgb}{0.85000,0.32500,0.09800}%
\definecolor{mycolor3}{rgb}{0.5,0.5,0.5}
\begin{tikzpicture}[scale=0.8]

\begin{axis}[%
width=1.5in,
height=1.5in,
at={(1.011in,0.642in)},
scale only axis,
xmin=0,
xmax=5,
xlabel style={font=\color{white!15!black}},
xlabel={Time},
%ymode=log,
ymin=0,
ymax=60,
yminorticks=true,
ylabel style={font=\color{white!15!black}},
ylabel={$m_c$},
axis background/.style={fill=white},
legend style={fill opacity=0.8, 
font=\tiny,
  draw opacity=1,
  text opacity=1,
  at={(0.02,0.87)},
  anchor=west}
  % draw=lightgray204}
]
\addplot [color=mycolor1, line width=1.2pt]
  table[row sep=crcr]{%
0	1\\
0.0025	0\\
0.005	0\\
0.0075	0\\
0.01	0\\
0.0125	0\\
0.015	0\\
0.0175	0\\
0.02	0\\
0.0225	0\\
0.025	0\\
0.0275	0\\
0.03	0\\
0.0325	0\\
0.035	0\\
0.0375	0\\
0.04	0\\
0.0425	0\\
0.045	0\\
0.0475	0\\
0.05	0\\
0.0525	0\\
0.055	0\\
0.0575	0\\
0.06	0\\
0.0625	0\\
0.065	0\\
0.0675	0\\
0.07	0\\
0.0725	0\\
0.075	0\\
0.0775	0\\
0.08	0\\
0.0825	0\\
0.085	0\\
0.0875	0\\
0.09	0\\
0.0925	0\\
0.095	0\\
0.0975	0\\
0.1	0\\
0.1025	0\\
0.105	0\\
0.1075	0\\
0.11	0\\
0.1125	0\\
0.115	0\\
0.1175	0\\
0.12	0\\
0.1225	0\\
0.125	0\\
0.1275	0\\
0.13	0\\
0.1325	0\\
0.135	0\\
0.1375	0\\
0.14	0\\
0.1425	0\\
0.145	0\\
0.1475	0\\
0.15	0\\
0.1525	0\\
0.155	0\\
0.1575	0\\
0.16	0\\
0.1625	0\\
0.165	0\\
0.1675	0\\
0.17	0\\
0.1725	0\\
0.175	0\\
0.1775	0\\
0.18	0\\
0.1825	0\\
0.185	0\\
0.1875	0\\
0.19	0\\
0.1925	0\\
0.195	0\\
0.1975	0\\
0.2	0\\
0.2025	0\\
0.205	1\\
0.2075	1\\
0.21	1\\
0.2125	1\\
0.215	1\\
0.2175	1\\
0.22	0\\
0.2225	0\\
0.225	0\\
0.2275	0\\
0.23	0\\
0.2325	0\\
0.235	0\\
0.2375	0\\
0.24	0\\
0.2425	0\\
0.245	0\\
0.2475	0\\
0.25	0\\
0.2525	0\\
0.255	0\\
0.2575	0\\
0.26	1\\
0.2625	1\\
0.265	1\\
0.2675	1\\
0.27	1\\
0.2725	1\\
0.275	0\\
0.2775	0\\
0.28	9\\
0.2825	1\\
0.285	1\\
0.2875	1\\
0.29	1\\
0.2925	1\\
0.295	1\\
0.2975	1\\
0.3	1\\
0.3025	1\\
0.305	1\\
0.3075	1\\
0.31	1\\
0.3125	1\\
0.315	1\\
0.3175	1\\
0.32	1\\
0.3225	1\\
0.325	1\\
0.3275	1\\
0.33	1\\
0.3325	1\\
0.335	1\\
0.3375	1\\
0.34	1\\
0.3425	1\\
0.345	1\\
0.3475	1\\
0.35	1\\
0.3525	1\\
0.355	1\\
0.3575	1\\
0.36	1\\
0.3625	1\\
0.365	1\\
0.3675	1\\
0.37	1\\
0.3725	1\\
0.375	1\\
0.3775	1\\
0.38	1\\
0.3825	1\\
0.385	1\\
0.3875	1\\
0.39	0\\
0.3925	0\\
0.395	0\\
0.3975	0\\
0.4	0\\
0.4025	0\\
0.405	0\\
0.4075	0\\
0.41	0\\
0.4125	0\\
0.415	0\\
0.4175	0\\
0.42	0\\
0.4225	0\\
0.425	0\\
0.4275	0\\
0.43	0\\
0.4325	0\\
0.435	0\\
0.4375	0\\
0.44	0\\
0.4425	0\\
0.445	0\\
0.4475	1\\
0.45	1\\
0.4525	1\\
0.455	1\\
0.4575	1\\
0.46	1\\
0.4625	1\\
0.465	1\\
0.4675	1\\
0.47	1\\
0.4725	1\\
0.475	1\\
0.4775	1\\
0.48	1\\
0.4825	1\\
0.485	1\\
0.4875	1\\
0.49	1\\
0.4925	1\\
0.495	1\\
0.4975	1\\
0.5	1\\
0.5025	1\\
0.505	1\\
0.5075	1\\
0.51	1\\
0.5125	1\\
0.515	0\\
0.5175	0\\
0.52	0\\
0.5225	0\\
0.525	0\\
0.5275	1\\
0.53	1\\
0.5325	1\\
0.535	1\\
0.5375	1\\
0.54	1\\
0.5425	1\\
0.545	1\\
0.5475	1\\
0.55	1\\
0.5525	1\\
0.555	0\\
0.5575	0\\
0.56	0\\
0.5625	0\\
0.565	1\\
0.5675	1\\
0.57	0\\
0.5725	1\\
0.575	8\\
0.5775	1\\
0.58	1\\
0.5825	1\\
0.585	1\\
0.5875	1\\
0.59	1\\
0.5925	1\\
0.595	1\\
0.5975	1\\
0.6	1\\
0.6025	1\\
0.605	1\\
0.6075	1\\
0.61	1\\
0.6125	1\\
0.615	1\\
0.6175	1\\
0.62	1\\
0.6225	1\\
0.625	1\\
0.6275	1\\
0.63	1\\
0.6325	1\\
0.635	1\\
0.6375	1\\
0.64	1\\
0.6425	1\\
0.645	1\\
0.6475	1\\
0.65	1\\
0.6525	1\\
0.655	1\\
0.6575	1\\
0.66	1\\
0.6625	1\\
0.665	1\\
0.6675	1\\
0.67	1\\
0.6725	1\\
0.675	1\\
0.6775	1\\
0.68	1\\
0.6825	1\\
0.685	1\\
0.6875	1\\
0.69	1\\
0.6925	1\\
0.695	1\\
0.6975	1\\
0.7	1\\
0.7025	1\\
0.705	1\\
0.7075	1\\
0.71	1\\
0.7125	1\\
0.715	1\\
0.7175	1\\
0.72	1\\
0.7225	1\\
0.725	1\\
0.7275	1\\
0.73	1\\
0.7325	1\\
0.735	1\\
0.7375	1\\
0.74	1\\
0.7425	1\\
0.745	1\\
0.7475	1\\
0.75	1\\
0.7525	1\\
0.755	1\\
0.7575	1\\
0.76	1\\
0.7625	1\\
0.765	1\\
0.7675	1\\
0.77	1\\
0.7725	1\\
0.775	1\\
0.7775	1\\
0.78	1\\
0.7825	1\\
0.785	1\\
0.7875	0\\
0.79	0\\
0.7925	0\\
0.795	0\\
0.7975	0\\
0.8	0\\
0.8025	0\\
0.805	0\\
0.8075	0\\
0.81	0\\
0.8125	0\\
0.815	0\\
0.8175	0\\
0.82	0\\
0.8225	0\\
0.825	0\\
0.8275	0\\
0.83	0\\
0.8325	0\\
0.835	0\\
0.8375	1\\
0.84	1\\
0.8425	0\\
0.845	0\\
0.8475	0\\
0.85	1\\
0.8525	1\\
0.855	0\\
0.8575	1\\
0.86	1\\
0.8625	1\\
0.865	1\\
0.8675	1\\
0.87	1\\
0.8725	1\\
0.875	1\\
0.8775	1\\
0.88	0\\
0.8825	0\\
0.885	11\\
0.8875	0\\
0.89	0\\
0.8925	0\\
0.895	0\\
0.8975	0\\
0.9	0\\
0.9025	0\\
0.905	0\\
0.9075	0\\
0.91	0\\
0.9125	0\\
0.915	0\\
0.9175	0\\
0.92	0\\
0.9225	0\\
0.925	0\\
0.9275	0\\
0.93	0\\
0.9325	0\\
0.935	0\\
0.9375	0\\
0.94	0\\
0.9425	0\\
0.945	0\\
0.9475	0\\
0.95	0\\
0.9525	0\\
0.955	0\\
0.9575	0\\
0.96	0\\
0.9625	0\\
0.965	0\\
0.9675	0\\
0.97	0\\
0.9725	0\\
0.975	0\\
0.9775	0\\
0.98	0\\
0.9825	1\\
0.985	1\\
0.9875	1\\
0.99	1\\
0.9925	1\\
0.995	1\\
0.9975	1\\
1	1\\
1.0025	1\\
1.005	1\\
1.0075	1\\
1.01	1\\
1.0125	1\\
1.015	1\\
1.0175	1\\
1.02	1\\
1.0225	1\\
1.025	1\\
1.0275	1\\
1.03	1\\
1.0325	1\\
1.035	0\\
1.0375	0\\
1.04	0\\
1.0425	0\\
1.045	0\\
1.0475	0\\
1.05	0\\
1.0525	0\\
1.055	0\\
1.0575	0\\
1.06	0\\
1.0625	0\\
1.065	0\\
1.0675	0\\
1.07	0\\
1.0725	1\\
1.075	1\\
1.0775	1\\
1.08	1\\
1.0825	1\\
1.085	1\\
1.0875	1\\
1.09	1\\
1.0925	1\\
1.095	1\\
1.0975	1\\
1.1	1\\
1.1025	1\\
1.105	1\\
1.1075	1\\
1.11	1\\
1.1125	1\\
1.115	1\\
1.1175	1\\
1.12	1\\
1.1225	1\\
1.125	1\\
1.1275	1\\
1.13	1\\
1.1325	1\\
1.135	0\\
1.1375	0\\
1.14	0\\
1.1425	0\\
1.145	0\\
1.1475	0\\
1.15	0\\
1.1525	0\\
1.155	0\\
1.1575	0\\
1.16	0\\
1.1625	0\\
1.165	0\\
1.1675	0\\
1.17	0\\
1.1725	0\\
1.175	0\\
1.1775	0\\
1.18	0\\
1.1825	0\\
1.185	0\\
1.1875	0\\
1.19	0\\
1.1925	0\\
1.195	0\\
1.1975	0\\
1.2	0\\
1.2025	0\\
1.205	0\\
1.2075	0\\
1.21	10\\
1.2125	0\\
1.215	0\\
1.2175	0\\
1.22	0\\
1.2225	0\\
1.225	0\\
1.2275	0\\
1.23	0\\
1.2325	0\\
1.235	1\\
1.2375	0\\
1.24	0\\
1.2425	0\\
1.245	1\\
1.2475	1\\
1.25	1\\
1.2525	0\\
1.255	0\\
1.2575	0\\
1.26	0\\
1.2625	1\\
1.265	1\\
1.2675	1\\
1.27	1\\
1.2725	1\\
1.275	1\\
1.2775	1\\
1.28	1\\
1.2825	1\\
1.285	0\\
1.2875	0\\
1.29	0\\
1.2925	0\\
1.295	0\\
1.2975	0\\
1.3	0\\
1.3025	0\\
1.305	0\\
1.3075	1\\
1.31	1\\
1.3125	1\\
1.315	1\\
1.3175	1\\
1.32	1\\
1.3225	1\\
1.325	1\\
1.3275	1\\
1.33	1\\
1.3325	1\\
1.335	1\\
1.3375	1\\
1.34	0\\
1.3425	1\\
1.345	1\\
1.3475	1\\
1.35	0\\
1.3525	0\\
1.355	1\\
1.3575	1\\
1.36	1\\
1.3625	1\\
1.365	1\\
1.3675	1\\
1.37	1\\
1.3725	1\\
1.375	1\\
1.3775	1\\
1.38	1\\
1.3825	1\\
1.385	0\\
1.3875	0\\
1.39	0\\
1.3925	0\\
1.395	0\\
1.3975	0\\
1.4	0\\
1.4025	0\\
1.405	0\\
1.4075	1\\
1.41	1\\
1.4125	1\\
1.415	1\\
1.4175	1\\
1.42	1\\
1.4225	1\\
1.425	1\\
1.4275	1\\
1.43	1\\
1.4325	1\\
1.435	1\\
1.4375	1\\
1.44	1\\
1.4425	1\\
1.445	1\\
1.4475	1\\
1.45	1\\
1.4525	1\\
1.455	1\\
1.4575	1\\
1.46	1\\
1.4625	1\\
1.465	1\\
1.4675	1\\
1.47	1\\
1.4725	1\\
1.475	1\\
1.4775	1\\
1.48	0\\
1.4825	0\\
1.485	0\\
1.4875	0\\
1.49	0\\
1.4925	0\\
1.495	0\\
1.4975	0\\
1.5	0\\
1.5025	0\\
1.505	0\\
1.5075	0\\
1.51	0\\
1.5125	0\\
1.515	0\\
1.5175	0\\
1.52	0\\
1.5225	0\\
1.525	0\\
1.5275	0\\
1.53	0\\
1.5325	0\\
1.535	0\\
1.5375	0\\
1.54	0\\
1.5425	0\\
1.545	0\\
1.5475	0\\
1.55	0\\
1.5525	1\\
1.555	1\\
1.5575	1\\
1.56	1\\
1.5625	1\\
1.565	0\\
1.5675	0\\
1.57	0\\
1.5725	0\\
1.575	0\\
1.5775	0\\
1.58	0\\
1.5825	0\\
1.585	0\\
1.5875	0\\
1.59	0\\
1.5925	0\\
1.595	0\\
1.5975	0\\
1.6	0\\
1.6025	0\\
1.605	0\\
1.6075	0\\
1.61	0\\
1.6125	0\\
1.615	0\\
1.6175	0\\
1.62	0\\
1.6225	0\\
1.625	0\\
1.6275	0\\
1.63	0\\
1.6325	0\\
1.635	0\\
1.6375	0\\
1.64	0\\
1.6425	0\\
1.645	0\\
1.6475	0\\
1.65	0\\
1.6525	0\\
1.655	0\\
1.6575	0\\
1.66	0\\
1.6625	0\\
1.665	0\\
1.6675	0\\
1.67	0\\
1.6725	0\\
1.675	0\\
1.6775	0\\
1.68	0\\
1.6825	0\\
1.685	0\\
1.6875	0\\
1.69	0\\
1.6925	1\\
1.695	1\\
1.6975	1\\
1.7	1\\
1.7025	1\\
1.705	1\\
1.7075	1\\
1.71	1\\
1.7125	1\\
1.715	1\\
1.7175	1\\
1.72	1\\
1.7225	1\\
1.725	1\\
1.7275	1\\
1.73	1\\
1.7325	1\\
1.735	1\\
1.7375	1\\
1.74	1\\
1.7425	1\\
1.745	1\\
1.7475	1\\
1.75	1\\
1.7525	1\\
1.755	1\\
1.7575	1\\
1.76	1\\
1.7625	1\\
1.765	1\\
1.7675	1\\
1.77	1\\
1.7725	1\\
1.775	4\\
1.7775	20\\
1.78	1\\
1.7825	1\\
1.785	1\\
1.7875	1\\
1.79	1\\
1.7925	1\\
1.795	1\\
1.7975	1\\
1.8	1\\
1.8025	1\\
1.805	1\\
1.8075	1\\
1.81	1\\
1.8125	1\\
1.815	1\\
1.8175	0\\
1.82	0\\
1.8225	0\\
1.825	0\\
1.8275	0\\
1.83	0\\
1.8325	0\\
1.835	1\\
1.8375	1\\
1.84	1\\
1.8425	1\\
1.845	1\\
1.8475	0\\
1.85	0\\
1.8525	0\\
1.855	0\\
1.8575	0\\
1.86	0\\
1.8625	0\\
1.865	0\\
1.8675	0\\
1.87	0\\
1.8725	0\\
1.875	0\\
1.8775	1\\
1.88	1\\
1.8825	1\\
1.885	1\\
1.8875	1\\
1.89	1\\
1.8925	1\\
1.895	1\\
1.8975	1\\
1.9	1\\
1.9025	1\\
1.905	1\\
1.9075	1\\
1.91	1\\
1.9125	1\\
1.915	1\\
1.9175	1\\
1.92	1\\
1.9225	1\\
1.925	1\\
1.9275	1\\
1.93	1\\
1.9325	1\\
1.935	1\\
1.9375	1\\
1.94	1\\
1.9425	1\\
1.945	1\\
1.9475	1\\
1.95	1\\
1.9525	1\\
1.955	1\\
1.9575	1\\
1.96	1\\
1.9625	1\\
1.965	1\\
1.9675	1\\
1.97	1\\
1.9725	1\\
1.975	1\\
1.9775	1\\
1.98	1\\
1.9825	1\\
1.985	1\\
1.9875	1\\
1.99	1\\
1.9925	1\\
1.995	1\\
1.9975	1\\
2	1\\
2.0025	1\\
2.005	1\\
2.0075	1\\
2.01	1\\
2.0125	1\\
2.015	1\\
2.0175	1\\
2.02	1\\
2.0225	1\\
2.025	1\\
2.0275	1\\
2.03	1\\
2.0325	1\\
2.035	1\\
2.0375	1\\
2.04	1\\
2.0425	1\\
2.045	1\\
2.0475	1\\
2.05	1\\
2.0525	1\\
2.055	1\\
2.0575	1\\
2.06	1\\
2.0625	1\\
2.065	1\\
2.0675	1\\
2.07	1\\
2.0725	1\\
2.075	1\\
2.0775	1\\
2.08	1\\
2.0825	1\\
2.085	1\\
2.0875	1\\
2.09	1\\
2.0925	1\\
2.095	1\\
2.0975	1\\
2.1	1\\
2.1025	1\\
2.105	1\\
2.1075	1\\
2.11	1\\
2.1125	1\\
2.115	1\\
2.1175	1\\
2.12	1\\
2.1225	1\\
2.125	1\\
2.1275	1\\
2.13	1\\
2.1325	1\\
2.135	1\\
2.1375	1\\
2.14	1\\
2.1425	1\\
2.145	1\\
2.1475	1\\
2.15	1\\
2.1525	1\\
2.155	1\\
2.1575	1\\
2.16	1\\
2.1625	1\\
2.165	1\\
2.1675	1\\
2.17	1\\
2.1725	1\\
2.175	1\\
2.1775	1\\
2.18	1\\
2.1825	1\\
2.185	1\\
2.1875	1\\
2.19	1\\
2.1925	1\\
2.195	1\\
2.1975	1\\
2.2	1\\
2.2025	1\\
2.205	1\\
2.2075	1\\
2.21	1\\
2.2125	1\\
2.215	1\\
2.2175	1\\
2.22	0\\
2.2225	0\\
2.225	0\\
2.2275	0\\
2.23	0\\
2.2325	0\\
2.235	0\\
2.2375	0\\
2.24	0\\
2.2425	0\\
2.245	0\\
2.2475	0\\
2.25	0\\
2.2525	0\\
2.255	0\\
2.2575	0\\
2.26	0\\
2.2625	1\\
2.265	1\\
2.2675	1\\
2.27	1\\
2.2725	1\\
2.275	1\\
2.2775	1\\
2.28	0\\
2.2825	0\\
2.285	0\\
2.2875	0\\
2.29	0\\
2.2925	0\\
2.295	0\\
2.2975	0\\
2.3	1\\
2.3025	1\\
2.305	1\\
2.3075	1\\
2.31	1\\
2.3125	1\\
2.315	1\\
2.3175	1\\
2.32	1\\
2.3225	1\\
2.325	1\\
2.3275	1\\
2.33	1\\
2.3325	1\\
2.335	1\\
2.3375	1\\
2.34	1\\
2.3425	1\\
2.345	1\\
2.3475	1\\
2.35	1\\
2.3525	1\\
2.355	1\\
2.3575	1\\
2.36	1\\
2.3625	0\\
2.365	0\\
2.3675	0\\
2.37	0\\
2.3725	0\\
2.375	0\\
2.3775	0\\
2.38	0\\
2.3825	0\\
2.385	0\\
2.3875	0\\
2.39	0\\
2.3925	0\\
2.395	0\\
2.3975	0\\
2.4	0\\
2.4025	0\\
2.405	0\\
2.4075	0\\
2.41	0\\
2.4125	0\\
2.415	0\\
2.4175	0\\
2.42	0\\
2.4225	0\\
2.425	0\\
2.4275	0\\
2.43	0\\
2.4325	0\\
2.435	0\\
2.4375	0\\
2.44	0\\
2.4425	20\\
2.445	0\\
2.4475	0\\
2.45	0\\
2.4525	0\\
2.455	0\\
2.4575	0\\
2.46	0\\
2.4625	1\\
2.465	1\\
2.4675	1\\
2.47	1\\
2.4725	1\\
2.475	1\\
2.4775	1\\
2.48	1\\
2.4825	1\\
2.485	1\\
2.4875	1\\
2.49	1\\
2.4925	1\\
2.495	1\\
2.4975	1\\
2.5	1\\
2.5025	1\\
2.505	1\\
2.5075	1\\
2.51	1\\
2.5125	1\\
2.515	1\\
2.5175	0\\
2.52	0\\
2.5225	0\\
2.525	0\\
2.5275	0\\
2.53	0\\
2.5325	0\\
2.535	0\\
2.5375	0\\
2.54	0\\
2.5425	1\\
2.545	1\\
2.5475	1\\
2.55	1\\
2.5525	1\\
2.555	1\\
2.5575	1\\
2.56	1\\
2.5625	1\\
2.565	1\\
2.5675	0\\
2.57	1\\
2.5725	1\\
2.575	1\\
2.5775	1\\
2.58	1\\
2.5825	1\\
2.585	1\\
2.5875	1\\
2.59	1\\
2.5925	1\\
2.595	1\\
2.5975	1\\
2.6	1\\
2.6025	1\\
2.605	1\\
2.6075	1\\
2.61	1\\
2.6125	1\\
2.615	1\\
2.6175	1\\
2.62	1\\
2.6225	1\\
2.625	1\\
2.6275	1\\
2.63	1\\
2.6325	1\\
2.635	1\\
2.6375	1\\
2.64	1\\
2.6425	1\\
2.645	1\\
2.6475	1\\
2.65	1\\
2.6525	1\\
2.655	1\\
2.6575	1\\
2.66	1\\
2.6625	1\\
2.665	1\\
2.6675	1\\
2.67	1\\
2.6725	1\\
2.675	1\\
2.6775	1\\
2.68	1\\
2.6825	1\\
2.685	1\\
2.6875	1\\
2.69	1\\
2.6925	1\\
2.695	1\\
2.6975	0\\
2.7	1\\
2.7025	1\\
2.705	1\\
2.7075	1\\
2.71	1\\
2.7125	1\\
2.715	1\\
2.7175	1\\
2.72	1\\
2.7225	1\\
2.725	1\\
2.7275	1\\
2.73	1\\
2.7325	1\\
2.735	1\\
2.7375	1\\
2.74	1\\
2.7425	1\\
2.745	1\\
2.7475	1\\
2.75	1\\
2.7525	1\\
2.755	1\\
2.7575	1\\
2.76	1\\
2.7625	1\\
2.765	1\\
2.7675	1\\
2.77	1\\
2.7725	1\\
2.775	1\\
2.7775	1\\
2.78	1\\
2.7825	1\\
2.785	1\\
2.7875	1\\
2.79	1\\
2.7925	1\\
2.795	1\\
2.7975	1\\
2.8	1\\
2.8025	1\\
2.805	1\\
2.8075	1\\
2.81	0\\
2.8125	0\\
2.815	0\\
2.8175	0\\
2.82	0\\
2.8225	0\\
2.825	0\\
2.8275	0\\
2.83	0\\
2.8325	0\\
2.835	0\\
2.8375	0\\
2.84	1\\
2.8425	1\\
2.845	1\\
2.8475	1\\
2.85	1\\
2.8525	1\\
2.855	1\\
2.8575	1\\
2.86	1\\
2.8625	1\\
2.865	1\\
2.8675	1\\
2.87	1\\
2.8725	1\\
2.875	1\\
2.8775	1\\
2.88	1\\
2.8825	1\\
2.885	1\\
2.8875	1\\
2.89	1\\
2.8925	1\\
2.895	1\\
2.8975	0\\
2.9	1\\
2.9025	1\\
2.905	1\\
2.9075	1\\
2.91	1\\
2.9125	1\\
2.915	1\\
2.9175	1\\
2.92	1\\
2.9225	1\\
2.925	1\\
2.9275	1\\
2.93	1\\
2.9325	1\\
2.935	1\\
2.9375	1\\
2.94	1\\
2.9425	1\\
2.945	1\\
2.9475	1\\
2.95	1\\
2.9525	1\\
2.955	1\\
2.9575	1\\
2.96	1\\
2.9625	1\\
2.965	1\\
2.9675	1\\
2.97	1\\
2.9725	1\\
2.975	1\\
2.9775	1\\
2.98	1\\
2.9825	1\\
2.985	1\\
2.9875	1\\
2.99	1\\
2.9925	1\\
2.995	1\\
2.9975	1\\
3	1\\
3.0025	1\\
3.005	1\\
3.0075	1\\
3.01	1\\
3.0125	1\\
3.015	1\\
3.0175	1\\
3.02	1\\
3.0225	1\\
3.025	1\\
3.0275	1\\
3.03	1\\
3.0325	1\\
3.035	1\\
3.0375	1\\
3.04	1\\
3.0425	1\\
3.045	1\\
3.0475	1\\
3.05	1\\
3.0525	1\\
3.055	1\\
3.0575	1\\
3.06	1\\
3.0625	1\\
3.065	1\\
3.0675	0\\
3.07	0\\
3.0725	0\\
3.075	0\\
3.0775	0\\
3.08	0\\
3.0825	0\\
3.085	0\\
3.0875	0\\
3.09	0\\
3.0925	0\\
3.095	0\\
3.0975	0\\
3.1	0\\
3.1025	0\\
3.105	0\\
3.1075	0\\
3.11	0\\
3.1125	0\\
3.115	0\\
3.1175	0\\
3.12	0\\
3.1225	0\\
3.125	0\\
3.1275	0\\
3.13	0\\
3.1325	0\\
3.135	0\\
3.1375	0\\
3.14	0\\
3.1425	0\\
3.145	0\\
3.1475	44\\
3.15	1\\
3.1525	1\\
3.155	35\\
3.1575	1\\
3.16	1\\
3.1625	1\\
3.165	4\\
3.1675	48\\
3.17	3\\
3.1725	1\\
3.175	1\\
3.1775	29\\
3.18	1\\
3.1825	1\\
3.185	19\\
3.1875	6\\
3.19	1\\
3.1925	1\\
3.195	48\\
3.1975	1\\
3.2	1\\
3.2025	1\\
3.205	22\\
3.2075	5\\
3.21	1\\
3.2125	1\\
3.215	1\\
3.2175	0\\
3.22	23\\
3.2225	0\\
3.225	40\\
3.2275	0\\
3.23	1\\
3.2325	0\\
3.235	1\\
3.2375	1\\
3.24	0\\
3.2425	0\\
3.245	0\\
3.2475	1\\
3.25	29\\
3.2525	1\\
3.255	1\\
3.2575	0\\
3.26	0\\
3.2625	0\\
3.265	0\\
3.2675	1\\
3.27	0\\
3.2725	0\\
3.275	0\\
3.2775	0\\
3.28	8\\
3.2825	10\\
3.285	33\\
3.2875	25\\
3.29	34\\
3.2925	1\\
3.295	1\\
3.2975	1\\
3.3	1\\
3.3025	1\\
3.305	1\\
3.3075	1\\
3.31	1\\
3.3125	1\\
3.315	1\\
3.3175	1\\
3.32	1\\
3.3225	1\\
3.325	1\\
3.3275	1\\
3.33	1\\
3.3325	1\\
3.335	1\\
3.3375	1\\
3.34	1\\
3.3425	1\\
3.345	1\\
3.3475	1\\
3.35	13\\
3.3525	1\\
3.355	1\\
3.3575	1\\
3.36	5\\
3.3625	1\\
3.365	1\\
3.3675	1\\
3.37	1\\
3.3725	1\\
3.375	1\\
3.3775	2\\
3.38	1\\
3.3825	1\\
3.385	1\\
3.3875	1\\
3.39	1\\
3.3925	1\\
3.395	1\\
3.3975	1\\
3.4	1\\
3.4025	1\\
3.405	1\\
3.4075	1\\
3.41	1\\
3.4125	1\\
3.415	1\\
3.4175	1\\
3.42	1\\
3.4225	1\\
3.425	1\\
3.4275	1\\
3.43	0\\
3.4325	0\\
3.435	0\\
3.4375	0\\
3.44	0\\
3.4425	0\\
3.445	0\\
3.4475	0\\
3.45	0\\
3.4525	0\\
3.455	0\\
3.4575	0\\
3.46	0\\
3.4625	0\\
3.465	0\\
3.4675	0\\
3.47	0\\
3.4725	0\\
3.475	0\\
3.4775	0\\
3.48	0\\
3.4825	0\\
3.485	0\\
3.4875	0\\
3.49	0\\
3.4925	0\\
3.495	0\\
3.4975	0\\
3.5	0\\
3.5025	0\\
3.505	0\\
3.5075	0\\
3.51	0\\
3.5125	0\\
3.515	0\\
3.5175	0\\
3.52	0\\
3.5225	0\\
3.525	0\\
3.5275	0\\
3.53	0\\
3.5325	0\\
3.535	0\\
3.5375	0\\
3.54	0\\
3.5425	0\\
3.545	0\\
3.5475	0\\
3.55	0\\
3.5525	0\\
3.555	0\\
3.5575	0\\
3.56	0\\
3.5625	0\\
3.565	0\\
3.5675	0\\
3.57	0\\
3.5725	0\\
3.575	0\\
3.5775	0\\
3.58	0\\
3.5825	0\\
3.585	0\\
3.5875	0\\
3.59	0\\
3.5925	0\\
3.595	0\\
3.5975	0\\
3.6	0\\
3.6025	4\\
3.605	1\\
3.6075	1\\
3.61	1\\
3.6125	1\\
3.615	1\\
3.6175	9\\
3.62	1\\
3.6225	20\\
3.625	1\\
3.6275	1\\
3.63	9\\
3.6325	1\\
3.635	12\\
3.6375	1\\
3.64	1\\
3.6425	1\\
3.645	1\\
3.6475	1\\
3.65	1\\
3.6525	1\\
3.655	1\\
3.6575	1\\
3.66	1\\
3.6625	1\\
3.665	1\\
3.6675	0\\
3.67	1\\
3.6725	1\\
3.675	1\\
3.6775	1\\
3.68	1\\
3.6825	1\\
3.685	1\\
3.6875	1\\
3.69	1\\
3.6925	17\\
3.695	1\\
3.6975	1\\
3.7	1\\
3.7025	1\\
3.705	1\\
3.7075	1\\
3.71	1\\
3.7125	1\\
3.715	1\\
3.7175	1\\
3.72	1\\
3.7225	1\\
3.725	1\\
3.7275	1\\
3.73	1\\
3.7325	1\\
3.735	1\\
3.7375	1\\
3.74	1\\
3.7425	1\\
3.745	1\\
3.7475	1\\
3.75	1\\
3.7525	1\\
3.755	1\\
3.7575	1\\
3.76	1\\
3.7625	1\\
3.765	1\\
3.7675	1\\
3.77	1\\
3.7725	1\\
3.775	1\\
3.7775	1\\
3.78	1\\
3.7825	1\\
3.785	1\\
3.7875	1\\
3.79	1\\
3.7925	1\\
3.795	1\\
3.7975	1\\
3.8	1\\
3.8025	1\\
3.805	1\\
3.8075	1\\
3.81	1\\
3.8125	1\\
3.815	1\\
3.8175	1\\
3.82	1\\
3.8225	1\\
3.825	1\\
3.8275	1\\
3.83	1\\
3.8325	1\\
3.835	1\\
3.8375	1\\
3.84	1\\
3.8425	1\\
3.845	1\\
3.8475	1\\
3.85	1\\
3.8525	1\\
3.855	1\\
3.8575	1\\
3.86	1\\
3.8625	1\\
3.865	1\\
3.8675	1\\
3.87	1\\
3.8725	1\\
3.875	1\\
3.8775	1\\
3.88	1\\
3.8825	1\\
3.885	1\\
3.8875	1\\
3.89	1\\
3.8925	1\\
3.895	1\\
3.8975	1\\
3.9	1\\
3.9025	1\\
3.905	1\\
3.9075	1\\
3.91	1\\
3.9125	1\\
3.915	1\\
3.9175	1\\
3.92	1\\
3.9225	1\\
3.925	1\\
3.9275	1\\
3.93	1\\
3.9325	1\\
3.935	1\\
3.9375	1\\
3.94	1\\
3.9425	1\\
3.945	1\\
3.9475	1\\
3.95	1\\
3.9525	1\\
3.955	1\\
3.9575	1\\
3.96	1\\
3.9625	1\\
3.965	1\\
3.9675	1\\
3.97	1\\
3.9725	1\\
3.975	1\\
3.9775	1\\
3.98	1\\
3.9825	1\\
3.985	1\\
3.9875	1\\
3.99	1\\
3.9925	1\\
3.995	1\\
3.9975	1\\
4	1\\
4.0025	1\\
4.005	1\\
4.0075	1\\
4.01	1\\
4.0125	1\\
4.015	1\\
4.0175	1\\
4.02	1\\
4.0225	1\\
4.025	1\\
4.0275	1\\
4.03	1\\
4.0325	1\\
4.035	1\\
4.0375	1\\
4.04	1\\
4.0425	1\\
4.045	1\\
4.0475	1\\
4.05	1\\
4.0525	1\\
4.055	1\\
4.0575	1\\
4.06	1\\
4.0625	1\\
4.065	1\\
4.0675	1\\
4.07	1\\
4.0725	1\\
4.075	1\\
4.0775	1\\
4.08	1\\
4.0825	1\\
4.085	1\\
4.0875	1\\
4.09	1\\
4.0925	1\\
4.095	1\\
4.0975	1\\
4.1	0\\
4.1025	0\\
4.105	0\\
4.1075	1\\
4.11	1\\
4.1125	1\\
4.115	1\\
4.1175	1\\
4.12	1\\
4.1225	1\\
4.125	1\\
4.1275	1\\
4.13	1\\
4.1325	1\\
4.135	1\\
4.1375	1\\
4.14	1\\
4.1425	1\\
4.145	1\\
4.1475	1\\
4.15	1\\
4.1525	1\\
4.155	1\\
4.1575	1\\
4.16	1\\
4.1625	1\\
4.165	1\\
4.1675	1\\
4.17	1\\
4.1725	1\\
4.175	1\\
4.1775	26\\
4.18	1\\
4.1825	0\\
4.185	1\\
4.1875	1\\
4.19	1\\
4.1925	1\\
4.195	0\\
4.1975	1\\
4.2	1\\
4.2025	0\\
4.205	0\\
4.2075	0\\
4.21	0\\
4.2125	0\\
4.215	0\\
4.2175	0\\
4.22	0\\
4.2225	0\\
4.225	0\\
4.2275	0\\
4.23	0\\
4.2325	0\\
4.235	0\\
4.2375	0\\
4.24	0\\
4.2425	0\\
4.245	0\\
4.2475	0\\
4.25	0\\
4.2525	0\\
4.255	0\\
4.2575	0\\
4.26	0\\
4.2625	0\\
4.265	0\\
4.2675	0\\
4.27	0\\
4.2725	0\\
4.275	0\\
4.2775	0\\
4.28	0\\
4.2825	0\\
4.285	0\\
4.2875	0\\
4.29	0\\
4.2925	0\\
4.295	0\\
4.2975	0\\
4.3	1\\
4.3025	1\\
4.305	1\\
4.3075	1\\
4.31	1\\
4.3125	1\\
4.315	0\\
4.3175	0\\
4.32	0\\
4.3225	0\\
4.325	0\\
4.3275	0\\
4.33	0\\
4.3325	0\\
4.335	0\\
4.3375	0\\
4.34	1\\
4.3425	1\\
4.345	1\\
4.3475	1\\
4.35	1\\
4.3525	1\\
4.355	1\\
4.3575	1\\
4.36	1\\
4.3625	1\\
4.365	1\\
4.3675	1\\
4.37	1\\
4.3725	1\\
4.375	1\\
4.3775	12\\
4.38	0\\
4.3825	0\\
4.385	0\\
4.3875	1\\
4.39	1\\
4.3925	14\\
4.395	1\\
4.3975	0\\
4.4	4\\
4.4025	0\\
4.405	1\\
4.4075	1\\
4.41	1\\
4.4125	1\\
4.415	1\\
4.4175	1\\
4.42	1\\
4.4225	1\\
4.425	1\\
4.4275	1\\
4.43	1\\
4.4325	1\\
4.435	1\\
4.4375	1\\
4.44	1\\
4.4425	1\\
4.445	1\\
4.4475	1\\
4.45	1\\
4.4525	1\\
4.455	1\\
4.4575	1\\
4.46	1\\
4.4625	1\\
4.465	1\\
4.4675	1\\
4.47	1\\
4.4725	1\\
4.475	1\\
4.4775	1\\
4.48	1\\
4.4825	1\\
4.485	1\\
4.4875	1\\
4.49	1\\
4.4925	1\\
4.495	1\\
4.4975	1\\
4.5	1\\
4.5025	1\\
4.505	1\\
4.5075	1\\
4.51	0\\
4.5125	1\\
4.515	5\\
4.5175	0\\
4.52	0\\
4.5225	0\\
4.525	0\\
4.5275	0\\
4.53	0\\
4.5325	1\\
4.535	19\\
4.5375	19\\
4.54	1\\
4.5425	1\\
4.545	28\\
4.5475	9\\
4.55	1\\
4.5525	11\\
4.555	1\\
4.5575	1\\
4.56	1\\
4.5625	1\\
4.565	1\\
4.5675	1\\
4.57	1\\
4.5725	1\\
4.575	1\\
4.5775	8\\
4.58	4\\
4.5825	9\\
4.585	4\\
4.5875	7\\
4.59	1\\
4.5925	1\\
4.595	1\\
4.5975	1\\
4.6	13\\
4.6025	25\\
4.605	1\\
4.6075	1\\
4.61	1\\
4.6125	1\\
4.615	1\\
4.6175	1\\
4.62	1\\
4.6225	1\\
4.625	1\\
4.6275	1\\
4.63	1\\
4.6325	5\\
4.635	1\\
4.6375	1\\
4.64	54\\
4.6425	1\\
4.645	1\\
4.6475	5\\
4.65	1\\
4.6525	1\\
4.655	3\\
4.6575	1\\
4.66	1\\
4.6625	23\\
4.665	1\\
4.6675	1\\
4.67	1\\
4.6725	10\\
4.675	0\\
4.6775	1\\
4.68	10\\
4.6825	0\\
4.685	1\\
4.6875	4\\
4.69	0\\
4.6925	1\\
4.695	1\\
4.6975	19\\
4.7	0\\
4.7025	1\\
4.705	4\\
4.7075	1\\
4.71	1\\
4.7125	1\\
4.715	1\\
4.7175	1\\
4.72	1\\
4.7225	1\\
4.725	1\\
4.7275	1\\
4.73	1\\
4.7325	1\\
4.735	1\\
4.7375	1\\
4.74	1\\
4.7425	1\\
4.745	1\\
4.7475	1\\
4.75	1\\
4.7525	1\\
4.755	1\\
4.7575	1\\
4.76	1\\
4.7625	1\\
4.765	1\\
4.7675	1\\
4.77	1\\
4.7725	0\\
4.775	0\\
4.7775	0\\
4.78	0\\
4.7825	0\\
4.785	0\\
4.7875	0\\
4.79	0\\
4.7925	0\\
4.795	0\\
4.7975	0\\
4.8	0\\
4.8025	0\\
4.805	0\\
4.8075	0\\
4.81	0\\
4.8125	0\\
4.815	0\\
4.8175	0\\
4.82	0\\
4.8225	0\\
4.825	0\\
4.8275	0\\
4.83	1\\
4.8325	0\\
4.835	0\\
4.8375	1\\
4.84	0\\
4.8425	1\\
4.845	1\\
4.8475	0\\
4.85	0\\
4.8525	0\\
4.855	0\\
4.8575	0\\
4.86	1\\
4.8625	1\\
4.865	1\\
4.8675	1\\
4.87	1\\
4.8725	1\\
4.875	1\\
4.8775	1\\
4.88	1\\
4.8825	1\\
4.885	1\\
4.8875	1\\
4.89	1\\
4.8925	1\\
4.895	1\\
4.8975	1\\
4.9	1\\
4.9025	1\\
4.905	1\\
4.9075	1\\
4.91	1\\
4.9125	1\\
4.915	1\\
4.9175	1\\
4.92	1\\
4.9225	1\\
4.925	1\\
4.9275	1\\
4.93	1\\
4.9325	1\\
4.935	1\\
4.9375	1\\
4.94	1\\
4.9425	1\\
4.945	1\\
4.9475	1\\
4.95	1\\
4.9525	1\\
4.955	1\\
4.9575	1\\
4.96	1\\
4.9625	1\\
4.965	1\\
4.9675	1\\
4.97	1\\
4.9725	0\\
4.975	0\\
4.9775	1\\
4.98	1\\
4.9825	1\\
4.985	1\\
4.9875	1\\
4.99	1\\
4.9925	1\\
4.995	1\\
4.9975	1\\
};
\addlegendentry{$\epsilon_u = 10^{-8}$}

\addplot [color=mycolor2, line width=1.2pt,opacity=.5]
  table[row sep=crcr]{%
0	1\\
0.0025	1\\
0.005	0\\
0.0075	0\\
0.01	0\\
0.0125	1\\
0.015	1\\
0.0175	1\\
0.02	0\\
0.0225	0\\
0.025	0\\
0.0275	1\\
0.03	1\\
0.0325	0\\
0.035	0\\
0.0375	1\\
0.04	1\\
0.0425	1\\
0.045	1\\
0.0475	1\\
0.05	0\\
0.0525	0\\
0.055	0\\
0.0575	0\\
0.06	0\\
0.0625	0\\
0.065	0\\
0.0675	0\\
0.07	0\\
0.0725	0\\
0.075	0\\
0.0775	0\\
0.08	0\\
0.0825	0\\
0.085	0\\
0.0875	0\\
0.09	0\\
0.0925	0\\
0.095	0\\
0.0975	0\\
0.1	0\\
0.1025	0\\
0.105	0\\
0.1075	0\\
0.11	0\\
0.1125	0\\
0.115	1\\
0.1175	0\\
0.12	0\\
0.1225	0\\
0.125	0\\
0.1275	0\\
0.13	0\\
0.1325	0\\
0.135	0\\
0.1375	0\\
0.14	0\\
0.1425	0\\
0.145	0\\
0.1475	0\\
0.15	0\\
0.1525	1\\
0.155	1\\
0.1575	1\\
0.16	1\\
0.1625	0\\
0.165	0\\
0.1675	0\\
0.17	0\\
0.1725	1\\
0.175	1\\
0.1775	1\\
0.18	1\\
0.1825	1\\
0.185	1\\
0.1875	1\\
0.19	1\\
0.1925	0\\
0.195	0\\
0.1975	0\\
0.2	0\\
0.2025	0\\
0.205	0\\
0.2075	0\\
0.21	0\\
0.2125	1\\
0.215	1\\
0.2175	1\\
0.22	1\\
0.2225	1\\
0.225	1\\
0.2275	0\\
0.23	0\\
0.2325	0\\
0.235	0\\
0.2375	0\\
0.24	0\\
0.2425	0\\
0.245	0\\
0.2475	0\\
0.25	0\\
0.2525	0\\
0.255	1\\
0.2575	1\\
0.26	1\\
0.2625	1\\
0.265	0\\
0.2675	1\\
0.27	0\\
0.2725	1\\
0.275	1\\
0.2775	1\\
0.28	0\\
0.2825	0\\
0.285	0\\
0.2875	0\\
0.29	0\\
0.2925	0\\
0.295	0\\
0.2975	1\\
0.3	1\\
0.3025	1\\
0.305	1\\
0.3075	0\\
0.31	0\\
0.3125	0\\
0.315	0\\
0.3175	0\\
0.32	0\\
0.3225	0\\
0.325	0\\
0.3275	1\\
0.33	0\\
0.3325	0\\
0.335	0\\
0.3375	0\\
0.34	0\\
0.3425	1\\
0.345	0\\
0.3475	0\\
0.35	0\\
0.3525	1\\
0.355	1\\
0.3575	0\\
0.36	0\\
0.3625	0\\
0.365	1\\
0.3675	1\\
0.37	1\\
0.3725	0\\
0.375	0\\
0.3775	0\\
0.38	0\\
0.3825	0\\
0.385	0\\
0.3875	0\\
0.39	0\\
0.3925	0\\
0.395	0\\
0.3975	0\\
0.4	0\\
0.4025	0\\
0.405	0\\
0.4075	0\\
0.41	0\\
0.4125	0\\
0.415	0\\
0.4175	0\\
0.42	1\\
0.4225	0\\
0.425	0\\
0.4275	0\\
0.43	0\\
0.4325	0\\
0.435	0\\
0.4375	0\\
0.44	0\\
0.4425	0\\
0.445	0\\
0.4475	0\\
0.45	0\\
0.4525	0\\
0.455	0\\
0.4575	0\\
0.46	0\\
0.4625	0\\
0.465	0\\
0.4675	0\\
0.47	0\\
0.4725	0\\
0.475	0\\
0.4775	0\\
0.48	0\\
0.4825	0\\
0.485	0\\
0.4875	0\\
0.49	0\\
0.4925	0\\
0.495	0\\
0.4975	0\\
0.5	0\\
0.5025	0\\
0.505	0\\
0.5075	0\\
0.51	0\\
0.5125	0\\
0.515	0\\
0.5175	0\\
0.52	0\\
0.5225	0\\
0.525	0\\
0.5275	0\\
0.53	0\\
0.5325	0\\
0.535	0\\
0.5375	0\\
0.54	0\\
0.5425	0\\
0.545	0\\
0.5475	0\\
0.55	0\\
0.5525	0\\
0.555	0\\
0.5575	0\\
0.56	0\\
0.5625	0\\
0.565	0\\
0.5675	0\\
0.57	0\\
0.5725	0\\
0.575	0\\
0.5775	0\\
0.58	1\\
0.5825	1\\
0.585	0\\
0.5875	0\\
0.59	0\\
0.5925	0\\
0.595	0\\
0.5975	9\\
0.6	0\\
0.6025	1\\
0.605	1\\
0.6075	1\\
0.61	0\\
0.6125	0\\
0.615	0\\
0.6175	0\\
0.62	0\\
0.6225	0\\
0.625	0\\
0.6275	1\\
0.63	1\\
0.6325	1\\
0.635	1\\
0.6375	1\\
0.64	1\\
0.6425	1\\
0.645	1\\
0.6475	0\\
0.65	0\\
0.6525	1\\
0.655	1\\
0.6575	1\\
0.66	1\\
0.6625	1\\
0.665	0\\
0.6675	1\\
0.67	0\\
0.6725	0\\
0.675	0\\
0.6775	1\\
0.68	1\\
0.6825	1\\
0.685	0\\
0.6875	0\\
0.69	0\\
0.6925	0\\
0.695	0\\
0.6975	0\\
0.7	1\\
0.7025	1\\
0.705	0\\
0.7075	0\\
0.71	0\\
0.7125	0\\
0.715	0\\
0.7175	0\\
0.72	0\\
0.7225	0\\
0.725	1\\
0.7275	0\\
0.73	0\\
0.7325	0\\
0.735	0\\
0.7375	0\\
0.74	1\\
0.7425	0\\
0.745	0\\
0.7475	1\\
0.75	1\\
0.7525	1\\
0.755	1\\
0.7575	1\\
0.76	1\\
0.7625	1\\
0.765	1\\
0.7675	1\\
0.77	0\\
0.7725	0\\
0.775	0\\
0.7775	0\\
0.78	0\\
0.7825	0\\
0.785	0\\
0.7875	0\\
0.79	1\\
0.7925	1\\
0.795	0\\
0.7975	0\\
0.8	1\\
0.8025	1\\
0.805	1\\
0.8075	1\\
0.81	1\\
0.8125	1\\
0.815	1\\
0.8175	1\\
0.82	1\\
0.8225	0\\
0.825	0\\
0.8275	0\\
0.83	0\\
0.8325	1\\
0.835	1\\
0.8375	1\\
0.84	1\\
0.8425	1\\
0.845	1\\
0.8475	1\\
0.85	1\\
0.8525	0\\
0.855	0\\
0.8575	0\\
0.86	0\\
0.8625	1\\
0.865	1\\
0.8675	1\\
0.87	1\\
0.8725	0\\
0.875	1\\
0.8775	1\\
0.88	0\\
0.8825	0\\
0.885	0\\
0.8875	0\\
0.89	0\\
0.8925	1\\
0.895	1\\
0.8975	1\\
0.9	0\\
0.9025	0\\
0.905	1\\
0.9075	1\\
0.91	0\\
0.9125	0\\
0.915	0\\
0.9175	0\\
0.92	0\\
0.9225	0\\
0.925	0\\
0.9275	0\\
0.93	0\\
0.9325	0\\
0.935	0\\
0.9375	1\\
0.94	1\\
0.9425	1\\
0.945	1\\
0.9475	0\\
0.95	0\\
0.9525	0\\
0.955	0\\
0.9575	1\\
0.96	1\\
0.9625	0\\
0.965	0\\
0.9675	0\\
0.97	1\\
0.9725	1\\
0.975	0\\
0.9775	0\\
0.98	0\\
0.9825	1\\
0.985	1\\
0.9875	0\\
0.99	0\\
0.9925	0\\
0.995	1\\
0.9975	1\\
1	0\\
1.0025	1\\
1.005	1\\
1.0075	1\\
1.01	1\\
1.0125	1\\
1.015	0\\
1.0175	0\\
1.02	0\\
1.0225	0\\
1.025	0\\
1.0275	0\\
1.03	1\\
1.0325	1\\
1.035	1\\
1.0375	1\\
1.04	0\\
1.0425	0\\
1.045	0\\
1.0475	0\\
1.05	0\\
1.0525	0\\
1.055	0\\
1.0575	0\\
1.06	0\\
1.0625	0\\
1.065	0\\
1.0675	0\\
1.07	0\\
1.0725	0\\
1.075	1\\
1.0775	1\\
1.08	1\\
1.0825	0\\
1.085	1\\
1.0875	1\\
1.09	1\\
1.0925	1\\
1.095	1\\
1.0975	0\\
1.1	0\\
1.1025	1\\
1.105	1\\
1.1075	0\\
1.11	0\\
1.1125	1\\
1.115	1\\
1.1175	1\\
1.12	1\\
1.1225	1\\
1.125	1\\
1.1275	1\\
1.13	0\\
1.1325	0\\
1.135	1\\
1.1375	1\\
1.14	1\\
1.1425	1\\
1.145	1\\
1.1475	1\\
1.15	1\\
1.1525	0\\
1.155	0\\
1.1575	0\\
1.16	0\\
1.1625	0\\
1.165	1\\
1.1675	0\\
1.17	1\\
1.1725	0\\
1.175	1\\
1.1775	0\\
1.18	0\\
1.1825	0\\
1.185	0\\
1.1875	0\\
1.19	0\\
1.1925	0\\
1.195	0\\
1.1975	0\\
1.2	0\\
1.2025	0\\
1.205	0\\
1.2075	0\\
1.21	0\\
1.2125	0\\
1.215	1\\
1.2175	1\\
1.22	0\\
1.2225	0\\
1.225	0\\
1.2275	0\\
1.23	0\\
1.2325	0\\
1.235	1\\
1.2375	0\\
1.24	1\\
1.2425	0\\
1.245	0\\
1.2475	1\\
1.25	1\\
1.2525	1\\
1.255	1\\
1.2575	0\\
1.26	0\\
1.2625	0\\
1.265	1\\
1.2675	0\\
1.27	1\\
1.2725	0\\
1.275	0\\
1.2775	0\\
1.28	0\\
1.2825	0\\
1.285	0\\
1.2875	0\\
1.29	0\\
1.2925	1\\
1.295	1\\
1.2975	1\\
1.3	1\\
1.3025	1\\
1.305	1\\
1.3075	1\\
1.31	1\\
1.3125	1\\
1.315	1\\
1.3175	1\\
1.32	1\\
1.3225	1\\
1.325	0\\
1.3275	0\\
1.33	0\\
1.3325	0\\
1.335	0\\
1.3375	0\\
1.34	1\\
1.3425	1\\
1.345	1\\
1.3475	1\\
1.35	0\\
1.3525	0\\
1.355	0\\
1.3575	0\\
1.36	1\\
1.3625	1\\
1.365	0\\
1.3675	0\\
1.37	0\\
1.3725	0\\
1.375	1\\
1.3775	0\\
1.38	1\\
1.3825	0\\
1.385	0\\
1.3875	0\\
1.39	0\\
1.3925	0\\
1.395	0\\
1.3975	0\\
1.4	0\\
1.4025	1\\
1.405	0\\
1.4075	0\\
1.41	0\\
1.4125	0\\
1.415	0\\
1.4175	0\\
1.42	1\\
1.4225	0\\
1.425	1\\
1.4275	1\\
1.43	0\\
1.4325	1\\
1.435	0\\
1.4375	0\\
1.44	0\\
1.4425	0\\
1.445	0\\
1.4475	1\\
1.45	1\\
1.4525	1\\
1.455	1\\
1.4575	1\\
1.46	1\\
1.4625	1\\
1.465	1\\
1.4675	0\\
1.47	1\\
1.4725	0\\
1.475	0\\
1.4775	0\\
1.48	0\\
1.4825	0\\
1.485	0\\
1.4875	0\\
1.49	0\\
1.4925	0\\
1.495	0\\
1.4975	1\\
1.5	1\\
1.5025	1\\
1.505	1\\
1.5075	1\\
1.51	1\\
1.5125	0\\
1.515	0\\
1.5175	0\\
1.52	0\\
1.5225	0\\
1.525	0\\
1.5275	0\\
1.53	0\\
1.5325	1\\
1.535	1\\
1.5375	1\\
1.54	1\\
1.5425	1\\
1.545	1\\
1.5475	1\\
1.55	1\\
1.5525	0\\
1.555	0\\
1.5575	0\\
1.56	0\\
1.5625	0\\
1.565	0\\
1.5675	0\\
1.57	0\\
1.5725	1\\
1.575	0\\
1.5775	0\\
1.58	0\\
1.5825	1\\
1.585	0\\
1.5875	0\\
1.59	0\\
1.5925	0\\
1.595	0\\
1.5975	0\\
1.6	0\\
1.6025	0\\
1.605	0\\
1.6075	0\\
1.61	0\\
1.6125	0\\
1.615	0\\
1.6175	0\\
1.62	0\\
1.6225	0\\
1.625	0\\
1.6275	0\\
1.63	0\\
1.6325	1\\
1.635	1\\
1.6375	1\\
1.64	1\\
1.6425	0\\
1.645	0\\
1.6475	1\\
1.65	1\\
1.6525	2\\
1.655	2\\
1.6575	2\\
1.66	2\\
1.6625	2\\
1.665	0\\
1.6675	0\\
1.67	1\\
1.6725	1\\
1.675	1\\
1.6775	1\\
1.68	1\\
1.6825	1\\
1.685	1\\
1.6875	0\\
1.69	0\\
1.6925	1\\
1.695	1\\
1.6975	0\\
1.7	0\\
1.7025	0\\
1.705	0\\
1.7075	0\\
1.71	0\\
1.7125	1\\
1.715	1\\
1.7175	1\\
1.72	1\\
1.7225	1\\
1.725	0\\
1.7275	1\\
1.73	1\\
1.7325	0\\
1.735	2\\
1.7375	1\\
1.74	0\\
1.7425	0\\
1.745	0\\
1.7475	0\\
1.75	0\\
1.7525	0\\
1.755	1\\
1.7575	1\\
1.76	1\\
1.7625	1\\
1.765	0\\
1.7675	1\\
1.77	1\\
1.7725	1\\
1.775	1\\
1.7775	1\\
1.78	1\\
1.7825	1\\
1.785	1\\
1.7875	1\\
1.79	1\\
1.7925	1\\
1.795	0\\
1.7975	0\\
1.8	1\\
1.8025	0\\
1.805	0\\
1.8075	0\\
1.81	0\\
1.8125	0\\
1.815	1\\
1.8175	1\\
1.82	1\\
1.8225	0\\
1.825	0\\
1.8275	0\\
1.83	0\\
1.8325	0\\
1.835	0\\
1.8375	0\\
1.84	0\\
1.8425	0\\
1.845	0\\
1.8475	0\\
1.85	0\\
1.8525	0\\
1.855	0\\
1.8575	0\\
1.86	0\\
1.8625	0\\
1.865	0\\
1.8675	0\\
1.87	0\\
1.8725	0\\
1.875	0\\
1.8775	0\\
1.88	0\\
1.8825	0\\
1.885	0\\
1.8875	0\\
1.89	0\\
1.8925	0\\
1.895	0\\
1.8975	0\\
1.9	0\\
1.9025	0\\
1.905	0\\
1.9075	0\\
1.91	1\\
1.9125	1\\
1.915	1\\
1.9175	1\\
1.92	1\\
1.9225	1\\
1.925	1\\
1.9275	1\\
1.93	1\\
1.9325	1\\
1.935	0\\
1.9375	0\\
1.94	0\\
1.9425	0\\
1.945	0\\
1.9475	0\\
1.95	0\\
1.9525	0\\
1.955	0\\
1.9575	0\\
1.96	0\\
1.9625	1\\
1.965	1\\
1.9675	1\\
1.97	1\\
1.9725	1\\
1.975	0\\
1.9775	0\\
1.98	0\\
1.9825	0\\
1.985	0\\
1.9875	0\\
1.99	1\\
1.9925	1\\
1.995	1\\
1.9975	0\\
2	0\\
2.0025	0\\
2.005	0\\
2.0075	0\\
2.01	0\\
2.0125	0\\
2.015	0\\
2.0175	0\\
2.02	0\\
2.0225	0\\
2.025	0\\
2.0275	0\\
2.03	0\\
2.0325	0\\
2.035	1\\
2.0375	0\\
2.04	0\\
2.0425	0\\
2.045	0\\
2.0475	0\\
2.05	0\\
2.0525	1\\
2.055	1\\
2.0575	1\\
2.06	1\\
2.0625	1\\
2.065	1\\
2.0675	1\\
2.07	1\\
2.0725	1\\
2.075	1\\
2.0775	1\\
2.08	1\\
2.0825	0\\
2.085	0\\
2.0875	0\\
2.09	0\\
2.0925	0\\
2.095	0\\
2.0975	0\\
2.1	1\\
2.1025	1\\
2.105	1\\
2.1075	0\\
2.11	0\\
2.1125	0\\
2.115	1\\
2.1175	1\\
2.12	0\\
2.1225	1\\
2.125	1\\
2.1275	1\\
2.13	1\\
2.1325	1\\
2.135	1\\
2.1375	1\\
2.14	1\\
2.1425	1\\
2.145	1\\
2.1475	1\\
2.15	1\\
2.1525	0\\
2.155	0\\
2.1575	0\\
2.16	0\\
2.1625	0\\
2.165	0\\
2.1675	0\\
2.17	0\\
2.1725	0\\
2.175	0\\
2.1775	0\\
2.18	0\\
2.1825	1\\
2.185	1\\
2.1875	1\\
2.19	1\\
2.1925	1\\
2.195	1\\
2.1975	1\\
2.2	0\\
2.2025	0\\
2.205	0\\
2.2075	0\\
2.21	0\\
2.2125	0\\
2.215	0\\
2.2175	0\\
2.22	0\\
2.2225	0\\
2.225	0\\
2.2275	1\\
2.23	1\\
2.2325	0\\
2.235	1\\
2.2375	1\\
2.24	1\\
2.2425	0\\
2.245	0\\
2.2475	0\\
2.25	0\\
2.2525	0\\
2.255	0\\
2.2575	0\\
2.26	0\\
2.2625	0\\
2.265	0\\
2.2675	0\\
2.27	0\\
2.2725	0\\
2.275	0\\
2.2775	0\\
2.28	0\\
2.2825	0\\
2.285	0\\
2.2875	0\\
2.29	0\\
2.2925	0\\
2.295	0\\
2.2975	0\\
2.3	1\\
2.3025	0\\
2.305	0\\
2.3075	0\\
2.31	0\\
2.3125	0\\
2.315	1\\
2.3175	1\\
2.32	1\\
2.3225	0\\
2.325	0\\
2.3275	0\\
2.33	0\\
2.3325	0\\
2.335	0\\
2.3375	0\\
2.34	0\\
2.3425	0\\
2.345	1\\
2.3475	1\\
2.35	1\\
2.3525	0\\
2.355	0\\
2.3575	0\\
2.36	0\\
2.3625	0\\
2.365	0\\
2.3675	0\\
2.37	0\\
2.3725	1\\
2.375	0\\
2.3775	0\\
2.38	0\\
2.3825	0\\
2.385	0\\
2.3875	0\\
2.39	0\\
2.3925	0\\
2.395	0\\
2.3975	1\\
2.4	1\\
2.4025	1\\
2.405	1\\
2.4075	1\\
2.41	1\\
2.4125	1\\
2.415	1\\
2.4175	0\\
2.42	0\\
2.4225	0\\
2.425	0\\
2.4275	1\\
2.43	0\\
2.4325	1\\
2.435	1\\
2.4375	0\\
2.44	1\\
2.4425	0\\
2.445	1\\
2.4475	0\\
2.45	1\\
2.4525	0\\
2.455	0\\
2.4575	1\\
2.46	1\\
2.4625	1\\
2.465	1\\
2.4675	0\\
2.47	0\\
2.4725	0\\
2.475	0\\
2.4775	0\\
2.48	0\\
2.4825	0\\
2.485	0\\
2.4875	0\\
2.49	0\\
2.4925	1\\
2.495	1\\
2.4975	1\\
2.5	1\\
2.5025	1\\
2.505	1\\
2.5075	0\\
2.51	0\\
2.5125	1\\
2.515	1\\
2.5175	1\\
2.52	1\\
2.5225	0\\
2.525	1\\
2.5275	1\\
2.53	1\\
2.5325	1\\
2.535	1\\
2.5375	1\\
2.54	1\\
2.5425	1\\
2.545	1\\
2.5475	0\\
2.55	0\\
2.5525	1\\
2.555	0\\
2.5575	0\\
2.56	0\\
2.5625	1\\
2.565	1\\
2.5675	1\\
2.57	1\\
2.5725	0\\
2.575	0\\
2.5775	0\\
2.58	0\\
2.5825	0\\
2.585	0\\
2.5875	0\\
2.59	0\\
2.5925	0\\
2.595	1\\
2.5975	0\\
2.6	0\\
2.6025	0\\
2.605	1\\
2.6075	1\\
2.61	1\\
2.6125	2\\
2.615	2\\
2.6175	1\\
2.62	1\\
2.6225	1\\
2.625	1\\
2.6275	1\\
2.63	1\\
2.6325	1\\
2.635	0\\
2.6375	0\\
2.64	0\\
2.6425	0\\
2.645	0\\
2.6475	0\\
2.65	0\\
2.6525	0\\
2.655	1\\
2.6575	1\\
2.66	1\\
2.6625	1\\
2.665	1\\
2.6675	1\\
2.67	1\\
2.6725	1\\
2.675	1\\
2.6775	1\\
2.68	1\\
2.6825	1\\
2.685	1\\
2.6875	1\\
2.69	1\\
2.6925	0\\
2.695	0\\
2.6975	0\\
2.7	2\\
2.7025	1\\
2.705	1\\
2.7075	1\\
2.71	1\\
2.7125	1\\
2.715	0\\
2.7175	0\\
2.72	1\\
2.7225	0\\
2.725	1\\
2.7275	0\\
2.73	0\\
2.7325	1\\
2.735	0\\
2.7375	0\\
2.74	0\\
2.7425	0\\
2.745	0\\
2.7475	0\\
2.75	0\\
2.7525	0\\
2.755	1\\
2.7575	1\\
2.76	1\\
2.7625	1\\
2.765	1\\
2.7675	1\\
2.77	1\\
2.7725	1\\
2.775	1\\
2.7775	0\\
2.78	0\\
2.7825	1\\
2.785	0\\
2.7875	0\\
2.79	0\\
2.7925	1\\
2.795	1\\
2.7975	0\\
2.8	1\\
2.8025	0\\
2.805	1\\
2.8075	1\\
2.81	1\\
2.8125	1\\
2.815	1\\
2.8175	1\\
2.82	1\\
2.8225	1\\
2.825	1\\
2.8275	1\\
2.83	0\\
2.8325	0\\
2.835	1\\
2.8375	1\\
2.84	2\\
2.8425	1\\
2.845	1\\
2.8475	1\\
2.85	0\\
2.8525	0\\
2.855	0\\
2.8575	1\\
2.86	0\\
2.8625	0\\
2.865	0\\
2.8675	1\\
2.87	1\\
2.8725	1\\
2.875	1\\
2.8775	1\\
2.88	1\\
2.8825	1\\
2.885	1\\
2.8875	1\\
2.89	1\\
2.8925	0\\
2.895	0\\
2.8975	0\\
2.9	1\\
2.9025	1\\
2.905	1\\
2.9075	1\\
2.91	1\\
2.9125	1\\
2.915	1\\
2.9175	1\\
2.92	1\\
2.9225	1\\
2.925	1\\
2.9275	1\\
2.93	1\\
2.9325	1\\
2.935	1\\
2.9375	1\\
2.94	1\\
2.9425	1\\
2.945	1\\
2.9475	1\\
2.95	1\\
2.9525	1\\
2.955	1\\
2.9575	1\\
2.96	1\\
2.9625	1\\
2.965	1\\
2.9675	1\\
2.97	1\\
2.9725	1\\
2.975	1\\
2.9775	1\\
2.98	1\\
2.9825	0\\
2.985	0\\
2.9875	0\\
2.99	0\\
2.9925	1\\
2.995	0\\
2.9975	1\\
3	1\\
3.0025	1\\
3.005	1\\
3.0075	1\\
3.01	1\\
3.0125	1\\
3.015	0\\
3.0175	0\\
3.02	0\\
3.0225	1\\
3.025	1\\
3.0275	1\\
3.03	1\\
3.0325	1\\
3.035	1\\
3.0375	0\\
3.04	0\\
3.0425	1\\
3.045	1\\
3.0475	1\\
3.05	1\\
3.0525	1\\
3.055	1\\
3.0575	1\\
3.06	1\\
3.0625	1\\
3.065	1\\
3.0675	0\\
3.07	0\\
3.0725	0\\
3.075	1\\
3.0775	1\\
3.08	1\\
3.0825	1\\
3.085	1\\
3.0875	1\\
3.09	1\\
3.0925	0\\
3.095	0\\
3.0975	0\\
3.1	0\\
3.1025	0\\
3.105	0\\
3.1075	1\\
3.11	1\\
3.1125	0\\
3.115	0\\
3.1175	0\\
3.12	0\\
3.1225	0\\
3.125	0\\
3.1275	0\\
3.13	0\\
3.1325	0\\
3.135	0\\
3.1375	0\\
3.14	0\\
3.1425	0\\
3.145	0\\
3.1475	0\\
3.15	0\\
3.1525	0\\
3.155	0\\
3.1575	0\\
3.16	0\\
3.1625	0\\
3.165	0\\
3.1675	0\\
3.17	0\\
3.1725	1\\
3.175	1\\
3.1775	1\\
3.18	2\\
3.1825	2\\
3.185	2\\
3.1875	2\\
3.19	2\\
3.1925	2\\
3.195	2\\
3.1975	2\\
3.2	2\\
3.2025	2\\
3.205	2\\
3.2075	2\\
3.21	2\\
3.2125	2\\
3.215	0\\
3.2175	0\\
3.22	0\\
3.2225	0\\
3.225	1\\
3.2275	1\\
3.23	1\\
3.2325	1\\
3.235	1\\
3.2375	1\\
3.24	1\\
3.2425	0\\
3.245	0\\
3.2475	0\\
3.25	0\\
3.2525	0\\
3.255	1\\
3.2575	1\\
3.26	1\\
3.2625	1\\
3.265	1\\
3.2675	1\\
3.27	0\\
3.2725	0\\
3.275	0\\
3.2775	0\\
3.28	0\\
3.2825	0\\
3.285	0\\
3.2875	0\\
3.29	0\\
3.2925	0\\
3.295	0\\
3.2975	0\\
3.3	0\\
3.3025	0\\
3.305	0\\
3.3075	0\\
3.31	0\\
3.3125	0\\
3.315	0\\
3.3175	0\\
3.32	0\\
3.3225	0\\
3.325	0\\
3.3275	0\\
3.33	0\\
3.3325	0\\
3.335	1\\
3.3375	0\\
3.34	0\\
3.3425	0\\
3.345	0\\
3.3475	0\\
3.35	0\\
3.3525	0\\
3.355	0\\
3.3575	0\\
3.36	1\\
3.3625	0\\
3.365	0\\
3.3675	0\\
3.37	0\\
3.3725	1\\
3.375	1\\
3.3775	1\\
3.38	1\\
3.3825	1\\
3.385	1\\
3.3875	1\\
3.39	1\\
3.3925	1\\
3.395	1\\
3.3975	1\\
3.4	1\\
3.4025	1\\
3.405	1\\
3.4075	1\\
3.41	1\\
3.4125	1\\
3.415	1\\
3.4175	1\\
3.42	0\\
3.4225	0\\
3.425	0\\
3.4275	0\\
3.43	0\\
3.4325	0\\
3.435	0\\
3.4375	0\\
3.44	0\\
3.4425	0\\
3.445	0\\
3.4475	1\\
3.45	1\\
3.4525	1\\
3.455	0\\
3.4575	0\\
3.46	0\\
3.4625	0\\
3.465	0\\
3.4675	0\\
3.47	0\\
3.4725	0\\
3.475	1\\
3.4775	0\\
3.48	0\\
3.4825	0\\
3.485	0\\
3.4875	0\\
3.49	0\\
3.4925	0\\
3.495	0\\
3.4975	0\\
3.5	0\\
3.5025	0\\
3.505	0\\
3.5075	1\\
3.51	1\\
3.5125	1\\
3.515	0\\
3.5175	1\\
3.52	0\\
3.5225	1\\
3.525	0\\
3.5275	1\\
3.53	1\\
3.5325	0\\
3.535	0\\
3.5375	1\\
3.54	1\\
3.5425	1\\
3.545	1\\
3.5475	1\\
3.55	1\\
3.5525	0\\
3.555	1\\
3.5575	1\\
3.56	0\\
3.5625	0\\
3.565	0\\
3.5675	0\\
3.57	1\\
3.5725	1\\
3.575	1\\
3.5775	1\\
3.58	0\\
3.5825	1\\
3.585	1\\
3.5875	1\\
3.59	1\\
3.5925	1\\
3.595	1\\
3.5975	1\\
3.6	0\\
3.6025	1\\
3.605	1\\
3.6075	0\\
3.61	1\\
3.6125	1\\
3.615	1\\
3.6175	1\\
3.62	0\\
3.6225	2\\
3.625	0\\
3.6275	0\\
3.63	0\\
3.6325	1\\
3.635	0\\
3.6375	0\\
3.64	1\\
3.6425	1\\
3.645	0\\
3.6475	1\\
3.65	1\\
3.6525	1\\
3.655	1\\
3.6575	1\\
3.66	0\\
3.6625	1\\
3.665	1\\
3.6675	1\\
3.67	1\\
3.6725	1\\
3.675	1\\
3.6775	1\\
3.68	0\\
3.6825	1\\
3.685	1\\
3.6875	1\\
3.69	0\\
3.6925	1\\
3.695	1\\
3.6975	1\\
3.7	0\\
3.7025	0\\
3.705	0\\
3.7075	0\\
3.71	0\\
3.7125	0\\
3.715	0\\
3.7175	0\\
3.72	0\\
3.7225	0\\
3.725	0\\
3.7275	0\\
3.73	1\\
3.7325	0\\
3.735	0\\
3.7375	1\\
3.74	1\\
3.7425	1\\
3.745	0\\
3.7475	0\\
3.75	0\\
3.7525	0\\
3.755	0\\
3.7575	0\\
3.76	0\\
3.7625	0\\
3.765	0\\
3.7675	0\\
3.77	1\\
3.7725	0\\
3.775	0\\
3.7775	0\\
3.78	1\\
3.7825	0\\
3.785	0\\
3.7875	0\\
3.79	0\\
3.7925	0\\
3.795	0\\
3.7975	0\\
3.8	1\\
3.8025	1\\
3.805	1\\
3.8075	1\\
3.81	1\\
3.8125	1\\
3.815	1\\
3.8175	1\\
3.82	1\\
3.8225	1\\
3.825	1\\
3.8275	1\\
3.83	2\\
3.8325	1\\
3.835	2\\
3.8375	1\\
3.84	1\\
3.8425	1\\
3.845	1\\
3.8475	1\\
3.85	1\\
3.8525	1\\
3.855	1\\
3.8575	1\\
3.86	1\\
3.8625	1\\
3.865	1\\
3.8675	1\\
3.87	1\\
3.8725	1\\
3.875	1\\
3.8775	1\\
3.88	1\\
3.8825	1\\
3.885	1\\
3.8875	1\\
3.89	1\\
3.8925	1\\
3.895	1\\
3.8975	1\\
3.9	1\\
3.9025	1\\
3.905	1\\
3.9075	1\\
3.91	1\\
3.9125	1\\
3.915	1\\
3.9175	1\\
3.92	1\\
3.9225	1\\
3.925	1\\
3.9275	1\\
3.93	1\\
3.9325	1\\
3.935	1\\
3.9375	1\\
3.94	1\\
3.9425	1\\
3.945	1\\
3.9475	1\\
3.95	1\\
3.9525	1\\
3.955	1\\
3.9575	1\\
3.96	1\\
3.9625	1\\
3.965	1\\
3.9675	1\\
3.97	1\\
3.9725	1\\
3.975	1\\
3.9775	1\\
3.98	1\\
3.9825	1\\
3.985	1\\
3.9875	1\\
3.99	1\\
3.9925	1\\
3.995	0\\
3.9975	0\\
4	0\\
4.0025	0\\
4.005	0\\
4.0075	1\\
4.01	1\\
4.0125	1\\
4.015	1\\
4.0175	1\\
4.02	1\\
4.0225	1\\
4.025	1\\
4.0275	1\\
4.03	1\\
4.0325	1\\
4.035	1\\
4.0375	1\\
4.04	1\\
4.0425	1\\
4.045	1\\
4.0475	1\\
4.05	1\\
4.0525	1\\
4.055	1\\
4.0575	1\\
4.06	1\\
4.0625	1\\
4.065	1\\
4.0675	1\\
4.07	1\\
4.0725	0\\
4.075	1\\
4.0775	0\\
4.08	0\\
4.0825	0\\
4.085	1\\
4.0875	0\\
4.09	0\\
4.0925	0\\
4.095	0\\
4.0975	0\\
4.1	0\\
4.1025	0\\
4.105	0\\
4.1075	0\\
4.11	0\\
4.1125	0\\
4.115	0\\
4.1175	1\\
4.12	1\\
4.1225	0\\
4.125	1\\
4.1275	0\\
4.13	0\\
4.1325	1\\
4.135	1\\
4.1375	1\\
4.14	0\\
4.1425	1\\
4.145	1\\
4.1475	1\\
4.15	1\\
4.1525	1\\
4.155	1\\
4.1575	1\\
4.16	1\\
4.1625	1\\
4.165	1\\
4.1675	1\\
4.17	1\\
4.1725	1\\
4.175	1\\
4.1775	1\\
4.18	0\\
4.1825	0\\
4.185	0\\
4.1875	1\\
4.19	1\\
4.1925	1\\
4.195	1\\
4.1975	1\\
4.2	1\\
4.2025	1\\
4.205	1\\
4.2075	1\\
4.21	1\\
4.2125	1\\
4.215	1\\
4.2175	1\\
4.22	0\\
4.2225	1\\
4.225	1\\
4.2275	1\\
4.23	0\\
4.2325	0\\
4.235	1\\
4.2375	1\\
4.24	1\\
4.2425	0\\
4.245	0\\
4.2475	1\\
4.25	0\\
4.2525	0\\
4.255	2\\
4.2575	1\\
4.26	1\\
4.2625	0\\
4.265	0\\
4.2675	0\\
4.27	0\\
4.2725	0\\
4.275	0\\
4.2775	0\\
4.28	1\\
4.2825	1\\
4.285	1\\
4.2875	1\\
4.29	1\\
4.2925	0\\
4.295	0\\
4.2975	0\\
4.3	0\\
4.3025	0\\
4.305	0\\
4.3075	0\\
4.31	0\\
4.3125	0\\
4.315	0\\
4.3175	1\\
4.32	1\\
4.3225	0\\
4.325	0\\
4.3275	0\\
4.33	0\\
4.3325	0\\
4.335	0\\
4.3375	1\\
4.34	1\\
4.3425	1\\
4.345	1\\
4.3475	1\\
4.35	0\\
4.3525	0\\
4.355	0\\
4.3575	0\\
4.36	0\\
4.3625	0\\
4.365	0\\
4.3675	1\\
4.37	1\\
4.3725	1\\
4.375	1\\
4.3775	1\\
4.38	1\\
4.3825	1\\
4.385	1\\
4.3875	1\\
4.39	1\\
4.3925	0\\
4.395	0\\
4.3975	1\\
4.4	1\\
4.4025	0\\
4.405	0\\
4.4075	0\\
4.41	0\\
4.4125	1\\
4.415	0\\
4.4175	1\\
4.42	0\\
4.4225	0\\
4.425	0\\
4.4275	0\\
4.43	0\\
4.4325	0\\
4.435	0\\
4.4375	0\\
4.44	0\\
4.4425	1\\
4.445	1\\
4.4475	1\\
4.45	1\\
4.4525	1\\
4.455	1\\
4.4575	1\\
4.46	1\\
4.4625	1\\
4.465	1\\
4.4675	1\\
4.47	0\\
4.4725	1\\
4.475	1\\
4.4775	1\\
4.48	1\\
4.4825	0\\
4.485	1\\
4.4875	1\\
4.49	1\\
4.4925	1\\
4.495	1\\
4.4975	0\\
4.5	0\\
4.5025	0\\
4.505	0\\
4.5075	1\\
4.51	1\\
4.5125	0\\
4.515	0\\
4.5175	1\\
4.52	1\\
4.5225	1\\
4.525	0\\
4.5275	0\\
4.53	0\\
4.5325	0\\
4.535	0\\
4.5375	1\\
4.54	1\\
4.5425	1\\
4.545	1\\
4.5475	0\\
4.55	0\\
4.5525	0\\
4.555	0\\
4.5575	0\\
4.56	0\\
4.5625	1\\
4.565	1\\
4.5675	1\\
4.57	1\\
4.5725	1\\
4.575	1\\
4.5775	1\\
4.58	1\\
4.5825	1\\
4.585	1\\
4.5875	1\\
4.59	0\\
4.5925	1\\
4.595	1\\
4.5975	1\\
4.6	1\\
4.6025	0\\
4.605	1\\
4.6075	1\\
4.61	1\\
4.6125	1\\
4.615	1\\
4.6175	1\\
4.62	1\\
4.6225	1\\
4.625	1\\
4.6275	0\\
4.63	0\\
4.6325	1\\
4.635	1\\
4.6375	0\\
4.64	0\\
4.6425	1\\
4.645	1\\
4.6475	1\\
4.65	1\\
4.6525	1\\
4.655	1\\
4.6575	1\\
4.66	1\\
4.6625	0\\
4.665	0\\
4.6675	0\\
4.67	1\\
4.6725	1\\
4.675	1\\
4.6775	1\\
4.68	1\\
4.6825	1\\
4.685	1\\
4.6875	1\\
4.69	1\\
4.6925	1\\
4.695	1\\
4.6975	1\\
4.7	1\\
4.7025	0\\
4.705	1\\
4.7075	1\\
4.71	0\\
4.7125	0\\
4.715	1\\
4.7175	1\\
4.72	0\\
4.7225	0\\
4.725	1\\
4.7275	1\\
4.73	1\\
4.7325	1\\
4.735	1\\
4.7375	1\\
4.74	1\\
4.7425	1\\
4.745	1\\
4.7475	0\\
4.75	1\\
4.7525	1\\
4.755	1\\
4.7575	1\\
4.76	0\\
4.7625	1\\
4.765	1\\
4.7675	1\\
4.77	1\\
4.7725	1\\
4.775	1\\
4.7775	1\\
4.78	1\\
4.7825	1\\
4.785	1\\
4.7875	1\\
4.79	1\\
4.7925	1\\
4.795	1\\
4.7975	1\\
4.8	1\\
4.8025	1\\
4.805	1\\
4.8075	1\\
4.81	1\\
4.8125	1\\
4.815	1\\
4.8175	0\\
4.82	0\\
4.8225	1\\
4.825	1\\
4.8275	1\\
4.83	0\\
4.8325	0\\
4.835	0\\
4.8375	0\\
4.84	0\\
4.8425	0\\
4.845	0\\
4.8475	0\\
4.85	1\\
4.8525	1\\
4.855	0\\
4.8575	0\\
4.86	0\\
4.8625	0\\
4.865	0\\
4.8675	0\\
4.87	1\\
4.8725	1\\
4.875	1\\
4.8775	1\\
4.88	1\\
4.8825	1\\
4.885	0\\
4.8875	0\\
4.89	0\\
4.8925	0\\
4.895	1\\
4.8975	0\\
4.9	0\\
4.9025	0\\
4.905	0\\
4.9075	0\\
4.91	1\\
4.9125	1\\
4.915	1\\
4.9175	1\\
4.92	1\\
4.9225	0\\
4.925	0\\
4.9275	0\\
4.93	0\\
4.9325	0\\
4.935	0\\
4.9375	0\\
4.94	0\\
4.9425	0\\
4.945	0\\
4.9475	0\\
4.95	0\\
4.9525	0\\
4.955	0\\
4.9575	0\\
4.96	0\\
4.9625	0\\
4.965	0\\
4.9675	0\\
4.97	0\\
4.9725	0\\
4.975	0\\
4.9775	0\\
4.98	0\\
4.9825	0\\
4.985	0\\
4.9875	1\\
4.99	1\\
4.9925	0\\
4.995	0\\
4.9975	0\\
};
\addlegendentry{$\epsilon_u = 10^{-11}$}

\end{axis}

\end{tikzpicture}%

%% file: figures/mc_AC.tex
% This file was created with tikzplotlib v0.10.1.
\definecolor{mycolor1}{rgb}{0.00000,0.44700,0.74100}%
\definecolor{mycolor2}{rgb}{0.85000,0.32500,0.09800}%
\definecolor{mycolor3}{rgb}{0.5,0.5,0.5}
\begin{tikzpicture}[scale=0.8]

\begin{axis}[%
width=1.5in,
height=1.5in,
at={(1.011in,0.642in)},
scale only axis,
xmin=0,
xmax=200,
xlabel style={font=\color{white!15!black}},
xlabel={Time},
%ymode=log,
ymin=0,
ymax=60,
yminorticks=true,
ylabel style={font=\color{white!15!black}},
axis background/.style={fill=white},
legend style={fill opacity=0.8, 
font=\tiny,
  draw opacity=1,
  text opacity=1,
  at={(0.02,0.87)},
  anchor=west}
  % draw=lightgray204}
]

\addplot [line width=1.2pt,  mycolor1]
table {%
0 7
0.5 0
1 1
1.5 1
2 1
2.5 1
3 0
3.5 1
4 1
4.5 1
5 1
5.5 1
6 1
6.5 1
7 1
7.5 0
8 2
8.5 0
9 3
9.5 1
10 3
10.5 1
11 0
11.5 0
12 1
12.5 0
13 0
13.5 0
14 1
14.5 1
15 13
15.5 0
16 7
16.5 2
17 2
17.5 1
18 0
18.5 13
19 11
19.5 0
20 0
20.5 8
21 1
21.5 1
22 2
22.5 0
23 1
23.5 1
24 5
24.5 1
25 7
25.5 4
26 2
26.5 1
27 1
27.5 2
28 2
28.5 2
29 1
29.5 1
30 0
30.5 6
31 2
31.5 10
32 1
32.5 6
33 0
33.5 1
34 7
34.5 2
35 5
35.5 46
36 1
36.5 1
37 1
37.5 13
38 5
38.5 10
39 1
39.5 1
40 2
40.5 4
41 12
41.5 8
42 7
42.5 1
43 1
43.5 8
44 0
44.5 1
45 1
45.5 1
46 6
46.5 0
47 4
47.5 0
48 7
48.5 4
49 18
49.5 6
50 9
50.5 11
51 2
51.5 5
52 8
52.5 15
53 1
53.5 0
54 6
54.5 1
55 0
55.5 8
56 0
56.5 6
57 2
57.5 0
58 0
58.5 11
59 1
59.5 9
60 1
60.5 0
61 1
61.5 6
62 1
62.5 1
63 2
63.5 4
64 2
64.5 2
65 2
65.5 8
66 11
66.5 5
67 5
67.5 5
68 2
68.5 6
69 0
69.5 0
70 1
70.5 1
71 1
71.5 1
72 1
72.5 0
73 0
73.5 1
74 0
74.5 1
75 0
75.5 0
76 5
76.5 1
77 1
77.5 2
78 1
78.5 4
79 9
79.5 0
80 1
80.5 0
81 0
81.5 2
82 1
82.5 1
83 0
83.5 2
84 0
84.5 1
85 1
85.5 1
86 2
86.5 0
87 6
87.5 4
88 0
88.5 0
89 2
89.5 1
90 1
90.5 8
91 0
91.5 0
92 5
92.5 5
93 2
93.5 1
94 7
94.5 0
95 0
95.5 0
96 0
96.5 2
97 0
97.5 0
98 1
98.5 2
99 1
99.5 1
100 0
100.5 0
101 0
101.5 1
102 1
102.5 14
103 0
103.5 1
104 1
104.5 11
105 2
105.5 10
106 1
106.5 1
107 0
107.5 3
108 0
108.5 1
109 1
109.5 0
110 12
110.5 1
111 0
111.5 1
112 9
112.5 0
113 1
113.5 0
114 0
114.5 0
115 0
115.5 11
116 11
116.5 11
117 0
117.5 0
118 0
118.5 0
119 1
119.5 1
120 0
120.5 4
121 2
121.5 0
122 1
122.5 9
123 1
123.5 0
124 6
124.5 1
125 6
125.5 0
126 0
126.5 1
127 1
127.5 2
128 0
128.5 1
129 4
129.5 5
130 1
130.5 1
131 0
131.5 14
132 4
132.5 1
133 9
133.5 0
134 0
134.5 0
135 11
135.5 0
136 0
136.5 1
137 4
137.5 1
138 0
138.5 0
139 1
139.5 4
140 0
140.5 1
141 1
141.5 1
142 1
142.5 1
143 8
143.5 2
144 3
144.5 2
145 0
145.5 2
146 4
146.5 3
147 0
147.5 1
148 0
148.5 2
149 12
149.5 0
150 1
150.5 1
151 1
151.5 2
152 2
152.5 0
153 0
153.5 2
154 2
154.5 1
155 2
155.5 4
156 1
156.5 11
157 1
157.5 2
158 0
158.5 7
159 1
159.5 0
160 1
160.5 6
161 1
161.5 2
162 1
162.5 8
163 1
163.5 2
164 1
164.5 0
165 1
165.5 1
166 1
166.5 5
167 4
167.5 9
168 1
168.5 0
169 1
169.5 1
170 7
170.5 1
171 1
171.5 1
172 1
172.5 2
173 6
173.5 0
174 2
174.5 5
175 2
175.5 5
176 2
176.5 0
177 1
177.5 1
178 1
178.5 1
179 0
179.5 1
180 1
180.5 1
181 1
181.5 1
182 1
182.5 1
183 1
183.5 1
184 0
184.5 1
185 4
185.5 9
186 0
186.5 6
187 4
187.5 0
188 0
188.5 1
189 8
189.5 6
190 1
190.5 4
191 5
191.5 6
192 1
192.5 1
193 1
193.5 10
194 1
194.5 1
195 1
195.5 7
196 0
196.5 6
197 1
197.5 5
198 0
198.5 8
199 2
199.5 3
};
%\addlegendentry{$\epsilon_u = 10^{-8}$}

\addplot [line width=1.2pt,  mycolor2,opacity=0.5]
table {%
0 10
0.5 1
1 1
1.5 1
2 1
2.5 0
3 1
3.5 1
4 1
4.5 1
5 1
5.5 0
6 0
6.5 0
7 0
7.5 0
8 0
8.5 1
9 1
9.5 2
10 0
10.5 1
11 8
11.5 1
12 1
12.5 3
13 0
13.5 6
14 1
14.5 1
15 1
15.5 1
16 2
16.5 0
17 2
17.5 1
18 0
18.5 0
19 1
19.5 1
20 0
20.5 11
21 2
21.5 4
22 1
22.5 0
23 2
23.5 2
24 1
24.5 6
25 1
25.5 8
26 0
26.5 7
27 1
27.5 0
28 0
28.5 0
29 0
29.5 0
30 0
30.5 1
31 1
31.5 2
32 1
32.5 2
33 2
33.5 1
34 8
34.5 2
35 1
35.5 1
36 3
36.5 1
37 1
37.5 0
38 0
38.5 0
39 1
39.5 0
40 0
40.5 2
41 0
41.5 1
42 16
42.5 1
43 7
43.5 2
44 1
44.5 1
45 2
45.5 2
46 1
46.5 1
47 1
47.5 1
48 1
48.5 2
49 1
49.5 5
50 0
50.5 8
51 1
51.5 11
52 0
52.5 0
53 4
53.5 1
54 0
54.5 2
55 2
55.5 3
56 4
56.5 1
57 1
57.5 2
58 0
58.5 6
59 0
59.5 2
60 12
60.5 2
61 1
61.5 0
62 0
62.5 1
63 0
63.5 0
64 0
64.5 8
65 1
65.5 1
66 2
66.5 2
67 0
67.5 2
68 0
68.5 1
69 0
69.5 3
70 4
70.5 1
71 0
71.5 0
72 13
72.5 1
73 0
73.5 0
74 5
74.5 4
75 1
75.5 0
76 1
76.5 0
77 4
77.5 8
78 0
78.5 7
79 7
79.5 4
80 1
80.5 0
81 1
81.5 0
82 2
82.5 0
83 2
83.5 3
84 0
84.5 2
85 0
85.5 0
86 1
86.5 3
87 0
87.5 0
88 1
88.5 1
89 0
89.5 0
90 1
90.5 0
91 0
91.5 0
92 0
92.5 0
93 1
93.5 0
94 15
94.5 0
95 1
95.5 1
96 0
96.5 1
97 1
97.5 0
98 1
98.5 0
99 1
99.5 4
100 1
100.5 8
101 0
101.5 2
102 1
102.5 1
103 0
103.5 6
104 15
104.5 0
105 0
105.5 0
106 0
106.5 11
107 0
107.5 0
108 0
108.5 5
109 2
109.5 0
110 1
110.5 0
111 0
111.5 1
112 0
112.5 1
113 1
113.5 0
114 9
114.5 0
115 1
115.5 3
116 7
116.5 6
117 1
117.5 1
118 1
118.5 5
119 1
119.5 1
120 0
120.5 1
121 0
121.5 0
122 1
122.5 0
123 1
123.5 1
124 0
124.5 1
125 2
125.5 1
126 1
126.5 0
127 1
127.5 2
128 1
128.5 0
129 0
129.5 0
130 1
130.5 1
131 0
131.5 9
132 1
132.5 10
133 1
133.5 3
134 4
134.5 1
135 1
135.5 0
136 1
136.5 6
137 7
137.5 0
138 1
138.5 0
139 1
139.5 0
140 1
140.5 1
141 1
141.5 1
142 4
142.5 0
143 8
143.5 8
144 8
144.5 5
145 1
145.5 0
146 1
146.5 0
147 10
147.5 1
148 4
148.5 1
149 8
149.5 5
150 9
150.5 1
151 9
151.5 6
152 0
152.5 0
153 8
153.5 10
154 7
154.5 7
155 1
155.5 1
156 0
156.5 0
157 1
157.5 1
158 1
158.5 1
159 7
159.5 7
160 1
160.5 4
161 5
161.5 4
162 2
162.5 1
163 0
163.5 8
164 8
164.5 1
165 0
165.5 5
166 0
166.5 1
167 1
167.5 1
168 1
168.5 5
169 2
169.5 9
170 8
170.5 1
171 0
171.5 1
172 9
172.5 13
173 1
173.5 6
174 2
174.5 1
175 1
175.5 10
176 1
176.5 10
177 0
177.5 5
178 1
178.5 0
179 1
179.5 1
180 4
180.5 10
181 5
181.5 8
182 7
182.5 1
183 5
183.5 1
184 0
184.5 4
185 1
185.5 0
186 1
186.5 1
187 0
187.5 0
188 0
188.5 0
189 0
189.5 2
190 0
190.5 3
191 1
191.5 20
192 1
192.5 7
193 11
193.5 8
194 6
194.5 1
195 1
195.5 0
196 1
196.5 1
197 1
197.5 4
198 1
198.5 0
199 0
199.5 1
};
%\addlegendentry{$\epsilon_u = 10^{-11}$}
\end{axis}

\end{tikzpicture}

%% file: figures/mc_KdV.tex
% This file was created with tikzplotlib v0.10.1.
\definecolor{mycolor1}{rgb}{0.00000,0.44700,0.74100}%
\definecolor{mycolor2}{rgb}{0.85000,0.32500,0.09800}%
\definecolor{mycolor3}{rgb}{0.5,0.5,0.5}
\begin{tikzpicture}[scale=0.8]

\begin{axis}[%
width=1.5in,
height=1.5in,
at={(1.011in,0.642in)},
scale only axis,
xmin=0,
xmax=1,
xlabel style={font=\color{white!15!black}},
xlabel={Time},
%ymode=log,
ymin=0,
ymax=60,
yminorticks=true,
ylabel style={font=\color{white!15!black}},
axis background/.style={fill=white},
legend style={fill opacity=0.8, 
font=\tiny,
  draw opacity=1,
  text opacity=1,
  at={(0.02,0.87)},
  anchor=west}
  % draw=lightgray204}
]

\addplot [line width=1.2pt,  mycolor1]
table {%
0 1
0.0002 1
0.0004 1
0.0006 0
0.0016 0
0.0018 1
0.0054 1
0.0056 0
0.0058 1
0.0132 1
0.0134 1
0.0136 0
0.0186 0
0.0188 1
0.019 1
0.0344 1
0.0346 0
0.0348 0
0.051 0
0.0512 2
0.0514 0
0.0516 0
0.0598 0
0.06 1
0.0602 0
0.0604 0
0.075 0
0.0752 1
0.0754 1
0.0968 1
0.097 1
0.0972 0
0.0974 0
0.18 0
0.1802 2
0.1804 0
0.1806 0
0.1808 0
0.181 0
0.1812 0
0.1814 1
0.1816 1
0.1818 1
0.2016 1
0.2018 1
0.202 0
0.2022 0
0.498 0
0.4982 5
0.4984 0
0.4986 0
0.4988 0
0.499 1
0.4992 1
0.4994 1
0.552 1
0.5522 0
0.5524 0
0.567 0
0.5672 1
0.5674 1
0.5676 1
0.5678 1
0.568 1
0.5682 1
0.5684 0
0.5686 0
0.5688 0
0.569 0
0.6366 0
0.6368 0
0.637 1
0.6372 1
0.6374 1
0.668 1
0.6682 1
0.6684 1
0.6686 1
0.6688 1
0.669 1
0.6692 0
0.6694 0
0.68 0
0.6802 0
0.6804 0
0.6806 1
0.6808 1
0.681 1
0.6812 1
0.6814 1
0.6816 0
0.6818 0
0.6838 0
0.684 0
0.6842 1
0.6844 1
0.6846 1
0.6848 1
0.685 1
0.6852 0
0.6854 0
0.6856 0
0.6858 0
0.686 0
0.6862 0
0.6864 0
0.6866 0
0.6868 0
0.687 0
0.6872 1
0.6874 1
0.6876 1
0.6878 1
0.688 1
0.6882 1
0.6884 1
0.6886 1
0.6888 1
0.689 1
0.6892 1
0.6894 1
0.6896 1
0.6898 1
0.69 1
0.6902 0
0.6904 0
0.6906 0
0.6908 0
0.691 0
0.6912 0
0.6914 5
0.6916 1
0.6918 1
0.692 1
0.7238 1
0.724 1
0.7242 1
0.7244 1
0.7482 1
0.7484 1
0.7486 55
0.7488 1
0.749 1
0.7492 1
0.7494 1
0.7496 1
0.7498 1
0.75 0
0.7502 0
0.7504 0
0.7618 0
0.762 0
0.7622 1
0.7624 1
0.7626 1
0.7642 1
0.7644 1
0.7646 1
0.7648 1
0.765 1
0.7652 1
0.7654 0
0.7656 0
0.7658 0
0.766 0
0.7932 0
0.7934 0
0.7936 1
0.7938 0
0.794 0
0.8388 0
0.839 0
0.8392 0
0.8394 1
0.8396 1
0.8398 1
0.84 1
0.8402 1
0.8404 1
0.8406 1
0.8408 1
0.841 1
0.8412 1
0.8414 1
0.8416 1
0.8418 6
0.842 1
0.8422 1
0.8424 1
0.87 1
0.8702 0
0.8704 0
0.8706 0
0.8708 5
0.871 1
0.8712 1
0.8714 1
0.8716 1
0.933 1
0.9332 4
0.9334 1
0.9336 1
0.9338 1
0.9676 1
0.9678 12
0.968 1
0.9682 1
0.9684 1
0.9882 1
0.9884 0
0.9886 0
0.9888 0
0.989 0
0.9892 1
0.9894 1
0.9896 1
0.9898 1
0.9922 1
0.9924 1
0.9926 1
0.9928 1
0.993 1
0.9932 1
0.9934 0
0.9936 0
0.9938 0
0.994 0
0.9942 0
0.9944 0
0.9946 0
0.9948 0
0.9998 0
};
%\addlegendentry{$\epsilon_u = 10^{-8}$}

\addplot [line width=1.2pt,  mycolor2,opacity=0.5]
table {%
0 1
0.0085 1
0.009 1
0.0095 0
0.01 0
0.0105 0
0.011 0
0.0115 0
0.012 1
0.0125 1
0.013 1
0.0135 1
0.014 1
0.0145 1
0.015 1
0.0155 1
0.016 1
0.0165 1
0.017 0
0.0175 0
0.018 0
0.0185 0
0.019 0
0.0195 0
0.02 0
0.0205 0
0.021 0
0.0215 0
0.022 0
0.0225 0
0.023 0
0.0235 0
0.024 0
0.0245 0
0.025 0
0.0255 1
0.026 1
0.0265 1
0.027 1
0.0275 1
0.028 1
0.0285 1
0.0395 1
0.04 2
0.0465 2
0.047 1
0.0475 1
0.048 1
0.1145 1
0.115 0
0.1155 0
0.116 0
0.1165 1
0.117 1
0.1175 0
0.118 0
0.1185 0
0.119 0
0.1195 0
0.12 0
0.1205 0
0.121 0
0.1215 0
0.122 1
0.1225 1
0.123 1
0.1235 1
0.124 0
0.1245 0
0.131 0
0.1315 0
0.132 1
0.1325 1
0.133 1
0.1565 1
0.157 0
0.1575 0
0.158 0
0.1625 0
0.163 1
0.17 1
0.1705 1
0.171 1
0.1715 0
0.172 0
0.1725 0
0.173 0
0.1735 0
0.174 0
0.1745 0
0.175 0
0.1755 1
0.176 1
0.1765 1
0.177 1
0.1775 0
0.198 0
0.1985 1
0.199 0
0.1995 1
0.2035 1
0.204 0
0.2045 1
0.205 1
0.2125 1
0.213 0
0.2135 0
0.214 0
0.2145 1
0.218 1
0.2185 0
0.219 0
0.2195 0
0.22 0
0.2205 1
0.417 1
0.6 1
0.6005 1
0.601 1
0.6015 1
0.602 1
0.6025 1
0.603 1
0.6035 1
0.604 1
0.7 1
0.7005 1
0.701 1
0.7015 1
0.702 1
0.7025 0
0.703 0
0.7035 0
0.704 1
0.7045 1
0.705 1
0.73 1
0.7305 0
0.731 0
0.7315 1
0.732 1
0.7575 1
0.758 1
0.7585 1
0.759 0
0.7595 0
0.76 0
0.7605 0
0.761 0
0.7615 0
0.762 1
0.7625 1
0.763 1
0.7755 1
0.776 1
0.7765 2
0.777 2
0.7775 2
0.778 2
0.7785 2
0.779 2
0.7795 1
0.78 1
0.7805 1
0.781 1
0.813 1
0.8135 1
0.814 2
0.8145 2
0.815 1
0.8155 1
0.816 1
0.8475 1
0.848 2
0.8485 2
0.849 2
0.8495 2
0.85 2
0.8505 2
0.851 2
0.8515 2
0.852 1
0.8525 1
0.853 1
0.8535 1
0.854 1
0.8545 1
0.855 2
0.8555 1
0.856 1
0.8565 1
0.857 1
0.8575 1
0.858 1
0.8585 1
0.859 1
0.8595 1
0.86 1
0.8605 1
0.861 2
0.8615 1
0.862 1
0.8625 1
0.863 1
0.8635 1
0.864 1
0.911 1
0.9115 0
0.912 0
0.9125 0
0.913 0
0.9135 0
0.914 1
0.9145 1
0.915 1
0.9155 1
0.916 1
0.9165 1
0.917 1
0.9175 1
0.918 1
0.9185 2
0.919 1
0.9195 1
0.92 0
0.9205 1
0.921 1
0.9215 1
0.922 0
0.9225 0
0.923 1
0.9235 1
0.924 1
0.9245 1
0.925 1
0.9255 1
0.926 1
0.9265 1
0.927 1
0.9275 1
0.928 1
0.9285 1
0.929 1
0.9295 1
0.93 0
0.9305 0
0.931 2
0.9315 1
0.932 0
0.9325 2
0.933 1
0.9335 1
0.934 1
0.9345 1
0.935 1
0.9355 1
0.936 1
0.9365 1
0.937 1
0.9375 1
0.938 1
0.9385 1
0.939 1
0.9395 1
0.94 1
0.9405 1
0.955 1
0.9555 1
0.956 12
0.9565 1
0.957 1
0.973 1
0.9735 0
0.974 0
0.9745 0
0.975 0
0.9755 0
0.976 1
0.9765 1
0.977 1
0.9775 0
0.978 1
0.9785 1
0.979 1
0.9795 1
0.98 0
0.9805 6
0.981 0
0.9815 0
0.982 0
0.9825 0
0.983 1
0.9835 1
0.984 1
0.9845 1
0.985 1
0.9855 1
0.986 1
0.9865 1
0.987 0
0.9875 0
0.988 0
0.9885 0
0.989 0
0.9895 1
0.99 1
0.9905 1
0.991 1
0.9915 1
0.992 1
0.9925 1
0.993 1
0.9935 1
0.994 1
0.9945 1
0.995 1
0.9955 1
0.996 1
0.9965 1
0.997 2
0.9975 2
0.998 1
0.9985 1
0.999 1
0.9995 1
};
%\addlegendentry{$\epsilon_u = 10^{-11}$}
\end{axis}

\end{tikzpicture}